\documentclass[11pt,a4paper,english, reqno]{amsart}
\usepackage[english]{babel}
\usepackage{amsthm}
\usepackage{amsmath}
\usepackage{amssymb}
\usepackage{graphicx}
\usepackage{dsfont}
\usepackage[bf]{caption}
\usepackage{parskip}
\usepackage{mathrsfs} %pour mathscr

\newtheorem{thm}{Theorem}[section]
\newtheorem{cor}[thm]{Corollary}
\newtheorem{lem}[thm]{Lemma}
\newtheorem{prop}[thm]{Proposition}

\newtheorem{remark}[thm]{Remark}

\newtheorem{deff}[thm]{Definition}

\newcommand{\bb} {\textsc{b}}
\newcommand{\e} {\varepsilon}
\newcommand{\E} { \textnormal{\textsf{E}}}
\newcommand{\EE} { \mathbb{E} }
\newcommand{\G}{{\mathcal G}}
\newcommand{\I}{{\mathcal I}}
\newcommand{\N} { \mathbb{N} }
\def\o{\omega}

\newcommand{\PP} { \mathbb{P} }
\newcommand{\p} {\textnormal{\textsf{P}}}

\newcommand{\R} { \mathbb{R} }
\newcommand{\SSS}{{\mathcal S}}

\newcommand{\Z} { \mathbb{Z} }
\newcommand{\VVV} {\mathscr{V}}

\def\po{P_\omega}
\def\eo{\E_\omega}
\def\dd{\textnormal{d}}
\def\un{{\bf 1}}

\allowdisplaybreaks

\begin{document}

\title[Annealed local limit theorem for Sinai's random walk]
{Annealed local limit theorem for Sinai's random walk in random environment}
\date{\today}

\author{Alexis Devulder}
\address{Universit\'e Paris-Saclay, UVSQ, CNRS, Laboratoire de Math\'ematiques de Versailles, 78000, Versailles, France.}
\email{devulder@math.uvsq.fr }

\subjclass[2010]{60K37, 60F15, 60G50, 60K05, 60K50}
\keywords{Sinai's walk, random walk, random environment,
local limit theorem, localization,
path decomposition, random walks conditioned to stay positive, renewal theory}

\begin{abstract}
We consider Sinai's random walk in random environment $(S_n)_{n\in\N}$.
We prove a local limit theorem for $(S_n)_{n\in\N}$ under the annealed law $\PP$.
% that is,
%we provide an
As a consequence, we get an equivalent for the annealed probability $\PP(S_n=z_n)$ as $n$ goes to infinity, when $z_n=O\big((\log n)^2\big)$.
%We also prove a local limit theorem for the favorite location of Sinai's walk at time $n$, generally denoted by $b_{\log n}$.
To this aim, we develop a path decomposition for the potential of Sinai's walk, that is, for some random walks with i.i.d. increments.
The proof also relies on
%a path decomposition for the potential of Sinai's walk,
renewal theory, a coupling argument,
a very careful analysis
%description
of the environments and trajectories of Sinai's walk satisfying $S_n=z_n$,
%and of the different kinds of environments in which such trajectories have large or negligible quenched probability,
and on precise estimates for random walks conditioned to stay positive or nonnegative.
\end{abstract}

\maketitle

{\bf Table of contents}
\tableofcontents

\section{Introduction and statement of the main results}\label{Sect_Intro}

\subsection{Presentation of the model}
%\subsection{Presentation of the model and main results}

We consider a collection $\omega:=(\omega_x)_{x\in\Z}$  of i.i.d. random variables,
taking values in the interval $]0,1[$,
with joint law $\p$. A realization of $\omega$ is called an {\it environment}.
A random walk $(S_k)_{k\in\N}$ in the environment $\omega$ is defined as follows.
Conditionally on $\omega$, $(S_k)_{k\in\N}$ is a Markov chain
starting at $S_0=0$ and such that for every $k\in\N:=\{0,1,2,\dots\}$, $x\in\Z$ and $y\in\Z$,
\begin{equation}\label{eq_Def_Sn}
    \po\big(S_{k+1}=y|S_k=x\big)
=
\left\{
\begin{array}{ll}
\o_x & \text{ if } y=x+1,\\
1-\o_x & \text{ if } y=x-1,\\
0      & \text{ otherwise}.
\end{array}
\right.
\end{equation}
We call $P_\omega$ the {\it quenched law}, and
$S:=(S_k)_k$ is a {\it random walk in random environment (RWRE)}.
The {\it annealed law} is defined as follows:
$$
    \mathbb P[\cdot]:=\int P_\omega[\cdot]\p(\text{d}\omega).
$$
Notice that $\PP$ is not Markovian.
The expectations with respect to $\PP$, $ P_\omega$ and $\p$
are denoted respectively by $\mathbb E$, $E_\omega$ and $\E$.

One dimensional RWRE have many unusual properties, and have attracted much interest from mathematicians and physicists.
For applications in physics and in biology, see e.g. Cocco et al. \cite{CoccoMonasson}, Hughes \cite{Hug}
and more recently the introduction of Padash et al. \cite{Padash}.
Also, (one dimensional) RWRE are used to define or study some other mathematical models, see e.g.
Kochler \cite{These_Jochler} (chapter 3) for random walks in oriented lattices with random environments,
Zindy \cite{Zindy} for random walks in random environments with random scenery.
Aurzada et al. \cite{AurzadaDevulderGuillotinPene} for branching processes in random environments,
and Devulder \cite{Devulder_SPL} for branching random walks in random environments.
%and e.g. Avena et al. \cite{Avena_Chino_daCosta_denHollander}
%for random walks in dynamic random environments.
We refer to R{\'e}v{\'e}sz \cite{Revesz} and Zeitouni \cite{Z01}
for a general account on results on RWRE proved before 2005.
For a statistical point of view, see e.g. Diel et al. \cite{Diel_Stat} and references therein.

We assume that there exists $\e_0\in]0,1/2[$ such that
\begin{equation}\label{eqEllipticity}
    \p[\varepsilon_0\leq \omega_0\leq 1-\varepsilon_0]=1.
\end{equation}
This classical condition is known as the {\it ellipticity} condition.
We introduce $\rho_x:=\frac{1-\omega_x}{\omega_x}$, $x\in\Z$.
Solomon \cite{S75} proved that $(S_k)_k$ is recurrent for almost every environment $\o$ if
\begin{equation}\label{eqRecurrence}
    \E[\log \rho_0]=0,
\end{equation}
and transient for almost every $\o$ otherwise.
Throughout the paper, $\log$ denotes the natural logarithm.
We only consider the recurrent case \eqref{eqRecurrence} in the present paper.
Also, in order to avoid the degenerate case of simple random walks, we assume that
\begin{equation}\label{eq_def_sigma}
    \sigma
:=
    \big(\E\big[(\log \rho_0)^2\big]\big)^{1/2}
>0.
\end{equation}
The asymptotic behaviour of $S$ in the very delicate recurrent case was first analyzed in a celebrated paper of Sinai \cite{S82}.
Indeed, Sinai \cite{S82} showed that under Hypotheses \eqref{eqEllipticity}, \eqref{eqRecurrence} and \eqref{eq_def_sigma},
$S_n$ is localized at time $n$, with large annealed probability, in the neighborhood
of some random quantity $b_{\log n}'$, which depends only on the environment.
More precisely, he proved that  for every $\e>0$,
$$
    \PP\big[|S_n-b_{\log n}'|\leq \e (\log n)^2\big]
\to_{n\to+\infty}
    1.
$$
He also proved that $\sigma^2 b_{\log n}'/(\log n)^2$ converges in law, as $n\to+\infty$, to
some random variable $b_\infty$, which is non degenerate and non gaussian.
As a consequence, Sinai obtained the following convergence in law under the annealed law $\PP$:
\begin{equation*}
%\label{eq_cv_loi_Sinai}
    \frac{\sigma^2}{(\log n)^2}S_n
\to_{n\to+\infty}
    b_\infty.
\end{equation*}
It was proved independently by Kesten \cite{Kesten} and Golosov \cite{Golosov86} that
$\p[b_\infty \in \text{d}x]=\varphi_\infty(x)\text{d}x$, where
\begin{equation}\label{eq_def_phi_infini}
    \varphi_{\infty}(x)
:=
    \frac{2}{\pi}\sum_{k=0}^\infty \frac{(-1)^k}{2k+1}\exp\bigg(-\frac{(2k+1)^2\pi^2}{8}|x|\bigg),
\qquad
    x\in\R.
\end{equation}
This very slow movement of $(S_k)_{k\in\N}$, of order $(\log n)^2$ instead of $\sqrt{n}$ for simple random walks,
is due to the presence of some traps which slow down the walk.
Due to this result proved by Sinai, a random walk in random environment $(S_k)_{k\in\N}$ satisfying
Hypotheses \eqref{eqEllipticity}, \eqref{eqRecurrence} and \eqref{eq_def_sigma} is often called a {\it Sinai walk}.
Some other unusual properties of Sinai's walk are proved e.g. in
Dembo et al. \cite{Dembo_Gantert_Peres_Shi},
Gantert et al. \cite{GPS_2010}, \cite{Gantert_Sshi_2002},
Hu et al. \cite{HuShiLimits}, \cite{Hu_Shi_Problem} and Shi \cite{Shi_98}. See also Shi \cite{S2}
for a general account about Sinai's walk before 2001.

\subsection{Main results}
Throughout the paper, for sequences $(d_n)$ and $(m_n)$ with $m_n\neq 0$ for large $n$,
we write $d_n \sim_{n\to+\infty} m_n$
if $d_n$/$m_n\to 1$ as $n \to +\infty$,
$d_n = o(m_n)$
if $d_n$/$m_n\to 0$ as $n \to +\infty$,
and
$d_n = O(m_n)$
if $\limsup_{n\to+\infty} |d_n/m_n|<\infty$.

Our main result is the following local limit theorem for Sinai's walk $(S_n)_{n\in\N}$ under the annealed law $\PP$:

\begin{thm}\label{Th_Local_Limit_Sinai}
Assume \eqref{eqEllipticity}, \eqref{eqRecurrence} and \eqref{eq_def_sigma}.
%Let $A>0$.
As $n\to+\infty$,
$$
%    \sup_{z\in[-A(\log n)^2, A(\log n)^2]\cap(n+2\Z)}
    \sup_{z\in(2\Z+n)}
    \bigg|
        \PP\big(S_n=z\big)
        -
        \frac{2\sigma^2}{(\log n)^2}\, \varphi_\infty\bigg(\frac{\sigma^2 z}{(\log n)^2}\bigg)
    \bigg|
=
    o\left(\frac{1}{(\log n)^2}\right),
$$
where $2\Z+n$ denotes the set of integers having the same parity as $n$.
\end{thm}

Notice that $S:=(S_k)_{k\in\N}$ only makes $\pm 1$ jumps and starts from $0$ under $\PP$,
so $\PP(S_n=z)=0$ if $n$ and $z$ have different parity.
Since $\varphi_\infty>0$ and is continuous on $\R$, we get in particular:

\begin{cor}
Assume \eqref{eqEllipticity}, \eqref{eqRecurrence} and \eqref{eq_def_sigma}.
Let $(z_n)_{n\in\N}$ be a sequence of integers such that  $z_n=O\big((\log n)^2\big)$ as $n\to+\infty$,
and such that $z_n$ and $n$ have the same parity for every $n\in\N$. Then,
$$
    \PP\big(S_n=z_n\big)
\sim_{n\to+\infty}
    \frac{2\sigma^2}{(\log n)^2}\, \varphi_\infty\bigg(\frac{\sigma^2 z_n}{(\log n)^2}\bigg).
$$
\end{cor}
Also $\sum_{k=0}^\infty \frac{(-1)^k}{2k+1}=\arctan(1)=\pi/4$,
hence $\varphi_\infty(0)=1/2$, so this leads to:

\begin{cor}
Assume \eqref{eqEllipticity}, \eqref{eqRecurrence} and \eqref{eq_def_sigma}. We have,
$$
    \PP\big(S_{2n}=0\big)
\sim_{n\to+\infty}
    \frac{\sigma^2}{(\log n)^2},
$$
and more generally
$
    \PP(S_{2n}=2x)
\sim_{n\to+\infty}
    \frac{\sigma^2}{(\log n)^2}
$
for every fixed $x\in\Z$ since $\varphi_\infty$ is continuous on $\R$.
Also, for every fixed $x\in\Z$,
$$
    \PP\big(S_{2n}=2\big\lfloor (x/2)(\log n)^2\big\rfloor\big)
\sim_{n\to+\infty}
    \frac{2\sigma^2 \varphi_\infty(\sigma^2 x)}{(\log n)^2},
$$
where for $y\in\R$, $\lfloor y\rfloor$ denotes the integer part of $y$.
\end{cor}
In order to prove Theorem \ref{Th_Local_Limit_Sinai},
we introduce in Section \ref{Sect_Potential} (see \eqref{eqDefbh}) a random quantity $b_h$, $h>0$,
depending only on the environment.
%It is defined precisely in Section \ref{Sect_Potential} (see \eqref{eqDefbh}).
It is defined differently from the localization point $b_h'$ introduced by Sinai, but plays a similar role.
%(we believe that $b_h=b_h'$ with large probability, but we do not prove this result which we do not need).
%For the definition of our $b_h$, we use left $h$-extrema, also defined in Section \ref{Sect_Potential}.
Our $b_h$ is defined in terms of left $h$-extrema,
which are also introduced in Section \ref{Sect_Potential} (see Definition \ref{def_left_extrema}).
In order to prove our Theorem \ref{Th_Local_Limit_Sinai}, we first prove a local limit theorem for $b_h$:

\begin{thm}\label{Th_Local_Limit_b_h}
We have as $h\to+\infty$,
$$
    \sup_{x\in\Z}\bigg|\p\big(b_h=x\big)-\frac{\sigma^2}{h^2}\varphi_\infty\bigg(\frac{\sigma^2 x}{h^2}\bigg)\bigg|
=
    o\bigg(\frac{1}{h^2}\bigg).
$$
%For every $A>0$, we have as $h\to+\infty$,
%$$
%    \sup_{x\in[-Ah^2, Ah^2]\cap\Z}\bigg|\p\big(b_h=x\big)-\frac{\sigma^2}{h^2}\varphi_\infty\bigg(\frac{\sigma^2 x}{h^2}\bigg)\bigg|
%=
%    o\bigg(\frac{1}{h^2}\bigg).
%$$
\end{thm}
Even though Theorem \ref{Th_Local_Limit_b_h} looks, at first sight, very similar to Theorem \ref{Th_Local_Limit_Sinai},
Theorem \ref{Th_Local_Limit_Sinai} is not a direct consequence of Theorem \ref{Th_Local_Limit_b_h},
because, loosely speaking, the event $\{S_n=z\}$ can be decomposed into a union of events
$\{S_n=z\}\cap\{b_{\log n}=y\}$, and we will see that each one has a non-negligible probability for $y$ "close" to $z$.
Also, estimating the annealed probabilities of these events for $y$ close to $z$,
as well as proving that such probabilities are negligible for $y$ "far" from $z$, is not immediate,
since we have to decompose each of these events into many different cases,
corresponding to different kinds of environments and trajectories.
%Also for this reason, we will see that it is important to have a uniform convergence in Theorem \ref{Th_Local_Limit_b_h}

The probability $\PP(S_n=z_n)$ for Sinai's walk seems to have been first studied
in a physics paper in 1985 by Nauenberg \cite{Nauenberg}, by heuristic arguments in some particular cases
and numerical simulations.
However the function he obtained instead of our $\varphi_\infty$ is
%$\frac{1}{\sqrt{2}} \exp(-\sqrt{2}|x|)$,
$x\mapsto (C/2)\exp(-C|x|)$ for some $C>0$,
which is not correct.
This function was also claimed in Nauenberg \cite{Nauenberg} to be the density of the limit law of
$\frac{\sigma^2}{(\log n)^2}S_n$, and Kesten \cite{Kesten} already noticed that this is not the correct function,
although $\varphi_\infty(x)$ is equivalent to some exponential as $x\to+\infty$.

%For a recent survey about local limit theorems, we refer to Szewczak et al. \cite{Szewczak_Weber}.

There have been many papers dealing with local limit theorems for different models of random walks in random environments recently.
For example, Dolgopyat and Goldsheid \cite{Dolgo_gold_13}, \cite{Dolgo_gold_19},
Leskela and Stenlund \cite{Leskela_Stenlund} and Berger et al. \cite{Berger_CR}
prove local limit theorems for transient RWRE respectively on $\Z$ and on a strip, both in the diffusive regime,
 on $\Z$ with only $0$ or $1$ jumps, and for some ballistic multidimensional RWRE.
See also Dolgopyat et al. \cite{Dolgo_gold_21} for diffusive recurrent RWRE on a strip,
Takenami \cite{Takenami} for random walks on periodic environments,
Chiarini et al. \cite{Chiarini} for some diffusions in random environment,
and Andres et al. \cite{Andres_Taylor} for the random conductance model.
We refer to the first two sections of Dolgopyat et al.  \cite{Dolgo_gold_19} for a recent review of this subject.
However, the previously cited papers consider transient or diffusive random walks or diffusions,
whereas we consider Sinai's walk which is recurrent and subdiffusive.
Also, we obtain probabilities of order $(\log n)^{-2}$ with a non gaussian limit law,
instead of $n^{-1/2}$ with a gaussian limit law in their cases.
Therefore, to the extent of our knowledge, our Theorem \ref{Th_Local_Limit_Sinai} is the first local limit theorem for
(recurrent) subdiffusive RWRE.

Also, a similar local limit theorem for the quenched probability, replacing $\PP$ by $\po$, does not hold.
Indeed, $\po(S_n=0)$ almost surely takes very small values compared to $(\log n)^{-2}$  as $n\to+\infty$,
since for $\eta\in]0,1[$, $\p$-almost surely
$\po(S_n=0)= O\big( \exp(-(\log n)^{1-\eta})\big)$ as $n\to+\infty$
%$\po(S_n=0)\leq C \exp(-(\log n)^{1-\eta})$
(see Devulder et al. \cite{DGP_Collision_Sinai}, last inequality of page 6).
See also Gantert et al. (\cite{NMFa}, Theorem 1.1) for previous results,
Comets et al. (\cite{CometsPopov}, Theorem 2.1 and Corollary 2.1)
for estimates for a related model in continuous time,
 and Gantert et al. \cite{Gantert_Peterson} for transient RWRE.
So, contrarily to some of the previously cited papers on local limit theorems for RWRE,
our annealed local limit theorem, Theorem \ref{Th_Local_Limit_Sinai},
cannot be the consequence of a corresponding quenched local limit theorem.

We also mention that some estimates of $\PP(S_n=z_n)$ when $z_n$ is large, more precisely when $n=O(z_n)$,
are given by Comets et al. \cite{CometsGantertZeitouni}.
%, and in other articles about large deviations for RWRE.
For an overview of the vast literature about large deviations for RWRE, see e.g.
Gantert et al. \cite{Gantert_Zeitouni} and more recently
Buraczewski et al. \cite{Buraczewski_Dyszewski}.

Finally, we think that the tools and technics developed in the present paper, in particular the ones of Section \ref{Sect_Potential},
will be useful for future research projects, including \cite{Devulder_Rates_CV}, which
will study the rates of convergence in Sinai and Golosov localization theorems for Sinai's walk.

%Hu et al. \cite{HuShiModerate} in the moderate deviation case, and in

\noindent{\bf Acknowledgement:} I am thankful to Yueyun Hu for asking,
after a talk in a conference in Landela (France) in 2016, if I could give an estimate of $\PP(S_{2n}=0)$ as
$n\to +\infty$, which made me aware that this question was still open.
I also thank Fran{\c c}oise P\`ene for organizing this conference.
Part of this work was done during a six months sabbatical "d\'el\'egation CNRS".

\subsection{Organization of the proof and of the paper}\label{Sub_Sect_Sketch}

In Section \ref{Sect_Potential}, we recall the definition and use of the potential $V$. We also
define left and right $h$-extrema for $V$, for $h>0$. This allows us to introduce two path decompositions of the potential $V$,
one with left $h$-extrema and one with right $h$-extrema.
We can then define our localization point $b_h$.
We describe the law of the potential $V$
between two consecutive left (or right) $h$-extrema $x_i$ and $x_{i+1}$ when $0\notin [x_i, x_{i+1}]$,
which uses in particular the law of $V$ or $-V$ conditioned to stay
positive, or nonnegative, up to some hitting time (see Theorem \ref{Lemma_Law_of_Slopes}).
The law of $V$ between the two left $h$-extrema surrounding $0$
is given by a renewal theorem (see Theorem \ref{Lemma_Central_Slope}),
and some independence is provided by Theorem \ref{Lemma_Independence_h_extrema}.
A first application of this renewal theorem is that we can give a simple formula for the law of $b_h$,
that is, for $\p(b_h=x)$, $x\in\Z$ (in Lemma \ref{Lemma_Proba_bh_egal}),
which is an important tool in the proof of Theorem \ref{Th_Local_Limit_b_h}.

Section \ref{Sect_Proof_Th_LLT_b_h} is devoted to the proof of Theorem \ref{Th_Local_Limit_b_h}.

In Section \ref{Sect_Coupling}, we first define an event $E_C^{(n)}(z)$, depending only on the environment and on $z$.
On this event, we use a coupling argument, which helps us approximate
the quenched probability $\po(S_n=z)$ by $\widehat \nu_n(z)$,
where $\widehat \nu_n$ is an invariant probability measure.
This enables us to give an upper bound for the annealed probability that $S_n=z$
on $E_C^{(n)}(z)$ (see Proposition \ref{Ineg_Upper_Bound_ELT}),
giving the main contribution in the upper bound of Theorem  \ref{Th_Local_Limit_Sinai}.
To this aim, loosely speaking, we express the expectation of $\widehat \nu_n(z)$
on each event $\{b_{\log n}=k+z\}\cap E_C^{(n)}(z)$
with quantities depending only on the laws of the potential $V$ between
consecutive left or right $(\log n)$-extrema;
summing this over $k$ makes appear, after some inequalities and computations using the tools developed in Section \ref{Sect_Potential},
 a formula equal to $\p(b_{\log n}=z)$ by Lemma \ref{Lemma_Proba_bh_egal}.
 We conclude by applying Theorem \ref{Th_Local_Limit_b_h}.

In Section \ref{Sect_Negligible}, we prove that the environments and trajectories such that $S_n=z$
which were not considered in Section \ref{Sect_Coupling} have a negligible annealed probability.
%Section \ref{Sect_Negligible} is devoted to the study of negligible environments or trajectories,
%or of some particular trajectories in some particular environments, satisfying $S_n=z$,
%the ones which are not included on $E_C^{(n)}(z)$.
This covers many different cases, which often combine conditions on both environments and trajectories of $(S_k)_k$.
For example, $z$ can be far from $b_{\log n}$, or the origin $0$ can be very close
to the maximum of the potential between two valleys (defined before \eqref{eqDef_tau_1_V}), or some of the valleys around the origin can
have a height just slightly larger than $\log n$, or the central valley of height at least $\log n$ can
include one or several subvalleys of height slightly less than $\log n$.
The potentials for some of these cases are represented in Figures
\ref{figure_Proba_E3c} page \pageref{figure_Proba_E3c},
\ref{figure_Lemma_5_5} page \pageref{figure_Lemma_5_5},
\ref{figure_Ik_cas_1_et_2} page \pageref{figure_Ik_cas_1_et_2}
and \ref{figure_Ik} page \pageref{figure_Ik}.
In this section, we prove that all these cases, and some others,  with $S_n=z$ have a negligible annealed probability (compared to $(\log n)^{-2}$).
Combining this with the previous subsection,
%and using Theorem \ref{Th_Local_Limit_b_h},
we get (uniformly on $z$) an upper bound of $\PP(S_n=z)$,
which completes the proof of the upper bound in Theorem  \ref{Th_Local_Limit_Sinai}.
Even if this section mainly consider negligible events, it is maybe the most delicate of the paper.

Section \ref{Sect_lower_Bonud} is devoted to the proof of the lower bound in Theorem \ref{Th_Local_Limit_Sinai},
that is, we give (uniformly on $z$) a minoration of $\PP(S_n=z)$.
The proof is divided into three cases, depending on $z$ being negative and far from $0$, positive and far from $0$, or
$z$ being close to $0$.
This uses results of all the other   sections.

Finally, Section \ref{Sect_Technical} is devoted to some important technical lemmas and their proofs.
These lemmas mainly deal with the potential $V$, and with $V$ conditioned to stay positive or nonnegative.

Outlines or sketches of proofs of several lemmas or theorems are also provided throughout the paper.

%%%%%%%%%%%%%%%%%%%%%%%%%%%%%%%%%%%%%%%%%%%%%%%%%%%%%%%%%%%%%%%%%%%%%%%%%

\section{Potential, path decomposition and renewal theorem}\label{Sect_Potential}

\subsection{Definition and applications of the potential}
The {\it potential} $(V(x),\ x\in\Z)$,
which was first introduced by Sinai \cite{S82}, is an important quantity which depends only on the environment $\o$.
It is defined as follows:
\begin{equation}\label{eqDefPotentialV}
    V(x)
:=
    \left\{
    \begin{array}{lr}
        \sum_{i=1}^x \log\frac{1-\omega_i}{\omega_i}
    &
        \textnormal{if } x>0,
    \\
        0
    &
        \textnormal{if }x=0,
    \\
        -\sum_{i=x+1}^0 \log\frac{1-\omega_i}{\omega_i}
    &
        \textnormal{if } x<0.
    \end{array}
    \right.
\end{equation}
We denote by $\po^x$ the quenched probability for the RWRE $(S_k)_k$ starting at $x\in\Z$ instead of $0$,
and by $E_\o^x$ the expectation with respect to $\po^x$.
Also, let
$$
    \tau(y)
:=
    \inf\{k\ge 0\ :\ S_k=y\},
\qquad
    \tau^*(y)
:=
    \inf\{k\ge 1\ :\ S_k=y\},
\qquad
\qquad
    y\in\Z,
$$
where by convention, $\inf\emptyset=+\infty$.
In words, $\tau(y)$ (resp. $\tau^*(y)$)
is the hitting time of (resp. return time to) the site $y$ by the RWRE $(S_k)_k$.
We also define for
%$j\in\{1,\dots, d\}$,
$x\in\Z$ and $y\in\Z$,
\begin{equation*}
%\label{tautaujxy}
    \tau(x,y)
:=
    \inf\{k\in \N\ :\ S_{\tau(x)+k}=y\}.
\end{equation*}
We now recall some classical estimates, which explain why the potential is very useful.
These formulas will be used throughout the paper.
%The potential is very useful, due to the following properties.
First, we have  (see e.g. \cite[(2.1.4)]{Z01},
%and \cite[Lemma 2.2]{Devulder_Persistence} coming from \cite[p. 250]{Z01}),
\begin{equation}
\label{probaatteinte}
    P_\omega^b[\tau(c)< \tau(a)]
=
    \bigg(\sum_{j=a}^{b-1} e^{V(j)}\bigg)\bigg(\sum_{j=a}^{c-1} e^{V(j)}\bigg)^{-1},
\qquad
    a<b<c.
\end{equation}
Furthermore (see e.g. \cite {Devulder_Persistence} Lem. 2.2 coming from Zeitouni \cite{Z01} p. 250),
if $g<h<i$,
\begin{equation}\label{InegEsperance1}
    \eo^{h} [\tau(g)    \wedge\tau(i)]
\leq
    \sum_{k={h}}^{i-1}\sum_{\ell=g}^{k}\frac{\exp[V(k)-V(\ell)]}{\o_{\ell}}
\leq
    \e_0^{-1}(i-g)^2\exp\left[\max_{g\leq\ell\leq k\leq i-1, k\geq h}(V(k)-V(\ell))\right],
\end{equation}
%\begin{eqnarray}
%\label{InegEsperance1}
%    E_\omega^b[\tau(a)\wedge \tau(c)]
%& \leq &
%    \varepsilon_0^{-1}(c-a)^2
%    \exp\Big[\max_{a\leq \ell \leq k \leq c-1; k\ge b}\big(V(k)-V(\ell)\big)\Big],
%\qquad
%    a<b<c.
%\end{eqnarray}
%where $E_\omega^b$ denotes the expectation with respect to $P_\omega^b$.
where we used ellipticity \eqref{eqEllipticity} in the last inequality
and with $x\wedge y :=\min(x,y)$.
For symmetry reasons, we also have
\begin{equation}
\label{InegEsperance2}
    E_\omega^b[\tau(a)\wedge \tau(c)]
\leq
    \varepsilon_0^{-1}(c-a)^2
    \exp\Big[\max_{a\leq \ell \leq k \leq c-1,\ \ell\le b-1}\big(V(\ell)-V(k)\big)\Big],
\qquad
    a<b<c\, .~~~~
\end{equation}
Moreover, we have (see Golosov  \cite{Golosov84}, Lemma 7, proved for a RWRE on $\N$ but still true for a RWRE on $\Z$),
\begin{equation}
\label{InegProba1}
    P_\omega^b[\tau(c)<k]
\leq
     k \exp\left(\min_{\ell \in [b,c-1]}V( \ell)-V(c-1)\right),\qquad     b<c\, .
\end{equation}
Also by symmetry,
%replacing each $\omega_x$ by $1-\omega_{-x}$,
we get  (similarly as in Shi and Zindy \cite{ShiZindy}, eq. (2.5) but with some slight differences for the values of $\ell$)
\begin{eqnarray}
\label{InegProba2}
    P_\omega^b[\tau(a)<k]
& \leq & k \exp\left(\min_{\ell \in [a,b-1]}V(\ell) -V(a)\right),\qquad     a<b\, .
\end{eqnarray}
Moreover, we have by Devulder et al. (\cite{DGP_Collision_Sinai}, Lemma 4.10), if $a\neq b$,
\begin{equation}\label{ineg_Proba_Atteinte_egalite_DGP}
    \forall k\in\N,
\qquad
    P_\omega^b[\tau(a)=k]
\leq
    P_\omega^b[\tau(a)<\tau^*(b)].
\end{equation}
%\smallskip
Finally, we recall that, given $\omega$, the Markov chain $S$ is an electrical network
where, for every $x\in\Z$,  the conductance of the unoriented bond $(x,x+1)$ is $C_{(x,x+1)}=e^{-V(x)}$
(in the sense of Doyle and Snell \cite{Doyle_Snell}) (see also Levin et al. \cite{Levin_Peres}).
In particular, its reversible measure $\mu_\omega$
(unique up to a multiplication by a constant) is given by
\begin{equation}\label{reversiblemeas}
    \mu_\omega(x)
:=
    e^{-V(x)}+e^{-V(x-1)},
\qquad
    z\in\Z,
\end{equation}
where, for the sake of simplicity, we write $\mu_\omega(x)$ instead of $\mu_\omega(\{x\})$.
%{\bf (rappeler l'utilisation ?)}
For any process $Y$, we define
\begin{eqnarray}
\label{eq_def_TY}
    T_Y(A)
& := &
    \inf\{x\geq 0, \ Y(x)\in A\},
\qquad
    A\subset \R,
\\
    T_Y^*(A)
& := &
    \inf\{x> 0, \ Y(x)\in A\},
\qquad
    A\subset \R.
\label{eq_def_TY*}
\end{eqnarray}
We sometimes write $T_Y(a):=T_Y([a,+\infty[)$ when $a>0$
and $T_Y(a):=T_Y(]-\infty, a])$ when $a<0$.
Due to the ellipticity \eqref{eqEllipticity}, we have
\begin{equation}\label{eq_ellipticity_for_V}
    \forall x\in\Z,
\qquad
    \big|V(x)-V(x-1)\big|
\leq
    \log\bigg(\frac{1-\e_0}{\e_0}\bigg)
=:C_0.
\end{equation}
In particular, thanks to \eqref{eqRecurrence} and \eqref{eq_ellipticity_for_V},
the following fact follows from the optimal stopping theorem
applied to the martingale $(V(k), \ k\geq 0)$ at time $T_V([z, +\infty[)  \wedge  T_V(]-\infty,x]$:
%(see e.g. Zindy \cite{Zindy}, Lem. 2.1 applied to $-V$)
%\begin{equation}\label{eqOptimalStopping1}
%    \frac{z-y}{z-x+C_0}
%\leq
%    \p^y\big[ T_V(]-\infty,x]) < T_V([z, +\infty[) \big]
%\leq
%    \frac{z-y+C_0}{z-x+C_0},
%\qquad
%    x<y<z,
%\end{equation}
\begin{equation}\label{eqOptimalStopping2}
    \frac{y-x}{z-x+C_0}
\leq
    \p^y\big[T_V([z, +\infty[)  <  T_V(]-\infty,x])\big]
\leq
    \frac{y-x+C_0}{z-x+C_0},
\qquad
    x<y<z,
\end{equation}
where $\p^y$ denotes the law of $V$ starting from $y$ instead of $0$.
Moreover, these inequalities remain valid if we replace $]-\infty,x]$ and/or $[z,+\infty[$
by the corresponding open interval $]-\infty,x[$ and/or $]z,+\infty[$.
Also, there exist constants $c_1>0$ and $c_1^*>0$ such that
(see e.g. Lemma \ref{Lem_Proba_An_Ordre2}),
\begin{equation}\label{eq_Proba_Atteinte_logn_avant0}
    \p[T_V(h)<T_V^*(\R_-)]
\sim_{h\to+\infty}
    c_1^* h^{-1},
    %c h^{-1},
\quad
    \p[T_V(h)<T_V(\R_-^*)]
\sim_{h\to+\infty}
    c_1 h^{-1}.
%    c^* h^{-1}.
\end{equation}
%\begin{equation}\label{eq_Proba_Atteinte_logn_avant0}
%    \p[T_V(\log n)<T_V^*(\R_-)]
%\sim_{n\to+\infty}
%    c (\log n)^{-1},
%\quad
%    \p[T_V(\log n)<T_V(\R_-^*)]
%\sim_{n\to+\infty}
%    c^* (\log n)^{-1}.
%\end{equation}

%%%%%%%%%%%%%%%%%%%%%%%%%%%%%%%%%%%%%%%%%%%%%%%%%%%%%%%%%%%%%%

\subsection{Definition and properties of left and right $h$-extrema}
The point of view of $h$-extrema has been used recently in some papers for RWRE or diffusions in a random potential,
either to prove localization results, see e.g. \cite{AndreolettiDevulder}, \cite{Bovier_Faggionato},  \cite{DGP_Collision_Sinai} and \cite{Freire},
or to use localization techniques, see e.g.  \cite{AndreolettiDevulderVechambre}, \cite{Cheliotis_Favorite}, \cite{Devulder_Persistence} and \cite{CometsPopov} (where they are called $e^h$-stable points).

However, these studies use $h$-extrema of a (maybe drifted) two-sided Brownian motion $W$, and sometimes transfer results about
$W$ to the potential $V$ by
Koml\`os, Major and Tusn\'ady strong approximation theorem \cite{KMT}. This is not precise enough to prove our theorems, so we introduce and study variants of $h$-extrema directly for our potential $V$.

Let $h>0$, and $v$ be a function from $\Z$ to $\R$.
Following Neveu and Pitman \cite{NP}, we say that $y$ is an {\it $h$-minimum} for $v$ if
there exist integers $\alpha<y<\beta$ such that $v(y)=\min_{[\alpha,\beta]}v$,
$v(\alpha)\geq v(y)+h$ and $v(\beta)\geq v(y)+h$.
We say that $y$ is an {\it $h$-maximum} for $v$ if it is an $h$-minimum for $-v$.
In both cases, we say that $y$ is an {\it $h$-extremum} for $v$.

\begin{figure}[htbp]
\includegraphics[width=16.0cm,height=6.3cm]{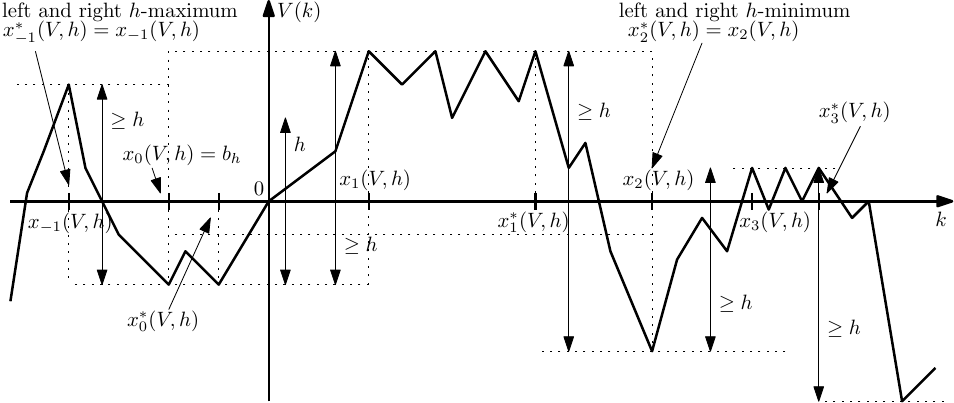}
\caption{Schema of the potential $V$ with left $h$-extrema $x_i(V,h)$ (defined before \eqref{eqDefbh})
and right $h$-extrema $x_i^*(V,h)$ (defined before \eqref{eq_relation_xi_xi*}).}
\label{figure_h_extrema}
\end{figure}

One of the main differences with  $h$-extrema of Brownian motion
is that unfortunately, in the general case,
$h$-maxima and $h$-minima for $V$   do not necessarily alternate.
For this reason, we introduce the following definitions (see Figure \ref{figure_h_extrema}).
%having properties similar to that of $h$-extrema of Brownian motion
%(which are provided by \cite{Cheliotis}).

\begin{deff}\label{def_left_extrema}
Let $h>0$ and $v$ be a function from $\Z$ to $\R$.
We say that $y\in\Z$ is a {\it left $h$-minimum} (resp. {\it right $h$-minimum}) for $v$ if
there exist $\alpha<y<\beta$ such that
\begin{itemize}
\item $\min_{[\alpha,y-1]}v>v(y)$ (resp $\geq$),
\item $\min_{[y+1,\beta]}v\geq v(y)$ (resp. $>$),
\item $v(\alpha)\geq v(y)+h$,
\item $v(\beta)\geq v(y)+h$.
\end{itemize}
We say that $y$ is a {\it left $h$-maximum} (resp. {\it right $h$-maximum}) for $v$
if it is a left $h$-minimum (resp. right $h$-minimum) for $-v$.
In both cases, we say that $y$ is a {\it left $h$-extremum} (resp.  {\it right $h$-extremum}) for $v$.
\end{deff}

With these definitions, left $h$-minima and left $h$-maxima for $v$ alternate, and similarly
right $h$-minima and right $h$-maxima for $v$ alternate. The elementary proof is
given in Lemma \ref{Lemma_Alternate}.
Also, between two consecutive left $h$-maxima
$y_1$ and $y_2$, more precisely in $[y_1,y_2[\cap\Z$,
there are one or several $h$-minima, among which the smallest one is
the only left $h$-minimum, which is $y_1$, and the largest one is the only right $h$-minimum,
which we will not use in the present paper.

Left and right $h$-extrema of $V$ have the disadvantage of not being stopping times.
However, we will see that they allow a very simple definition of the localization point $b_h$ (see \eqref{eqDefbh} below,
which can be compared e.g. to \eqref{eqDef_bK_h}),
%see also Lemma \ref{Lemma_Comparaison_b_bK}),
that they have nice independence properties,
that the properties of the law of trajectories of $V$ between consecutive left or right $h$-extrema are convenient,
and that we can use renewal theory,
which enables for example to prove very useful formulas such as the law of $b_h$
(see Lemma \ref{Lemma_Proba_bh_egal}).

We now focus on left $h$-extrema.
Let $\VVV$ be the set of functions $v$ from $\Z$ to $\R$, such that
$\liminf_{\pm \infty}v=-\infty$ and $\limsup_{\pm \infty}v=+\infty$.
If $v\in\VVV$ and $h>0$, then the set of left $h$-minima of $v$ is unbounded from above and below,
and so is the set of left $h$-maxima of $v$.
Consequently, for $v\in\VVV$ for every $h>0$, the set of left $h$-extrema of $v$
can be denoted by $\{x_k(v,h), \ k\in\Z\}$, such that
$k\mapsto x_k(v,h)$ is strictly increasing
and $x_0(v,h)\leq 0<x_1(v,h)$.
%and the left $h$-minima and left $h$ maxima of $V$ alternate.
%(we need to define the $x_i$ a bit differently so that this is true in the general case without the previous assumption).
And also, $\lim_{k\to \pm\infty} x_k(v,h)=\pm \infty$.
Notice that due to our hypotheses \eqref{eqRecurrence} and \eqref{eq_def_sigma},
%$\liminf_{\pm \infty}V=-\infty$ and $\limsup_{\pm \infty}V=+\infty$
$V\in\VVV$ almost surely.
%If $v\notin\VVV$, we define $x_k(v,h)=0$

Similarly as in the continuous case (see Cheliotis \cite{Cheliotis}),
we can now define for $h>0$,
\begin{equation}\label{eqDefbh}
    b_h
:=
    \left\{
    \begin{array}{ll}
        x_0(V,h)
    &
        \textnormal{if } x_0(V,h) \textnormal{ is a left $h$-minimum for }V,
    \\
        x_1(V,h)
    &
        \textnormal{otherwise}.
    \end{array}
    \right.
\end{equation}
As already mentioned,
the definition of the localization point $b_h'$ given by Sinai \cite{S82} is not the same.
%since it is defined by a sequence of refinements of valleys.
%, and the two can differ,
%however we think that they are equal with quite large probability.
%We will not prove this in the present paper since this is not needed.

Similarly as in the continuous case for $h$-slopes, we introduce for each function $v\in\VVV$ and
for each $i\in\Z$ and $h>0$     the  {\it left $h$-slope}
$T_i(v,h):=(v(j)-v[x_i(v,h)],\ x_i(v,h)\leq j\leq x_{i+1}(v,h))$.
Its {\it height} and its {\it excess height} are defined respectively as
 $$
    H[T_i(v, h)]
:=
    \big|v[x_{i+1}(v,h)]-v[x_{i}(v,h)]\big|
\geq
    h,
\qquad
    e[T_i(v, h)]=H[T_i(v, h)]-h\geq 0.
$$
If $x_i(v,h)$ is a left $h$-minimum (resp. maximum), then $T_i(v,h)$ is a nonnegative (resp. nonpositive) function,
it is said to be an {\it upward slope} (resp. a {\it downward slope})
and its maximum (resp. minimum) is attained at $x_{i+1}(v,h)$, with
$
    \sup_{[x_i(v,h),x_{i+1}(v,h)[}v
<
    v[x_{i+1}(v,h)]
$
(resp.
$
    \inf_{[x_i(v,h),x_{i+1}(v,h)[}v
>
    v[x_{i+1}(v,h)]
$
).

Similarly, if $y_i$ and $y_{i+1}$ are two consecutive right $h$-extrema of $v$,
we say that
%$T_i^{(r)}=
$(v(j)-v(y_i),\ y_i\leq j\leq y_{i+1})$ is a {\it right $h$-slope} of $v$
(see Subsection \ref{Sub_sec_right_extrema} for some properties of right $h$-slopes and extrema).
More generally, we call a {\it slope} each T=$(T(j), \ \alpha\leq j \leq \beta)\in\R^{\beta-\alpha+1}$,
with $\alpha\in\Z$, $\beta\in\Z\cap]\alpha,+\infty[$, such that either
$T(\alpha)=0=\min_{[\alpha, \beta]\cap \Z}T <\max_{[\alpha, \beta]\cap \Z}T=T(\beta)$
or
$T(\beta)=\min_{[\alpha, \beta]\cap \Z}T <\max_{[\alpha, \beta]\cap \Z}T=T(\alpha)=0$.
Also, for each slope $T=(T(j),\ \alpha\leq j\leq \beta)$, we define its {\it length} $\ell(T):=\beta-\alpha$,
its {\it height} $H(T)=|T(\beta)-T(\alpha)|$,
and the {\it translated slope} $\theta(T):=(T(j+\alpha),\ 0\leq j\leq \beta-\alpha)$.
%and similarly as before, its {\it height} $h(T)=|V(\lpha)-V(\beta)$
%and its {\it excess height} $e(T)=.

%For each interval $[a,b]$, we say that $[a,b]$ is a {\it valley of height at least $h$} for $V$
%if and only if $a<b$ are two consecutive left $h$-maxima for $V$.

We call {\it valleys of height at least $h$} of $V$ the intervals
$[x_i(V,h), x_{i+2}(V,h)]$, $i\in\Z$, such that $x_i(V,h)$ and $x_{i+2}(V,h)$ are (consecutive) left $h$-maxima.
The {\it bottom} of such a valley is the left $h$-minimum $x_{i+1}(V,h)$.
If its bottom is $b_h$,
that is, if $b_h=x_{i+1}(V,h)$,
%$0\in]x_i(V,h), x_{i+2}(V,h)[$,
then it is called the {\it central valley of height at least $h$} of $V$.

%For each $i\in\Z$,
%we will say that the  interval $[x_i(V,h), x_{i+2}(V,h)]$ is a {\it valley of height at least $h$}
%if $x_i(V,h)$ and $x_{i+2}(V,h)$ are  left $h$-maxima.
%Its {\it bottom} is the left $h$-minimum $x_{i+1}(V,h)$.
%If its bottom is $b_h$,
%that is, if $0\in[x_i(V,h), x_{i+2}(V,h)]$,
%then it is called the {\it central valley of height at least $h$}.

Knowing, for some $h>0$, $\theta[T_i(V,h)]$ for each $i\in\Z^*$ and $(\theta[T_0(V,h)], \ x_0(V,h))$
allows us to reconstitute totally the process $V$ since $V(0)=0$.
The two following subsections will provide their laws and independence properties.

\subsection{Definition and law of $\mathcal T_{V,h}^\uparrow$ and $\mathcal T_{V,h}^\downarrow$}

\begin{figure}[htbp]
\includegraphics[width=16.0cm,height=6.73cm]{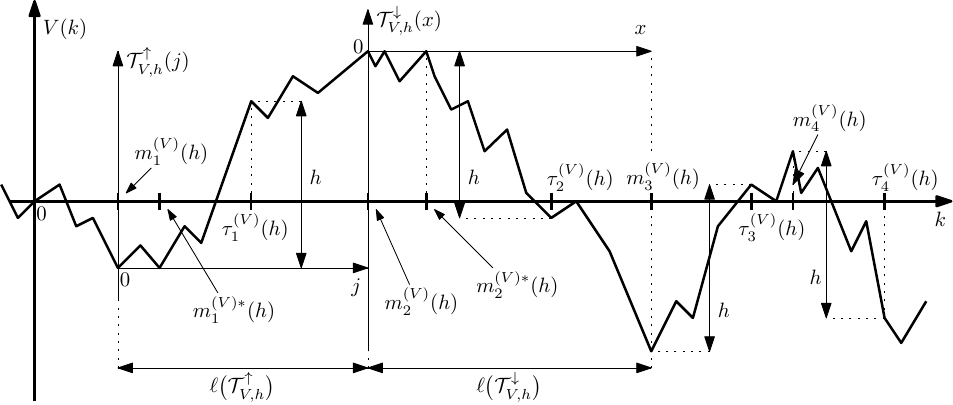}
\caption{Schema of the potential $V$ with the $\tau^{(V)}_{i}(h)$, $m^{(V)}_{i}(h)$,
$\mathcal T_{V,h}^\uparrow$, $\mathcal T_{V,h}^\downarrow$ and $m^{(V)*}_{i}(h)$
%(the last one being defined before \eqref{eq_def_slope_etoile_1}).
(defined between \eqref{eqDef_tau_1_V} and \eqref{eqDef_m_2_V}, in Definition \ref{deff_law_slopes}
and before \eqref{eq_def_slope_etoile_1}).
}
\label{figure_tau_i_m_i_slopes}
\end{figure}

Let $h>0$.
We define by induction the following notation. Let $\tau_0^{(V)}(h):=0$ and for $i\geq 0$ (see Figure \ref{figure_tau_i_m_i_slopes}),
\begin{eqnarray}
    \tau^{(V)}_{2i+1}(h)
& := &
    \min\Big\{k\geq \tau_{2i}^{(V)}(h),\ V(k)-\min\nolimits_{[\tau_{2i}^{(V)}(h),k]} V\geq h\Big\},
\label{eqDef_tau_1_V}
\\
    m^{(V)}_{2i+1}(h)
& := &
    \min\Big\{k\geq \tau_{2i}^{(V)}(h),\ V(k)=\min\nolimits_{[\tau_{2i}^{(V)}(h),\tau_{2i+1}^{(V)}(h)]} V\Big\},
\label{eqDef_m_1_V}
\\
    \tau^{(V)}_{2i+2}(h)
& := &
    \min\Big\{k\geq \tau_{2i+1}^{(V)}(h),\ \max\nolimits_{[\tau_{2i+1}^{(V)}(h),k]} V-V(k)\geq h\Big\},
\label{eqDef_tau_2_V}
\\
    m^{(V)}_{2i+2}(h)
& := &
    \min\Big\{k\geq \tau_{2i+1}^{(V)}(h),\ V(k)=\max\nolimits_{[\tau_{2i+1}^{(V)}(h),\tau_{2i+2}^{(V)}(h)]} V\Big\}.
\label{eqDef_m_2_V}
\end{eqnarray}
Notice that
that $\tau^{(V)}_{i}(h)<\infty$ $\p$-a.s. for $i\geq 0$
since $V\in\VVV$ $\p$-a.s. due to \eqref{eqRecurrence} and \eqref{eq_def_sigma},
and that the $\tau^{(V)}_{i}(h)$, $i\geq 0$, are stopping times for the natural filtration of $(V(\ell), \ \ell\geq 0)$.

Let $\bigsqcup$ denote the disjoint union.
Notice that, with a slight abuse of notation,
each translated (left $h$-) slope $T=(T(0), T(1), \dots,$
$T(\ell(T)))$
belongs to $\R^{\ell(T)+1}$.
So, we can consider our translated slopes (and $\mathcal T_{V,h}^\uparrow$ and $\mathcal T_{V,h}^\downarrow$ defined below) as random variables taking values into
$\bigsqcup_{t\in\mathbb N^*}\R^t$, equipped with the $\sigma$-algebra
$\{\bigsqcup_{t\in\mathbb N^*}A_t\, :\, \forall t\in\mathbb N^*,\ A_t\in\mathcal B(\mathbb R^t)\}$,
where $\mathcal B(\mathbb R^t)$ is the Borel $\sigma$-algebra of $\mathbb R^t$.
The following notations are useful to express the law of left $h$-slopes in the next subsection:

\begin{deff} \label{deff_law_slopes}
Let $h>0$. We introduce (see Figure \ref{figure_tau_i_m_i_slopes}),
\begin{eqnarray*}
    \mathcal T_{V,h}^\uparrow
& := &
    \big(V\big[m_{1}^{(V)}(h)+x\big]-V\big[m_{1}^{(V)}(h)\big], \ 0\leq x \leq m_{2}^{(V)}(h)-m_{1}^{(V)}(h)\big),
\\
    \mathcal T_{V,h}^\downarrow
& := &
    \big(V\big[m_{2}^{(V)}(h)+x\big]-V\big[m_{2}^{(V)}(h)\big], \ 0\leq x \leq m_{3}^{(V)}(h)-m_{2}^{(V)}(h)\big).
\end{eqnarray*}
\end{deff}

In particular,
$
    \ell\big(\mathcal T_{V,h}^\uparrow\big)
=
    m_{2}^{(V)}(h)-m_{1}^{(V)}(h)
$
and
$
    \ell\big(\mathcal T_{V,h}^\downarrow\big)
=
    m_{3}^{(V)}(h)-m_{2}^{(V)}(h)
$.
We sometimes write $\mathcal T_V^\uparrow$ and $\mathcal T_V^\downarrow$ instead of
$\mathcal T_{V,h}^\uparrow$ and $\mathcal T_{V,h}^\downarrow$
to simplify the notation
when no confusion is possible for the value of $h$.
The laws of $\mathcal T_{V,h}^\uparrow$ and $\mathcal T_{V,h}^\downarrow$ are given in the following theorem.

\begin{thm}\label{Lemma_Law_of_Slopes}
Assume \eqref{eqEllipticity}, \eqref{eqRecurrence} and \eqref{eq_def_sigma}. Let $h>0$.
\\
{\bf (i)} The process $\mathcal T_{V,h}^\uparrow$ up to its first hitting time $T_{\mathcal T_{V,h}^\uparrow}([h,+\infty[)$
of $[h,+\infty[$,
that is, $\big(\mathcal T_{V,h}^\uparrow(k),$
$0\leq k \leq T_{\mathcal T_{V,h}^\uparrow}([h,+\infty[)\big)$,
is equal in law to
$
    \big(V(k), \  0\leq k \leq T_V([h,+\infty[)
    \big)
$
conditioned on $\{ T_V([h,+\infty[)$
$ < T_V(]-\infty, 0[ )\}$.
Moreover, it is independent of
$
    \big(
    \mathcal T_{V,h}^\uparrow\big(T_{\mathcal T_{V,h}^\uparrow}([h,+\infty[)+k\big)
    -
    \mathcal T_{V,h}^\uparrow\big(T_{\mathcal T_{V,h}^\uparrow}([h,+\infty[)\big),\,
    0\leq k \leq \ell\big(\mathcal T_{V,h}^\uparrow\big)-T_{\mathcal T_{V,h}^\uparrow}([h,+\infty[)
    \big)
$,
%which has the same law as $V$ up to
%$\min\{k\in\N, \, V(k)=\max_{[0, \tilde\tau_1]}V\}$,
which has the same law as
%$V$ between $0$ and
$(V(k),\, 0\leq k \leq M^\sharp_h)$, with
$M^\sharp_h:=\min\{k\in\N, \, V(k)=\max_{[0, \tilde\tau_1(h)]}V\}$,
where
$
    \tilde\tau_1(h)
:=
    \min\{k\in \N,\, \max_{[0,k]}V-V(k)\geq h\}
$.

{\bf (ii)}
$\mathcal T_{-V,h}^\uparrow=_{law}-\mathcal T_{V,h}^\downarrow$
and
$\mathcal T_{-V,h}^\downarrow=_{law}-\mathcal T_{V,h}^\uparrow$.

{\bf (iii)} Also,
$
    \E\big(\ell\big(\mathcal T_{V,h}^\uparrow\big)\big)
<
    \infty
$
and
$
    \E\big(\ell\big(\mathcal T_{V,h}^\downarrow\big)\big)
<
    \infty
$.
\end{thm}

Before proving Theorem \ref{Lemma_Law_of_Slopes}, we introduce some notation.
For a slope $(T(i),\ 0\leq i \leq \ell(T))$ (recall that $T(0)=0$),
we define the slope
\begin{equation}\label{eq_def_zeta}
    \zeta(T)
:=
\big(
    T[\ell(T)-i]
-
    T[\ell(T)],
\
    0\leq i \leq \ell(T)
\big),
\end{equation}
with $\zeta\circ\zeta$ being identity (since $T(0)=0$ when $T$ is a slope).

\noindent{\bf Proof of Theorem \ref{Lemma_Law_of_Slopes}:}
Let $h>0$.
Applying (\cite{DGP_Collision_Transient}, Proposition 5.2, {\bf (ii)}),
$    \big(V\big[m_{1}^{(V)}(h)+x\big]-V\big[m_{1}^{(V)}(h)\big],
    \ 0\leq x \leq \tau_1^{(V)}(h)-m_{1}^{(V)}(h)\big),
$
is equal in law to
$
    \big(V(k), \  0\leq k \leq T_V([h,+\infty[)
    \big)
$
conditioned on $\{ T_V([h,+\infty[) < T_V(]-\infty, 0[ \}$,
which proves the first part of {\bf (i)}.
The second one follows from the strong Markov property applied to $(V(k), \ k\geq 0)$ at stopping time
$
    \tau_1^{(V)}(h)
$,
which is equal to
$
    m_{1}^{(V)}(h)
    +
    T_{\mathcal T_{V,h}^\uparrow}([h,+\infty[)
$.

We now prove some more general results, which will also be useful later.
Due to Lemma \ref{Lemma_Only_h_extrema}, the $m^{(V)}_{2i+1}(h)$, $i\geq 1$, are left $h$-minima,
 the $m^{(V)}_{2i+2}(h)$, $i\geq 0$, are left $h$-maxima,
and the $m^{(V)}_i(h)$, $i\geq 2$, are the only left $h$-extrema in $\big[\tau^{(V)}_1(h),+\infty\big[$.
However,  $m^{(V)}_{1}(h)$ is not necessarily a left $h$-minimum, depending on
the values taken by $(V(k),\  k\leq 0)$.

For $k\geq 1$, let
$
    \widehat \theta_{k,h}^{(\ell)}(V)
:=
    \widehat \theta_{k,h}^{(\ell)}
:=
    \big(V\big(m^{(V)}_{k}(h)-x\big)-V\big(m^{(V)}_{k}(h)\big),
    0\leq x \leq m^{(V)}_k(h)- \tau^{(V)}_{k-1}(h)\big)
$
and
$
    \widehat \theta_{k,h}^{(r)}(V)
:=
    \widehat \theta_{k,h}^{(r)}
:=
    \big(V\big(m^{(V)}_k(h)+x\big)-V\big(m^{(V)}_k(h)\big),\
    0\leq x \leq \tau^{(V)}_k(h)-m^{(V)}_k(h)\big)
$.
According to (\cite{DGP_Collision_Transient}, Proposition 5.2, {\bf (i)}),
the processes $\widehat \theta_{1,h}^{(\ell)}(V)$ and $\widehat \theta_{1,h}^{(r)}(V)$
are independent.
Also,
$
    \widehat \theta_{2,h}^{(\ell)}(V)
=
    -\widehat \theta_{1,h}^{(\ell)}\big(V\big(\tau_1^{(V)}(h)\big)-V\big(\tau_1^{(V)}(h)+.\big)\big)
$
and
$
    \widehat \theta_{2,h}^{(r)}(V)
=
    -\widehat \theta_{1,h}^{(r)}\big(V\big(\tau_1^{(V)}(h)\big)-V\big(\tau_1^{(V)}(h)+.\big)\big)
$,
so it follows from the previous result and from the strong Markov property applied at stopping time $\tau^{(V)}_1(h)$ that
$\widehat \theta_{2,h}^{(\ell)}(V)$ and $\widehat \theta_{2,h}^{(r)}(V)$ are independent
%. This and the strong Markov property applied at stopping time $\tau^{(V)}_1(h)$, show
and more precisely that all the trajectories
$\widehat \theta_{k,h}^{(\ell)}(V)$ and $\widehat \theta_{k,h}^{(r)}(V)$, $k\in\{1,2\}$, are independent.
Applying the same procedure by induction,
with the strong Markov property applied successively at stopping times
$\tau^{(V)}_{k}(h)$, $k\geq 1$,
proves  that all the trajectories
$\widehat \theta_{k,h}^{(\ell)}$ and $\widehat \theta_{k,h}^{(r)}$, $k\geq 1$, are independent.

In what follows we will "glue" trajectories.
For two trajectories $(f(i), \ a\leq i \leq b)$ and $(g(i), \ c\leq i \leq d)$,
by gluing $g$ to the right of $f$, we mean defining a new function $j$
: $\{a, \dots, b+d-c\}\to\R$ such that
\begin{equation}\label{eq_def_glue}
    j(i)
=
    \text{Glue}(f,g)(i)
:=
\left\{
\begin{array}{ll}
f(i) & \text{ if } a\leq i \leq b,\\
f(b)+g(i-b+c)-g(c) & \text{ if } b\leq i \leq b+d-c.
\end{array}
\right.
\end{equation}
Thanks to the previous paragraph, the trajectories
%{\bf (ici le $(V)$ est en exposant contrairement aux $\theta^{(r,\ell)}$)}
\begin{equation}\label{eq_theta_kl_hat}
    \widehat \theta_{k,h}^{(V)}
:=
    \big(V\big(x+m^{(V)}_{k}(h)\big)-V\big(m^{(V)}_{k}(h)\big),\
    0 \leq x \leq  m^{(V)}_{k+1}(h)-m^{(V)}_{k}(h)\big),
\qquad
    k\in\N^*
\end{equation}
are independent, since the $k$-th one is obtained by gluing
$\widehat \theta_{k,h}^{(r)}$ and, to its right,
$
    \big(V\big(\tau^{(V)}_{k}(h)+x\big)-V\big(\tau^{(V)}_{k}(h)\big),\
    0\leq x \leq m^{(V)}_{k+1}(h)- \tau^{(V)}_k(h)\big)
=
    \zeta\big(\widehat \theta_{k+1,h}^{(\ell)}\big)
$
(with $\zeta$ defined in \eqref{eq_def_zeta}),
that is,
$
   \widehat \theta_{k,h}^{(V)}
=
    \text{Glue}\big[\widehat \theta_{k,h}^{(r)}, \zeta\big(\widehat \theta_{k+1,h}^{(\ell)}\big)\big]
$.

%which is obtained from $\widehat \theta_{k+1,h}^{(\ell)}$ by time inversion as in \eqref{eq_def_zeta},

Also by the strong Markov property applied at stopping time
%$\tau^{(V)}_{2k+i-1}(h)$,
$\tau^{(V)}_{2k}(h)$,
$
    \widehat \theta_{2k+i,h}^{(\ell)}(V)
=_{law}
    \widehat \theta_{i,h}^{(\ell)}(V)
$
and
$
    \widehat \theta_{2k+i,h}^{(r)}(V)
=_{law}
    \widehat \theta_{i,h}^{(r)}(V)
$
for every $k\geq 1$ and $i\in\{1,2\}$.
Consequently, using the previous paragraph,
$
    \widehat \theta_{2k+1,h}^{(V)}
=
    \text{Glue}\big[\widehat \theta_{2k+1,h}^{(r)}, \zeta\big(\widehat \theta_{2k+2,h}^{(\ell)}\big)\big]
=_{law}
    \text{Glue}\big[\widehat \theta_{1,h}^{(r)}, \zeta\big(\widehat \theta_{2,h}^{(\ell)}\big)\big]
=
    \widehat \theta_{1,h}^{(V)}
=
    \mathcal T_{V,h}^\uparrow
$
and
$
    \widehat \theta_{2k+2,h}^{(V)}
=
    \text{Glue}\big[\widehat \theta_{2k+2,h}^{(r)}, \zeta\big(\widehat \theta_{2k+3,h}^{(\ell)}\big)\big]
=_{law}
    \text{Glue}\big[\widehat \theta_{2,h}^{(r)}, \zeta\big(\widehat \theta_{3,h}^{(\ell)}\big)\big]
=
    \widehat \theta_{2,h}^{(V)}
=
    \mathcal T_{V,h}^\downarrow
$
for every $k\in\N$.

Finally, by the strong Markov property applied at time $\tau^{(-V)}_1(h)$,
$\big(\widehat \theta_{2,h}^{(\ell)}(-V), \widehat \theta_{2,h}^{(r)}(-V)\big)$
is equal in law to
$\big(-\widehat \theta_{1,h}^{(\ell)}(V), -\widehat \theta_{1,h}^{(r)}(V)\big)$.
Similarly,
$\big(\widehat \theta_{3,h}^{(\ell)}(-V), \widehat \theta_{3,h}^{(r)}(-V)\big)
=_{law}
\big(-\widehat \theta_{2,h}^{(\ell)}(V), -\widehat \theta_{2,h}^{(r)}(V)\big)$.
As a consequence,
$
    \mathcal T_{-V,h}^\downarrow
=
    \text{Glue} \big[\widehat \theta_{2,h}^{(r)}(-V), \zeta\big(\widehat \theta_{3,h}^{(\ell)}(-V)\big)\big]
=_{law}
    \text{Glue}\big[-\widehat \theta_{1,h}^{(r)}(V), \zeta\big(-\widehat \theta_{2,h}^{(\ell)}(V)\big)\big]
$
$
=
    -\text{Glue}\big[\widehat \theta_{1,h}^{(r)}(V), \zeta\big(\widehat \theta_{2,h}^{(\ell)}(V)\big)\big]
=
    -    \mathcal T_{V,h}^\uparrow
$.
%$\mathcal T_{-V,h}^\downarrow=\widehat \theta_{2,h}^{(-V)}$, obtained as previously from
%$\widehat \theta_{2,h}^{(r)}(-V)$ and $\widehat \theta_{3,h}^{(\ell)}(-V)$ by gluing and time reversal,
%is equal in law to $-\widehat \theta_{1,h}^{(V)}=-\mathcal T_{V,h}^\uparrow$,
%obtained from $-\widehat \theta_{1,h}^{(r)}(V)$ and $-\widehat \theta_{2,h}^{(\ell)}(V)$.
Also, applying this to $-V$ instead of $V$ gives
$\mathcal T_{-V,h}^\uparrow=_{law}-\mathcal T_{V,h}^\downarrow$, which ends the proof of
{\bf (ii)}.

We now prove {\bf (iii)}. Due to \eqref{eqRecurrence} and \eqref{eq_def_sigma},
there exist $a>0$ such that $\p[V(1)\geq a]=:b>0$. Let $d:=\lfloor h/a\rfloor +1$.
Now, notice that
$
    \tau^{(V)}_{1}(h)
\leq
    d(N_d+1)
$,
where
$
    N_d
:=
    \min\{i\in\N,\ \forall 0\leq k\leq d, \, V(i d +k) -V(i d) \geq a k\}
$.
Hence,
$
    \E\big(\tau^{(V)}_{1}(h)\big)
\leq
    d(\E(N_d)+1)
<
    \infty
$
since $N_d$ is a geometric r.v. with parameter
$\p[\forall 0\leq k \leq d,\, V(k)\geq a k]\geq b^d>0$.
Using the strong Markov property,
we get similarly
$
    \E\big(\tau^{(V)}_{2}(h)-\tau^{(V)}_{1}(h)\big)
<
\infty
$.
Consequently,
$
    \E\big(\ell\big(\mathcal T_{V,h}^\uparrow\big)\big)
=
    \E\big(m_{2}^{(V)}(h)-m_{1}^{(V)}(h)\big)
\leq
    \E\big(\tau^{(V)}_{2}(h)\big)
<
\infty
$.
Finally, applying this to $-V$, we get
$
    \E\big[\ell\big(\mathcal T_{V,h}^\downarrow\big)\big]
=
    \E\big[\ell\big(\mathcal T_{-V,h}^\uparrow\big)\big]
<
    \infty
$,
since
$\mathcal T_{-V,h}^\uparrow=_{law}-\mathcal T_{V,h}^\downarrow$
by {\bf(ii)}.
This proves {\bf(iii)}.
\hfill$\Box$

%%%%%%%%%%%%%%%%%%%%%%%%%%%%%%%%%%%%%%%%%%%%%%%%%%%%%%%%%%%%%%%%%%%%%%%%%%%%%%%%%%%%%%%%%%%%

%%%%%%%%%%%%%%%%%%%%%%%%%%%%%%%%%%%%%%%%%%%%%%%%%%%%%%%%%%%%%%%%%%%%%%%%%%%%%%%%%%%%%%%%%%%%

%\subsection{Left $h$-extrema and renewal theory}

\subsection{Independence and law of translated left $h$-slopes via renewal theory}

%Let $\mathcal T_V^\uparrow$ having the same law as $\theta[T_1(V,h)]$ when
%$b_h>0$
%%$b_0(V,h)>0$
%(law of non central upward slopes),
%and $\mathcal T_V^\downarrow$ having the same law as $\theta[T_2(V,h)]$ when
%$b_h>0$
%%$b_0(V,h)>0$
%(law of non central downward slopes).

Notice that the law of $V$ may be nonsymmetric, so $\mathcal T_{V,h}^\uparrow$
and $-\mathcal T_{V,h}^\downarrow=_{law}\mathcal T_{-V,h}^\uparrow$
may have a different law, contrarily to what happens for Brownian motion
(imagine for example that the jumps of $V$ belong to $[-2,-1]\cup [4,5]$).

The following theorem is proved simultaneously as the next one. It says that for $h>0$,
roughly speaking,
conditionally on the central left $h$-slope $T_0(V,h)$ being upward (or being downward),
the translated left $h$-slopes $\theta[T_i(V, h)]$, $i\in\Z^*$, are independent
and are independent of the (non translated) central left $h$-slope $T_0(V,h)$,
and that the translated left $h$-slopes $\theta[T_i(V,h)]$, $i\in\Z^*$,
have the same law as $\mathcal T_{V,h}^\uparrow$ (under $\p$) for the upward ones
(ie the ones with $i\in(2\Z)-\{0\}$ when $T_0(V,h) $ is upward,
the ones for $i\in(2\Z+1)$ when $T_0(V,h) $ is downward)
and the same law as $\mathcal T_{V,h}^\downarrow$ (under $\p$) for the downward ones (the other ones).

We denote by $\mathscr{L}\big(\mathcal T_{V,h}^\uparrow\big)$
\big(resp. $\mathscr{L}\big(\mathcal T_{V,h}^\downarrow\big)$\big)
the law of $\mathcal T_{V,h}^\uparrow$ \big(resp. $\mathcal T_{V,h}^\downarrow$\big)
under $\p$.

\begin{thm}\label{Lemma_Independence_h_extrema}
Let $h>0$.
{\bf (i)} Conditionally on $\big\{V(x_1(V,h))>V(x_0(V,h))\big\}$ (i.e. on the central left $h$-slope $T_0(V,h)$ being upward),
the $\theta[T_{2i+1}(V,h)]$, $i\in\Z$ have
%the same law as $\mathcal T_{V,h}^\downarrow$
the law $\mathscr{L}\big(\mathcal T_{V,h}^\downarrow\big)$
whereas the $\theta[T_{2i}(V,h)]$, $i\in\Z^*$
%have the same law as $\mathcal T_{V,h}^\uparrow$,
have the law $\mathscr{L}\big(\mathcal T_{V,h}^\uparrow\big)$,
and
$(\theta[T_0(V,h)], x_0(V,h),$ $x_1(V,h))$, $\theta[T_{i}(V,h)]$, $i\in\Z^*$ are independent.

{\bf (ii)} Conditionally on $\big\{V(x_1(V,h))<V(x_0(V,h))\big\}$ (i.e. on the central left $h$-slope $T_0(V,h)$ being downward),
the $\theta[T_{2i+1}(V,h)]$, $i\in\Z$
%have the same law as $\mathcal T_{V,h}^\uparrow$,
have the law $\mathscr{L}\big(\mathcal T_{V,h}^\uparrow\big)$,
whereas the $\theta[T_{2i}(V,h)]$, $i\in\Z^*$ have
%the same law as $\mathcal T_{V,h}^\downarrow$,
the law $\mathscr{L}\big(\mathcal T_{V,h}^\downarrow\big)$,
and
$(\theta[T_0(V,h)], x_0(V,h), x_1(V,h))$, $\theta[T_{i}(V,h)]$, $i\in\Z^*$ are independent.

%the $\theta[T_{i}(V,h)]$, $i\in\Z^*$ are independent
%and independent of $(\theta[T_0(V,h)], x_0(V,h), x_1(V,h))$.
\end{thm}

However the law of the central left $h$-slope $T_0(V,h)$ is different.
It is provided by the following renewal theorem.

\begin{thm}\label{Lemma_Central_Slope}
Let $h>0$, $\Delta_0\subset \Z$ and $\Delta_1\subset \Z$.
For $A\in\{\bigsqcup_{t\in\mathbb N^*}A_t\, :\, \forall t\in\mathbb N^*,\ A_t\in\mathcal B(\mathbb R_+^t)\}$
%$A\subset\bigsqcup_{i=1}^\infty (\{0\}\times \R_+^{i})$
%(so that $A$ contains only upward slopes), we have
(so that the only slopes in $A$  are upward slopes), we have
\begin{eqnarray}
&&
    \p[\theta(T_0(V,h))\in A, \, x_0(V,h)\in\Delta_0, \, x_1(V,h)\in\Delta_1]
\nonumber\\
& = &
    \frac{
        \E\big[
        \sharp\big\{0\leq i < \ell(\mathcal T_{V,h}^\uparrow), \, (-i)\in\Delta_0,
        \, (\ell(\mathcal T_{V,h}^\uparrow)-i)\in\Delta_1\big\}
        \un_{\{\mathcal T_{V,h}^\uparrow\in A\}}
        \big]
    }
    {\E[\ell(\mathcal T_{V,h}^\uparrow)+\ell(\mathcal T_{V,h}^\downarrow)]}.
\label{eq_Central_Slope_Upward}
\end{eqnarray}
%\begin{equation}\label{eq_Central_Slope_Upward}
%    \p\big[\theta(T_0(V,h))\in A\big]
%=
%    \frac{\E[\ell(\mathcal T_V^\uparrow){\bf 1}_{\{\mathcal T_V^\uparrow\in A\}}]}
%    {\E[\ell(\mathcal T_V^\uparrow)+\ell( \mathcal T_V^\downarrow)]}
%=
%    \sum_{k=0}^\infty \p\big[\mathcal T_V^\uparrow\in A  \mid \ell( \mathcal T_V^\uparrow)=k\big]
%    \frac{k \p\big[\ell( \mathcal T_V^\uparrow)=k\big]}{\E[\ell(\mathcal T_V^\uparrow)+\ell(\mathcal T_V^\downarrow)]}.
%\end{equation}
Moreover if
$A\in\{\bigsqcup_{t\in\mathbb N^*}A_t\, :\, \forall t\in\mathbb N^*,\ A_t\in\mathcal B(\mathbb R_-^t)\}$
%$A\subset\cup_{i=0}^\infty (\{0\}\times \R_-^{i})$
%(so that $A$ contains only downward slopes), then
(so that the only slopes in $A$  are downward slopes), then
\begin{eqnarray}
&&
    \p[\theta(T_0(V,h))\in A, \, x_0(V,h)\in\Delta_0, \, x_1(V,h)\in\Delta_1]
\nonumber\\
& = &
    \frac{
        \E\big[
        \sharp\big\{0\leq i < \ell\big(\mathcal T_{V,h}^\downarrow\big), \, (-i)\in\Delta_0, \,
        \big(\ell\big(\mathcal T_{V,h}^\downarrow\big)-i\big)\in\Delta_1\big\}
        \un_{\{\mathcal T_{V,h}^\downarrow\in A\}}
        \big]
    }
    {\E\big[\ell\big(\mathcal T_{V,h}^\uparrow\big)+\ell\big(\mathcal T_{V,h}^\downarrow\big)\big]}.
\label{eq_Central_Slope_Downward}
\end{eqnarray}
Finally, for all  nonnegative function, $\varphi$ :
$\bigsqcup_{t\in\mathbb N^*}\mathbb R^t\rightarrow \mathbb [0,+\infty[$,
measurable with respect to the $\sigma$-algebra
$\{\bigsqcup_{t\in\mathbb N^*}A_t\, :\, \forall t\in\mathbb N^*,\ A_t\in\mathcal B(\mathbb R^t)\}$,
\begin{eqnarray}
&&
    \E\Big[\varphi\big[\theta(T_0(V,h))\big]\un_{\{x_0(V,h)\in\Delta_0\}}\un_{\{x_1(V,h)\in\Delta_1\}}\Big]
\nonumber\\
& = &
    \frac{
        \E\big[
        \sharp\big\{0\leq i < \ell\big(\mathcal T_{V,h}^\uparrow\big), \, (-i)\in\Delta_0, \,
        \big(\ell\big(\mathcal T_{V,h}^\uparrow\big)-i\big)\in\Delta_1\big\}
        \varphi\big(\mathcal T_{V,h}^\uparrow\big)
        \big]
    }
    {\E[\ell(\mathcal T_{V,h}^\uparrow)+\ell(\mathcal T_{V,h}^\downarrow)]}
\nonumber\\
&&
+
    \frac{
        \E\big[
        \sharp\big\{0\leq i < \ell\big(\mathcal T_{V,h}^\downarrow\big), \, (-i)\in\Delta_0, \,
        \big(\ell\big(\mathcal T_{V,h}^\downarrow\big)-i\big)\in\Delta_1\big\}
        \varphi\big(\mathcal T_{V,h}^\downarrow\big)
        \big]
    }
    {\E\big[\ell\big(\mathcal T_{V,h}^\uparrow\big)+\ell\big(\mathcal T_{V,h}^\downarrow\big)\big]}.
\label{eq_Central_Slope_General_phi}
\end{eqnarray}
\end{thm}

\noindent{\bf Proof of Theorems \ref{Lemma_Independence_h_extrema} and \ref{Lemma_Central_Slope}:}
Let $h>0$, $\Delta_0\subset \Z$, $\Delta_1\subset \Z$,
$q\leq 0 \leq r$,
%$q\leq r$,
and $B_i\in\{\bigsqcup_{t\in\mathbb N^*}A_t\, :\, \forall t\in\mathbb N^*,\ A_t\in\mathcal B(\mathbb R^t)\}=:\G$,
for $q\leq i\leq r$.
We first assume that
$B_0\in\{\bigsqcup_{t\in\mathbb N^*}A_t\, :\, \forall t\in\mathbb N^*,\ A_t\in\mathcal B(\mathbb R_+^t)\}$,
%$A\in\bigsqcup_{i=0}^\infty (\{0\}\times \R_+^{i})$
so that $B_0$ contains only upward slopes.
%We can restrict ourselves to $\Delta_0\subset \Z_-$ and $\Delta_1\subset \N^*$.

For $t\in\Z$ and $v\in\VVV$, let, loosely speaking, $T_0(t, v, h)$
be the left $h$-slope around $t$ for $v$,
that is, the left $h$-slope whose of $v$ domain contains $t$,
and denote its domain as $[x_0(t, v, h), x_1(t, v,h)]$.
More precisely and more generally, for $j\in\Z$,
we define $T_j(t, v, h)=T_{i+j}(v,h)$  if and only if $x_i(v,h)\leq t <x_{i+1}(v,h)$,
and for this unique $i$,
$x_j(t, v, h):=x_{i+j}(v,h)$ for $j\in\Z$
(recall that the notations $x_{i+j}$, $T_{i+j}$ are defined before and after \eqref{eqDefbh}).
%For $t\in\Z$, let, lossely speaking, $T_0(t, V, h)$ be the left $h$-slope around $t$ for $V$,
%that is, the left $h$-slope whose domain contains $t$,
%and denote its domain as $[x_0(t, V, h), x_1(t, V,h)]$.
%More precisely and more generally, for $j\in\Z$,
%we define $T_j(t, V, h)=T_{i+j}(V,h)$  iff $x_i(V,h)\leq t <x_{i+1}(V,h)$,
%and for this unique $i$,
%$x_j(t, V, h):=x_{i+j}(V,h)$ for $j\in\Z$
%(recall that the notations $x_{i+j}$, $T_{i+j}$ are defined before and after \eqref{eqDefbh}).
%and
%$x_1(t, V,h):=x_{i+1}(V,h)$.
We also introduce $V_{-t}(k):=V(k-t)-V(-t)$ for $t\in\Z$ and $k\in\Z$.
Hence, for $t\in\N$,
\begin{eqnarray}
&&
   \p\bigg(\big\{x_0(V,h)\in\Delta_0, \, x_1(V,h)\in\Delta_1\big\}
   \cap \bigcap_{i=q}^r\big\{\theta[T_i(V,h)]\in B_i\big\}\bigg)
\label{eqIntermediaireCalcul_Independent_Slopes-1}
\\
& = &
   \p\bigg(\big\{(x_0(t, V_{-t},h)-t)\in\Delta_0, \, (x_1(t, V_{-t},h)-t)\in\Delta_1\big\}
   \cap \bigcap_{i=q}^r\big\{\theta[T_i(t, V_{-t},h)]\in B_i\big\}\bigg)
\nonumber\\
& = &
    \p[E_B(t)]
\nonumber\\
& = &
    \p\big[ E_B(t), \, m_{-q+3}^{(V)}(h) \leq t \big]
+
    \p\big[E_B(t), \, m_{-q+3}^{(V)}(h) > t\big],
\label{eqIntermediaireCalcul_Independent_Slopes0}
\end{eqnarray}
where
$
    E_B(t)
:=
    \{
    (x_0(t,V,h)-t)\in \Delta_0, \, (x_1(t,V,h)-t)\in \Delta_1
    \}
\cap
    \cap_{i=q}^r\big\{\theta[T_i(t, V,h)]\in B_i\big\}
$,
because
%$\theta(T_0(V,h))$ is the same as the image by $\theta$ of the slope around $t$ for
%$V_{-t}$,
%=(V(k-t)-V(-t),\ k\in\Z)$,
$x_j(V,h)=x_j(t, V_{-t}, h)-t$ for $j\in\Z$,
$\theta[T_i(V,h)]=\theta[T_i(t, V_{-t},h)]$ for $i\in\Z$,
and $V_{-t}$ has the same law as $V$.

Let $(Y_k)_{k\in\Z}$ be a sequence of independent left $h$-slopes, such that
$Y_{2k}=_{law} \mathcal T_{V,h}^\uparrow$
and
$Y_{2k+1}=_{law} \mathcal T_{V,h}^\downarrow$
for every $k\in\Z$.
We glue sequentially (see \eqref{eq_def_glue}) $Y_0, Y_1, \dots, Y_k,\dots$ to get a process
$(Y(i), \ i\in\N)$, starting from $0$ (i.e. $Y(i)=Y_0(i)$ for $0\leq i \leq \ell(Y_0)$).
This process $(Y(i), \ i\in\N)$
has the same law as
$\big(V\big[m_{1}^{(V)}(h)+x\big]-V\big[m_{1}^{(V)}(h)\big], \ x\geq 0\big)$.
%(because the $\widehat \theta_{k,h}^{(V)}$, $k\in\N^*$ are independent, and
%$
%    \widehat \theta_{2k+1,h}^{(V)}
%=_{law}
%%    \widehat \theta_{1,h}^{(V)}
%%=
%%    \mathcal T_{V,h}^\uparrow
%    Y_{2k}
%$
%and
%$
%    \widehat \theta_{2k+2,h}^{(V)}
%=_{law}
%%    \widehat \theta_{2,h}^{(V)}
%%=
%%    \mathcal T_{V,h}^\downarrow
%    Y_{2k+1}
%$
%for every $k\in\N$,
%see \eqref{eq_theta_kl_hat} and below).
Indeed, this last process can be obtained from gluing
$\widehat \theta_{1,h}^{(V)}, \widehat \theta_{2,h}^{(V)}, \dots, \widehat \theta_{k,h}^{(V)}, \dots$
(see \eqref{eq_theta_kl_hat}), which are independent and such that $\widehat \theta_{k,h}^{(V)}=_{law} Y_{k-1}$, $k\in\N^*$ by
definition of the $Y_k$ and the law of the $\widehat \theta_{k,h}^{(V)}$ (see after \eqref{eq_theta_kl_hat}).
We also glue sequentially the $Y_k$, $k<0$ in the same way to the left of $(Y(i),\ i\in\N)$,
so that $Y_k$ is followed by $Y_{k+1}$, $k\in\Z$.
%{\bf (detailler ?)}.
The resulting process is denoted by $(Y(i), \ i\in\Z)$, with $Y(0)=0$.
Notice that $x_0(Y,h)=0$, for $i\in\N^*$. We also have $x_i(Y,h)=\ell(Y_0)+\dots+\ell(Y_{i-1})$,
and for $i\in\Z_-^*$, we have $x_i(Y,h)=-\ell(Y_{-1})-\dots-\ell(Y_{i})$.

%For $t\in\Z$, we denote, similarly as before, by $\widetilde T_0(t, Y, h)$ the left $h$-slope of $Y$ around $t$.
%More precisely and more generally,
%we define, for $j\in\Z$,
%$\widetilde T_j(t, Y, h)=T_{i+j}(Y,h)$
%for the unique $i\in\Z$ such that $x_i(Y,h)\leq t <x_{i+1}(Y,h)$,
%Also we define $\widetilde x_j(t, Y, h):=x_{i+j}(Y,h)$ with the previous $i$,
%and so the support of $\widetilde T_j(t, Y, h)$
%is $\big[\widetilde x_j(t, Y, h), \widetilde x_{j+1}(t, Y, h)\big]$.

We can assume without loss of generality that $q\in(-2\N^*)$, so $T_q(V,h)$ is an upward slope when $\theta(T_0(V,h))\in B_0$.
Using, in the second equality, the fact that
$
    \big(V\big(x+m^{(V)}_3(h)\big)-V\big(m^{(V)}_3(h)\big),\
    x \geq 0\big)
$
has the same law as $(Y(x),\, x\geq 0)$
(see \eqref{eq_theta_kl_hat} and below),
and is independent of $\big(V(x),\, x\leq m_3^{(V)}(h)\big)$
(see the paragraph before \eqref{eq_def_glue}),
we have for $t\in\N$,
\begin{eqnarray}
    \p\big[E_B(t), \, m_{-q+3}^{(V)}(h) \leq t\big]
& = &
    \sum_{y=0}^{t}
    \p\big[E_B(t), \, m_3^{(V)}(h)=y,\, m_{-q+3}^{(V)}(h) \leq t\big]
\nonumber\\
& = &
    \sum_{y=0}^{t}
    \p\big[m_3^{(V)}(h)=y\big]h_B(t-y)
\label{eqIntermediaireCalculSlopeCentrale1}
\end{eqnarray}
where
$
%    h_{B,q}(p)
    h_B(p)
:=
    \p\big[\widetilde E_B(p)\cap \{\ell(Y_0)+\dots+\ell(Y_{-q-1}) \leq p\}\big]
$,
%where $\ell(Y_0)+\dots+\ell(Y_{-q-1})=0$ when $q=0$, and
with
$
    \widetilde E_B(p)
:=
    \{
    \big(x_0(p,Y,h)-p\big) \in \Delta_0,\, \big(x_1(p,Y,h)-p\big) \in \Delta_1
    \}
\cap
    \cap_{i=q}^r\big\{\theta\big[T_i(p, Y,h)\big]\in B_i\big\}
$,
$p\in\N$.
Indeed, on $\big\{m_{-q+3}^{(V)}(h) \leq t\big\}$,
we have
$x_0(t,V,h)\geq m_{-q+3}^{(V)}(h)$ thanks to Lemma \ref{Lemma_Only_h_extrema},
thus
$x_q(t,V,h)\geq m_{3}^{(V)}(h)$, so
$E_B(t)$ depends only on
$
    \big(V\big(x+m^{(V)}_3(h)\big)-V\big(m^{(V)}_3(h)\big),\
    x \geq 0\big)
=:(Y'(x),\ x\geq 0)$,
with $x_i(t,V,h)=x_i(t-y,Y',h)+y$  and
$
    \theta[T_i(t,V,h)]
=
    \theta[T_i(t-y,Y',h)]
$
for $i\geq q$ on $\big\{m_3^{(V)}(h)=y\big\}$
and $Y'=_{law}Y$.

We want to prove that $h_B(p)$ has a limit as $p\to+\infty$.
For $p\in\N$, let
\begin{eqnarray*}
    a_p'
& := &
    \p\bigg[
        \big\{(\ell(Y_0)+\dots+
        \ell(Y_{-q-1})-p)
        \in\Delta_0,\, \big(\ell(Y_0)+\dots+\ell(Y_{-q})-p\big)\in\Delta_1
        \big\}
\\
&&
\qquad
    \cap
        \big\{0\leq p-\ell(Y_0)-\dots-\ell(Y_{-q-1})    <\ell(Y_{-q})\big\}
        \cap
        \bigcap_{i=q}^r \{Y_{i-q}\in B_i\}
    \bigg].
\end{eqnarray*}
We have for $p\in\N$, since $q\in(-2\N^*)$,
%since $\Delta_1\subset \N^*$,
\begin{eqnarray}
&&
    h_B(p)
\nonumber\\
& = &
    \p\big[\widetilde E_B(p),\, \ell(Y_0)+\dots+\ell(Y_{-q-1})\leq  p<\ell(Y_0)+\dots+\ell(Y_{-q-1})+\ell(Y_{-q}) \big]
\label{eq_line_proba_a_commenter_1}
\\
&&
    +
    \p\big[\widetilde E_B(p),\, \ell(Y_0)+\dots+\ell(Y_{-q})\leq p <\ell(Y_0)+\dots+\ell(Y_{-q})+\ell(Y_{-q+1}) \big]
\label{eq_line_proba_a_commenter}
\\
&& +
    \sum_{y=0}^p
    \p\big[\ell(Y_0)+\ell(Y_1)=y, \widetilde E_B(p),
    \, \ell(Y_2)+\dots+\ell(Y_{-q+1})\leq p-y \big]
\label{eq_line_proba_a_commenter_2}
\\
& = &
    a_p'
    +
    0
+
    \sum_{y=0}^p
    \p\big[ \ell(Y_0)+\ell(Y_1)=y \big]
    \p\big[\widetilde E_B(p-y), \, \ell(Y_0)+\dots+\ell(Y_{-q-1})\leq  p-y\big]
\nonumber\\
& = &
    a_p'
    +
    \sum_{y=0}^p
    \p\big[\ell(Y_0)+\ell(Y_1)=y \big]h_B(p-y).
\nonumber
\end{eqnarray}
Indeed in the probability in line \eqref{eq_line_proba_a_commenter_1},
the image by $\theta$ of the slope
$T_0(p, Y, h)$
%$\widetilde T_0(p, Y, h)$
containing $p$ is $Y_{-q}$
and the $0$ in the second equality comes from the fact that
on the set inside the probability of line \eqref{eq_line_proba_a_commenter},
$\theta\big[T_0(p, Y, h)\big]=Y_{-q+1}$
%$\theta\big[\widetilde T_0(p, Y, h)\big]=Y_{-q+1}$
is a downward slope,
whereas $B_0$ contains only upward slopes,
and in the sets appearing in \eqref{eq_line_proba_a_commenter_2},
there exists $j\geq 0$ such that
$x_i(p-y, Y,h)\big]=\ell(Y_2)+\dots+\ell(Y_{-q+1+j+i})$
and
$\theta\big[T_i(p-y, Y,h)\big]=Y_{-q+2+j+i}$ for $i\geq q$
%$\theta\big[\widetilde T_0(p, Y,h)\big]=Y_{-q+2+j}$
and
$(Y_0, Y_1)$ is independent of $(Y_{i+2},\ i\geq 0)$, which has the same law as
$(Y_i,\ i\geq 0)$.

So, $h_B(p)$ is solution of the discrete time renewal equation
$h_p=a_p'+\sum_{k=0}^p f_k h_{p-k}$, $p\in\N$, with
$h_p=h_B(p)$ and $f_k=\p\big[\ell(Y_0)+\ell(Y_1)=k \big]$.
Notice that  $a_p'\geq 0$, $p\in\N$ and
$
    \sum_{p=0}^\infty a_p'
\leq
    \E\big[\ell(Y_0)+\dots+\ell(Y_{-q})\big]+a_0
\leq
    (|q|+1)\E\big[\ell\big(\mathcal T_{V,h}^\uparrow\big)+\ell\big(\mathcal T_{V,h}^\downarrow\big)\big]+1
<
    \infty
$
by our Theorem \ref{Lemma_Law_of_Slopes} {\bf (iii)}.
So, Theorem 2.2 of Barbu and Limnios \cite{Barbu_Limnios}
%{\bf (or Nummelin Esa ?)}
with its notation $X_n=\ell(Y_{2n-2})+\ell(Y_{2n-1})>0$, $n\geq 1$ so that $f_k=\p[X_1=k]$
and
$
    u_n
:=
    \sum_{m=0}^n \p[X_1+\dots + X_m=n]
=
    \sum_{m=0}^n \p[\ell(Y_0)+\ell(Y_1)+\dots + \ell(Y_{2m-2})+\ell(Y_{2m-1})=n]
$
with $X_1+\dots + X_0=0$ by convention,
give us that
this renewal equation has a unique solution, which is
\begin{equation*}\
%label{eq_Limite_Renewal}
    h_B(p)
=
    h_p
=
    (u*a')_p
=
    \sum_{k=0}^p u_{p-k}a_k',
\qquad
    p\in\N.
\end{equation*}
Let $n_1\in\N^*$ and $n_2\in\N^*$ be such that
$\p[T_V(h)=n_1 \mid T_V(h)<T_V(\R_-^*)]=:c_2>0$
and $\p[\ell(Y_1)=n_2]=:c_3>0$
and let $c_4:=\p[T_V(-h)<T_V(]0,\infty[)]>0$ due to \eqref{eqRecurrence} and \eqref{eq_def_sigma}.
Hence, using the law of $\mathcal T^\uparrow_{V,h}$ (see Theorem \ref{Lemma_Law_of_Slopes} {\bf (i)}),
 $\p[\ell(Y_0)=n_1]=\p\big[\ell\big(\mathcal T^\uparrow_{V,h}\big)=n_1\big]\geq c_2c_4>0$.
Also, $\p[\ell(Y_0)=n_1+1]\geq c_2\p[V(1)>0]c_4>0$.
Thus, $\p[\ell(Y_0)+\ell(Y_1)=n_1+n_2]>0$ and $\p[\ell(Y_0)+\ell(Y_1)=n_1+n_2+1]>0$,
and then the renewal chain $(X_1+\dots+X_n)_n$  is aperiodic.
It is also recurrent since $X_1<\infty$ a.s.,
e.g. because
 $\E(X_1)=\E\big[\ell\big(\mathcal T^\uparrow_{V,,h}\big)
+\ell\big(\mathcal T^\downarrow_{V,h}\big)\big]<\infty$
by Theorem \ref{Lemma_Law_of_Slopes} {\bf (iii)}.
So by the renewal theorem (see e.g. Barbu and Limnios \cite{Barbu_Limnios}, Theorem 2.6), we have
$
    u_p
\to_{p\to+\infty}
%    1/\E[\ell(Y_0)+\ell(Y_1)]
    1/\E(X_1)
=
    1/\E\big[\ell\big(\mathcal T^\uparrow_{V,h}\big)+\ell\big(\mathcal T^\downarrow_{V,h}\big)\big]
$.
Moreover since this renewal chain is recurrent and aperiodic and
$\sum_{p=0}^\infty |a_p'|<\infty$, we have
by the key renewal theorem (see e.g. Barbu and Limnios \cite{Barbu_Limnios}, Theorem 2.7),
\begin{equation}\label{eq_lim_gA}
    h_B(p)
=
    h_p
=
    \sum_{k=0}^p u_{p-k}a_k'
\to_{p\to+\infty}
    \frac{1}{\E[\ell(\mathcal T_{V,h}^\uparrow)+\ell(\mathcal T_{V,h}^\downarrow)]}
    \sum_{p=0}^\infty a_p'.
\end{equation}
Also, let
$
    A_{k_0,\dots, k_{r-q}}
:=
        \big\{p\in\N, \ (k_0+\dots+
        k_{-q-1}-p)
        \in\Delta_0,\,
        \big(k_0+\dots+k_{-q}-p\big)\in\Delta_1
        \big\}
\cap
        \big\{0\leq p-k_0-\dots-k_{-q-1}<k_{-q} \big\}
$
for $(k_0, \dots, k_{r-q})\in\N^{r-q+1}$. We have,
\begin{eqnarray}
&&
    \sum_{p=0}^\infty a_p'
\nonumber\\
& = &
    \sum_{p=0}^\infty \sum_{(k_0,\dots, k_{r-q})\in\N^{r-q+1}}
    \p\bigg[
        \bigcap_{j=0}^{r-q} \big\{\ell(Y_j)=k_j\big\}\cap \{Y_j\in B_{j+q}\}
    \cap \{p\in A_{k_0,\dots, k_{r-q}}\}
    \bigg]
\nonumber\\
%&&
%        \big\{(k_0+\dots+
%        k_{-q-1}-p)
%        \in\Delta_0,\,
%        \big(k_0+\dots+k_{-q}-p\big)\in\Delta_1
%        \big\}
%\nonumber\\
%&&
%\qquad
%    \cap
%        \big\{0\leq p-k_0-\dots-k_{-q-1}<k_{-q} \big\}
%    \bigg]
%\nonumber\\
& = &
    \sum_{(k_0,\dots, k_{r-q})\in\N^{r-q+1}}
        \bigg(
        \prod_{j=0}^{r-q}
        \p\Big[
          \big\{\ell(Y_j)=k_j\big\}\cap \{Y_j\in B_{j+q}\}
        \Big]
        \bigg)
    \sum_{p=0}^\infty {\bf 1}_{A_{k_0,\dots, k_{r-q}}}(p)
\nonumber\\
%&&
%\qquad
%\times
%    \sum_{p=0}^\infty
%        {\bf 1}_{\{(k_0+\dots+
%        k_{-q-1}-p)
%        \in\Delta_0,\,
%        \big(k_0+\dots+k_{-q}-p\big)\in\Delta_1
%        \big\}
%   \cap
%        \big\{0\leq p-k_0-\dots-k_{-q-1}<k_{-q} \big\}}
\nonumber\\
& = &
    \sum_{(k_0,\dots, k_{r-q})\in\N^{r-q+1}}
        \bigg(
        \prod_{0\leq j\leq r-q,\, j\neq -q}
        \p\Big[
          \big\{\ell(Y_j)=k_j\big\}\cap \{Y_j\in B_{j+q}\}
        \Big]
        \bigg)
\nonumber\\
&&
\times
    \E\bigg(
    {\bf 1}_{\{\ell(Y_{-q})=k_{-q}, \ Y_{-q}\in B_0\}}
    \sum_{p=0}^\infty
        {\bf 1}_{A_{k_0,\dots, k_{r-q}}}(p)
%        {\bf 1}_{\{(k_0+\dots+
%        k_{-q-1}-p)
%        \in\Delta_0,\,
%        \big(k_0+\dots+k_{-q}-p\big)\in\Delta_1
%        \big\}
%   \cap
%        \big\{0\leq p-k_0-\dots-k_{-q-1}<k_{-q} \big\}}
    \bigg)
\nonumber\\
& = &
    \sum_{(k_0,\dots, k_{r-q})\in\N^{r-q+1}}
        \bigg(
        \prod_{0\leq j\leq r-q,\, j\neq -q}
        \p\Big[
          \big\{\ell(Y_j)=k_j\big\}\cap \{Y_j\in B_{j+q}\}
        \Big]
        \bigg)
\nonumber\\
&&
\times
%    \sum_{k_{-q}=0}^\infty
    \E\bigg(
    {\bf 1}_{\{\ell(Y_{-q})=k_{-q}, \ Y_{-q}\in B_0\}}
    \sharp\{0\leq m <\ell(Y_{-q}), (-m)\in\Delta_0, (\ell(Y_{-q})-m)\in\Delta_1\}
    \bigg)
\nonumber\\
& = &
        \bigg(
        \prod_{0\leq j\leq r-q,\, j\neq -q}
        \p\Big[
          Y_j\in B_{j+q}
        \Big]
        \bigg)
\nonumber\\
&&
\times
    \E\bigg(
    {\bf 1}_{\{ Y_{-q}\in B_0\}}
    \sharp\{0\leq m <\ell(Y_{-q}), (-m)\in\Delta_0, (\ell(Y_{-q})-m)\in\Delta_1\}
    \bigg).
\label{eq_calcul_independence_slopes_Up}
\end{eqnarray}
Now, notice that by definition of $(Y_k)_{k\in\Z}$ and since $q\in(-2\N^*)$, the product in \eqref{eq_calcul_independence_slopes_Up} is equal to
\begin{equation}\label{eq_simplification_produit}
        \bigg(
        \prod_{q\leq i\leq r,\, i\neq 0,\, i\in(2\Z)}
        \p\Big[
         \mathcal T^\uparrow_{V,h} \in B_i
        \Big]
        \bigg)
\times
        \bigg(
        \prod_{q\leq i\leq r,\, i\in(2\Z+1)}
        \p\Big[
         \mathcal T^\downarrow_{V,h} \in B_i
        \Big]
        \bigg).
\end{equation}
The second probability in \eqref{eqIntermediaireCalcul_Independent_Slopes0}
is less than  $\p\big[m_{-q+3}^{(V)}(h) > t\big]$
and then it goes to $0$ as $t\to+\infty$ since
$m_{-q+3}^{(V)}(h)<\tau_{-q+3}^{(V)}(h)<\infty$ a.s.
since $V\in\VVV$ a.s.
%and
%because
%$
%    \E\big(\tau_{-q+3}^{(V)}(h)\big)
%\leq
%    (|q|+3)\E\big(\mathcal T^\uparrow_V+\mathcal T^\downarrow_V\big)
%<
%    \infty
%$ by Theorem \ref{Lemma_Law_of_Slopes}.

Combining this with \eqref{eqIntermediaireCalculSlopeCentrale1},
letting $t\to+\infty$  and applying the dominated convergence theorem
gives
$
\eqref{eqIntermediaireCalcul_Independent_Slopes-1}
%\eqref{eqIntermediaireCalcul_Independent_Slopes0}
=
    \lim_{p\to+\infty} h_B(p)
$
(since this limit exists by \eqref{eq_lim_gA}).
This, together with
%\eqref{eqIntermediaireCalculSlopeCentrale1},
\eqref{eq_lim_gA},
\eqref{eq_calcul_independence_slopes_Up},
 and \eqref{eq_simplification_produit} leads to
%and letting $t\to+\infty$ gives, thanks to the dominated convergence theorem,
\begin{align}
&
\eqref{eqIntermediaireCalcul_Independent_Slopes-1}
%\eqref{eqIntermediaireCalcul_Independent_Slopes0}
=
        \bigg(
        \prod_{q\leq i\leq r,\, i\neq 0,\, i\in(2\Z)}
        \p\Big[
         \mathcal T^\uparrow_{V,h} \in B_i
        \Big]
        \bigg)
\times
        \bigg(
        \prod_{q\leq i\leq r,\, i\in(2\Z+1)}
        \p\Big[
         \mathcal T^\downarrow_{V,h} \in B_i
        \Big]
        \bigg)
\nonumber\\
&
    \times
    \E\bigg(
    \frac{{\bf 1}_{\{ Y_{-q}\in B_0\}}}{\E[\ell(\mathcal T^\uparrow_{V,h})+\ell(\mathcal T^\downarrow_{V,h})]}
    \sharp\{0\leq m <\ell(Y_{-q}), (-m)\in\Delta_0, (\ell(Y_{-q})-m)\in\Delta_1\}
    \bigg).
\label{eqIntermediaireCalcul_Independent_Slopes_Presque_fini}
\end{align}
    Moreover, taking (only here) all the $B_i$ equal to $\bigsqcup_{t\in\mathbb N^*}\R^t$, except $B_0$ in
\eqref{eqIntermediaireCalcul_Independent_Slopes_Presque_fini},
we get
\begin{align}
&
   \p\big(x_0(V,h)\in\Delta_0, \, x_1(V,h)\in\Delta_1,\,
   \theta[T_0(V,h)]\in B_0
   \big)
\nonumber\\
& =
    \frac{
    \E\Big({\bf 1}_{\{ \mathcal T^\uparrow_{V,h} \in B_0\}}
    \sharp\{0\leq m <\ell(\mathcal T^\uparrow_{V,h}), (-m)\in\Delta_0, (\ell(\mathcal T^\uparrow_{V,h})-m)\in\Delta_1\}
    \Big)
    }{\E[\ell(\mathcal T^\uparrow_{V,h})+\ell(\mathcal T^\downarrow_{V,h})]},
\label{eq_Proba_Central_Slope_Preuve}
\end{align}
since $Y_{-q}$ has the same law as $\mathcal T^\uparrow_{V,h}$ because $q\in(2\Z)$.
This proves \eqref{eq_Central_Slope_Upward}.
Consequently, \eqref{eqIntermediaireCalcul_Independent_Slopes_Presque_fini}
becomes
\begin{align}
&
\eqref{eqIntermediaireCalcul_Independent_Slopes-1}
%\eqref{eqIntermediaireCalcul_Independent_Slopes0}
=
        \bigg(
        \prod_{q\leq i\leq r,\, i\neq 0,\, i\in(2\Z)}
        \p\Big[
         \mathcal T^\uparrow_{V,h} \in B_i
        \Big]
        \bigg)
\times
        \bigg(
        \prod_{q\leq i\leq r,\, i\in(2\Z+1)}
        \p\Big[
         \mathcal T^\downarrow_{V,h} \in B_i
        \Big]
        \bigg)
\nonumber\\
&
    \times
   \p\big(\big\{x_0(V,h)\in\Delta_0, \, x_1(V,h)\in\Delta_1\big\}
   \cap \big\{\theta[T_0(V,h)]\in B_0\big\}\big).
\label{eqIntermediaireCalcul_Independent_Slopes_fini_Up}
\end{align}
This proves Theorem \ref{Lemma_Independence_h_extrema} {\bf (i)}.

We now prove \eqref{eq_Central_Slope_Downward} and Theorem \ref{Lemma_Independence_h_extrema} {\bf (ii)}.
We assume that
$B_0\in\{\bigsqcup_{t\in\mathbb N^*}A_t\, :\, \forall t\in\mathbb N^*,\ A_t\in\mathcal B(\mathbb R_-^t)\}$,
%$A\in\bigsqcup_{i=0}^\infty (\{0\}\times \R_+^{i})$
so that $B_0$ contains only downward slopes.
Notice that $x_i(-V,h)=x_i(V,h)$ and
$
    \theta[T_i(-V,h)]
=
    -\theta[T_i(V,h)]
$ for $i\in\Z$.
Then,
$-B_0=\{-f, \     f\in B_0\}\in\{\bigsqcup_{t\in\mathbb N^*}A_t\, :\, \forall t\in\mathbb N^*,\ A_t\in\mathcal B(\mathbb R_+^t)\}$,
% contains only upward slopes,
and for each $q\leq i \leq r$,
$\theta(T_i(V,h))\in B_i$ iff $\theta(T_i(-V,h))\in (-B_i)$,
for which we can apply \eqref{eqIntermediaireCalcul_Independent_Slopes_fini_Up}
and \eqref{eq_Proba_Central_Slope_Preuve} as follows.
We get,
\begin{eqnarray}\label{eq_Proba_Independence_Downward_Preuve}
&&
   \p\bigg(\big\{x_0(V,h)\in\Delta_0, \, x_1(V,h)\in\Delta_1\big\}
   \cap \bigcap_{i=q}^r\big\{\theta[T_i(V,h)]\in B_i\big\}\bigg)
\\
& = &
   \p\bigg(\big\{x_0(-V,h)\in\Delta_0, \, x_1(-V,h)\in\Delta_1\big\}
   \cap \bigcap_{i=q}^r\big\{\theta[T_i(-V,h)]\in (-B_i)\big\}\bigg)
\nonumber\\
& = &
        \bigg(
        \prod_{q\leq i\leq r,\, i\neq 0,\, i\in(2\Z)}
        \p\Big[
         \mathcal T^\uparrow_{-V,h} \in (-B_i)
        \Big]
        \bigg)
\times
        \bigg(
        \prod_{q\leq i\leq r,\, i\in(2\Z+1)}
        \p\Big[
         \mathcal T^\downarrow_{-V,h} \in (-B_i)
        \Big]
        \bigg)
\nonumber\\
&&
    \times
    \frac{
    \E\Big({\bf 1}_{\{ \mathcal T^\uparrow_{-V,h} \in (-B_0)\}}
    \sharp\{0\leq m <\ell(\mathcal T^\uparrow_{-V,h}), (-m)\in\Delta_0, (\ell(\mathcal T^\uparrow_{-V,h})-m)\in\Delta_1\}
    \Big)
    }{\E[\ell(\mathcal T^\uparrow_{-V,h})+\ell(\mathcal T^\downarrow_{-V,h})]},
\nonumber\\
& = &
        \bigg(
        \prod_{q\leq i\leq r,\, i\neq 0,\, i\in(2\Z)}
        \p\Big[
         \mathcal T^\downarrow_{V,h} \in B_i
        \Big]
        \bigg)
\times
        \bigg(
        \prod_{q\leq i\leq r,\, i\in(2\Z+1)}
        \p\Big[
         \mathcal T^\uparrow_{V,h} \in B_i
        \Big]
        \bigg)
\nonumber\\
&&
    \times
    \frac{
    \E\Big({\bf 1}_{\{ \mathcal T^\downarrow_{V,h} \in B_0\}}
    \sharp\{0\leq m <\ell(\mathcal T^\downarrow_{V,h}), (-m)\in\Delta_0, (\ell(\mathcal T^\downarrow_{V,h})-m)\in\Delta_1\}
    \Big)
    }{\E[\ell(\mathcal T^\uparrow_{V,h})+\ell(\mathcal T^\downarrow_{V,h})]},
\nonumber
\end{eqnarray}
since
$\mathcal T_{-V,h}^\uparrow=_{law}-\mathcal T_{V,h}^\downarrow$
and
$
\E\big[\ell\big(\mathcal T_{-V,h}^\uparrow\big)+\ell\big(\mathcal T_{-V,h}^\downarrow\big)\big]
=
\E\big[\ell\big(\mathcal T_{V,h}^\downarrow\big)+\ell\big(\mathcal T_{V,h}^\uparrow\big)\big]
$
by Theorem \ref{Lemma_Law_of_Slopes} {\bf (ii)}.
Taking all the $B_i$, $i\neq 0$, equal to $\bigsqcup_{t\in\mathbb N^*}\R^t$, this proves \eqref{eq_Central_Slope_Downward}.
This, in turn, proves that
\begin{eqnarray*}
    \eqref{eq_Proba_Independence_Downward_Preuve}
& = &
        \bigg(
        \prod_{q\leq i\leq r,\, i\neq 0,\, i\in(2\Z)}
        \p\Big[
         \mathcal T^\downarrow_{V,h} \in B_i
        \Big]
        \bigg)
\times
        \bigg(
        \prod_{q\leq i\leq r,\, i\in(2\Z+1)}
        \p\Big[
         \mathcal T^\uparrow_{V,h} \in B_i
        \Big]
        \bigg)
\\
&&
    \times
    \p\big(x_0(V,h)\in\Delta_0, \, x_1(V,h)\in\Delta_1,\,
    \theta[T_0(V,h)]\in B_0\big),
\end{eqnarray*}
which proves Theorem \ref{Lemma_Independence_h_extrema} {\bf (ii)}.

In order to prove \eqref{eq_Central_Slope_General_phi}, we first show that
\eqref{eq_Central_Slope_General_phi} is true for $\varphi=\un_A$
for any $A\in\{\bigsqcup_{t\in\mathbb N^*}A_t\, :\, \forall t\in\mathbb N^*,\ A_t\in\mathcal B(\mathbb R^t)\}=\G$.
%$A\subset\cup_{i=0}^\infty (\{0\}\times \R^{i})$,
To this aim, let
%$B\in \{\bigsqcup_{t\in\mathbb N^*}A_t\, :\, \forall t\in\mathbb N^*,\ A_t\in\mathcal B(\mathbb R^t)\}$.
$A\in\G$.
We introduce $\SSS_\pm:=\bigsqcup_{t\in\mathbb N^*} \R_\pm^t$.
Applying \eqref{eq_Central_Slope_Upward}
to $A\cap \SSS_+$ (resp. \eqref{eq_Central_Slope_Downward} to $A\cap \SSS_-$)
proves
\eqref{eq_Central_Slope_General_phi} for $\varphi=\un_{A\cap \SSS_+}$
(resp. $\varphi=\un_{A\cap \SSS_-}$),
since the second (resp. first) expectation in \eqref{eq_Central_Slope_General_phi} is $0$
when for $\varphi=\un_{A\cap \SSS_+}$ (resp. $\varphi=\un_{A\cap \SSS_-}$),
because
$\mathcal T_{V,h}^\downarrow\notin \SSS_+$
(resp.
$\mathcal T_{V,h}^\uparrow\notin \SSS_-$).
%downward (resp. upward) slopes are nonpositive (resp. nonnegative) functions.
Also, \eqref{eq_Central_Slope_General_phi} is true for $\varphi=\un_{A\cap (\SSS_+\cup \SSS_-)^c}$
since every term is equal to $0$ in \eqref{eq_Central_Slope_General_phi} in this case,
since when $\theta[T_0(V,h)]$ is a downward (resp. upward) slope, it belongs to $\SSS_-$ (resp. $\SSS_+$)
and
$\mathcal T_{V,h}^\downarrow\in \SSS_-$
(resp.
$\mathcal T_{V,h}^\uparrow\in \SSS_+$).
Hence, adding \eqref{eq_Central_Slope_General_phi} in the three previous cases proves
that \eqref{eq_Central_Slope_General_phi} is true for $\varphi=\un_A$,
 for every $A\in\G$.

Then by linearity, \eqref{eq_Central_Slope_General_phi} is true for
every simple function $\sum_{i=1}^p \alpha_i {\bf 1}_{B_i}$ for $p\geq 1$, $\alpha_i\geq 0$ and $B_i\in\G$, $1\leq i \leq p$.
Finally, \eqref{eq_Central_Slope_General_phi} is true
for any nonnegative $\G$-measurable function by the monotone convergence theorem,
since every nonnegative $\G$-measurable function is the pointwise limit of
a nondecreasing sequence of nonnegative simple $\G$-measurable functions.
\hfill$\Box$

\subsection{A simple expression for $\p(b_h= x)$}
A first application of our renewal Theorem \ref{Lemma_Central_Slope}
is the following lemma, which contains key formulas to prove Theorem \ref{Th_Local_Limit_b_h}
and study the main contribution in Theorem \ref{Th_Local_Limit_Sinai}
(see e.g. \eqref{Ineg_Premiere}, %\eqref{Ineg_J1j_1} as well as
\eqref{eq_def_J6c} and \eqref{Ineg_J9}).

\begin{lem}\label{Lemma_Proba_bh_egal}
For $h>0$,
\begin{eqnarray}
    \forall x \geq 0,
\qquad
    \p(b_h= x)
& = &
    \frac{\p\big[\ell\big(\mathcal T_{V,h}^\downarrow\big)\geq x\big]}
    {\E\big[\ell\big(\mathcal T_{V,h}^\uparrow\big)+\ell\big(\mathcal T_{V,h}^\downarrow\big)\big]},
\label{eq_proba_bh_x_positif}
\\
    \forall x\leq 0,
\qquad
    \p(b_h= x)
& = &
    \frac{\p\big[\ell\big(\mathcal T_{V,h}^\uparrow\big)> -x\big]}
    {\E\big[\ell\big(\mathcal T_{V,h}^\uparrow\big)+\ell\big(\mathcal T_{V,h}^\downarrow\big)\big]}.
\label{eq_proba_bh_x_negatif}
\end{eqnarray}
\end{lem}

\noindent{\bf Proof:}
Let $h>0$ and $x\in\Z$.
If $x> 0$,
applying Theorem \ref{Lemma_Central_Slope} eq. \eqref{eq_Central_Slope_Downward},
%using \eqref{eq_Central_Slope_Upward}
with $A=\bigsqcup_{i=1}^\infty \R_-^i$, $\Delta_1=\{x\}$ and $\Delta_0=-\N$,
\begin{eqnarray*}
    \p(b_h= x)
& = &
    \p[b_h= x,\ x_1(V,h)=b_h>0]
=
    \p\bigg[x_1(V,h)=x,\ \theta(T_0(V,h))\in\bigsqcup_{i=1}^\infty \R_-^i\bigg]
\nonumber\\
& = &
    \frac{\E\big(\sharp\big\{0\leq i < \ell\big(\mathcal T_{V,h}^\downarrow\big),
        \ell\big(\mathcal T_{V,h}^\downarrow\big)-i=x\big\}\big)}
    {\E\big[\ell\big(\mathcal T_{V,h}^\uparrow\big)+\ell\big(\mathcal T_{V,h}^\downarrow\big)\big]}
=
    \frac{\p\big[\ell\big(\mathcal T_{V,h}^\downarrow\big)\geq x\big]}
    {\E\big[\ell\big(\mathcal T_{V,h}^\uparrow\big)+\ell\big(\mathcal T_{V,h}^\downarrow\big)\big]}.
\end{eqnarray*}
Similarly if $x\leq 0$, applying Theorem \ref{Lemma_Central_Slope} eq. \eqref{eq_Central_Slope_Upward}
with $A=\bigsqcup_{i=1}^\infty \R_+^i$, $\Delta_0=\{x\}$ and $\Delta_1=\N^*$,
\begin{equation*}
    \p(b_h= x)
%& = &
%    \p[b_h= x, x_0(V,h)=b_h\leq0]
%=
%    \p\bigg[x_0(V,h)=x,\ \theta(T_0(V,h))\in\bigsqcup_{i=1}^\infty \R_+^i\bigg]
%\nonumber\\
%& = &
=
    \frac{\E\big(\sharp\big\{0\leq i < \ell\big(\mathcal T_{V,h}^\uparrow\big), i=-x\big\}\big)}
    {\E\big[\ell\big(\mathcal T_{V,h}^\uparrow\big)+\ell\big(\mathcal T_{V,h}^\downarrow\big)\big]}
=
    \frac{\p\big[\ell\big(\mathcal T_{V,h}^\uparrow\big)> -x\big]}
    {\E\big[\ell\big(\mathcal T_{V,h}^\uparrow\big)+\ell\big(\mathcal T_{V,h}^\downarrow\big)\big]}.
\end{equation*}
In particular,
\begin{equation}
\label{eq_proba_bh_x_zero}
    \p(b_h= 0)
=
    \frac{\p\big[\ell\big(\mathcal T_{V,h}^\uparrow\big)> 0\big]}
    {\E\big[\ell\big(\mathcal T_{V,h}^\uparrow\big)+\ell\big(\mathcal T_{V,h}^\downarrow\big)\big]}
=
    \frac{1}
    {\E\big[\ell\big(\mathcal T_{V,h}^\uparrow\big)+\ell\big(\mathcal T_{V,h}^\downarrow\big)\big]}
=
    \frac{\p\big[\ell\big(\mathcal T_{V,h}^\downarrow\big)\geq 0\big]}
    {\E\big[\ell\big(\mathcal T_{V,h}^\uparrow\big)+\ell\big(\mathcal T_{V,h}^\downarrow\big)\big]},
\end{equation}
so both formulas of Lemma \ref{Lemma_Proba_bh_egal} are true for $x=0$.
%which is greater than or equal to any $\p(b_h=x)$, $x\in\Z$.
\hfill$\Box$

\subsection{About right $h$-extrema and right $h$-slopes}\label{Sub_sec_right_extrema}
We have detailed in the previous subsections, for $h>0$, a path decomposition of the potential $V$,
which we cut into different trajectories, called left $h$-slopes, between random times which are the left $h$-extrema.
We have also given the laws and independence properties of these left $h$-slopes,
in particular in Theorems \ref{Lemma_Law_of_Slopes}, \ref{Lemma_Independence_h_extrema} and \ref{Lemma_Central_Slope}.

We now focus on right $h$-extrema and provide a similar path decomposition of $V$ with right $h$-slopes and right $h$-extrema.
Similarly as for left $h$-minima,
%almost surely,
for $v\in\VVV$,
for every $h>0$,
the set of right $h$-extrema of $v$ can be denoted by $\{x_k^*(v,h), \ k\in\Z\}$, such that
$k\mapsto x_k^*(v,h)$ is strictly increasing
and $x_0^*(v,h)< 0\leq x_1^*(v,h)$ (see Figure \eqref{figure_h_extrema}),
the first inequality being strict and second one being large,
contrarily to inequalities for left $h$-extrema $x_i(v,h)$, $i\in\Z$,
in order to get relation \eqref{eq_relation_xi_xi*} below.
%and the left $h$-minima and left $h$ maxima of $V$ alternate.
%(we need to define the $x_i$ a bit differently so that this is true in the general case without the previous assumption).
Also, we prove below that the right $h$-extrema of $v$ can be obtained from the left $h$-extrema of
$v^-(.):=v_-(.):=v(-.)$
(and in particular, $V^-(.):=V_-(.):=V(-.)$; both notations $V^-$ and $V_-$ will be used throughout the paper, depending on which one is more convenient).
More precisely, we have:

\begin{lem} \label{Lemma_relation_xi_xi*}
Let $v\in\VVV$. For $h>0$,
\begin{equation}\label{eq_relation_xi_xi*}
    \forall i\in\Z,
\qquad
    x_i^*(v,h)
=
    -x_{1-i}(v^-,h).
\end{equation}
\end{lem}

\noindent{\bf Proof:}
Let $v\in\VVV$ and $h>0$.
First, notice that, applying Definition \ref{def_left_extrema},
$-x_j(v^-,h)$ is a right $h$-extremum for $v$ for each $j\in\Z$, so
$
    \{-x_j(v^-,h),\ j\in\Z\}
\subset
    \{x_i^*(v,h),\ i\in\Z\}
$.
Similarly, for $i\in\Z$, $-x_i^*(v,h)$ is a left $h$-extremum for $v^-$, so
$
    \{x_i^*(v,h),\ i\in\Z\}
\subset
    \{-x_j(v^-,h),\ j\in\Z\}
$, thus these two sets are equal.
Moreover,
$(x_i^*(v,h))_{i\in\Z}$
and
$(-x_{-j}(v^-,h))_{j\in\Z}$
are two strictly increasing sequences,
taking the same values, so there exists $k\in\Z$ such that
$x_i^*(v,h)=-x_{k-i}(v^-,h)$ for every $i\in\Z$.
Since $x_0^*(v,h)< 0\leq x_1^*(v,h)$ and
%$x_0(v,h)\leq 0<x_1(v,h)$,
%$-x_0(v,h)\geq 0>-x_1(v,h)$,
$-x_1(v^-,h)< 0 \leq -x_0(v^-,h)$,
we have $k=1$, which proves the lemma.
\hfill$\Box$

Let $h>0$.
Similarly as for left $h$-extrema, for $v\in\VVV$,
we introduce for each $i\in\Z$ the right {\it $h$-slope}
$T_i^*(v,h):=(v(j)-v[x_i^*(V,h)],\ x_i^*(v,h)\leq j\leq x_{i+1}^*(v,h))$.
If $x_i^*(v,h)$ is a right $h$-minimum (resp. maximum), then $\theta[T_i^*(v,h)]$ is strictly positive
(resp. strictly negative) on $\{1,\dots,\ell(T_i^*(v,h))\}$.
and its maximum (resp. minimum) is attained at $\ell(T_i^*(v,h))$.
The notation with a star for $x_i^*$ and $T_i^*$ corresponds to this fact that
the translated slopes $\theta[T_i^*(v,h)]$
are non-zero except at the origin.

Using the previous definition of $\tau_i^{(V)}(h)$ (see around \eqref{eqDef_tau_1_V}), we define for $i\geq 0$
(see Figure \ref{figure_tau_i_m_i_slopes}, in which
$m^{(V)*}_{3}(h)=m^{(V)}_{3}(h)$ and $m^{(V)*}_{4}(h)=m^{(V)}_{4}(h)$),
\begin{eqnarray}
%    \tau^{(V)}_{2i+1}(h)
%& := &
%    \min\Big\{k\geq \tau_{2i}^{(V)}(h),\ V(k)-\min\nolimits_{[\tau_{2i}^{(V)}(h),k]} V\geq h\Big\},
%\label{eqDef_tau_1_V}
%\\
    m^{(V)*}_{2i+1}(h)
& := &
    \max\Big\{k\in \big[\tau_{2i}^{(V)}(h), \tau_{2i+1}^{(V)}(h)\big]\cap \N,\ V(k)=\min\nolimits_{[\tau_{2i}^{(V)}(h),\tau_{2i+1}^{(V)}(h)]} V\Big\},
\nonumber\\
%    \tau^{(V)}_{2i+2}(h)
%& := &
%    \min\Big\{k\geq \tau_{2i+1}^{(V)}(h),\ \max\nolimits_{[\tau_{2i+1}^{(V)}(h),k]} V-V(k)\geq h\Big\},
%\nonumber\\
    m^{(V)*}_{2i+2}(h)
& := &
    \max\Big\{k\in \big[ \tau_{2i+1}^{(V)}(h), \tau_{2i+2}^{(V)}(h)\big]\cap\N,\ V(k)=\max\nolimits_{[\tau_{2i+1}^{(V)}(h),\tau_{2i+2}^{(V)}(h)]} V\Big\}.
\nonumber
\end{eqnarray}
Also, similarly as in Definition \ref{deff_law_slopes}, we introduce for $h>0$,
\begin{eqnarray}
\label{eq_def_slope_etoile_1}
    \mathcal T_{V,h}^{\uparrow*}
& := &
    \big(V\big[m_{1}^{(V)*}(h)+x\big]-V\big[m_{1}^{(V)*}(h)\big], \ 0\leq x \leq m_{2}^{(V)*}(h)-m_{1}^{(V)*}(h)\big),
\\
    \mathcal T_{V,h}^{\downarrow*}
& := &
    \big(V\big[m_{2}^{(V)*}(h)+x\big]-V\big[m_{2}^{(V)*}(h)\big], \ 0\leq x \leq m_{3}^{(V)*}(h)-m_{2}^{(V)*}(h)\big).
\label{eq_def_slope_etoile_2}
\end{eqnarray}

Recall $T_V$ and $T_V^*$ from \eqref{eq_def_TY} and \eqref{eq_def_TY*}.
The following proposition is similar to (\cite{DGP_Collision_Transient}, Proposition 5.2)
with $m^{(V)*}_{1}(h)$ instead of $m^{(V)}_{1}(h)$.
The other main difference is that in {\bf (ii)}, we condition by $\{T_V([h,+\infty[)<T_V^*(]-\infty,0])\}$, closed at $0$,
instead of $\{T_V([h,+\infty[)<T_V(]-\infty,0[)\}$.
Since we did not find this lemma in the literature (in which our stopping time $\tau_1^{(V)}(h)$
is generally replaced by a deterministic time, see e.g. \cite{Bertoin_Split}),
we give a detailed proof.
% for the sake of completeness.

\begin{prop}\label{Lemma_Loi_Autour_m1_Star}
Let $h>0$.
Let $V$ be a random walk given as in \eqref{eqDefPotentialV}
by a sequence of partial sums of i.i.d.
%non degenerate
r.v. $\log \rho_i$, $i\in\Z$, such that $\p[\log \rho_0>0]>0$ and $\p[\log \rho_0<0]>0$
(this result does not require Hypotheses \eqref{eqEllipticity}, \eqref{eqRecurrence} or \eqref{eq_def_sigma}).
If moreover $\liminf_{x\to+\infty}V(x)=-\infty$ a.s., then\\
\noindent{\bf (i)}
The processes
$\big(V\big[m^{(V)*}_{1}(h)-k\big]-V\big[m^{(V)*}_{1}(h)\big],\ 0\leq k \leq m^{(V)*}_{1}(h)\big)$
and
$\big(V\big[m^{(V)*}_{1}(h)+k\big]-V\big[m^{(V)*}_{1}(h)\big],\ 0\leq k \leq \tau_1^{(V)}(h)-m^{(V)*}_{1}(h)\big)$
are independent.\\
\noindent{\bf (ii)}
The process
$\big(V\big[m^{(V)*}_{1}(h)+k\big]-V\big[m^{(V)*}_{1}(h)\big],\ 0\leq k \leq \tau_1^{(V)}(h)-m^{(V)*}_{1}(h)\big)$
is equal in law to
$\big(V(k),\ 0\leq k\leq T_V([h,+\infty[) \big)$
conditioned on $\{T_V([h,+\infty[)<T_V^*(]-\infty,0])\}$.
%conditioned to stay nonnegative between times $0$ and $T_V([h,\infty[)$.
\end{prop}

\noindent{\bf Proof:}
We fix $h>0$, and consider $V$ satisfying the hypotheses.
Let $\psi_1$ and $\psi_2$ be two nonnegative functions,
$\bigsqcup_{t\in\mathbb N^*}\mathbb R^t\rightarrow \mathbb [0,+\infty[$,
measurable with respect to the $\sigma$-algebra $\{\bigsqcup_{t\in\mathbb N^*}A_t\, :\, \forall t\in\mathbb N^*,\ A_t\in\mathcal B(\mathbb R^t)\}$.
To simplify the notation, we set $m_1^*:=m^{(V)*}_{1}(h)$ and $\tau_1^*:=\tau_1^{(V)}(h)$.

We now define by induction, e.g. as in Enriquez et al. \cite{ESZ2} and \cite{ESZ3},
the weak descending ladder epochs for $V$ as
\begin{equation}\label{eq_def_ei}
    e_0
:=
    0,
\qquad
    e_{i}
:=
    \inf\{k>e_{i-1}\ :\ V(k)\leq V(e_{i-1})\},
\qquad
    i\geq 1,
\end{equation}
with $e_i<\infty$ a.s. for each $i\geq 1$ since $\liminf_{x\to+\infty}V(x)=-\infty$.
In particular, the excursions $(V(k+e_i)-V(e_i),\ 0\leq k\leq e_{i+1}-e_i)$, $i\geq 0$ are i.i.d.
by the Strong Markov property.
Also, the height $H_i$ of the excursion $[e_i,e_{i+1}]$ is defined as
\begin{equation}\label{eq_def_Hi}
    H_i
:=
    \max_{e_i\leq k\leq e_{i+1}}[V(k)-V(e_i)],
\qquad
    i\geq 0.
\end{equation}
Notice in particular that
$m_1^*=e_L$, where
$L:=\min\{\ell\geq 0,\ H_\ell\geq h\}<\infty$ a.s.
Hence, summing over the values of $L$, we get
\begin{align*}
&
    \E\left[\psi_1\big(V(m_1^*-k)-V(m_1^*),\ 0\leq k \leq m_1^*\big)
            \psi_2\big(V(m_1^*+k)-V(m_1^*),\ 0\leq k \leq \tau_1^*-m_1^*\big)\right]
\\
%& =
%    \sum_{\ell=0}^\infty
%    \E\big[\psi_1\big(V\big(e_\ell'-k\big)-V\big(e_\ell'\big),\ 0\leq k \leq e_\ell'\big)
%\\
%&
%    \qquad\qquad
%    \times
%    \psi_2\big(V\big(e_\ell'+k\big)-V\big(e_\ell'\big),\ 0\leq k \leq T{}^{\uparrow}-e_\ell'\big)\un_{\{L=\ell\}}\big]
%\\
& =
    \sum_{\ell=0}^\infty
    \E\big[\psi_1\big(V\big(e_\ell-k\big)-V\big(e_\ell\big),\ 0\leq k \leq e_\ell\big)\un_{\cap_{i=0}^{\ell-1}\{H_i<h\}}
    \un_{\{H_\ell\geq h\}}
\\
&
    \qquad\qquad
    \times
    \psi_2\big(V\big(e_\ell+k\big)-V\big(e_\ell\big),\ 0\leq k \leq T_{V(\cdot+e_\ell)-V(e_\ell)}([h,+\infty[)\big)\big]
\\
& =
    \Pi_1\Pi_2,
\end{align*}
due to the strong Markov property at stopping time $e_\ell$, where, applying it again on the second equality,
\begin{align*}
    \Pi_1
& :=
    \sum_{\ell=0}^\infty
    \E\big[\psi_1\big(V\big(e_\ell-k\big)-V\big(e_\ell\big),\ 0\leq k \leq e_\ell\big)\un_{\cap_{i=0}^{\ell-1}\{H_i<h\}}
    \big]
    \p[H_\ell\geq h]
\\
& =
    \sum_{\ell=0}^\infty
    \E\big[\psi_1\big(V\big(e_\ell-k\big)-V\big(e_\ell\big),\ 0\leq k \leq e_\ell\big)\un_{\{L=\ell\}}
    \big]
\\
& =
    \E\big[\psi_1\big(V(m_1^*-k)
    -V(m_1^*),\ 0\leq k \leq m_1^*\big)\big]
\end{align*}
and, since $\p[H_\ell\geq h]=\p[T_V([h,+\infty[)<T_V^*(]-\infty,0])]$,
\begin{align*}
    \Pi_2
& :=
    \E\big[\psi_2\big(V(k),\ 0\leq k \leq T_V([h,+\infty[))
     \mid T_V([h,+\infty[)<T_V^*(]-\infty,0])
     \big].
\end{align*}
Since this is true for all $\psi_1$ and $\psi_2$, this proves the proposition.
\hfill$\Box$

As a consequence, we get

\begin{thm}\label{Lemma_Law_of_Slopes_Right}
Assume \eqref{eqEllipticity}, \eqref{eqRecurrence} and \eqref{eq_def_sigma}. Let $h>0$.
\\
{\bf (i)} The process $\mathcal T_{V,h}^{\uparrow*}$
up to its first hitting time $T_{\mathcal T_{V,h}^{\uparrow*}}([h,+\infty[)$
of $[h,+\infty[$,
that is, $\big(\mathcal T_{V,h}^{\uparrow*}(k),$\ $0\leq k \leq T_{\mathcal T_{V,h}^{\uparrow*}}([h,+\infty[)\big)$,
is equal in law to
$
    \big(V(k), \  0\leq k \leq T_V([h,+\infty[)
    \big)
$
conditioned on $\{ T_V([h,+\infty[)< T_V^*(]-\infty, 0] )\}$.
Moreover, it is independent of
$
    \big(
    \mathcal T_{V,h}^{\uparrow*}\big(T_{\mathcal T_{V,h}^{\uparrow*}}([h,+\infty[)+k\big)
    -
    \mathcal T_{V,h}^{\uparrow*}\big(T_{\mathcal T_{V,h}^{\uparrow*}}([h,+\infty[)\big),\,
    0\leq k \leq \ell\big(\mathcal T_{V,h}^{\uparrow*}\big)-T_{\mathcal T_{V,h}^{\uparrow*}}([h,+\infty[)
    \big)
$,
which has the same law as
%$V$ between $0$ and
$\big(V(k),\, 0\leq k \leq \widetilde M^\sharp_h\big)$, with
$
    \widetilde M^\sharp_h
:=
    \max\{0\leq k\leq \tilde\tau_1(h), \, V(k)=\max_{[0, \tilde\tau_1(h)]}V\}$,
where
$
    \tilde\tau_1(h)
:=
    \min\{k\geq 0,\, \max_{[0,k]}V-V(k)\geq h\}
$.

{\bf (ii)}
$\mathcal T_{-V,h}^{\uparrow*}=_{law}-\mathcal T_{V,h}^{\downarrow*}$
and
$\mathcal T_{-V,h}^{\downarrow*}=_{law}-\mathcal T_{V,h}^{\uparrow*}$.

{\bf (iii)} Also,
$
    \E\big(\ell\big(\mathcal T_{V,h}^{\uparrow*}\big)\big)
<
    \infty
$
and
$
    \E\big(\ell\big(\mathcal T_{V,h}^{\downarrow*}\big)\big)
<
    \infty
$.
\end{thm}

\noindent{\bf Proof:}
The proof of this theorem is the same as the proof of Theorem \ref{Lemma_Law_of_Slopes},
with Proposition \ref{Lemma_Loi_Autour_m1_Star}, $\mathcal T_{V,h}^{\uparrow*}$, $m^{(V)*}_{i}(h)$
and right extrema instead of (\cite{DGP_Collision_Transient}, Proposition 5.2),
$\mathcal T_{V,h}^{\uparrow}$, $m^{(V)}_{i}(h)$ and  left extrema respectively.
\hfill$\Box$

The following lemma says that $\zeta$, defined in \eqref{eq_def_zeta},
transforms translated left (resp. right) $h$-slopes for $V$
into right (resp. left) ones for $V^-$ (see Lemma \ref{Lemma_zeta} below),
and  upward ones into downward ones.

\begin{lem}\label{Lemma_zeta}
For $i\in\Z$, $\zeta[\theta(T_i(V,h)))]=\theta[T_{-i}^*(V^-,h)]$.
\end{lem}

\noindent{\bf Proof:}
Recall that
$
    x_k^*(V,h)
=
    -x_{1-k}(V^-,h)
$
for $k\in\Z$ by Lemma \ref{Lemma_relation_xi_xi*}.
Hence for $i\in\Z$,
\begin{eqnarray}
&&
    \zeta[\theta(T_i(V,h)))]
=
    \zeta\big[\big(V[x_i(V,h)+j]-V[x_i(V,h)], \ 0\leq j\leq x_{i+1}(V,h)-x_i(V,h)\big)\big]
\nonumber\\
& = &
    (V^-[-x_{i+1}(V,h)+j] - V^-[-x_{i+1}(V,h)], \ 0\leq j\leq x_{i+1}(V,h)-x_i(V,h))
\nonumber\\
& = &
    \big(V^-[x_{-i}^*(V^-,h)+j] - V^-[x_{-i}^*(V^-,h)], \ 0\leq j\leq -x_{-i}^*(V^-,h))+x_{1-i}^*(V^-,h))\big)
\nonumber\\
& = &
    \theta[T_{-i}^*(V^-,h)].
\label{eq_egalite_slopes_droite_gauche}
\end{eqnarray}
This proves the lemma.

As a consequence, we get the following result.

\begin{thm}\label{Theorem_Etoile}
Theorems \ref{Lemma_Independence_h_extrema} and \ref{Lemma_Central_Slope} remain valid if we replace
"left" and each $x_k(V,h)$, $T_k(V,h)$,  $\mathcal T_{V,h}^{\uparrow}$ and $\mathcal T_{V,h}^{\downarrow}$
respectively by
"right", $x_k^*(V,h)$, $T_k^*(V,h)$,  $\mathcal T_{V,h}^{\uparrow*}$ and $\mathcal T_{V,h}^{\downarrow*}$,
and $<$ and $\leq$ respectively by $\leq$ and $<$ in Theorem \ref{Lemma_Central_Slope}.
\end{thm}

\noindent{\bf Proof:}
Indeed, their proofs remain valid
if we make these replacements and also replace
$m_k^{(V)}(h)$ by $m_k^{(V)*}(h)$, $k\in\Z$,
and $T_j(t,V,h)$ by $T_j^*(t,V,h)=T_{i+j}^*(V,h)$ if $x_i^*(V,h)<t\leq x_{i+1}^*(V,h)$,
and for this unique $i$,
$x_j^*(t, V, h):=x_{i+j}^*(V,h)$ for $j\in\Z$,
and as a consequence, replace $<$ and $\leq$ respectively by $\leq$ and $<$ throughout the proof.
\hfill$\Box$

%We can now prove the following proposition.
The following proposition, combined with some other results such as Theorem \ref{Lemma_Law_of_Slopes_Right}, will be useful to
obtain the law of $V$ on the left of $x_i(V,h)$ (for $i\in\Z$) conditionally on $b_h\leq 0$ or $b_h>0$,
in view of Theorems \ref{Lemma_Independence_h_extrema} and \ref{Lemma_Central_Slope}.

\begin{prop}\label{Prop_Egalite_Loi_Zeta_Slopes}
Let $h>0$. Then,
{\bf (i)}
$\zeta\big(\mathcal T_{V,h}^{\uparrow}\big)=_{law} \mathcal T_{V^-,h}^{\downarrow*}$
and
{\bf (ii)}
$\zeta\big(\mathcal T_{V,h}^{\downarrow}\big)=_{law} \mathcal T_{V^-,h}^{\uparrow*}$.
\end{prop}

\noindent{\bf Proof:}
% of Proposition \ref{Prop_Egalite_Loi_Zeta_Slopes}:}
We denote by $\mathscr{L}\big(\mathcal T_{V^-,h}^{\downarrow*}\big)$
the law of $\mathcal T_{V^-,h}^{\downarrow*}$
under $\p$.
Conditionally on
$
    \{V[x_0(V,h)]<V[x_1(V,h)]\}
%=
%    \{V^-[-x_0(V,h)]<V^-[-x_1(V,h)]\}
=
    \{V^-[x_1^*(V^-,h)]<V^-[x_0^*(V^-,h)]\}
$ (thanks to \eqref{eq_relation_xi_xi*}),
$\theta[T_2(V,h)]$ has the law $\mathcal{L}\big(\mathcal T_{V,h}^{\uparrow}\big)$
by Theorem \ref{Lemma_Independence_h_extrema} {\bf (i)},
whereas $\theta[T_{-2}^*(V^-,h)]$ has
the law $\mathscr{L}\big(\mathcal T_{V^-,h}^{\downarrow*}\big)$
%the same law as $\mathcal T_{V^-,h}^{\downarrow*}$
by the version of Theorem \ref{Lemma_Independence_h_extrema} {\bf (ii)}
with stars (see Theorem \ref{Theorem_Etoile}) applied to $V^-$.
This and \eqref{eq_egalite_slopes_droite_gauche} prove our {\bf (i)}.
Applying the same arguments to $\theta[T_1(V,h)]$ and $\theta[T_{-1}^*(V^-,h)]$ proves {\bf (ii)}.
\hfill$\Box$

%%%%%%%%%%%%%%%%%%%%%%%%%%%%%%%%%%%%%%%%%%%%%%%%%%%%%%%%%%%%%%%%%%%

%%%%%%%%%%%%%%%%%%%%%%%%%%%%%%%%%%%%%%%%%%%%%%%%%%%%%%%%%%%%%%%%%%%%%%%%%%%%%%%%%%%%%%%%

%%%%%%%%%%%%%%%%%%%%%%%%%%%%%%%%%%%%%%%%%%%%%%%%%%%%%%%%%%%%%%%%%%%%%%%%%%%%%%%%%%%%%%%%

\subsection{Relation with another localization point}\label{Subsection_Kestn_Loc_Def}
In this subsection, we recall another way to define a localization point denoted by $b_h^{(K)}$,
%so that $(S_n)_n$ will be localized close to it at time $\lfloor e^h\rfloor$ with large probability,
and we prove that $b_h^{(K]}$ is equal to $b_h$ (defined in \eqref{eqDefbh}) with large probability.
The localization point $b_h^{(K)}$ is useful because we will apply
the previous result of Kesten (\cite{Kesten}, Thm 1.2)
to the limit law of $b_h^{(K)}/h^2$
(in the proof of Theorem \ref{Th_Local_Limit_b_h}, see after \eqref{Ineg_THé_Uniforme_Loin}),
whereas our $b_h$ is convenient e.g. due to Lemma \ref{Lemma_Proba_bh_egal}
and to the law of the potential near $b_h$ (by Theorems \ref{Lemma_Law_of_Slopes}, \ref{Lemma_Independence_h_extrema},
and \ref{Lemma_Central_Slope}).

%{\bf (homogeneiser notations $V^-$ et $V_-$)}

\begin{figure}[htbp]
\includegraphics[width=16.0cm,height=6.22cm]{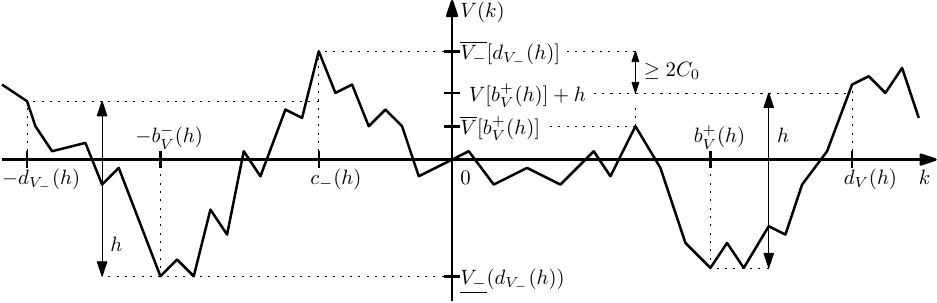}
\caption{Schema of the potential $V$ for the first case of the proof of Lemma \ref{Lemma_Comparaison_b_bK}
when $\overline{V}[b_V^+(h)]<V[b_V^+(h)] + h$.}
\label{figure_preuve_proba_b_egal_bK}
\end{figure}

To this aim, we define for any process $(Z(k),\ k\geq 0)$,
similarly as in Hu (\cite{HuLocal} from eq. (2.1) to eq. (2.6))
but for processes indexed by $\N$,
$$
    \overline{Z}(t):=\sup_{0\leq k\leq t} Z(k),
\qquad
    \underline{Z}(t):=\inf_{0\leq k\leq t} Z(k),
\qquad
    Z^\sharp(t):=\sup_{0\leq s \leq t}\big(Z(s)-\underline{Z}(s)\big)
\qquad
    t\geq 0,
$$
\begin{equation}\label{eq_def_dV_dZ}
    d_Z(h)
:=
    \inf\{t\geq 0, \ Z^\sharp(t)\geq h\},
\qquad
    h> 0.
\end{equation}
%\begin{equation}
%    b_Z(h)
%:=
%    \inf\big\{0\leq u\leq d_Z(h), \ Z(u)=\underline{Z}(d_Z(h))\big\},
%\qquad
%    h> 0,
%\label{eq_def_bV-}
%\end{equation}
Also, with $V_-(k):=V(-k)$ for $k\geq 0$ as before, we introduce (see Figure \ref{figure_preuve_proba_b_egal_bK})
\begin{eqnarray}
    b_V^+(h)
& := &
    \inf\{0\leq u\leq d_V(h), \ V(u)=\underline{V}(d_V(h))\},
\qquad
    h> 0,
\nonumber\\
    b_V^-(h)
& := &
%    \sup\big\{-d_{V_-}(h)\leq u\leq 0, \ V(u)=\underline{V_-}(d_{V_-}(h))\big\},
    \sup\{0\leq u\leq d_{V_-}(h), \ V_-(u)=\underline{V_-}(d_{V_-}(h))\},
\qquad
    h> 0.
\label{eq_def_bV-}
\end{eqnarray}
The $\sup$ instead of $\inf$ in the last line will be necessary so that in some cases,
$-b_V^-(h)$ is a left $h$-minimum for $V$ instead of a right one
(as in Figure \ref{figure_preuve_proba_b_egal_bK}).
%Notice that these last two definitions
%are symmetric.
% if the last $\inf$ was replaced by a $\sup$.
Finally, we introduce
\begin{equation}\label{eqDef_bK_h}
    b_h^{(K)}
:=
    \left\{
    \begin{array}{ll}
        b_V^+(h)
    &
        \textnormal{if } \overline{V}[d_V(h)] < \overline{V_-}[d_{V_-}(h)],
    \\
%        b_V^-(h)
        -b_V^-(h)
    &
        \textnormal{otherwise}.
    \end{array}
    \right.
\end{equation}
Let $(W(x),\, x\in\R)$ be a two-sided Brownian motion, and $W_-:=(W(-x), \, x\geq 0)$.
As in Hu (\cite{HuLocal} eq. (2.6)), for $w=W$ or $w=\sigma W$,
we define $b_h^{(K,w)}$
%{\bf (notation confuse sans le $K$ ?)}
by the same formula as in \eqref{eqDef_bK_h},
the previous notations of this Subsection \ref{Subsection_Kestn_Loc_Def} being the same, with
$V$ replaced by $w$, and the $\inf$ and $\sup$ being taken for real numbers instead of integers.
As already stated by Hu \big(\cite{HuLocal} after eq. (2.6), his $b(1)$ being a.s. equal to our $b_1^{(K,W)}$
since the $\sup$ in \eqref{eq_def_bV-} is a.s. a $\inf$ when $V$ is replaced by $W$\big),
the density of $b_1^{(K,W)}$ is
$\varphi_{\infty}$, defined in \eqref{eq_def_phi_infini}.
Indeed, it is easy to check that $b_1^{(K,W)}$ is a.s. equal to the r.v. $L$ of
Kesten (\cite{Kesten}, as expressed in the statement of his Lemma 2.1),
which has density $\varphi_{\infty}$ by (\cite{Kesten}, Thm 1.2).

For some choices of $\p$, we have $\p\big[b_h\neq b_h^{(K)}\big]>0$ for some $h>0$.
Indeed, for example, if $\p[V(1)=z]>0$ for every $z\in\{-2,-1,0,1,2\}$, we have for $h\in\N^*$, with non zero probability,
$V(-1)=V(0)=V(1)=0$, with $V(k)=k-1$ for $1\leq k \leq h+1$,
$V(k)=|k|-1$ for $-h\leq k \leq -1$ and $V(-h-1)=h+1$, and so $b_h^{(K)}=b_V^+(h)=0$ whereas
$b_h=-1\neq b_h^{(K)}$.
However, we prove that $b_h^{(K)}=b_h$ with large probability. More precisely, we have:

\begin{lem}\label{Lemma_Comparaison_b_bK}
There exists a constant $c_5>0$ such that, for large $h$,
$$
    \p\big[b_h^{(K)}\neq b_h\big]
\leq
    c_5 h^{-1}.
$$
\end{lem}

This lemma will be useful to prove Lemma \ref{Lemma_Proba_bh_positif} and Theorem \ref{Th_Local_Limit_b_h}.
Moreover, we think it will also be necessary in a work in progress \cite{Devulder_Rates_CV}.

\noindent{\bf Proof of Lemma \ref{Lemma_Comparaison_b_bK}:}
Let $h>0$.
{\bf First case:} we assume that
%In these two cases, i.e. \eqref{eq_Hypothese_Cas1} and \eqref{eq_Hypothese_Cas2}, we have
\begin{equation}\label{eq_Hypothese_Cas1_et_2}
    \max\Big[
    %\sup_{[0, b_V^+(h)]} V,
    \overline{V}[b_V^+(h)],
    V[b_V^+(h)] + h
    \Big]
\leq
    \overline{V_-}[d_{V_-}(h)] -2C_0.
\end{equation}
Let $c_-(h):=\sup\{k\leq 0,\  V(k)=\overline{V_-}[d_{V_-}(h)]\}$ (which may be $-d_{V_-}(h)$ or not).
First, by definition of $b_V^+(h)$, we have $V[b_V^+(h)]=\min_{[b_V^+(h), d_V(h)]}V$.
Also for the same reason,
$V[b_V^+(h)]<\min_{[0,b_V^+(h)-1]}V$,
with $\min\emptyset=+\infty$ by convention,
and since $-d_{V_-}(h)\leq c_-(h)\leq 0$, we have
$
    \min_{[c_-(h), 0]} V
\geq
    V(c_-(h))-h -C_0
=
    \overline{V_-}[d_{V_-}(h)] -h -C_0
>
    V[b_V^+(h)]
$
first by definition of $d_{V_-}(h)$ and ellipticity,
then by definition of $c_-(h)$ followed by \eqref{eq_Hypothese_Cas1_et_2}.
So,
$
    \min_{[c_-(h), b_V^+(h)-1]} V
>
    V[b_V^+(h)]
$.
Moreover, by definition,
$
    V[d_V(h)]
\geq
%    \overline{V}[d_V(h)]+h
%=
    V[b_V^+(h)]+h
$.
Finally,
$
    V[c_-(h)]
=
    \overline{V_-}[d_{V_-}(h)]
\geq
    V[b_V^+(h)]+h
$
first by definition, then by \eqref{eq_Hypothese_Cas1_et_2}.
Consequently, $b_V^+(h)$ is a left $h$-minimum.

Assume that there exists a left $h$-extremum in $[0, b_V^+(h)-1]$.
Since $b_V^+(h)$ is a left $h$-minimum, and left $h$-maxima and minima for $V$
alternate by Lemma \ref{Lemma_Alternate},
there would be at least one left $h$-maximum in this interval, which we denote by $\alpha\in[0,b_V^+(h)[$.
Now, denote by $\gamma$ the largest left $h$-minimum such that $\gamma<\alpha$,
so that $[\gamma,\alpha]$ is (the support of) an upward left $h$-slope of $V$.
In particular,
\begin{equation}\label{Ineg_alpha_gamma}
    V(\gamma)\leq V(\alpha)-h
\quad
    \text{and}
\quad
    \gamma<\alpha.
\end{equation}
Assume that $0\leq \gamma$.
So, $\gamma\in[0,\alpha]$, hence $V(\gamma)\geq \inf_{[0, \alpha]} V=\underline{V}(\alpha)$, and then using \eqref{Ineg_alpha_gamma},
$$
    V^\sharp(\alpha)
\geq
    V(\alpha)-\underline{V}(\alpha)
    %\inf_{[0, \alpha]} V
\geq
    V(\alpha)-V(\gamma)
\geq
    h.
$$
By definition \eqref{eq_def_dV_dZ} of $d_V$, this would give
$
    d_V(h)
\leq
    \alpha
$,
which contradicts
$
    \alpha
<
    b_V^+(h)
\leq
    d_V(h)
$.

Hence we would have $\gamma<0\leq \alpha$.
%As before,  let $c_-(h):=\sup\{k\leq 0,\  V(k)=\overline{V_-}[d_{V_-}(h)]\}$.
Using first the fact that $[\gamma,\alpha]$ is an upward left $h$-slope,
then $\alpha\in[0, b_V^+(h)]$ and finally \eqref{eq_Hypothese_Cas1_et_2}
would give
$$
    \sup\nolimits_{[\gamma,\alpha]} V
\leq
    V(\alpha)
\leq
    \overline{V}[b_V^+(h)]
<
    \overline{V_-}[d_{V_-}(h)]
=
    V[c_-(h)].
$$
So $c_-(h)\notin [\gamma,\alpha]$ and since $c_-(h)\leq 0\leq \alpha$ by definition,
this gives $c_-(h)<\gamma<0$.
Using ellipticity \eqref{eq_ellipticity_for_V}, then \eqref{eq_Hypothese_Cas1_et_2}, we get
$$
    V[c_-(h)+1]
    %-V(\gamma)
\geq
    V[c_-(h)]-C_0
    %-V(\gamma)
=
    \overline{V_-}[d_{V_-}(h)]-C_0
    %-V(\gamma)
\geq
    \overline{V}[b_V^+(h)]+C_0.
$$
Thus, using \eqref{Ineg_alpha_gamma} in the second inequality,
$$
    V[c_-(h)+1]-V(\gamma)
\geq
%    \sup_{[0, b_V^+(h)]} V
    \overline{V}[b_V^+(h)]
    +C_0-V(\gamma)
\geq
%    \sup_{[0, b_V^+(h)]} V
    \overline{V}[b_V^+(h)]
    +C_0+h-V(\alpha)
\geq
    C_0+h,
$$
since
$
%\sup_{[0, b_V^+(h)]} V
    \overline{V}[b_V^+(h)]
\geq V(\alpha)$ because $\alpha\in[0, b_V^+(h)]$.
So, $V_-(|c_-(h)|-1)-V_-(|\gamma|)>h$
with $0<|\gamma|\leq |c_-(h)|-1$, which gives
$
    d_{V_-}(h)
\leq
    |c_-(h)|-1
<
    |c_-(h)|
\leq
    d_{V_-}(h)
$,
which is not possible.

Hence there is no left $h$-extremum in $[0, b_V^+(h)-1]$.
Since $b_V^+(h)$ is a left $h$-minimum, this gives
$x_1(V,h)=b_V^+(h)$ if $b_V^+(h)\neq 0$
and $x_0(V,h)=b_V^+(h)$ if $b_V^+(h)= 0$,
and by definition \eqref{eqDefbh} of $b_h$, it follows that
$b_h=b_V^+(h)$.
Since
$
    \overline{V}[d_V(h)]
\leq
    \max\big[
    \overline{V}[b_V^+(h)],
    \max_{[b_V^+(h), d_V(h)]} V
    \big]
\leq
    \max\Big[
    \overline{V}[b_V^+(h)],
%    \sup_{[0, b_V^+(h)]} V,
    V[b_V^+(h)] + h+C_0
    \Big]
<
    \overline{V_-}[d_{V_-}(h)]
$
by ellipticity and \eqref{eq_Hypothese_Cas1_et_2}, we also have
$b_h^{(K)}= b_V^+(h)$ by \eqref{eqDef_bK_h}.
Hence, $b_h=b_h^{(K)}$
when \eqref{eq_Hypothese_Cas1_et_2} holds.

{\bf Second case:} we assume that
\begin{equation}\label{eq_Hypothese_Cas3_4}
    \max\Big[
%    \sup_{[b_V^-(h),0]} V,
    \overline{V_-}[b_V^-(h)],
    V_-[b_V^-(h)] + h
    \Big]
\leq
    \overline{V}[d_{V}(h)] -2C_0.
\end{equation}
This case is nearly the symmetric of the previous one, the only asymmetry being the $\sup$ in \eqref{eq_def_bV-}
(which is necessary for $-b_V^-(h)$ to be a left $h$-minimum instead of a right one).
So we prove similarly as in the first case that
$b_h=-b_V^-(h)=b_h^{(K)}$
when \eqref{eq_Hypothese_Cas3_4} holds.

{\bf Third step:}
Consequently, if $b_h\neq b_h^{(K)}$ then neither \eqref{eq_Hypothese_Cas1_et_2}
nor \eqref{eq_Hypothese_Cas3_4} hold, and so
\begin{eqnarray*}
    \overline{V_-}[d_{V_-}(h)] -2C_0
& < &
    \max\left[
%    \sup_{[0, b_V^+(h)]} V,
    \overline{V}[b_V^+(h)],
    V[b_V^+(h)] + h
    \right]
\leq
    \overline{V}[d_{V}(h)]
\\
& < &
    \max\left[
%    \sup_{[b_V^-(h),0]} V,
    \overline{V_-}[b_V^-(h)],
    V_-[b_V^-(h)] + h
    \right]
    +2C_0
\leq
    \overline{V_-}[d_{V_-}(h)] +2C_0,
\end{eqnarray*}
where we first used the negation of \eqref{eq_Hypothese_Cas1_et_2},
then the definitions of $d_{V}(h)$ and $b_V^+(h)$,
then the negation of \eqref{eq_Hypothese_Cas3_4}
and finally the definitions of $d_{V_-}(h)$ and $b_V^-(h)$.
In view of these inequalities, we define
\begin{eqnarray*}
    E_1
& := &
    \left\{-2C_0<
        \max\left[ \overline{V}[b_V^+(h)], V[b_V^+(h)] + h\right]
        -
        \overline{V_-}[d_{V_-}(h)]
        <2C_0
    \right\},
\\
    E_2
& := &
    \big\{
        V[b_V^+(h)] + h
        <
        \overline{V}[b_V^+(h)]
    \big\},
\end{eqnarray*}
so that $\p\big[b_h\neq b_h^{(K)}\big]\leq \p[E_1]$.

First, notice that on $E_1\cap E_2$,
writing here $\beta:=\overline{V_-}[d_{V_-}(h)]$ to simplify the notation,
we have $\beta-2C_0<\overline{V}[b_V^+(h)]<\beta+2C_0$, and so
$V[b_V^+(h)]<\beta+2C_0-h$
thanks to $E_2$.
Hence,
$T_V([\beta-2C_0,+\infty[) \leq b_V^+(h)$
and
$V[.+T_V([\beta-2C_0,+\infty[)]$ hits $V[b_V^+(h)]\in]-\infty, \beta+2C_0-h]$
before $[\beta+2C_0,+\infty[$.
Thus, since $V_-$ is independent of $(V(x),\, x\geq 0)$,
the strong Markov property, and then \eqref{eqOptimalStopping2} lead to, if $h>4C_0$,
\begin{eqnarray*}
&&
    \p[E_1\cap E_2\mid V_-]
\\
& \leq &
    \E\big(
        \p^{T_V([y-2C_0,+\infty[)}
            [
                T_V(]-\infty, y+2C_0-h])
                <
                T_V([y+2C_0,+\infty[)
            ]_{|y=\beta}
    \mid V_-\big)
\\
& \leq &
    5C_0(h+C_0)^{-1}.
\end{eqnarray*}
Consequently,
$
    \p[E_1\cap E_2]
\leq
    6C_0 h^{-1}
$
for large $h$.

Similarly, notice that on $E_1\cap E_2^c$, once more with the notation $\beta:=\overline{V_-}[d_{V_-}(h)]$,
we have $\beta-2C_0<V[b_V^+(h)]+h<\beta+2C_0$.
So, $T_V(]-\infty,\beta+2C_0-h])\leq b_V^+(h)$.
Also, $\min_{[0, d_V(h)]}V=V[b_V^+(h)]> \beta-2C_0-h$ and $V[d_V(h)]\geq V[b_V^+(h)]+h>\beta-2C_0$,
thus $V[.+T_V(]-\infty,\beta+2C_0-h])$ hits
$[\beta-2C_0,+\infty[$ before $]-\infty, \beta-2C_0-h]$. Hence as previously,
since $V_-$ is independent of $(V(x),\, x\geq 0)$,
by the strong Markov property, and then by \eqref{eqOptimalStopping2}, if $h>4C_0$,
\begin{eqnarray*}
&&
    \p[E_1\cap E_2^c\mid V_-]
\\
& \leq &
    \E\big(
        \p^{T_V(]-\infty, y+2C_0-h])}
            [
                T_V([y-2C_0,+\infty[)
                <
                T_V(]-\infty, y-2C_0-h])
            ]_{y=\beta}
    \mid V_-\big)
\\
& \leq &
    5C_0(h+C_0)^{-1}.
\end{eqnarray*}
Consequently,
$
    \p[E_1\cap E_2^c]
\leq
    6C_0 h^{-1}
$
for large $h$.
Finally,
$\p\big[b_h\neq b_h^{(K)}\big]\leq \p[E_1]\leq 12C_0 h^{-1}$ for large $h$, which proves the lemma.
\hfill$\Box$

\begin{lem}\label{Lemma_Proba_bh_positif}
There exists a constant $c_6>0$ such that
$$
    \p\big[b_h>0\big]\to_{h\to+\infty}1/2,
\qquad
    \p\big[b_h=0\big]\sim_{h\to+\infty}      c_6 h^{-2}.
$$
\end{lem}

\noindent{\bf Proof:}
For the equivalent, observe that by \eqref{eqDefbh}, $b_h=0$ if and only if
$0$ is a left $h$-minimum for $V$, that is if and only if
$V$ and $V(-.)=:V_-(.)$ hit $[h,+\infty[$ before going back to $]-\infty,0]$ for $V_-$,
and
before hitting $]-\infty,0[$ for $(V(k), \ k\geq 0)$.
So by independence of $(V(k), \, k\geq 0)$ and $V_-$ and \eqref{eq_Proba_Atteinte_logn_avant0} (or \eqref{eq_Proba_An}),
\begin{align}
    \p\big[b_h=0\big]
& =
    \p\big[T_{V_-}([h,+\infty[)<T_{V_-}^*(]-\infty,0])\big]
    \p\big[T_V([h,+\infty[)<T_V(]-\infty,0[)\big]
\nonumber\\
& \sim_{h\to+\infty}      c_6 h^{-2}
\label{equivalent_Proba_b_nul}
\end{align}
with $c_6>0$ being the product of
%$c$
$c_1^*$
(for the law of $V_-$)  and of
%$c^*$
$c_1$
(for the law of $V$)
with the notation of \eqref{eq_Proba_Atteinte_logn_avant0} (and \eqref{eq_Proba_An}).
This proves the second claim in Lemma \ref{Lemma_Proba_bh_positif}.
Notice that this constant $c_6$ depends on the law of $\omega_0$, that is, $c_6$ depends on $\p$.

For the first limit of the lemma,  notice that $\p[b_h>0]=\p\big[b_h^{(K)}>0\big]+O(1/h)$
as $h\to+\infty$
by Lemma \ref{Lemma_Comparaison_b_bK}, so we just have to prove that $\p\big[b_h^{(K)}>0\big]\to_{h\to+\infty}1/2$.
We now consider a two sided Brownian motion $(W(x),\ x\in\R)$, and consider
$W_-(x):=W(-x)$ for $x\geq 0$, and define
$\overline{W}$, $\overline{W_-}$, $d_W$, $d_{W_-}$, as explained after \eqref{eqDef_bK_h}.
%similarly as before \eqref{eq_def_bV-}.
By \eqref{eqDef_bK_h}, we have for $h>0$,
\begin{eqnarray}
    \p\big[b_h^{(K)}>0\big]
& = &
    \p\big[\overline{V}[d_V(h)] < \overline{V_-}[d_{V_-}(h)],\ b_V^+(h)\neq 0\big]
\nonumber\\
& = &
    \p\big[\overline{V}[d_V(h)] < \overline{V_-}[d_{V_-}(h)]\big]+O(1/h)
\label{eq_Proba_Pour_Donsker}
\end{eqnarray}
since
$\p\big[b_V^+(h)= 0\big]=\p\big[T_V([h,+\infty[)<T_V(]-\infty,0[)\big]=O(1/h)$
as $h\to+\infty$ similarly as in \eqref{equivalent_Proba_b_nul}.
By the theorem of Donsker, the limit of the probability in \eqref{eq_Proba_Pour_Donsker}
as $h\to+\infty$ is
$\p\big[\overline{\sigma W}[d_{\sigma W}(1)] < \overline{\sigma W_-}[d_{\sigma W_-}(1)]\big]$,
which is $1/2$ by symmetry
and because
$\p\big[\overline{\sigma W}[d_{\sigma W}(1)] = \overline{\sigma W_-}[d_{\sigma W_-}(1)]\big]=0$
since the r.v. $\overline{\sigma W}[d_{\sigma W}(1)]$ and $\overline{\sigma W_-}[d_{\sigma W_-}(1)]$ are independent
and have a density (by Hu \cite{HuLocal} Lemma 2.1 and by scaling).
Hence
$\p\big[b_h^{(K)}>0\big]\to_{h\to+\infty}1/2$
and so
$\p\big[b_h>0\big]\to_{h\to+\infty}1/2$.
\hfill$\Box$

\begin{lem}\label{Lemma_Esperance_Longueur_Slope}
There exists a constant $c_7:=(2c_6)^{-1}>0$ such that
$$
    \E[\ell(\mathcal T_{V,h}^\uparrow)]
\sim_{h\to+\infty}
    \E[\ell(\mathcal T_{V,h}^\downarrow)]
\sim_{h\to+\infty}
    c_7 h^2.
$$
\end{lem}

\noindent{\bf Proof:}
Applying \eqref{eq_proba_bh_x_zero}, and Theorem \ref{Lemma_Central_Slope}, using \eqref{eq_Central_Slope_Upward}
with $A=\SSS_+=\bigsqcup_{t=1}^\infty \R_+^t$,
 $\Delta_1=\N^*$ and
%first with
$\Delta_0=-\N$,
% and then
%with $\Delta_0=\{0\}$,
we have, since $b_h\leq 0$ if and only if $\theta(T_0(V,h))\in \SSS_+$ by \eqref{eqDefbh},
$$
    \p(b_h= 0)
=
    \frac{1}
    {\E[\ell(\mathcal T_{V,h}^\uparrow)+\ell(\mathcal T_{V,h}^\downarrow)]},
\qquad
    \p(b_h\leq 0)
=
    \frac{ \E\big[\ell(\mathcal T_{V,h}^\uparrow)\big] }
    {\E[\ell(\mathcal T_{V,h}^\uparrow)+\ell(\mathcal T_{V,h}^\downarrow)]}.
$$
Consequently,
$
    \E\big[\ell(\mathcal T_{V,h}^\uparrow)\big]
=
    \frac{\p(b_h\leq 0)}{\p(b_h= 0)}
\sim_{h\to+\infty}
    h^2/(2c_6)
%    \frac{h^2}{2c_1(V)c_1(V_-)}
$ by Lemma \ref{Lemma_Proba_bh_positif}.
Similarly, we obtain
$
    \E\big[\ell(\mathcal T_{V,h}^\downarrow)\big]
=
    \frac{\p(b_h> 0)}{\p(b_h= 0)}
\sim_{h\to+\infty}
    h^2/(2c_6)
%    \frac{h^2}{2c_1(V)c_1(V_-)}
$,
which proves the lemma.
\hfill$\Box$

%%%%%%%%%%%%%%%%%%%%%%%%%%%%%%%%%%%%%%%%%%%%%%%%%%%%%%%%%%%%%%%%%%%

\subsection{An inequality for the excess height of left $h$-slopes}

%The following lemma is not optimal but is enough for our needs.

\begin{lem}\label{Lemma_Loi_excess_height_h_extrema}
There exists a constant $c_8>0$ such that, for large $h$,
\begin{equation}\label{eqLemmaHeightSlopes}
%    \forall z>0,\,
%\forall \Delta \geq C_0,\,
%\forall h>C_0,\,
\forall i\in\Z,\,
\forall C_0<\Delta<h,
\qquad
    \p\Big(e[T_{i}(V, h)]\leq \Delta |
    %V(x_1(V,h))>V(x_0(V,h))
    b_h\leq 0
    \Big)
%    \p\Big(e[T_i(V, h)]\leq \Delta\Big)
\leq
    c_8 \frac{\Delta}{h}.
\end{equation}
This remains true if $b_h\leq 0$ is replaced by $b_h>0$.
\end{lem}

\noindent{\bf Proof:}
Let $h>0$ and $C_0<\Delta<h$.
%We first assume that $i\neq 0$.
Applying Theorem \ref{Lemma_Independence_h_extrema} {\bf (i)}
since
$
    \{V(x_1(V,h))>V(x_0(V,h))\}
=
    \{b_h\leq 0\}
$,
then Theorem \ref{Lemma_Law_of_Slopes} {\bf (i)}, and then \eqref{eqOptimalStopping2},
we have for $i\neq 0$,
since $C_0<\Delta<h$,
\begin{align}
&
    \p\big(e[T_{2i}(V, h)]\leq \Delta | b_h\leq 0\big)
=
    \p\Big(H[\theta(T_{2i}(V, h))]-h\leq \Delta | b_h\leq 0\Big)
 =
    \p\Big(H\big(\mathcal T_{V,h}^{\uparrow}\big)-h\leq \Delta\Big)
\nonumber\\
&
\leq
    \p\Big( T_V(-h+\Delta)\leq \tilde\tau_1(h) < T_V(]\Delta, +\infty[)\Big)
\leq
    \frac{\Delta+C_0}{h+C_0}
\leq
    \frac{2\Delta}{h}.
\label{ineg_Proba_Excess_1}
\end{align}
Similarly,
applying Theorem \ref{Lemma_Independence_h_extrema} {\bf (i)},
then $\mathcal T_{V,h}^{\downarrow}=_{law}-\mathcal T_{-V,h}^{\uparrow}$
by Theorem \ref{Lemma_Law_of_Slopes} {\bf (ii)},
%then Theorem \ref{Lemma_Law_of_Slopes} {\bf (i)},
%and then \eqref{eqOptimalStopping2},
\begin{equation*}
    \p\Big(e[T_{2i+1}(V, h)]\leq \Delta | b_h\leq 0\Big)
=
    \p\Big(H\big(\mathcal T_{V,h}^{\downarrow}\big)-h\leq \Delta\Big)
=
    \p\Big(H\big(\mathcal T_{-V,h}^{\uparrow}\big)-h\leq \Delta\Big)
%\\
%& \leq &
%    \p\big( T_{-V}(-h+\Delta) < T_{-V}(]\Delta,+\infty[)\big)
%\leq
%    \frac{\Delta+C_0}{h+C_0}
\leq
    \frac{2\Delta}{h}
\end{equation*}
similarly as before for $i\in\Z$ and $C_0<\Delta<h$.
This proves \eqref{eqLemmaHeightSlopes} for $i\neq 0$.

The proof is similar when conditioning by $b_h>0$,
applying Theorem \ref{Lemma_Independence_h_extrema} {\bf (ii)} instead of {\bf (i)}.

We now consider the case $i=0$.
We have, by Theorem \ref{Lemma_Central_Slope} eq. \eqref{eq_Central_Slope_Upward} applied with $\Delta_0=\Delta_1=\Z$,
\begin{eqnarray}
&&
    \p\Big(e[T_{0}(V, h)]\leq \Delta | V(x_1(V,h))>V(x_0(V,h))\Big)
%\nonumber\\
%& = &
%    \frac{\p\Big(H[T_{0}(V, h)]-h\leq \Delta , V(x_1(V,h))>V(x_0(V,h))\Big)
%    }{\p\big(V(x_1(V,h))>V(x_0(V,h))\big)}
%\nonumber\\
%& = &
%    \frac{
%        \E\big[
%        \sharp\big\{0\leq i < \ell(\mathcal T_V^\uparrow),
%        \, (-i)\in\Z, \, (\ell(\mathcal T_V^\uparrow)-i)\in\Z\big\}
%        \un_{\{H(\mathcal T_V^\uparrow)-h \leq \Delta\}}
%        \big]
%    }
%    {\E[\ell(\mathcal T_V^\uparrow)+\ell(\mathcal T_V^\downarrow)]
%    \p\big(V(x_1(V,h))>V(x_0(V,h))\big)}.
\nonumber\\
& = &
    \frac{
        \E\big[
        \ell(\mathcal T_{V,h}^\uparrow)
        \un_{\{H(\mathcal T_{V,h}^\uparrow)-h \leq \Delta\}}
        \big]
    }
    {\E[\ell(\mathcal T_{V,h}^\uparrow)+\ell(\mathcal T_{V,h}^\downarrow)]
    \p\big[V(x_1(V,h))>V(x_0(V,h))\big]}.
\label{Ineg_Intermediaire_Preuve_Lemme_Excess_c}
\end{eqnarray}
Notice  by Theorem \ref{Lemma_Law_of_Slopes} {\bf (i)}
and since
$
    H(\mathcal T_{V,h}^\uparrow)
=
    \mathcal T_{V,h}^\uparrow\big[\ell\big(\mathcal T_{V,h}^\uparrow\big)\big]
$,
$
    T_{\mathcal T_{V,h}^\uparrow}([h,+\infty[)
\leq
    \ell(\mathcal T_{V,h}^\uparrow)
$
and $\Delta-h<0$, and finally by \eqref{eqOptimalStopping2},
\begin{eqnarray}
&&
        \E\big[
        T_{\mathcal T_{V,h}^\uparrow}([h,+\infty[)
        \un_{\{H(\mathcal T_{V,h}^\uparrow)-h \leq \Delta\}}
        \big]
\nonumber\\
& \leq &
        \E\big[
        T_{\mathcal T_{V,h}^\uparrow}([h,+\infty[)
        \un_{\{H(\mathcal T_{V,h}^\uparrow)- T_{V,h}^\uparrow(T_{\mathcal T_{V,h}^\uparrow}([h,+\infty[))\leq \Delta\}}
        \big]
\nonumber\\
& = &
        \E\big[
        T_{\mathcal T_{V,h}^\uparrow}([h,+\infty[)
        \big]
        \p\big[H(\mathcal T_{V,h}^\uparrow)- T_{V,h}^\uparrow(T_{\mathcal T_{V,h}^\uparrow}([h,+\infty[))\leq \Delta\big]
\nonumber\\
& \leq &
        \E\big[\ell\big(\mathcal T_{V,h}^\uparrow)
        \big]
    \p\big( T_V(-h+\Delta) < T_V(]\Delta, +\infty[)\big)
\nonumber\\
& \leq &
        \E\big[\ell\big(\mathcal T_{V,h}^\uparrow)
        \big]
%    \frac{2\Delta}{h}.
    2\Delta h^{-1}.
\label{Ineg_Intermediaire_Preuve_Lemme_Excess_b}
\end{eqnarray}
Finally, once more by Theorem \ref{Lemma_Law_of_Slopes} {\bf (i)} with its notation,
\begin{eqnarray}
        \E\big[
        \big(\ell(\mathcal T_{V,h}^\uparrow)-T_{\mathcal T_{V,h}^\uparrow}([h,+\infty[)\big)
        \un_{\{H(\mathcal T_{V,h}^\uparrow)-h \leq \Delta\}}
        \big]
& \leq &
        \E\big[
        M_h^\sharp
%        \un_{\{ T_V(-h+\Delta) < T_V(]\Delta, +\infty[) \}}
        \un_{\{ \widetilde \tau_1(h) < T_V(]\Delta, +\infty[) \}}
        \big]
\nonumber\\
& \leq &
        \E\big[
        \widetilde \tau_1(h)
        \un_{\{ \widetilde \tau_1(h) < T_V(]\Delta, +\infty[) \}}
        \big].
\label{Ineg_Intermediaire_Preuve_Lemme_Excess}
\end{eqnarray}
Notice that $\widehat X_k:= (V(k))^2- \sigma ^2 k$, $k\in\N$ is a martingale for
%the natural filtration of $V$.
the filtration $\mathcal{F}_{V,k}:=\sigma(V(1),\dots,$ $V(k))$, $k\in\N$.
Moreover, the stopping time $\widetilde \tau_1(h) \wedge T_V(]\Delta, +\infty[)$
has finite expectation, since
$
    \E[\widetilde \tau_1(h) \wedge T_V(]\Delta, +\infty[)]
\leq
    \E[\widetilde \tau_1(h)]
=
    \E\big[\tau_2^{(V)}(h)-\tau_1^{(V)}(h)\big]
\leq
    \E[\ell(\mathcal T_{V,h}^\uparrow)+\ell(\mathcal T_{V,h}^\downarrow)]
<
    \infty
$
by \eqref{eqDef_tau_1_V}, \eqref{eqDef_tau_2_V}, Definition \ref{deff_law_slopes}
(see also Figure \ref{figure_tau_i_m_i_slopes})
and
Theorem \ref{Lemma_Law_of_Slopes} {\bf (iii)}.
Also, for every $k\in\N$,
$$
    \E\big[\big|\widehat X_{k+1}-\widehat X_k\big|
    %\mid \sigma(V(1),\dots, V(k))\big]
    \mid \mathcal{F}_{V,k} \big]
=
    \E\big[\big|(V(k+1))^2-(V(k))^2-\sigma^2\big|
    \mid \mathcal{F}_{V,k} \big]
\leq
    2C_0(\Delta + h+ C_0)
+
    \sigma^2
$$
a.s. on $\{k< \widetilde \tau_1(h) \wedge T_V(]\Delta, +\infty[)\}$,
since $V(k)$ and $V(k+1)$ belong to $[-h-C_0, \Delta+C_0]$ on this event
and $|V(k+1)-V(k)|\leq C_0$.
Hence by the optimal stopping time theorem
(see e.g. \cite{Geoffrey_Stirzaker_3}, (9) p. 492),
we have
$
    \E\big[\widehat X_{\widetilde \tau_1(h) \wedge T_V(]\Delta, +\infty[)}]
=
    \E\big[\widehat X_{0}\big]=0
$.
This gives
%$$
%    \E\big[\widehat X_{\widetilde \tau_1(h) \wedge T_V(]\Delta, +\infty[)}\un_{\{\widetilde \tau_1(h) < T_V(]\Delta, +\infty[)\}}]
%+
%    \E\big[\widehat X_{\widetilde \tau_1(h) \wedge T_V(]\Delta, +\infty[)}\un_{\{\widetilde \tau_1(h) > T_V(]\Delta, +\infty[)\}}]
%=0,
%$$
%so
$$
    \E\big[
    \big(
        [V(\widetilde \tau_1(h))]^2 - \sigma^2 \widetilde \tau_1(h)
    \big)
    \un_{\{\widetilde \tau_1(h) < T_V(]\Delta, +\infty[)\}}
    \big]
+
    \E\big[\widehat X_{T_V(]\Delta, +\infty[)}\un_{\{\widetilde \tau_1(h) > T_V(]\Delta, +\infty[)\}}]
=0,
$$
since $\widetilde \tau_1(h) \neq T_V(]\Delta, +\infty[)$ a.s.
Consequently, using $\widehat X_k\leq (V(k))^2$ and ellipticity \eqref{eq_ellipticity_for_V},
\begin{eqnarray*}
&&
    \sigma^2
    \E\big[
         \widetilde  \tau_1(h)
    \un_{\{\widetilde \tau_1(h) < T_V(]\Delta, +\infty[)\}}
    \big]
\\
& = &
    \E\big[
        [V(\widetilde \tau_1(h))]^2
    \un_{\{\widetilde \tau_1(h) < T_V(]\Delta, +\infty[)\}}]
+
    \E\big[\widehat X_{T_V(]\Delta, +\infty[)}\un_{\{\widetilde \tau_1(h) > T_V(]\Delta, +\infty[)\}}]
\\
& \leq &
    (h+C_0)^2
    \p[T_V(-h+\Delta) < T_V(]\Delta, +\infty[)]
+
    (\Delta+C_0)^2
\\
& \leq &
    (h+C_0)^2
%    \frac{2\Delta}{h}
    2\Delta h^{-1}
+
    (\Delta+C_0)^2,
\end{eqnarray*}
as before since $C_0<\Delta<h$.
This and \eqref{Ineg_Intermediaire_Preuve_Lemme_Excess} give for large $h$ for every $\Delta\in]C_0,h[$,
$$
        \E\big[
        \big(\ell(\mathcal T_{V,h}^\uparrow)-T_{\mathcal T_{V,h}^\uparrow}([h,+\infty[)\big)
        \un_{\{H(\mathcal T_{V,h}^\uparrow)-h \leq \Delta\}}
        \big]
\leq
    \sigma^{-2}(3\Delta h +3\Delta h).
$$
This together with \eqref{Ineg_Intermediaire_Preuve_Lemme_Excess_b} gives
\begin{equation}\label{Ineg_E_lTVup}
        \E\big[
        \ell(\mathcal T_{V,h}^\uparrow)
        \un_{\{H(\mathcal T_{V,h}^\uparrow)-h \leq \Delta\}}
        \big]
\leq
    6\sigma^{-2}\Delta h
+
        \E\big[\ell\big(\mathcal T_{V,h}^\uparrow)
        \big]
    %\frac{2\Delta}{h}.
    2\Delta h^{-1}.
\end{equation}
Moreover,
$\p\big[V(x_1(V,h))>V(x_0(V,h))\big]=\p(b_h\leq 0)\to 1/2$ as $h\to+\infty$ by Lemma \ref{Lemma_Proba_bh_positif},
so \eqref{Ineg_Intermediaire_Preuve_Lemme_Excess_c}, \eqref{Ineg_E_lTVup} and Lemma \ref{Lemma_Esperance_Longueur_Slope} give
for large $h$ for every $\Delta\in]C_0,h[$,
\begin{eqnarray*}
    \p\Big(e[T_{0}(V, h)]\leq \Delta
    | b_h\leq 0\Big)
%    | V(x_1(V,h))>V(x_0(V,h))\Big)
\leq
    \frac{6 \sigma^{-2}\Delta h }{2c_7h^2.1/3}
+
    5\frac{\Delta}{h}
\leq
    c_8\frac{\Delta}{h}
\end{eqnarray*}
with $c_8:=9\sigma^{-2}/c_7+5$.
The proof is similar if we replace $b_h \leq 0$ by $b_h>0$,
using Theorem \ref{Lemma_Central_Slope} eq. \eqref{eq_Central_Slope_Downward}
instead of eq. \eqref{eq_Central_Slope_Upward}
and since
$
    \mathcal T_{V,h}^\downarrow
=_{law}
    -\mathcal T_{-V,h}^\uparrow
$
by Theorem \ref{Lemma_Law_of_Slopes} {\bf (ii)}.
This proves \eqref{eqLemmaHeightSlopes} in the case $i=0$, which ends the proof of the lemma.
\hfill$\Box$

%%%%%%%%%%%%%%%%%%%%%%%%%%%%%%%%%%%%%%%%%%%%%%%%%%%%%%%%%%%%%%%%%%%%%%%%%%%%%%%%%%%%%%%

%%%%%%%%%%%%%%%%%%%%%%%%%%%%%%%%%%%%%%%%%%%%%%%%%%%%%%%%%%%%%%%%%%%%%%%%%%%%%%%%%%%%%%%

%%%%%%%%%%%%%%%%%%%%%%%%%%%%%%%%%%%%%%%%%%%%%%%%%%%%%%%%%%%%%%%%%%%%%%%%%%%%%%%%%%%%%%%%
%                                                                                      %
%                                   SECTION PROOF 1                                    %
%                                                                                      %
%%%%%%%%%%%%%%%%%%%%%%%%%%%%%%%%%%%%%%%%%%%%%%%%%%%%%%%%%%%%%%%%%%%%%%%%%%%%%%%%%%%%%%%%

\section{Proof of Theorem \ref{Th_Local_Limit_b_h}}\label{Sect_Proof_Th_LLT_b_h}

%In this proof, we use the notations of \cite{Devulder_Rates_CV}. (left $h$-extrema)

The proof relies mainly on the
expression of $\p(b_h= x)$ provided by Lemma \ref{Lemma_Proba_bh_egal},
the monotonicity of $x\mapsto \p(b_h= x)$ on $\N$ and $-\N$
due to Lemma \ref{Lemma_Proba_bh_egal}, the uniform continuity of $\varphi_ \infty$,
Donsker's theorem, Kesten \cite{Kesten}'s result
and some estimates on the laws of left $h$-slopes.
The proof is divided into three steps, depending on whether $x$ is far from $0$, close to $0$,
or in between.

\noindent
{\bf Proof of Theorem \ref{Th_Local_Limit_b_h}:} Let $0<\e<1/2$.

\noindent{\bf First step:}
Notice that by Lemma \ref{Lemma_Proba_bh_egal} and Markov inequality,
for $h>0$,
%$$
%    \forall x\neq 0,
%\quad
%    \p(b_h= x)
%\leq
%    \frac{\p\big[\ell\big(\mathcal T_{V,h}^\downarrow\big)\geq |x|\big]
%        +
%        \p\big[\ell\big(\mathcal T_{V,h}^\uparrow\big)> -|x|\big]}
%    {\E\big[\ell\big(\mathcal T_{V,h}^\uparrow\big)+\ell\big(\mathcal T_{V,h}^\downarrow\big)\big]}
%\leq
%    \frac{\E\big[\ell\big(\mathcal T_{V,h}^\downarrow\big)\big]}
%    {|x|\E\big[\ell\big(\mathcal T_{V,h}^\uparrow\big)+\ell\big(\mathcal T_{V,h}^\downarrow\big)\big]}
%\leq
%    \frac{2}{|x|}.
%$$
\begin{eqnarray*}
    \forall x> 0,
&&
\quad
    \p(b_h= x)
=
    \frac{\p\big[\ell\big(\mathcal T_{V,h}^\downarrow\big)\geq x\big]}
    {\E\big[\ell\big(\mathcal T_{V,h}^\uparrow\big)+\ell\big(\mathcal T_{V,h}^\downarrow\big)\big]}
\leq
    \frac{\E\big[\ell\big(\mathcal T_{V,h}^\downarrow\big)\big]}
    {x\E\big[\ell\big(\mathcal T_{V,h}^\uparrow\big)+\ell\big(\mathcal T_{V,h}^\downarrow\big)\big]}
\leq
    \frac{1}{x}
%\\
%%$$
%%Similarly,
%%$$
%    \forall x< 0,
%&& \quad
%    \p(b_h= x)
%=
%    \frac{\p\big[\ell\big(\mathcal T_{V,h}^\uparrow\big)> -x\big]}
%    {\E\big[\ell\big(\mathcal T_{V,h}^\uparrow\big)+\ell\big(\mathcal T_{V,h}^\downarrow\big)\big]}
%\leq
%    \frac{\E\big[\ell\big(\mathcal T_{V,h}^\uparrow\big)\big]}
%    {|x|\E\big[\ell\big(\mathcal T_{V,h}^\uparrow\big)+\ell\big(\mathcal T_{V,h}^\downarrow\big)\big]}
%\leq
%    \frac{1}{|x|}.
\end{eqnarray*}
and similarly
$
    \p(b_h= x)
\leq
    \frac{1}{|x|}
$
for all $x<0$.
Moreover, $\lim_{x\to\pm \infty} \varphi_\infty(x)=0$, so we can fix some $A>0$ such that,
for every $h>0$,
for all
$x\in \Z$ such that $|x|>Ah^2$, we have
\begin{eqnarray}
    \bigg|\p\big(b_h=x\big)-\frac{\sigma^2}{h^2}\varphi_\infty\bigg(\frac{\sigma^2 x}{h^2}\bigg)\bigg|
& \leq &
    \p\big(b_h=x\big)+\frac{\sigma^2}{h^2}\varphi_\infty\bigg(\frac{\sigma^2 x}{h^2}\bigg)
\leq
    \frac{1}{|x|}+\frac{\sigma^2}{h^2}\sup_{|y|\geq A\sigma^2}\varphi_\infty(y)
\nonumber\\
& \leq &
    \frac{1}{h^2}\left(\frac{1}{A}+\sigma^2\sup_{|y|\geq A\sigma^2}\varphi_\infty(y)\right)
\leq
    \frac{\e}{h^2}.
\label{Ineg_THé_Uniforme_Loin}
\end{eqnarray}

\noindent{\bf Second step:}
By Donsker's theorem,
$b_h^{(K)}/h^2$ converges in law as $h\to+\infty$ under $\p$
to $b_1^{(K,\sigma W)}$ (defined after \eqref{eqDef_bK_h}), which has the same law as $\sigma^{-2} b_1^{(K,W)}$ by scaling.
Also, the law of $b_1^{(K,W)}$ is $\varphi_\infty(x)\dd x$
by Kesten \cite{Kesten} as explained after our \eqref{eqDef_bK_h}.
Also, $\p\big[b_h\neq b_h^{(K)}\big]\to_{h\to+\infty} 0$ by Lemma \ref{Lemma_Comparaison_b_bK},
so $\sigma^2 b_h/h^2$ converges in law under $\p$ to $\varphi_\infty(x)\dd x$ as $h\to+\infty$.

%Let
%$h>0$and
%$A>0$.
Since $\varphi_\infty$ is continuous on $\R$ and $\lim_{x\to\pm \infty}\varphi_\infty(x)=0$,
$\varphi_\infty$ is uniformly continuous on $\R$.
%Let $0<\e<1/2$.
Hence, there exists $\eta>0$ such that
\begin{equation}\label{eq_uniform_continuity}
    \forall x\in\R,\ \forall y\in\R,
\qquad
    |x-y|\leq \eta
\Rightarrow
    |\varphi_\infty(x)-\varphi_\infty(y)|<\e,
\end{equation}
and we can choose
$\eta>0$ small enough so that $5\eta\sigma^{-2}\leq A$,
$5c_7^{-1} \exp[-90^{-1}\eta^{-1}] \leq \e\sigma^2$ and  $3\sqrt{5\eta}\leq 1$,
where $c_7>0$ is a constant introduced in Lemma \ref{Lemma_Esperance_Longueur_Slope}.
We can now fix $N_0\in\N$ such that $[-A,A]\subset [-N_0\eta\sigma^{-2}, N_0\eta \sigma^{-2}]$.
Since $\sigma^2 b_h/h^2$ converges in law under $\p$ to $\varphi_\infty(x)\dd x$ as $h\to+\infty$,
for all $j\in\{-N_0-3, \dots, N_0+3\}$,
$$
    \p\big(\sigma^2 b_h/h^2\in[j\eta, (j+1)\eta[\big)
\to_{h\to+\infty}
%    \p\big[b_\infty\in[j\eta, (j+1)\eta[\big]
%=
    \int_{j\eta}^{(j+1)\eta} \varphi_\infty(u)\dd u.
$$
Hence there exists $h_0>0$ such that
$\eta\sigma^{-2}h_0^2>2$,
$
    1
\leq
    [(1-\e)^{-1}-1]\eta\sigma^{-2}h_0^2
$,
$
    1
\leq
    [1-(1+\e)^{-1}]\eta\sigma^{-2}h_0^2
$
 and
$$
    \forall h\geq h_0,  \forall j\in\{-N_0-3, \dots, N_0+3\},
\quad
%    \int_{j\eta}^{(j+1)\eta} \varphi_\infty(u)\dd u -\e_2
%\leq
    \bigg|
        \p\bigg(\frac{\sigma^2 b_h}{h^2}\in[j\eta, (j+1)\eta[\bigg)
        -
        \int_{j\eta}^{(j+1)\eta} \varphi_\infty(u)\dd u
    \bigg|
\leq
    \eta\e.
%\leq
%    \int_{j\eta}^{(j+1)\eta} \varphi_\infty(u)\dd u +\e_2.
$$
%Set $\e_2=\eta\e$.
This, combined with \eqref{eq_uniform_continuity}, gives
for all $j\in\{-N_0-3, \dots, N_0+3\}$,
\begin{equation}\label{Ineg_Proba_bh_Intervalle}
    \forall h\geq h_0,
\quad
    \eta[\varphi_\infty(j\eta)-\e] -\eta \e
\leq
    \p\big(\sigma^2 b_h/h^2\in[j\eta, (j+1)\eta[\big)
\leq
    \eta[\varphi_\infty(j\eta)+\e] +\eta \e.
\end{equation}
We consider $h\geq h_0$.
Due to Lemma \ref{Lemma_Proba_bh_egal},
%\eqref{eq_proba_bh_x_positif} and \eqref{eq_proba_bh_x_negatif},
%for $h\geq h_0$,
$x\mapsto \p(b_h= x)$ is nonincreasing on $\N$, and nondecreasing on $-\N$.
Hence,
%Now,
for $0\leq j \leq N_0+3$,
% and $h\geq h_0$,
%since $x\mapsto \p[b_h=x]$ is nonincreasing on $\N$,
\begin{eqnarray*}
    \p\big(\sigma^2 b_h/h^2\in[j\eta, (j+1)\eta[\big)
%& = &
%    \p\big[b_h\in[j\eta \sigma^{-2}h^2, (j+1)\eta \sigma^{-2}h^2[\big]
%\\
& = &
    \sum_{i\in[j\eta \sigma^{-2}h^2,(j+1)\eta \sigma^{-2}h^2[\cap\N}
    \p(b_h=i)
%\\
%& \leq &
%    (\eta\sigma^{-2} h^2+1) \p\big[b_h=\big\lfloor j\eta \sigma^{-2}h^2 \big\rfloor\big]
\\
& \leq &
    (1-\e)^{-1}\eta\sigma^{-2} h^2 \p\big(b_h=\big\lfloor j\eta \sigma^{-2}h^2 \big\rfloor\big),
\end{eqnarray*}
due to the second inequality defining $h_0$.
This and \eqref{Ineg_Proba_bh_Intervalle} give for such $j$,
\begin{eqnarray}
%    \forall h\geq h_0,
%\qquad
    \p\big(b_h=\big\lfloor j\eta \sigma^{-2}h^2 \big\rfloor\big)
& \geq &
    \frac{\eta[\varphi_\infty(j\eta)-\e] -\eta \e}{(1-\e)^{-1}\eta\sigma^{-2} h^2 }
=
    \sigma^2 [\varphi_\infty(j\eta)-2\e] h^{-2}(1-\e)
\nonumber\\
& \geq &
    \sigma^2 [\varphi_\infty(j\eta)-3\e] h^{-2},
\label{Ineg_Proba_bh_geq}
\end{eqnarray}
since $\varphi_\infty(u)\in[0,2/\pi]$ for all $u\in\R$.
Similarly for such $j$,
% and $h\geq h_0$,
\begin{equation*}
    \p\big(\sigma^2 b_h/h^2\in[j\eta, (j+1)\eta[\big)
%& = &
%    \p\big[b_h\in[j\eta \sigma^{-2}h^2, (j+1)\eta \sigma^{-2}h^2[\big]
%\\
%& = &
%    \sum_{i\in[j\eta \sigma^{-2}h^2,(j+1)\eta \sigma^{-2}h^2[]\cap\N}
%    \p[b_h=i]
%\\
\geq
    \big[\eta\sigma^{-2} h^2-1\big]
    \p\big(b_h=\big\lfloor (j+1)\eta \sigma^{-2}h^2 \big\rfloor\big).
%    \p\big[b_h=\big\lceil (j+1)\eta \sigma^{-2}h^2 \big\rceil\big].
\end{equation*}
This and \eqref{Ineg_Proba_bh_Intervalle} give, using the third inequality in the definition of $h_0$,
\begin{eqnarray}
%    \p\big[b_h=\big\lceil (j+1)\eta \sigma^{-2}h^2 \big\rceil\big]
    \p\big(b_h=\big\lfloor (j+1)\eta \sigma^{-2}h^2 \big\rfloor\big)
& \leq &
    \frac{\eta[\varphi_\infty(j\eta)+\e] +\eta \e}{(1+\e)^{-1}\eta\sigma^{-2} h^2 }
=
    (1+\e)\sigma^2 [\varphi_\infty(j\eta)+2\e] h^{-2}
\nonumber\\
& \leq &
    \sigma^2 [\varphi_\infty(j\eta)+4\e] h^{-2},
\label{Ineg_Proba_bh_leq}
\end{eqnarray}
since $0<\e<1/2$ and $\varphi_\infty(u)\in[0,2/\pi]$ for all $u\in\R$.

%Now, let $j\in\{3,\dots, N_0\}$ and $x\in\N$ such that
Now, let $j\in\{2,\dots, N_0\}$ and $x\in\N$ such that
$j\eta \sigma^{-2}h^2\leq x < (j+1)\eta \sigma^{-2}h^2$.
%$j\eta \sigma^{-2}h^2\leq x \leq (j+1)\eta \sigma^{-2}h^2$.
We have since $\p(b_h=.)$ is nonincreasing on $\N$ and
%$x\leq \lceil (j+1)\eta \sigma^{-2}h^2 \rceil \leq \lfloor (j+2)\eta \sigma^{-2}h^2 \rfloor$
$x\leq \lfloor (j+1)\eta \sigma^{-2}h^2 \rfloor$,
%since $h$ is large enough so that $\eta\sigma^{-2}h^2\geq \eta\sigma^{-2}h_0^2>2$,
then by \eqref{Ineg_Proba_bh_geq} and finally by \eqref{eq_uniform_continuity},
\begin{eqnarray*}
    \p(b_h=x)
& \geq &
    \p(b_h=\lfloor (j+1)\eta \sigma^{-2}h^2 \rfloor)
%    \p[b_h=\lceil (j+1)\eta \sigma^{-2}h^2 \rceil]
%\geq
%    \p\big[b_h=\big\lfloor (j+2)\eta \sigma^{-2}h^2 \big\rfloor\big]
\geq
%    \sigma^2 [\varphi_\infty((j+2)\eta)-3\e] h^{-2}
    \sigma^2 [\varphi_\infty((j+1)\eta)-3\e] h^{-2}
\\
& \geq &
%    \sigma^2 \big[\varphi_\infty\big(x \sigma^2 h^{-2}\big)-5\e\big] h^{-2}.
    \sigma^2 \big[\varphi_\infty\big(x \sigma^2 h^{-2}\big)-4\e\big] h^{-2}.
\end{eqnarray*}
Similarly, using \eqref{Ineg_Proba_bh_leq} applied to $j-1\geq 1$ instead of \eqref{Ineg_Proba_bh_geq},
followed by \eqref{eq_uniform_continuity},
\begin{eqnarray*}
    \p(b_h=x)
& \leq &
    \p(b_h=\lfloor j\eta \sigma^{-2}h^2 \rfloor)
%\leq
%    \p\big[b_h=\big\lceil (j-1)\eta \sigma^{-2}h^2 \big\rceil\big]
\\
& \leq &
%    \sigma^2 [\varphi_\infty((j-2)\eta)+4\e] h^{-2}
    \sigma^2 [\varphi_\infty((j-1)\eta)+4\e] h^{-2}
\leq
%    \sigma^2 \big[\varphi_\infty\big(x \sigma^2 h^{-2}\big)+7\e\big] h^{-2}.
    \sigma^2 \big[\varphi_\infty\big(x \sigma^2 h^{-2}\big)+6\e\big] h^{-2}.
\end{eqnarray*}
Since this is true for all $h\geq h_0$, every $j\in\{2,\dots,N_0\}$
and for every $x\in\N$ such that $j\eta \sigma^{-2}h^2\leq x < (j+1)\eta \sigma^{-2}h^2$
%for $j\in\{3,\dots, N_0\}$
for such $j$, and $A\leq N_0\eta \sigma^{-2}$, this gives
\begin{equation}\label{Ineg_Comparaison_Probas_bh_phi}
    \forall h\geq h_0,
\qquad
    \max_{x\in[2\eta \sigma^{-2}h^2, A h^2]\cap\Z}
    \big| \p(b_h=x) - \varphi_\infty\big(x \sigma^2 h^{-2}\big)\sigma^2h^{-2}\big|
%\leq
%    \max_{x\in[3\eta \sigma^{-2}h^2, N_0\eta \sigma^{-2}h^2]\cap\Z}
%   \big| \p[b_h=x] - \varphi_\infty\big(x \sigma^2 h^{-2}\big)\sigma^2h^{-2}\big|
\leq
%    7\e \sigma^2 h^{-2}.
    6\e \sigma^2 h^{-2}.
\end{equation}
We get similarly
\begin{equation}\label{Ineg_Comparaison_Probas_bh_phi_neg}
    \forall h\geq h_0,
\qquad
    \max_{x\in[-A h^2, -2\eta \sigma^{-2}h^2]\cap\Z}
    \big| \p(b_h=x) - \varphi_\infty\big(x \sigma^2 h^{-2}\big)\sigma^2h^{-2}\big|
%\leq
%    \max_{x\in[3\eta \sigma^{-2}h^2, N_0\eta \sigma^{-2}h^2]\cap\Z}
%   \big| \p[b_h=x] - \varphi_\infty\big(x \sigma^2 h^{-2}\big)\sigma^2h^{-2}\big|
\leq
%    7\e \sigma^2 h^{-2}.
    6\e \sigma^2 h^{-2}.
\end{equation}

\noindent{\bf Third step:}
Now, for $-5\eta \sigma^{-2}h^2\leq x\leq 0$, we have by \eqref{eq_proba_bh_x_negatif}
%\eqref{eq_proba_bh_x_positif}
and \eqref{eq_proba_bh_x_zero},
\begin{equation}\label{Ineg_Difference_Proba_bh_0}
    \big| \p(b_h=x)-\p(b_h=0)\big|
=
    \frac{\p\big[\ell\big(\mathcal T_{V,h}^\uparrow\big)\leq - x\big]}
    {\E\big[\ell\big(\mathcal T_{V,h}^\uparrow\big)+\ell\big(\mathcal T_{V,h}^\downarrow\big)\big]}
\leq
    \frac{\p\big[\ell\big(\mathcal T_{V,h}^\uparrow\big)\leq 5\eta \sigma^{-2}h^2 \big]}
    {c_7 h^2}
\end{equation}
(uniformly) for all $-5\eta \sigma^{-2}h^2\leq x\leq 0$ for large $h$, since
$
    \E\big[\ell\big(\mathcal T_{V,h}^\uparrow\big)\big]
\sim_{h\to+\infty}
    \E\big[\ell\big(\mathcal T_{V,h}^\downarrow\big)\big]
$
$
\sim_{h\to+\infty}
    c_7 h^2
$
by Lemma \ref{Lemma_Esperance_Longueur_Slope}.

We know from Theorem \ref{Lemma_Law_of_Slopes} {\bf (i)} that
up to its first hitting time of $[h,+\infty[)$, $\mathcal T_{V,h}^\uparrow$ has the same law as
$(V(k),\ 0\leq k \leq T_V([h,+\infty[)$ conditioned by $\{T_V([h,+\infty[)<T_V(]-\infty,0[)\}$.
Thus for $\alpha>0$, applying the strong Markov property in the last equality,
and ellipticity \eqref{eq_ellipticity_for_V} in the last line (for $h$ large enough so that $C_0<h/6$),
\begin{eqnarray*}
&&
    \p\big[\ell\big(\mathcal T_{V,h}^\uparrow\big)\leq \alpha h^2 \big]
\leq
%    \p\big[T_{\mathcal T_V^\uparrow}([h,+\infty[)\leq \alpha h^2 \big]
    \p\big[T_{\mathcal T_{V,h}^\uparrow}([h,+\infty[)
        -T_{\mathcal T_{V,h}^\uparrow}([h/2,+\infty[)
    \leq \alpha h^2 \big]
\\
& = &
%    \p\big[T_V([h,+\infty[)\leq \alpha h^2 \mid T_V([h,+\infty[)<T_V(]-\infty,0[)\big]
%\\
%& \leq &
%    \p\big[T_V([h,+\infty[)-T_V([h/2,+\infty[)\leq \alpha h^2
%    \mid
%    T_V([h,+\infty[)<T_V(]-\infty,0[)\big]
%\\
%& = &
    \frac{\p\big[T_V([h,+\infty[)-T_V([h/2,+\infty[)\leq \alpha h^2 , T_V([h,+\infty[)<T_V(]-\infty,0[)\big]}
    {\p\big[T_V([h,+\infty[)<T_V(]-\infty,0[)\big]}
\\
& = &
    \frac{\E\big[\un_{\{T_V([h/2,+\infty[)<T_V(]-\infty,0[)\}}
                 \p^{V(T_V([h/2,+\infty[))}\big[T_V([h,+\infty[)\leq (\alpha h^2)\wedge T_V(]-\infty,0[)\big]
          \big]}
    {\p\big[T_V([h,+\infty[)<T_V(]-\infty,0[)\big]}
\\
& \leq &
    \frac{\p\big[T_V([h/2,+\infty[)<T_V(]-\infty,0[)\big]}{\p\big[T_V([h,+\infty[)<T_V(]-\infty,0[)\big]}
                 \p\big[T_V([h/3,+\infty[)
                    %\wedge T_V(]-\infty,-h/3])
                    \leq \alpha h^2\big].
\end{eqnarray*}
Using $\p\big[T_V([h,+\infty[)<T_V(]-\infty,0[)\big]\sim_{h\to+\infty}
%c^*
c_1
h^{-1}$
(see \eqref{eq_Proba_Atteinte_logn_avant0})
and Donsker's theorem, the last line is equivalent, as $h\to+\infty$, to
\begin{eqnarray*}
    2 \p\big[T_{\sigma W}([1/3,+\infty[))
    %\wedge T_{\sigma W}(]-\infty,-h/3])
        \leq \alpha\big]
& = &
    2\p\Big[\sup_{[0, \alpha ]}(\sigma W)\geq 1/3\Big]
=
    2\p\big[\sigma |W(\alpha )|\geq 1/3\big]
\\
& = &
    2\p\big[|W(1)|\geq
    (3\sigma \sqrt{\alpha})^{-1}\big]
\leq
    4\exp[-(3\sigma \sqrt{\alpha})^{-2}/2]
\end{eqnarray*}
if $3\sigma \sqrt{\alpha}\leq 1$, where $(W(x),\ x\in\R)$ is a two-sided Brownian motion as before.
Since $3\sqrt{5\eta}\leq 1$,
this and \eqref{Ineg_Difference_Proba_bh_0} give for large $h$,
\begin{equation*}
    \max_{-5\eta \sigma^{-2}h^2\leq x\leq 0} \big| \p(b_h=x)-\p(b_h=0)\big|
\leq
%    5c_7^{-1}h^{-2} \exp\big[-\big(3\sigma \sqrt{ 5\eta \sigma^{-2} }\big)^{-2}/2\big]
    5c_7^{-1}h^{-2} \exp\big[-\big(3 \sqrt{ 5\eta  }\big)^{-2}/2\big]
%\\
%& = &
%    5c_7^{-1}h^{-2}     \exp[-(3.5^{1/2} \eta^{1/2})^{-2}/2]
%=
%    5c_7^{-1} \exp[-90^{-1}\eta^{-1}] h^{-2}
\leq
    \e\sigma^2 h^{-2}
\end{equation*}
by the second inequality after \eqref{eq_uniform_continuity}.
Since we have a similar result for $0\leq x \leq 5\eta \sigma^{-2}h^2$,
using \eqref{eq_proba_bh_x_positif} instead of \eqref{eq_proba_bh_x_negatif}
and e.g. $\ell\big(\mathcal T_{V,h}^\downarrow\big)=_{law}\ell\big(\mathcal T_{-V,h}^\uparrow\big)$
(see Theorem \ref{Lemma_Law_of_Slopes} {\bf (ii)}),
there exists $h_1>h_0$ such that
\begin{equation}\label{Ineg_Difference_Proba_bh_x_0}
    \forall h\geq h_1,
\qquad
    \max_{-5\eta \sigma^{-2}h^2\leq x\leq 5\eta \sigma^{-2}h^2} \big| \p(b_h=x)-\p(b_h=0)\big|
\leq
    \e\sigma^2 h^{-2}.
\end{equation}
We already know, from \eqref{Ineg_Proba_bh_geq}, that
$
    \forall h\geq h_1\geq h_0,
    \p(b_h=0)
\geq
    \sigma^2 [\varphi_\infty(0)-3\e] h^{-2}.
$
Moreover, using \eqref{Ineg_Difference_Proba_bh_x_0}, \eqref{Ineg_Comparaison_Probas_bh_phi} and then \eqref{eq_uniform_continuity},
\begin{eqnarray*}
    \p(b_h=0)
& \leq &
    \p(b_h=\lfloor 4\eta \sigma^{-2}h^2\rfloor) + \e\sigma^2 h^{-2}
\leq
    [\phi_\infty(\lfloor 4\eta \sigma^{-2}h^2\rfloor \sigma^2 h^{-2})+7\e]\sigma^2 h^{-2}
\\
& \leq &
    [\phi_\infty(0)+11\e]\sigma^2 h^{-2}
\end{eqnarray*}
for $h\geq h_1$. So,
\begin{equation}\label{Ineg_varphi_infini}
    \forall h\geq h_1,
\qquad
    \big|
    \p(b_h=0)-\sigma^2 \varphi_\infty(0) h^{-2}
    \big|
\leq
    11\e\sigma^2 h^{-2}.
\end{equation}
Finally, once more by \eqref{eq_uniform_continuity},
$\big|\varphi_\infty\big(x \sigma ^2  h^{-2}\big)-\varphi_\infty(0)\big|\leq 5\e$
for $x\in\Z$ such that  $|x| \leq 5\eta \sigma^{-2}h^2$. This, combined with \eqref{Ineg_Difference_Proba_bh_x_0} and \eqref{Ineg_varphi_infini} and the triangular inequality yields to
\begin{equation}\label{Ineg_Difference_Proba_bh_varphi}
    \forall h\geq h_1,
\qquad
   \max_{x\in[-5\eta \sigma^{-2}h^2, 5\eta \sigma^{-2}h^2]\cap \Z}
   \big| \p(b_h=x)- \varphi_\infty\big(x \sigma ^2 h^{-2}\big)\sigma^2  h^{-2}\big|
\leq
    17\e\sigma^2 h^{-2}.
\end{equation}
This, together with \eqref{Ineg_Comparaison_Probas_bh_phi} and \eqref{Ineg_Comparaison_Probas_bh_phi_neg} leads to
\begin{equation}\label{Ineg_Comparaison_Probas_bh_phi_AA}
    \forall h\geq h_1,
\qquad
    \max_{x\in[-A h^2, A h^2]\cap\Z}
    \big| \p(b_h=x) - \varphi_\infty\big(x \sigma^2 h^{-2}\big)\sigma^2h^{-2}\big|
\leq
    17\e \sigma^2 h^{-2}.
\end{equation}
This, combined with \eqref{Ineg_THé_Uniforme_Loin}, proves Theorem \ref{Th_Local_Limit_b_h}.
\hfill$\Box$

%%%%%%%%%%%%%%%%%%%%%%%%%%%%%%%%%%%%%%%%%%%%%%%%%%%%%%%%%%%%%%%%%%%%%%%%%%%%%%%%%%%%%%%%
%                                                                                      %
%                                   COUPLING                                           %
%                                                                                      %
%%%%%%%%%%%%%%%%%%%%%%%%%%%%%%%%%%%%%%%%%%%%%%%%%%%%%%%%%%%%%%%%%%%%%%%%%%%%%%%%%%%%%%%%

\section{Coupling argument when $b_{\log n}$ is close to $z$}\label{Sect_Coupling}

In this section, we use a coupling argument, in order to approximate the quenched probability $\po[S_n=z]$
by the invariant probability measure at $z$ of a RWRE reflected inside the central valley of the potential. In order to make this approximation, we require some conditions,  mainly for the environment.

\subsection{An inequality related to hitting times of $(S_k)_k$}\label{Subsec_Inequality_Hitting_Times}
\begin{figure}[htbp]
\includegraphics[width=16.0cm,height=6.88cm]{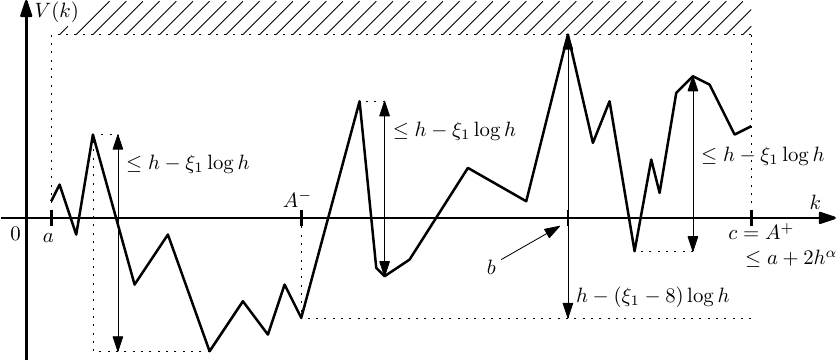}
\caption{Schema of the potential $V$ for Lemma \ref{Lemma_Proba_Descente} between $a$ and $c$ when $c=A^+$.}
\label{figure_Lemma_Proba_Descente}
\end{figure}

Before dealing with the coupling argument, we prove a useful inequality about hitting times.
This lemma is in the same spirit as (\cite{DGP_Collision_Sinai}, Lemma 4.7), but is more general.
We will use this lemma with different values of $\xi_1$.
See Figure \ref{figure_Lemma_Proba_Descente} for the schema of
the potential $V$ under the hypotheses of this lemma.

\begin{lem}\label{Lemma_Proba_Descente}
Assume \eqref{eqEllipticity}.
Let $\xi_1>0$, $\xi_2>0$ and $\alpha>0$.
There exists $\widehat h_2=\widehat h_2(\xi_1,\xi_2)>1$ such that,
for almost every environment $\omega$,
for every
$a<b<c$ and $h\geq \widehat h_2$ such that {\bf (i)} $V(b)=\max_{[a,c]}V$,
{\bf (ii)}
$
    \max_{b\leq \ell\leq k\leq c-1}\big(V(k)-V(\ell)\big)
\leq
    h-\xi_1\log h
$,
{\bf (iii)}
$
    \max_{a\leq \ell\leq k\leq b-1}\big(V(\ell)-V(k)\big)
\leq
    h-\xi_1\log h
$
and
{\bf (iv)}
$|c-a|\leq 2 h^{\alpha}$,
and for every $a\leq x \leq c$, we have
\begin{equation}\label{Ineg_Lemma_Proba_Descente}
    \po^x\big[\tau(a)\wedge\tau(c)\geq \xi_2 e^h\big]
\leq
    24\xi_2^{-1}\e_0^{-2} h^{2\alpha-\xi_1+8}
+
    4\e_0^{-1}h^{\alpha-8},
\end{equation}
and is, in particular, uniformly less than $h^{-4}$ for all $h\geq \widehat h_2$ if $\alpha=3$ and $\xi_1>19$.
\end{lem}

\noindent{\bf Proof:}
We cannot apply directly \eqref{InegEsperance1} or \eqref{InegEsperance2} to $\eo[\tau(a)\wedge \tau(c)]$,
because the $\max(\dots)$
which appear in these inequalities can be much too large, since they can be
respectively nearly     as large as $V(b)-V(a)$ or $V(b)-V(c)$, which can be much larger than our $h$.
Consider $\widehat h_2>1$ such that $h-(\xi_1-8)\log h>0$ for every $h\geq \widehat h_2$.
We fix $h\geq \widehat h_2$, and assume that the hypotheses of the lemma are satisfied for this $h$.
We define
(see Figure \ref{figure_Lemma_Proba_Descente}), with $x\vee y:=\max(x,y)$,
\begin{eqnarray*}
    A^-
& := &
    a\vee \big(\max\{y\leq b, \ V(b)-V(y)\geq h-(\xi_1-8)\log h\}\big),
\\
    A^+
& := &
    c \wedge \big(\min\{y\geq b, \ V(b)-V(y)\geq h-(\xi_1-8)\log h\}\big).
\end{eqnarray*}

{\bf First case:} we assume that $a\leq x \leq A^-$.
We start with the sub-case $a<x\leq A^-$, which implies that $A^-=\max\{\dots\}\neq a$ in the definition of $A^-$.
Then, by Markov inequality, \eqref{InegEsperance2} and Hypotheses {\bf (iii)} and {\bf (iv)},
\begin{eqnarray}
    \po^x[\tau(a)\wedge \tau(b)\geq \xi_2 e^h/2]
& \leq &
    2\xi_2^{-1}e^{-h}\e_0^{-1}(b-a)^2
    \exp\Big[\max_{a\leq \ell\leq k\leq b-1}\big(V(\ell)-V(k)\big)\Big]
\nonumber\\
& \leq &
    8\xi_2^{-1}e^{-h}\e_0^{-1} h^{2\alpha} \exp(h-\xi_1\log h)
    =
    8\xi_2^{-1}\e_0^{-1} h^{2\alpha-\xi_1}.
~~~~
\label{Ineg_Proba_Atteinte_ab_x}
\end{eqnarray}
Also, notice that since $a< A^-< b$, using Hypothesis {\bf (iii)},
\begin{eqnarray}
    \max_{[a,A^-]}V
& \leq &
    V(A^-)
    +
    \max_{a\leq \ell\leq k\leq b-1}\big(V(\ell)-V(k)\big)
\nonumber\\
& \leq &
    V(b)-(h-(\xi_1-8)\log h)+(h-\xi_1\log h)
\leq
    V(b)-8\log h.
\label{Ineg_Max_a_A-}
\end{eqnarray}
Hence using \eqref{probaatteinte}, $a< x \leq A^-<b$,
then Hypothesis {\bf (iv)}, ellipticity \eqref{eqEllipticity} and \eqref{Ineg_Max_a_A-},
\begin{eqnarray*}
    \po^x[\tau(b) < \tau(a)]
& \leq &
    (x-a)\exp[\max\nolimits_{[a,x-1]}V -V(b-1)]
\leq
    2h^\alpha
    \e_0^{-1}h^{-8}.
\end{eqnarray*}
Consequently, this and \eqref{Ineg_Proba_Atteinte_ab_x} lead to
\begin{eqnarray}
&&
    \po^x[\tau(a)\wedge \tau(c)\geq \xi_2 e^h/2]
\nonumber\\
& \leq &
    \po^x[\tau(b) < \tau(a)]
+
    \po^x[\tau(a)\wedge \tau(c)\geq \xi_2 e^h/2, \tau(a)<\tau(b)<\tau(c)]
\nonumber\\
& \leq &
    2\e_0^{-1}h^{\alpha-8}
+
    \po^x[\tau(a)\wedge \tau(b)\geq \xi_2 e^h/2]
\leq
    2\e_0^{-1}h^{\alpha-8}
+
    8\xi_2^{-1} \e_0^{-1} h^{2\alpha-\xi_1}.
\label{Ineg_Proba_Sortir_Cas_1}
\end{eqnarray}
This remains true if $x=a$, whether $a=A^-$ or $a\neq A^-$, and so for every $a\leq x \leq A^-$.
This already proves \eqref{Ineg_Lemma_Proba_Descente} in this case.

{\bf Second case:} we now assume that $A^+\leq  x \leq c $. This case is similar as the first one, so we get by symmetry,
using {\bf (ii)} instead of {\bf (iii)} and \eqref{InegEsperance1} instead of \eqref{InegEsperance2},
\begin{equation}
    \po^x[\tau(a)\wedge \tau(c)\geq \xi_2 e^h/2]
\leq
    2\e_0^{-1}h^{\alpha-8}
+
    8\xi_2^{-1} \e_0^{-1} h^{2\alpha-\xi_1}.
\label{Ineg_Proba_Sortir_Cas_2}
\end{equation}
This already proves \eqref{Ineg_Lemma_Proba_Descente} in the this case.

{\bf Third case:} We now assume that $A^-<x<A^+$.
Using Markov inequality, \eqref{InegEsperance1} and Hypothesis {\bf (iv)}
and $a\leq A^-< A^+\leq c\leq a+2h^\alpha$ in the first line,
then
$
    \max_{[A^-, A^+]} V
%=\max_{[a,c]} V
=V(b)
$
(due to Hypothesis {\bf (i)} and $b\in[A^-,A^+]\subset[a,c]$)
%, \eqref{Ineg_Max_a_A-} and the corresponding inequality obtained for $[A^+,c]$ by symmetry),
and $\min_{[A^-, A^+]} V\geq  V(b)-(h-(\xi_1-8)\log h)-\log\e_0^{-1}$
(by definition of $A^\pm$ and ellipticity \eqref{eq_ellipticity_for_V}), we have
\begin{eqnarray*}
    \po^x[\tau(A^-)\wedge \tau(A^+)\geq \xi_2 e^h/2]
& \leq &
    2\xi_2^{-1}e^{-h}
    \e_0^{-1} (2h^{\alpha})^2 \exp\Big[\max_{[A^-,A^+]}V-\min_{[A^-,A^+]}V\Big]
\\
& \leq &
    8\xi_2^{-1}e^{-h}
    \e_0^{-1} h^{2\alpha} \exp\Big[h-(\xi_1-8)\log h+\log\e_0^{-1}\Big]
\\
& = &
    8\xi_2^{-1}\e_0^{-2} h^{2\alpha-\xi_1+8}.
\end{eqnarray*}
Consequently, we have by the strong Markov property applied at time $\tau(A^-)\wedge \tau(A^+)$,
\begin{eqnarray*}
    \po^x\big[\tau(a)\wedge \tau(c)\geq \xi_2 e^h\big]
& \leq &
    \po^x\big[\tau(A^-)\wedge \tau(A^+)\geq \xi_2 e^h/2\big]
+
    \po^{A^-}\big[\tau(a)\wedge \tau(c)\geq \xi_2 e^h/2\big]
\\
&&
+
    \po^{A^+}\big[\tau(a)\wedge \tau(c)\geq \xi_2 e^h/2\big]
\\
& \leq &
    8\xi_2^{-1}\e_0^{-2} h^{2\alpha-\xi_1+8}
+
    2\big(2\e_0^{-1}h^{\alpha-8}+8\xi_2^{-1} \e_0^{-1} h^{2\alpha-\xi_1}\big)
\\
& \leq &
    24\xi_2^{-1}\e_0^{-2} h^{2\alpha-\xi_1+8}
+
    4\e_0^{-1}h^{\alpha-8}
\end{eqnarray*}
by \eqref{Ineg_Proba_Sortir_Cas_1} and \eqref{Ineg_Proba_Sortir_Cas_2} applied respectively at $A^-$ and $A^+$.
%{\bf (for example, if $\alpha=3$, $h=\log n$, $\gamma=8$, $\xi_1>18$, this is $\leq (\log n)^{-4}$ for large $n$)}.
This proves \eqref{Ineg_Lemma_Proba_Descente} in this third case,
so \eqref{Ineg_Lemma_Proba_Descente} is proved
in every case for every $h$ larger than some constant $\widehat h_2>1$.
Finally, when $\alpha=3$ and $\xi_1>19$,
we have
$
    24\xi_2^{-1}\e_0^{-2} h^{2\alpha-\xi_1+8}
+
    4\e_0^{-1}h^{\alpha-8}
\leq
    (24\xi_2^{-1}\e_0^{-2}
+
    4\e_0^{-1})h^{-5}
$
which is $o(h^{-4})$ as $h\to+\infty$,
so, up to a change of $\widehat h_2$, the right hand side of \eqref{Ineg_Lemma_Proba_Descente}
is less than $h^{-4}$ for all $h>\widehat h_2$, which ends the proof of the lemma.
\hfill$\Box$

\subsection{Some events useful for the coupling argument}
In order to evaluate the probability $\PP(S_n=z)$, we decompose the event $\{S_n=z\}$ into smaller ones, and to this aim we introduce some conditions on the environment $\o$.
First,
we fix $C_1>20$, $C_2>9$, and $\delta_1\in]0,2/3[$.
%$c_{834}>0$,
For $n\geq 3$, we introduce
%{\bf (pourquoi pas $C_2=C_1$ ?)}
$$
    h_n
:=
    \log n -C_1\log_2 n,
\qquad
    \widetilde h_n
:=
    h_n-C_1\log_2 n,
\qquad
    \Gamma_n
:=
    \big\lfloor (\log n)^{4/3+\delta_1}\big\rfloor,
$$
where for $x>1$, $\log_2 x:=\log\log x$.
We also fix an integer $n_3\geq 3$ such that, for all $n\geq n_3$,
$\log_2 n>C_0+1$,
 $\log n>\max\big[2\e_0^{-1}, \widehat h_2(2C_1,1/10), \widehat h_2(C_1,1/10), \widehat h_2(2C_1,1), p_5\big]$,
$h_n-C_1\log_2 n >\max\{3C_0+10\log_2 n, (\log n)/2+(2C_1+C_2+2)\log_2 n\}$,
%$h_n-C_1\log_2 n >4\log_2 n >4C_0$,
$(\log_2 n)^6\leq \log n$,
$n\geq (\log n)^{C_1+4}$
and
$\Gamma_n\geq p_4$,
with $p_4$ and $p_5$ defined in Proposition \ref{Lemma_Laplace_V_Conditionne}
and $\widehat h_2$ in Lemma \ref{Lemma_Proba_Descente}.
We also define for $n\geq n_3$ and $z\in\Z$,
\begin{eqnarray}
\label{eq_def_E1}
    E_{-}^{(n)}
& := &
    \{b_{\log n}\leq 0\}
=
    \{b_{\log n}=x_0(V,\log n)\},
\\
    E_{+}^{(n)}
& := &
    \{b_{\log n}> 0\}
=
    \{b_{\log n}=x_1(V,\log n)\}
=
    \big(     E_{-}^{(n)}\big)^c,
\nonumber\\
    E_3^{(n)}
& := &
    \cap_{i=-10}^{10}
%    \cap_{i=-9}^{9}
    \big\{
        H[T_i(V,h_n-C_1\log_2 n)]\geq \log n+C_2 \log_2 n
    \big\},
%\\
%    E_3^{(n)}
%& := &
%    \big\{
%        \sharp\{i\in\Z,\ -10\leq i \leq 10,\ H(T_i(V,h_n-C_1\log_2 n))<\log n+C_2 \log_2 n \}=0
%    \big\},
\nonumber\\
    E_4^{(n)}(z)
& := &
    \{V(z)-V(b_{\log n})\geq 5\log_2 n\}
\nonumber\\
&&
    \cup
    \Big(E_-^{(n)}\cap \Big\{\max_{[b_{\log n}, 0]}V<V[x_1(V,\log n)]-9\log_2 n\Big\}\Big)
\nonumber\\
&&
    \cup
    \Big(E_+^{(n)}\cap \Big\{\max_{[0,b_{\log n}]}V<V[x_0(V,\log n)]-9\log_2 n\Big\}\Big),
\nonumber\\
    E_5^{(n)}
& := &
    \big\{-(\log n)^{2+\delta_1} \leq x_{-12}(V, \log n)\leq x_{12}(V,\log n)\leq (\log n)^{2+\delta_1}\big\},
\nonumber\\
    E_6^{(n)}
& := &
    \big\{\max\{V(b_{\log n}+i)-V(b_{\log n}), |i|\leq \Gamma_n\}< \log n\big\},
\label{eq_def_E5}
\nonumber\\
    E_7^{(n)}(z)
& := &
    \{|b_{\log n}-z|\leq \Gamma_n\}.
%    \{|b_{\log n}-z|\leq (\log n)^{4/3+\delta_1}\}.
\end{eqnarray}
%{\bf (attention $E_5$ change avec $12$ au lieu de $10$)}
%{\bf (ou remplacer $\leq \log n$ par $\widetilde h_n$ dans $E_6$ ?? eventuellement pour lower bound, pas indispensable pour cela)}
Finally, let
\begin{equation}\label{eq_def_EL}
    E_C^{(n)}(z)
:=
    %E_{\pm}^{(n)}\cap
    E_3^{(n)} \cap E_4^{(n)}(z)\cap E_5^{(n)}\cap E_6^{(n)}
    \cap
    E_7^{(n)}(z).
\end{equation}
%and
%$
%    E_{C,\pm}^{(n)}(z)
%:=
%    E_C^{(n)}(z)
%\cap
%    E_{\pm}^{(n)}(z)
%$.

%The reason why we consider environments $\o$ belonging to $E_-^{(n)}$, instead of $E_+^{(n)}$
%is that this event contains $\{b_{\log n}=0\}$ (as was the case in continuous setting with the definition of \cite{NP}
%and we preferred keeping this tradition),
%and this last event is not negligible, since its probability is approximatively $(\log n)^{-2}$.

\begin{remark}\label{Remark_x_i}
For $\omega\in E_3^{(n)}$, for every $-9\leq i \leq 10$,
$H(T_{i-1}(V, h_n-C_1\log_2 n)) \geq \log n$ and
$H(T_i(V, h_n-C_1\log_2 n)) \geq \log n$, so $x_i(V,h_n-C_1\log_2 n)$  is also a left $(\log n)$-extremum.
So, $x_i(V,\log n)=x_i(V,h_n-C_1\log_2 n)$ for every $-9\leq i \leq 10$, and as a consequence,
$H[T_i(V,\log n)]=H[T_i(V,h_n-C_1\log_2 n)]$ for every $-9\leq i \leq 9$.
\end{remark}

The previous events depend only on the environment $\o$ and on $z$.
They are useful for the coupling argument used in this section.
More precisely, we saw in Remark \ref{Remark_x_i} that $E_3^{(n)}$ ensures that
$x_i(V, h_n-C_1\log_2 n)=x_i(V,\log n)$ for $|i|\leq 9$, and
as a consequence, there is no subvalley of height slightly less than $\log n$
in the $(\log n)$-central valley (defined after \eqref{eq_def_Mpm}),
so $(S_k)_k$ is not trapped a long time in such subvalleys,
which helps $(S_k)_k$ to go quickly to $b_{\log n}$  with large quenched probability.

Also,
$E_4^{(n)}(z)$ is useful to prove a technical lemma, Lemma \ref{Lemma_Produit_Nu_Proba}.
$E_5^{(n)}$ says that the
%$|x_i(V, h_n-C_1\log_2 n)|$
$|x_i(V, \log n)|$
are quite small,
which will often be useful in applying inequalities
such as \eqref{probaatteinte}, \dots, \eqref{InegProba2} to prove that some events are negligible.
Finally, $E_6^{(n)}$ and $E_7^{(n)}(z)$ will imply
in particular
that $z$ is inside the $(\log n)$-central valley (see \eqref{Ineg_zn_Mplus}).

We will use, in the proof of Theorem \ref{Th_Local_Limit_Sinai}, left $h$-extrema of $V$ for three different values of $h$.
In particular,
left $(\log n)$-extrema are useful to define $b_{\log n}$,
left $\widetilde h_n$-extrema are useful e.g. to use $E_3^{(n)}$ as explained previously,
and the proof of Lemma \ref{Lem_Proba_zn_Proche_Fond} uses
left $h$-extrema with two different values strictly less than $\log n$, which are
$h_n$ and $\widetilde h_n$; left $h_n$-extrema are used in Lemma \ref{Lemma_Valle_k_E2c} (in view of \eqref{eq_def_Ik}
and Lemma \ref{Lemma_Nombre_Vallees_Visitees}),
whereas left $\widetilde h_n$-extrema are also used
in Lemma \ref{Lem_Fond_Velle_k_Droite} and in
the proof of the lower bound of Theorem \ref{Th_Local_Limit_Sinai}
(see Section \ref{Sect_lower_Bonud}).

In the rest of the paper, the $n_i$, $3\leq i \leq 19$, denote some integers with $n_i\leq n_{i+1}$ for $3\leq i \leq 18$,
which are useful to get the uniformity in Theorem \ref{Th_Local_Limit_Sinai}
($n_3$ being defined before \eqref{eq_def_E1}).

\subsection{Definition of the coupling}
\begin{figure}[htbp]
\includegraphics[width=15.98cm,height=7.02cm]{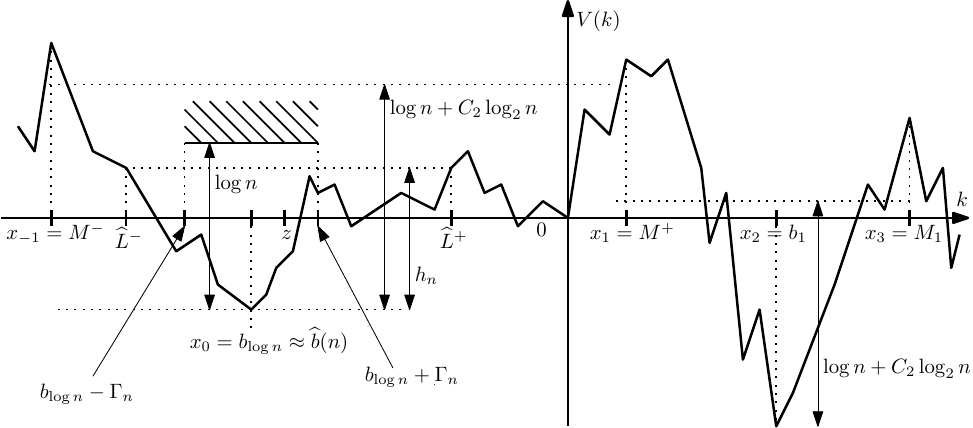}
\caption{Schema of the potential $V$ for $\omega \in E_C^{(n)}(z)$ in the case $b_{\log n}\leq 0$.}
\label{figure_Good_Env_Coupling}
\end{figure}

We fix an integer
%$n\geq 3$,
$n\geq n_3$,
$z\in\Z$, and an environment $\omega\in E_C^{(n)}(z)$.
In all the remaining of Section \ref{Sect_Coupling}, we set $x_i:=x_i(V,\log n)$, $i\in\Z$ (defined before \eqref{eqDefbh}),  to simplify the notation.
Notice that, since $\omega\in E_3^{(n)}$,
%for every $-9\leq i \leq 10$,
%$H(T_{i-1}(V, h_n-C_1\log_2 n)) \geq \log n$ and
%$H(T_i(V, h_n-C_1\log_2 n)) \geq \log n$, so $x_i(V,h_n-C_1\log_2 n)$  is also a left $(\log n)$-extremum.
%Hence, we also have
$x_i=x_i(V,\log n)=x_i(V,h_n-C_1\log_2 n)$ for every $-9\leq i \leq 10$
by Remark \ref{Remark_x_i}.
We also introduce
\begin{equation}\label{eq_def_b_chapeau}
    \widehat b(n)
:=
    2\lfloor b_{\log n}/2\rfloor
    +\un_{2\N+1}(n),
\end{equation}
which belongs  to $\{b_{\log n}-1,\, b_{\log n},\, b_{\log n}+1\}$ and has the same parity as $n$.
We define
\begin{equation}\label{eq_def_Mpm}
    M^-
:=
    \left\{
    \begin{array}{ll}
    x_{-1} & \text{ if } b_{\log n}\leq 0,\\
    x_0      &\text{ if } b_{\log n}> 0,
    \end{array}
     \right.
\qquad
    M^+
:=
    \left\{
    \begin{array}{ll}
    x_1 & \text{ if } b_{\log n}\leq 0,\\
    x_2      &\text{ if } b_{\log n}> 0.
    \end{array}
    \right.
\end{equation}
Since $b_{\log n}=x_0$ when $b_{\log n}\leq 0$
and $b_{\log n}=x_1$ when $b_{\log n}> 0$,
$M^-$ and $M^+$ are the two left $(\log n)$-maxima surrounding $b_{\log n}$, respectively on its left and on its right.
For this reason, $[M^-,M^+]$ is called the $(\log n)$-{\it central valley}
(see Figure \ref{figure_Good_Env_Coupling});
also $0\in [M^-,M^+]$.

Similarly as in Brox \cite{Brox} and Andreoletti et al. \cite{AndreolettiDevulder} for diffusions in a random environment,
and as in Devulder et al. \cite{DGP_Collision_Sinai} and \cite{DGP_Collision_Transient} for RWRE,
but with some adaptations,
we use a coupling between $S=(S_k)_k$ \big(under $P_\omega^{\widehat b(n)}$\big) and a reflected RWRE $\widehat S$ defined below.
To this aim, we define, for fixed $n$, $\big(\widehat \omega_x\big)_{x\in\Z}$ as follows:
$$
    \widehat \omega_{M^-}:=1,
\qquad
    \widehat \omega_x:=\omega_x \textnormal{ if }x\notin\{M^-,M^+\},
\qquad
    \widehat \omega_{M^+}:=0.
$$
We can now introduce, for fixed $\omega$ and $n$, a random walk $\widehat S:=\big(\widehat S_k\big)_{k\in\N}$
in the environment $\widehat \omega:=\big(\widehat \omega_x\big)_{M^-\leq x \leq M^+}$,
starting from $y\in[M^-,M^+]$,
and denote its law by $P_{\widehat \omega}^y$.
So, $\widehat S$ satisfies \eqref{eq_Def_Sn} with $\omega$ and $S$ replaced respectively by $\widehat \omega$ and $\widehat S$.
In words, $\widehat S$ is a random walk in the environment $\omega$, starting from $y\in[M^-,M^+]$,
and reflected at $M^-$ and $M^+$.
We also define the measure $\widehat \mu_n$ on $\Z$ by
\begin{equation*}
    \widehat \mu_n (M^-)
:=
    e^{-V(M^-)},
\qquad
    \widehat \mu_n (M^+)
:=
    e^{-V(M^+-1)},
\quad
\end{equation*}
\begin{equation}
\label{eq_def_mi_chapeau}
    \widehat \mu_n (x)
:=
    e^{-V(x)}+e^{-V(x-1)},
\quad
    M^-<x<M^+,
%\qquad
%    \widehat \mu_n(x)=0,
%\quad
%    x\notin [M^-,M^+].
\end{equation}
and $\widehat \mu_n(x):=0$ for $x\notin [M^-,M^+]$
(where $\widehat \mu_n(x)$ denotes $\widehat \mu_n(\{x\})$ for simplicity).

Observe that for fixed $n$ and $\omega$, $\widehat \mu_n(.)/\widehat \mu_n(\Z)$ is an invariant probability measure for
$\widehat S$.

Consequently, similarly as in
%(\cite{DGP_Collision_Sinai} after the proof of Lemma 5.7),
(\cite{DGP_Collision_Transient} eq. (55)),
for every fixed $n$ and $\omega$,
the measure $\widehat \nu_n$ defined by
\begin{equation}\label{eq_def_nu_hat}
    \widehat \nu(x)
:=
    \widehat \nu_n(x)
:=
\left\{
\begin{array}{ll}
    \widehat \mu_n(x){\bf 1}_{2\Z}(x)/\widehat\mu_n(2\Z) & \text{ if } n\in(2\N),\\
    \widehat \mu_n(x){\bf 1}_{2\Z+1}(x)/\widehat\mu_n(2\Z+1) & \text{ if } n\in(2\N+1),\\
\end{array}
\right.
\qquad
    x\in\Z,
\end{equation}
is an invariant probability measure for $\big(\widehat S_{2k}\big)_{k\in\N}$.
%for fixed $\widehat \omega$.
This means that
$
    P_{\widehat \omega}^{\widehat \nu}\big(\widehat S_{2k}=x\big)
=
    \widehat \nu(x)
$
for all $x\in\Z$ and $k\in\N$, where
$
    P_{\widehat \omega}^{\widehat \nu}\big(.\big)
:=
    \sum_{y\in\Z}\widehat \nu(y) P_{\widehat \omega}^{y}(.)
$.
Observe that $\widehat \omega$, $\widehat S$, $\widehat \mu_n$, $\widehat \nu_n$ and some other notation of this subsection defined below,
depend on $M^-$ and $M^+$ and so on $n$ and $\omega$,
but we often do not write the subscript $n$ in the following to simplify the notation.

We now have all the ingredients to build, for fixed $n$ and $\omega$, our coupling $Q_\omega$ of $S$ and $\widehat S$ as follows
and similarly as in
%\cite{DGP_Collision_Sinai}:
(\cite{DGP_Collision_Transient} around eq. (56)):
\begin{equation}\label{eq_def_Q_omega}
    Q_\omega\big(\widehat S\in.\big)
=
    P_{\widehat \omega}^{\widehat \nu}\big(\widehat S\in.\big),
\qquad
    Q_\omega\big(S\in.\big)
=
    P_\omega^{\widehat b(n)}\big(S\in.\big),
\end{equation}
so that under $Q_\omega$, the two Markov chains $\widehat S$ and $S$ move independently until
$$
    \tau_{\widehat S=S}
:=
    \inf\big\{\ell\geq 0,\ \widehat S_\ell=S_\ell\big\},
$$
which is their first meeting time,
then $\widehat S_k=S_k$ for all $\tau_{\widehat S=S}\leq k< \tau_{\text{exit}}$,
where
$$
    \tau_{\text{exit}}
:=
    \inf\big\{\ell>\tau_{\widehat S=S},\ S_\ell\notin[M^-,M^+]\big\}
$$
is the first exit time of $S$ from the central valley $[M^-,M^+]$ after the meeting time $\tau_{\widehat S=S}$,
and then $\widehat S$ and $S$ move independently again after $ \tau_{\text{exit}}$.

\subsection{Approximation of the quenched probability measure}

The next step is to prove that, under $Q_\omega$, $\widehat S$ and $S$ meet quickly,
and more precisely that
$\tau_{\widehat S=S}\leq n/10$ with large probability.
For this purpose, we define, for $n\geq n_3$,
in view of $E_3^{(n)}$,
\begin{eqnarray}
\label{eq_def_L-_hat}
    \widehat L^-
& := &
    \max\{k\leq b_{\log n},\ V(k)-V(b_{\log n}) \geq h_n\},
\\
    \widehat L^+
& := &
    \min\{k\geq b_{\log n},\ V(k)-V(b_{\log n}) \geq h_n\}.
\label{eq_def_L+_hat}
\end{eqnarray}
Loosely speaking, $\widehat L^-$ and $\widehat L^+$ are useful because $V\big(\widehat L^\pm\big)-V(b_{\log n})$
%the increase of potential of these points compared to $b_{\log n}$,
is approximatively $h_n$ and then
is quite lower than $\log n$, so $\widehat L^-$ and $\widehat L^+$ will  be hit quickly by $S$ under $Q_\omega$
(see Lemma \ref{Lemma_Meeting_S_and_Shat} below),
but $V\big(\widehat L^\pm\big)-V(b_{\log n})$ is also chosen quite large
because the invariant measure $\widehat \nu$ outside of
$\big[\widehat L^-, \widehat L^+\big]$ needs to be small (see Lemma \ref{Lemma_nu_loin_du_centre}).
We introduce the notation $u \vee v:=\max(u,v)$.
We prove the three following lemmas, which are uniform on $z$ since they do not depend on $z$.

\begin{lem}\label{Lemma_Meeting_S_and_Shat}
%There exists $n_4\geq n_3$ such that,
%If $n$ is large enough, we have
%For $n\geq n_3$,
We have,
with $\tau(.)$ denoting the hitting times by $S$ as before,
\begin{equation*}
%\label{}
%    \forall n\geq n_4,\
    \forall n\geq n_3,\
   \forall \omega\in
   %E_3^{(n)}\capE_{-}^{(n)}\cap
   E_3^{(n)}\cap E_5^{(n)},
\qquad
    Q_\omega\big[\tau\big(\widehat L^-\big) \vee \tau\big(\widehat L^+\big) > n/10\big]
\leq
    (\log n)^{-3}.
\end{equation*}
\end{lem}

\noindent{\bf Proof:}
Assume that $n\geq n_3$ and  $\omega\in
%E_3^{(n)}\capE_{-}^{(n)} \cap
E_3^{(n)}\cap E_5^{(n)}$.
Since
$
    V(M^\pm)-V(b_{\log n})
\geq
    \log n+C_2\log_2 n
>
    h_n+C_0
\geq
    V\big(\widehat L^+\big)-V(b_{\log n})
$
by $E_3^{(n)}$ (see also Remark \ref{Remark_x_i})
and using ellipticity \eqref{eq_ellipticity_for_V},
and since $h_n>0$ by definition of $n_3$,
we have $b_{\log n}< \widehat L^+< M^+$.
Moreover,
$$
    \max_{M^-\leq \ell\leq k\leq \widehat L^+, k\geq \widehat b(n)}[V(k)-V(\ell)]
\leq
    \max_{[\widehat b(n), \widehat L^+]} V
    -
    \min_{[M^-, \widehat L^+]}V
\leq
    V\big(\widehat L^+\big)-V\big(b_{\log n}\big)
\leq
    h_n+\log(\e_0^{-1})
$$
by ellipticity, i.e. by
%\eqref{eqEllipticity}
\eqref{eq_ellipticity_for_V},
and because $[M^-, M^+]$ is the $(\log n)$-central valley, its bottom being $b_{\log n}$.
Consequently, using \eqref{InegEsperance1} and Markov's inequality
since $M^-< \widehat b(n)<\widehat L^+$
because
$
    V\big(M^-\big)
>
    V\big(\widehat L^+\big)
\geq
    V\big(b_{\log n}\big)+h_n
\geq
    V\big(b_{\log n}\big)+3C_0
>
    V\big(\widehat b(n)\big)
$,
then
$
    \big[M^-, \widehat L^+\big]
\subset
    [x_{-1}, x_2[
\subset
    [-(\log n)^3,(\log n)^3[
$
because $\omega\in E_5^{(n)}$ and $\delta_1\in]0, 2/3[$,
this leads to
\begin{eqnarray*}
    P_\omega^{\widehat b(n)}\big[\tau(M^-)\wedge \tau\big(\widehat L^+\big)>n/10\big]
& \leq &
    10 n^{-1} \e_0^{-1}(2(\log n)^3)^2\e_0^{-1}e^{h_n}
\\
& = &
    40 \e_0^{-2}(\log n)^{6-C_1}
\leq
   (\log n)^{-3}/4,
\end{eqnarray*}
%uniformly on $E_3^{(n)}\cap E_5^{(n)}$ for large $n$,
since $n\geq n_3$ and $C_1>20$.
Moreover, applying \eqref{probaatteinte}, then $\omega\in E_5^{(n)}$ and the definition of
$\widehat L^+$ and finally using $V(M^\pm)-V(b_{\log n})\geq \log n+C_2\log_2 n$ on $E_3^{(n)}$ as before,
\begin{eqnarray*}
    P_\omega^{\widehat b(n)}\big[\tau(M^-)<\tau\big(\widehat L^+\big)\big]
& \leq &
    \big[\widehat L^+-\widehat b(n)\big]
    \exp\Big[\max_{[\widehat b(n),\widehat L^+-1]}V-V(M^-)\Big]
\\
& \leq &
    2(\log n)^3
    \exp\big[V(b_{\log n})+h_n-(V(b_{\log n})+\log n)\big]
\\
& \leq &
    2(\log n)^{3-C_1}
\leq
%    2(\log n)^{-6},
   (\log n)^{-3}/4,
\end{eqnarray*}
%uniformly on $E_3^{(n)}\cap E_5^{(n)}$ for large $n$,
since $n\geq n_3$ and $C_1>20$.
As a consequence, using \eqref{eq_def_Q_omega},
\begin{eqnarray}
    Q_\omega\big[\tau\big(\widehat L^+\big)>n/10\big]
& = &
    P_\omega^{\widehat b(n)}\big[\tau\big(\widehat L^+\big)>n/10\big]
\nonumber\\
& \leq &
    P_\omega^{\widehat b(n)}\big[\tau(M^-)<\tau\big(\widehat L^+\big)\big]
    +
    P_\omega^{\widehat b(n)}\big[\tau(M^-)\wedge \tau\big(\widehat L^+\big)>n/10\big]
\nonumber\\
& \leq &
    (\log n)^{-3}/2.
\label{Ineg_Proba_Q_Atteinte_L+}
\end{eqnarray}
%once more because $n\geq n_3$ and $C_1>20$.
%uniformly on $E_3^{(n)}\cap E_5^{(n)}$ for large $n$.
We prove similarly that
$
    Q_\omega\big[\tau\big(\widehat L^-\big)>n/10\big]
\leq
    (\log n)^{-3}/2
$
for all $n\geq n_3$ and $\omega \in E_3^{(n)}\cap E_5^{(n)}$,
using \eqref{InegEsperance2} instead of \eqref{InegEsperance1}.
% and $M^+$ instead of $M^-$.
This, together with \eqref{Ineg_Proba_Q_Atteinte_L+}, proves Lemma \ref{Lemma_Meeting_S_and_Shat}.
\hfill$\Box$

We now prove that the invariant measure outside $\big]\widehat L^-, \widehat L^+\big[$ is small for $n\geq n_3$.

\begin{lem}\label{Lemma_nu_loin_du_centre}
%For large $n$,
%There exists $n_5\geq n_3$ such that,
We have,
\begin{equation}\label{eq_lemma_nu_loin_du_centre}
%    \forall n\geq n_5, \
    \forall n\geq n_3, \
   \forall \omega\in E_3^{(n)}\cap E_5^{(n)},
\qquad
    \widehat \nu\big(\big[M^-,\widehat L^-\big]\big)
    +
    \widehat \nu\big(\big[\widehat L^+,M^+\big]\big)
\leq
    (\log n)^{-4}.
\end{equation}
\end{lem}

\noindent{\bf Proof:}
Let $n\geq n_3$ and $\omega\in E_3^{(n)}\cap E_5^{(n)}$.
As explained in Remark \ref{Remark_x_i},
%before \eqref{eq_def_b_chapeau},
due to $E_3^{(n)}$,
%, for $i\in\{-1,0,1\}$, we have
%$H(T_{i-1}(V, h_n-C_1\log_2 n) \geq \log n$ and
%$H(T_i(V, h_n-C_1\log_2 n) \geq \log n$, so
%$x_i(V,h_n-C_1\log_2 n)$ is a $(\log n)$-left extremum.
$x_i(V,h_n-C_1\log_2 n)=x_i(V,\log n)=x_i$ for every $i\in\{-1,0,1,2\}$.
So when $b_{\log n}\leq 0$, there is no left $(h_n-C_1\log_2 n)$-extremum in
$]x_0(V,h_n-C_1\log_2 n), x_1(V,h_n-C_1\log_2 n )[=]x_0,x_1[=]b_{\log n},M^+[$.
Similarly, when $b_{\log n}>0$, there is no left  $(h_n-C_1\log_2 n)$-extremum in $]x_1,x_2[=]b_{\log n},M^+[$.

%there are at maximum two $(h_n -C_1\log_2 n)$-left extrema
%in $]b_{\log n},x_1[=]x_0,x_1[$ (we could even say at least one).
%Indeed, there would otherwise be at least four, since $h$-left minima and $h$-left maxima alternate,
%and then five $(h_n -C_1\log_2 n)$-left slopes included in $[x_1,x_2]$,
%among which there would be, due to $E_0$ at least two consecutive $(h_n+\xi_1\log_2 n)$-left slopes,
%so there would be at least one $(h_n+\xi_1\log_2 n)$-extremum in $]x_1,x_2[$, which is not possible.

We first prove that
\begin{equation}\label{eq_Minimum_V_Couplage_Meeting}
    \min_{[\widehat L^+,M^+]}V
\geq
    V(b_{\log n})+C_1\log_2 n.
\end{equation}
Assume that $\min_{[\widehat L^+,M^+]}V<V(b_{\log n})+C_1\log_2 n$,
and let $u\in\big[\widehat L^+,M^+\big]$ be such that $V(u)=\min_{[\widehat L^+,M^+]}V$,
and $y:=\min\{\ell\in[b_{\log n}, u], V(\ell)=\max_{[b_{\log n}, u]}V\}$, so $y\geq \widehat L^+$.
Notice that  $V(y)\geq V\big(\widehat L^+\big)\geq V(b_{\log n})+h_n$
and $V(y)\geq V(b_{\log n})+h_n\geq V(u)-C_1\log_2 n+h_n$,
so $y$ would be a left $(h_n-C_1\log_2 n)$-maximum for $V$.
Since $b_{\log n}<y<u\leq M^+$, this contradicts the remark before \eqref{eq_Minimum_V_Couplage_Meeting}.
So, \eqref{eq_Minimum_V_Couplage_Meeting} is true.
We prove similarly that
\begin{equation}\label{eq_Minimum_V_Couplage_Meeting_2}
    \min_{[M^-, \widehat L^-]}V
\geq
    V(b_{\log n})+C_1\log_2 n.
\end{equation}
We have by \eqref{eq_Minimum_V_Couplage_Meeting} and since $\omega\in E_5^{(n)}$
and $\widehat\mu_n(2\Z)=\widehat\mu_n(2\Z+1)=\sum_{i=M^-}^{M^+-1}e^{-V(i)}\geq e^{-V(b_{\log n})}$,
\begin{eqnarray}
    \widehat \nu\big(\big[\widehat L^+,M^+\big]\big)
& \leq &
    \big[M^+-\widehat L^++1\big]\max_{x\in[\widehat L^+,M^+]}\big(e^{-V(x)}+e^{-V(x-1)}\big)
    e^{V(b_{\log n})}
\nonumber\\
& \leq &
    3(\log n)^3
    \big(1+\e_0^{-1}\big)
    (\log n)^{-C_1}
\leq
    (\log n)^{-4}/2
\label{Ineg_Mesure_Inv_Droite}
\end{eqnarray}
%uniformly on $E_3^{(n)}\cap E_5^{(n)}$ for large $n$
since $n\geq n_3$ and $C_1>20$,
and where we used
$-V(x-1)\leq -V(x)+\log(\e_0^{-1})$, $x\in\Z$
%$|V(x)-V(x-1)|\leq \log(\e_0^{-1})$, $x\in\Z$
by \eqref{eq_ellipticity_for_V}.
We prove similarly that
$
    \widehat \nu\big(\big[M^-,\widehat L^-\big]\big)
\leq
    (\log n)^{-4}/2
$
%uniformly on $E_3^{(n)}\cap E_5^{(n)}$ for large $n$
for all $n\geq n_3$ and $\omega \in E_3^{(n)}\cap E_5^{(n)}$
thanks to \eqref{eq_Minimum_V_Couplage_Meeting_2}.
This, together with \eqref{Ineg_Mesure_Inv_Droite}
proves \eqref{eq_lemma_nu_loin_du_centre}.
\hfill$\Box$

We can now prove that, with large enough probability, the coupling (i.e. $\widehat S=S$)
occurs quickly, and lasts at least until time $n$.

\begin{lem}\label{Lemma_Time_Meeting_Couplage}
%There exists $n_{60}\geq n_5$ such that
We have,
\begin{equation}\label{eq_Time_meeting_1}
    \forall n\geq n_3,\,
    \forall \omega\in E_3^{(n)}\cap E_5^{(n)},
\qquad
    Q_\omega\big[\tau_{\widehat S=S}>n/10\big]
\leq
    2(\log n)^{-3},
\end{equation}
and
\begin{equation}
    \forall n\geq n_3,\ \forall \omega\in E_3^{(n)},
\qquad
    Q_\omega\big[\tau_{\textnormal{exit}}\leq n\big]
\leq
    (\log n)^{-3}.
%    (\log n)^{-7}.
\label{eq_Time_meeting_exit}
\end{equation}
\end{lem}

\noindent{\bf Proof:}
Let $n\geq n_3$, and $\omega\in E_3^{(n)}\cap E_5^{(n)}$.
We have by Lemma \ref{Lemma_Meeting_S_and_Shat},
\begin{eqnarray*}
&&
    Q_\omega\big[\tau_{\widehat S=S}>n/10\big]
\\
& \leq &
    Q_\omega\big[\tau\big(\widehat L^-\big) \vee \tau\big(\widehat L^+\big) <\tau_{\widehat S=S}\big]
    +
    Q_\omega\big[\tau\big(\widehat L^-\big) \vee \tau\big(\widehat L^+\big) >n/10\big]
\\
& \leq &
    Q_\omega\big[\tau\big(\widehat L^-\big)  <\tau_{\widehat S=S},\ \widehat S_0<\widehat b(n)\big]
+
    Q_\omega\big[\tau\big(\widehat L^+\big)  <\tau_{\widehat S=S},\ \widehat S_0 \geq \widehat b(n)\big]
+
    (\log n)^{-3}.
\end{eqnarray*}
Now, observe that a.s. under $Q_\omega$, $S_0=\widehat b(n)$  by \eqref{eq_def_Q_omega} and has the same parity as $n$ by \eqref{eq_def_b_chapeau},
and $\widehat S_0$ also has the same parity as $n$ by \eqref{eq_def_Q_omega} and \eqref{eq_def_nu_hat}.
Hence the process $\big(\widehat S_k-S_k\big)_{k\in\N}$ starts at $\big(\widehat S_0-\widehat b(n)\big)\in(2\Z)$,
and it only makes jumps belonging to $\{-2,0,2\}$, so up to time $\tau_{\widehat S=S}-1$ it is $<0$
(resp. $>0$)  on  $\big\{\widehat S_0<\widehat b(n)\big\}$
\big(resp. on  $\big\{\widehat S_0>\widehat b(n)\big\}$\big),
and in particular at time $\tau\big(\widehat L^-\big)$
on $\big\{\tau\big(\widehat L^-\big)<\tau_{\widehat S=S}, \ \widehat S_0<\widehat b(n)\big\}$
\big(resp. at time $\tau\big(\widehat L^+\big)$
on
$
    \big\{\tau\big(\widehat L^+\big)<\tau_{\widehat S=S}, \ \widehat S_0>\widehat b(n)\big\}
=
    \big\{\tau\big(\widehat L^+\big)<\tau_{\widehat S=S}, \ \widehat S_0\geq\widehat b(n)\big\}
$; for
the last equality, notice that $\tau_{\widehat S=S}=0$ on $\big\{\widehat S_0=\widehat b(n)\big\}$
\big).
So,
\begin{eqnarray}
&&
    Q_\omega\big[\tau_{\widehat S=S}>n/10\big]
\nonumber\\
& \leq &
    Q_\omega\big[\tau\big(\widehat L^-\big)  <\tau_{\widehat S=S},\ \widehat S_{\tau(\widehat L^-)}<\widehat L^-\big]
+
    Q_\omega\big[\tau\big(\widehat L^+\big)  <\tau_{\widehat S=S},\ \widehat S_{\tau(\widehat L^+)} >  \widehat L^+\big]
+
    (\log n)^{-3}
\nonumber\\
& \leq &
    Q_\omega\big[\tau\big(\widehat L^-\big)  <\tau_{\widehat S=S},\ \widehat  S_{2\lfloor \tau(\widehat L^-)/2\rfloor}\leq \widehat L^-\big]
+
    Q_\omega\big[\tau\big(\widehat L^+\big)  <\tau_{\widehat S=S},\ \widehat S_{2\lfloor \tau(\widehat L^+)/2\rfloor} \geq  \widehat L^+\big]
\nonumber\\
&&
+
    (\log n)^{-3}
\nonumber\\
& \leq &
        \widehat \nu\big(\big[M^-,\widehat L^-\big]\big)
    +
    \widehat \nu\big(\big[\widehat L^+,M^+\big]\big)
+
    (\log n)^{-3}.
\label{Ineg_Q_avec_nu_hat}
\end{eqnarray}
Indeed, the last inequality is a consequence of the fact that
$
    Q_\omega\big(\widehat S_{2k}=x\big)
=
    P_{\widehat \omega}^{\widehat \nu}\big(\widehat S_{2k}=x\big)
=
    \widehat\nu(x)
$
for all $x\in\Z$ and all (deterministic) $k\in\N$
(see \eqref{eq_def_Q_omega} and the explanations after \eqref{eq_def_nu_hat}),
and from the independence of $\widehat S$ with $S$ (and its hitting times $\tau(.)$) up to time $\tau_{\widehat S=S}$.
Hence, \eqref{Ineg_Q_avec_nu_hat} together with Lemma \ref{Lemma_nu_loin_du_centre} prove \eqref{eq_Time_meeting_1}.

Finally, \eqref{eq_def_Q_omega}, followed by \eqref{InegProba1} and \eqref{InegProba2}
give for every $n\geq n_3$ and $\omega\in E_3^{(n)}$,
\begin{align}
    Q_\omega\big[\tau_{\textnormal{exit}}\leq n\big]
& \leq
    Q_\omega[\tau(M^-)\wedge\tau(M^+)\leq n]
=
    P_\omega^{\widehat b(n)}[\tau(M^-)\wedge\tau(M^+)\leq n]
%    ~~~~~~~~~~~~~~~~~~~~~~~~
\nonumber\\
& \leq
    P_\omega^{\widehat b(n)}[\tau(M^-)\leq n]
    +
    P_\omega^{\widehat b(n)}[\tau(M^+)\leq n]
\nonumber\\
& \leq
    2(n+1)\e_0^{-2}  \exp[-(\log n+C_2 \log_2 n)]
\leq
    4\e_0^{-2}(\log n)^{-C_2}
\leq
    (\log n)^{-3},
\label{Ineg_Sortie_M_pm}
\end{align}
%uniformly for large $n$ on $E_3^{(n)}$,
since $\min_{[M^-, b_{\log n}]}V=\min_{[b_{\log n}, M^+]}V=V(b_{\log n})$,
$V(M^{\pm})-V(b_{\log n})\geq \log n+C_2 \log_2 n$ on $E_3^{(n)}$,
$|b_{\log n}-\widehat b(n)|\leq 1$,
 $|V(u)-V(u-1)|\leq \log(\e_0^{-1})$ for $u\in\Z$
by \eqref{eq_ellipticity_for_V}, $\log n>2\e_0^{-1}$ since $n\geq n_3\geq 3$, and $C_2>9$.
This proves \eqref{eq_Time_meeting_exit}.
%for some $n_{60}\geq n_5$.
\hfill$\Box$

\medskip

%Also, the following lemma says, roughly speaking, that
%either the invariant measure $\widehat \nu_n(z)$ of $z$ is small,
%either, with large probability, $(S_k)_k$ goes quickly to $\widehat b(n)$, which is nearly the bottom of the $(\log n)$-central %valley.
%This lemma will be useful to prove Lemma \ref{Lemma_Approx_PQuenched_Nu}
%(see \eqref{Ineg_Proba_Quenched_SinA_1}).

Also, the following lemma
will be useful to prove Lemma \ref{Lemma_Approx_PQuenched_Nu}
(see \eqref{Ineg_Proba_Quenched_SinA_1}).

\begin{lem}\label{Lemma_Produit_Nu_Proba}
%There exists $n_6\geq n_3$ such that
We have,
\begin{equation}\label{Ineg_Produit_Nu_Proba_Atteinte}
    \forall n\geq n_3, \
    \forall z\in\Z, \
    \forall \omega\in E_3^{(n)}\cap E_4^{(n)}(z)\cap E_5^{(n)},
\qquad
    \widehat \nu_n(z)
    \po\big[\tau(\widehat b(n))\geq n/10\big]
\leq
    (\log n)^{-3}.
\end{equation}
\end{lem}

\noindent{\bf Proof:}
Let $n\geq n_3$, $z\in\Z$ and $\omega\in E_3^{(n)}\cap E_4^{(n)}(z)\cap E_5^{(n)}$.
We treat separately the three different cases defining $E_4^{(n)}(z)$.

\noindent{\bf First case:} if in addition $\omega\in\{V(z)-V(b_{\log n})\geq 5\log_2 n\}$, we have
by ellipticity,
$$
    \widehat \nu_n(z)
\leq
    \big(e^{-V(z)}+e^{-V(z-1)}\big)e^{V(b_{\log n})}
\leq
    (1+\e_0^{-1})e^{-[V(z)-V(b_{\log n})]}
\leq
    2\e_0^{-1}(\log n)^{-5}
\leq
    (\log n)^{-3}
$$
%uniformly on  $z$ and on such $\omega$ for large $n$,
since $n\geq n_3$,
which proves \eqref{Ineg_Produit_Nu_Proba_Atteinte} in this case.

{\bf Second case:} if
$\omega\in E_-^{(n)}\cap \big\{\max_{[b_{\log n}, 0]}V<V(x_1)-9\log_2 n\big\}$,
we have $b_{\log n}=x_0\leq 0$
and
either $\widehat b(n)=1$, or
$-(\log n)^3-1\leq \widehat b(n)\leq 0< x_1<x_2$ since $\omega\in E_5^{(n)}$ and
$\big|\widehat b(n)-b_{\log n}\big|\leq 1$.
We start with  this second sub-case $\widehat b(n)\leq 0$. We have by \eqref{probaatteinte},
\begin{eqnarray}
    \po\big[\tau(x_2)<\tau\big(\widehat b(n)\big)\big]
& \leq &
    \big(\big|\widehat b(n)\big|+1\big)\exp\big[\max\nolimits_{[\widehat b(n), 0]}V-V(x_1)\big]
\nonumber\\
& \leq &
    2(\log n)^3\exp\big[-9\log_2 n+\log\e_0^{-1}\big]
\leq
    (\log n)^{-4}
\label{Ineg_Proba_Atteinte_b_avant_x2}
\end{eqnarray}
%uniformly for
%for every $\omega$ of this second subcase
%for every $n$ larger than some $n_6'\geq n_3$
%large $n$
by ellipticity since $\big|\widehat b(n)- b_{\log n}\big|\leq 1$
and because $n\geq n_3$ so $\log n\geq 2\e_0^{-1}\geq 2$.

Also, by Lemma \ref{Lemma_Proba_Descente} applied with
$\xi_1=2C_1>19$, $\xi_2=1/10$, $\alpha=3$,
$a=\widehat b(n)\leq 0<b=x_1<c=x_2$,
$h=\log n>\widehat h_2(2C_1,1/10)$ because $n\geq n_3$,
and $x=0$,
since
its hypothesis {\bf (i)} is satisfied because $b_{\log n}\leq 0$ and so $x_1$ is a left $(\log n)$-maximum,
and there is no
left $(h-\xi_1\log h)=(h_n-C_1\log_2 n)$-extremum in $]x_1,x_2[$ nor in $]x_0,x_1[$ by $E_3^{(n)}$
(as explained after \eqref{eq_lemma_nu_loin_du_centre}
since $x_0, x_1$ and $x_2$ are consecutive left $(h_n-C_1\log_2 n)$-extrema)
and so hypotheses {\bf (ii)}
and {\bf (iii)} of this lemma are satisfied (e.g. if {\bf (ii)} was not satisfied,
there would be a left $(h_n-C_1\log_2 n)$-maximum in $]x_1,x_2[$),
and hypothesis {\bf (iv)} is satisfied with $\alpha=3$ thanks to $E_5^{(n)}$ and $\delta_1<1$,
so
%there exists $n_6''\geq n_3$ such that for every $n\geq n_6''$,
\begin{equation}\label{Ineg_Proba_Atteinte_b_x2}
    \forall \omega\in E_3^{(n)}\cap E_5^{(n)},
\quad
    (b_{\log n}\leq 0  \text{ and } \widehat b(n)\leq 0)
\Rightarrow
    \po\big[\tau\big(\widehat b(n)\big)\wedge \tau(x_2)\geq n/10\big]
\leq
    (\log n)^{-4}.
\end{equation}
This and \eqref{Ineg_Proba_Atteinte_b_avant_x2} lead to
$
    \po\big[\tau\big(\widehat b(n)\big)\geq n/10\big]
\leq
    2(\log n)^{-4}
\leq
    (\log n)^{-3}
$
for every $\omega$ of this second subcase
since $n\geq n_3$.
%for $n\geq \max(n_6', n_6'')$.

%This inequality is true for any $\omega\in E_3^{(n)}\cap E_5^{(n)}$
%such that $b_{\log n}\leq 0$ and $\widehat b(n)\leq 0$

We now turn to the other subcase, that is, we assume that $\widehat b(n)=1$.
Then, $b_{\log n}=x_0=0$ since $\omega\in E_-^{(n)}$
and $\big|\widehat b(n)-b_{\log n}\big|\leq 1$.
In this subcase we have, using \eqref{probaatteinte}, Markov inequality and \eqref{InegEsperance1}
in the second inequality,
\begin{eqnarray}
&&
    \po\big[\tau\big(\widehat b(n)\big)\geq n/10\big]
=
    \po\big[\tau(1)\geq n/10\big]
\nonumber\\
& \leq &
    \po\big[\tau(x_{-1})< \tau(1)\big]
+
    \po\big[\tau(x_{-1})\wedge \tau(1)\geq n/10\big]
%    \po\big[\tau(x_{-1})< \tau\big(\widehat b(n)\big)\big]
%+
%    \po\big[\tau(x_{-1})\wedge \tau\big(\widehat b(n)\big)\geq n/10\big]
\nonumber\\
& \leq &
    \exp[V(x_0)-V(x_{-1})]
+
    10 n^{-1}\e_0^{-1}(1-x_{-1})^2   \exp[V(0)-\min\nolimits_{[x_{-1},0]}V]
\nonumber\\
& \leq &
    n^{-1}
    %(\log n)^{-C_2}
+
    40\e_0^{-1}n^{-1} (\log n)^{6}
\leq
    (\log n)^{-4}
\label{Ineg_Proba_Atteinte_b_chapeau_si_egal_a_1}
\end{eqnarray}
for every $\omega$ of this subcase
since $n\geq n_3$,
%for every $n$ large than some $n_6'''\geq n_3$,
%uniformly for large $n$
    $H[T_0(V,\log n)]=V(x_{-1})-V(x_0)\geq \log n$,
%+C_2\log_2 n$ since $\omega\in E_3^{(n)}$,
$|x_{-1}|\leq (\log n)^3$ since $\omega\in E_5^{(n)}$ and
$\min_{[x_{-1},0]}V=\min_{[x_{-1},x_0]}V=V(x_0)=V(0)=0$.
%, and $C_1>20$.
So, \eqref{Ineg_Produit_Nu_Proba_Atteinte} is proved in this second case (whenever $\widehat b(n)=1$ or not), since $\widehat \nu_n(z)\leq 1$,
for all $n\geq n_3$,
$z\in\Z$
and
$
    \omega
\in
    E_-^{(n)}\cap \big\{\max_{[b_{\log n}, 0]}V<V(x_1)-9\log_2 n\big\}
    \cap E_3^{(n)}
    %\cap E_4^{(n)}(z_n)
    \cap E_5^{(n)}
$.
%for every $n\geq \max(n_6', n_6'', n_6''')$.
%in particular, we choose $n_6\geq \max(n_6', n_6'', n_6''')$.

{\bf Third case:} finally, the proof is similar when $\omega\in E_+^{(n)}\cap \big\{\max_{[0,b_{\log n}]}V<V(x_0)-9\log_2 n\big\}$ with $x_{-1}$ instead of $x_2$ and
$x_1$ exchanged with $x_0$,
%and \eqref{InegEsperance2} instead of \eqref{InegEsperance1},
which ends the proof of the lemma.
\hfill$\Box$

We now have all the ingredients to approximate the quenched probability
$P_\omega(S_n=z)$ by the invariant probability measure $\widehat \nu_n(z)$
for $\omega\in E_C^{(n)}(z)$ (defined in \eqref{eq_def_EL}), uniformly for $n\geq n_3$
(recall that $P_\omega(S_n=z)$ and $\widehat \nu_n(z)$
are equal to $0$ if $z$ and $n$ do not have the same parity by \eqref{eq_def_nu_hat}).

\begin{lem}\label{Lemma_Approx_PQuenched_Nu}
We have,
\begin{equation}\label{Ineg_Comparaison_Pquenched_et_nuhat_NEW}
    \forall n\ge n_3,\,
    \forall z\in\Z,\,
    \forall \omega\in E_C^{(n)}(z),
\qquad
    \big|P_\omega(S_n=z)-\widehat \nu_n(z)\big|
\leq
    5(\log n)^{-3}.
\end{equation}
\end{lem}

\noindent{\bf Proof:}
Let $n\geq n_3$, $z\in\Z$ and $\omega\in E_C^{(n)}(z)$.
For $u\in\Z$, we define $V_u=V_u^+$ and $V_u^-$ by
$V_u(.):=V(u+.)-V(u)$ and $V_u^{\pm}(.):=V(u\pm.)-V(u)$.
Since $\omega\in E_6^{(n)}$,
$
    T_{V_{b_{\log n}}^{\pm}}(\log n)
>
    \Gamma_n$.
Also, $|b_{\log n}-z|\leq \Gamma_n$ because  $\omega\in E_7^{(n)}(z)$, so ($M^{\pm}$ being defined in \eqref{eq_def_Mpm}),
\begin{equation}\label{Ineg_zn_Mplus}
    z
\leq
    b_{\log n}+\Gamma_n
<
    b_{\log n}+T_{V_{b_{\log n}}^+}(\log n)
\leq
    M^+.
\end{equation}
Thus $z<M^+$, and similarly, $z>M^-$, and so $z\in]M^-, M^+[$.

Observe that for $k\in[n/10,n]\cap(2\N)$,
\begin{eqnarray}
    P_\omega^{\widehat b(n)}[S_k=z]
%& = &
=
    Q_\omega[S_k=z]
%\nonumber\\
& \geq &
    Q_\omega\big[S_k=z, \  \tau_{\widehat S=S}\leq n/10\leq k\leq n<\tau_{\textnormal{exit}}\big]
\nonumber\\
& = &
    Q_\omega\big[\widehat S_k=z, \  \tau_{\widehat S=S}\leq n/10\leq k\leq n<\tau_{\textnormal{exit}}\big]
\nonumber\\
& \geq  &
    Q_\omega\big[\widehat S_k=z\big] -Q_\omega\big[\tau_{\widehat S=S}> n/10\big]-Q_\omega\big[\tau_{\textnormal{exit}}\leq n\big]
\nonumber\\
& \geq  &
%    \widehat \nu(B) -2(\log n)^{-3}-Q_\omega\big[\tau_{exit}\leq n\big],
    \widehat \nu(z) -3(\log n)^{-3},
\label{Ineg_Pquenched_S_dans_A_1}
\end{eqnarray}
where we used \eqref{eq_def_Q_omega} in the first equality,
$S_k=\widehat S_k$ for $k\in\big[\tau_{\widehat S=S}, \tau_{\textnormal{exit}}\big[$  in the second one,
and
$
    Q_\omega\big[\widehat S_k=x\big]
=
    P_{\widehat \omega}^{\widehat \nu}\big[\widehat S_k=x\big]
=
    \widehat\nu(x)
$
for all $x\in\Z$ since $k$ is even (see \eqref{eq_def_Q_omega} and the remark after \eqref{eq_def_nu_hat}),
and Lemma \ref{Lemma_Time_Meeting_Couplage} in the last line
since $n\geq n_3$ and $\omega\in E_C^{(n)}(z)$.

Similarly, for every $k\in[n/10,n]\cap(2\N)$,
\begin{align}
%&&
    P_\omega^{\widehat b(n)}[S_k=z]
% =
%    Q_\omega[S_k=z]
%\nonumber\\
 \leq &
    Q_\omega\big[S_k=z, \  \tau_{\widehat S=S}\leq n/10, \ \tau_{\textnormal{exit}}> n\big]
    +Q_\omega\big[\tau_{\widehat S=S}> n/10\big]
    +Q_\omega\big[\tau_{\textnormal{exit}}\leq n\big]
%\nonumber\\
%& \leq &
%    Q_\omega\big[\widehat S_k=z, \  \tau_{\widehat S=S}\leq n/10, \ \tau_{\textnormal{exit}}> n\big]
%     +3(\log n)^{-3}
\nonumber\\
 \leq  &
    Q_\omega\big[\widehat S_k=z\big]
%    +Q_\omega\big[\tau_{\widehat S=S}> n/10\big]
%    +Q_\omega\big[\tau_{exit}\leq n\big]
     +3(\log n)^{-3}
=
    \widehat \nu(z) +3(\log n)^{-3}.
\label{Ineg_Pquenched_S_dans_A_2}
\end{align}
We have, applying the strong Markov property in the second line, %since $\omega\in E_4^{(n)}'z)$,
\begin{eqnarray}
    P_\omega[S_n=z]
& \geq &
    P_\omega\big[S_n=z,\ \tau\big(\widehat b(n)\big)<n/10\big]
\nonumber\\
& = &
    E_\omega\big[{\bf 1}_{\{ \tau(\widehat b(n))<n/10 \}} P_\omega^{\widehat b(n)} [S_k=z]_{\mid k=n-\tau(\widehat b(n))}\big]
\nonumber\\
& \geq &
    E_\omega\big[{\bf 1}_{\{ \tau(\widehat b(n))<n/10 \}} \big(\widehat \nu(z) -3(\log n)^{-3}\big)\big]
\nonumber\\
& \geq &
    \widehat \nu(z)
    -\widehat \nu(z)\po\big[\tau(\widehat b(n))\geq n/10\big]
    -3(\log n)^{-3}
\nonumber\\
& \geq &
    \widehat \nu(z)
    -4(\log n)^{-3}
\label{Ineg_Proba_Quenched_SinA_1}
\end{eqnarray}
where we used \eqref{Ineg_Pquenched_S_dans_A_1} in the second inequality
since $\big(n-\tau\big(\widehat b(n)\big)\big)\in[9n/10,n]\cap (2\N)$
because $\widehat b(n)$, and then $\tau\big(\widehat b(n)\big)$, has the same parity as $n$ by \eqref{eq_def_b_chapeau},
and Lemma \ref{Lemma_Produit_Nu_Proba} in the last inequality, since $n\geq n_3$ and $\omega\in E_C^{(n)}$.

Similarly, using \eqref{Ineg_Pquenched_S_dans_A_2} instead of \eqref{Ineg_Pquenched_S_dans_A_1}, we get
\begin{eqnarray}
    P_\omega\big[S_n=z,\ \tau\big(\widehat b(n)\big)<n/10\big]
%& = &
%    E_\omega\big[{\bf 1}_{\{ \tau(\widehat b(n))<n/10 \}} P_\omega^{\widehat b(n)} [S_k=z]_{\mid k=n-\tau(\widehat b(n))}\big]
%\nonumber\\
& \leq &
    E_\omega\big[{\bf 1}_{\{ \tau(\widehat b(n))<n/10 \}}
%    \big(\widehat \nu(B) +2(\log n)^{-3}+Q_\omega\big[\tau_{exit}\leq n\big]\big)
    \big(\widehat \nu(z) +3(\log n)^{-3}\big)
    \big]
\nonumber\\
& \leq &
%    \widehat \nu(B) +2(\log n)^{-3}+Q_\omega\big[\tau_{exit}\leq n\big]
    \widehat \nu(z) +3(\log n)^{-3}.
\label{Ineg_Proba_Quenched_SinA_2}
\end{eqnarray}

We now assume that $b_{\log n}\leq 0$, and so $b_{\log n}=x_0$ and $M^+=x_1$.
Also, we have once more
\begin{equation}\label{Ineg_Proba_Sortie_b0b1}
    \po\big[\tau\big(\widehat b(n)\big)\wedge\tau(x_2)\geq n/10\big]
\leq
    (\log n)^{-4}.
\end{equation}
Indeed this is proved in \eqref{Ineg_Proba_Atteinte_b_x2} when $\widehat b(n)\neq 1$
since $n\geq n_3$,
%as stated after \eqref{Ineg_Proba_Atteinte_b_chapeau_si_egal_a_1},
whereas when $\widehat b(n)=1$,
the left hand side of \eqref{Ineg_Proba_Sortie_b0b1} is equal to
$\po\big[\tau\big(\widehat b(n)\big)\geq n/10\big]$, which is $\leq (\log n)^{-4}$
by \eqref{Ineg_Proba_Atteinte_b_chapeau_si_egal_a_1} since
$\widehat b(n)=1<x_2$ in this case.

Moreover for $0\leq k\leq n$, using $z< M^+=x_1< x_2$ (see \eqref{Ineg_zn_Mplus}),
we have by \eqref{InegProba2} and ellipticity \eqref{eq_ellipticity_for_V}, and since
$
    V(x_1)-\min_{[x_1,x_2]}V
=
    H[T_1(V,\log n)]
=
    H[T_1(V,h_n-C_1 \log_2 n)]
\geq
    \log n+C_2\log_2 n
$
because $\omega\in E_3^{(n)}$ (see also Remark \ref{Remark_x_i}),
\begin{eqnarray*}
    \po^{x_2}(S_k=z)
& \leq &
    \po^{x_2}[\tau(x_1)< k]
\leq
    n\exp[-H(T_1[V,\log n)]+\log\e_0^{-1}]
\\
& \leq &
    \e_0^{-1}(\log n)^{-C_2}
\leq
    (\log n)^{-3}
\end{eqnarray*}
%uniformly \big(on $E_C^{(n)}$\big) for large $n$
since $C_2>9$ and $\log n>\e_0^{-1}$ because $n\geq n_3$.
Hence by the strong Markov property,
\begin{equation}\label{Ineg_Proba_Sortie_b0b1_b}
    \p_\omega\big[S_n=z,\ \tau(x_2)<n/10\big]
=
    \E_\omega\big[{\bf 1}_{\{\tau(x_2)<n/10\}}\po^{x_2}(S_k=z)_{|k=n-\tau(x_2)}\big]
\leq
    (\log n)^{-3}.
\end{equation}
Finally, \eqref{Ineg_Proba_Quenched_SinA_2}, \eqref{Ineg_Proba_Sortie_b0b1} and \eqref{Ineg_Proba_Sortie_b0b1_b} give
\begin{eqnarray}
    \po(S_n=z)
& \leq &
    \po\big[\tau\big(\widehat b(n)\big)\wedge\tau(x_2)\geq n/10\big]
+
    P_\omega\big[S_n=z,\ \tau\big(\widehat b(n)\big)<n/10\big]
\nonumber\\
&&
\qquad\qquad
+
    \p_\omega\big[S_n=z,\ \tau(x_2)<n/10\big]
\\
& \leq &
    \widehat \nu(z) +5(\log n)^{-3}.
\label{eq_105_bis}
\end{eqnarray}
We prove similarly this inequality $\po(S_n=z)\leq \widehat \nu(z) +5(\log n)^{-3}$
by symmetry when $b_{\log n}>0$, exchanging
$x_0$ and $x_1$ and replacing $x_2$ by $x_{-1}$
in
\eqref{Ineg_Proba_Sortie_b0b1} and \eqref{Ineg_Proba_Sortie_b0b1_b}
since $z> M^-=x_0>x_{-1}$ in this case, and using \eqref{InegProba1} instead of \eqref{InegProba2}.

Combining this with \eqref{eq_105_bis} and \eqref{Ineg_Proba_Quenched_SinA_1} proves \eqref{Ineg_Comparaison_Pquenched_et_nuhat_NEW}.
\hfill$\Box$

%%%%%%%%%%%%%%%%%%%%%%%%%%%%%%%%%%%%%%%%%%%%%%%%%%%%%%%%%%%%%%%%%%%%

\subsection{Upper bound of the annealed probability: main contribution}
The aim of this subsection is to give an upper bound of the annealed probability of $\{S_n=z\}$
on the event for which we used the coupling,
%when the coupling holds and $b_{\log n}$ is close to $z$,
that is, on $E_C^{(n)}(z)$.
More precisely, we prove the following estimate.

\begin{prop}\label{Ineg_Upper_Bound_ELT}
We have, under the hypotheses of Theorem \ref{Th_Local_Limit_Sinai}, as $n\to+\infty$,
\begin{equation}\label{Ineg_Minor_Sur_EL}
    \sup_{z\in(2\Z+n)}
    \bigg[
    \PP\big(S_n=z, E_C^{(n)}(z)\big)
-
    \frac{2\sigma^2}{(\log n)^2}\varphi_\infty\bigg(\frac{\sigma^2 z}{(\log n)^2}\bigg)
    \bigg]
\leq
    o\big((\log n)^{-2}\big).
\end{equation}
%\begin{equation}\label{Ineg_Minor_Sur_EL}
%    \PP\big(S_n=z, E_C^{(n)}\big)
%\leq
%    \frac{2\sigma^2}{(\log n)^2}\varphi_\infty\bigg(\frac{\sigma^2 z}{(\log n)^2}\bigg)
%    +o\big((\log n)^{-2}\big).
%\end{equation}
\end{prop}

The strategy of the proof is to use Lemma \ref{Lemma_Approx_PQuenched_Nu} to dominate $\PP\big(S_n=z, E_C^{(n)}(z)\big)$
by some quantity expressed in terms of left $(\log n)$-slopes $T_i(V,\log n)$ for $-1\leq i \leq 1$
(see e.g. \eqref{Ineg_Proba_J1_0_et_1_Plus_ET_Moins}, \eqref{eq_def_J3},
\eqref{eq_def_phi_v_t} and \eqref{eq_J2_Slopes}),
then use our Theorems \ref{Lemma_Independence_h_extrema} and \ref{Lemma_Central_Slope}
to obtain an expression with $\mathcal T_{V,\log n}^\uparrow$ and $\mathcal T_{V,\log n}^\downarrow$,
then Lemma \ref{Lemma_Proba_bh_egal}
to make appear the quantity $\p\big(b_{\log n}= z_n^+\big)$ for some $z_n^+\approx z$,
which, in turn, can be approximated by the
%left hand side of
expression with $\varphi_\infty$ in
\eqref{Ineg_Minor_Sur_EL} thanks to Theorem \ref{Th_Local_Limit_b_h}.

\noindent{\bf Proof:}
We assume that the hypotheses of Theorem \ref{Th_Local_Limit_Sinai} are satisfied.
%In particular, $S$ is a Sinai walk, and $n$ and $z_n$ have the same parity.
Let $n\geq n_3$ and $z\in(2\Z+n)$.
Using Lemma \ref{Lemma_Approx_PQuenched_Nu}
in the last line, $E_C^{(n)}(z)$ being defined in \eqref{eq_def_EL}, we have
%writing $E_i^{(n)}(z)$ instead of $E_i^{(n)}$ even for $i\neq 3$  for simplicity,
%thanks to \eqref{Ineg_zn_dans_vallee_centrale},
\begin{eqnarray}
    \PP\big(S_n=z, E_C^{(n)}(z)\big)
%& = &
%    \PP\big(S_n=z,
%    \cap_{i=2}^5 E_i^{(n)}(z)
%    %E_3^{(n)} \cap E_4^{(n)}(z) \cap E_5^{(n)}  \cap E_6^{(n)}
%    \cap \big\{|b_{\log n}-z|\leq \Gamma_n\big\}\big)
%\nonumber\\
%& = &
%    \sum_{k=-\lfloor (\log n)^{4/3+\delta_1}\rfloor}^{\lfloor (\log n)^{4/3+\delta_1}\rfloor}
%    \PP\big(S_n=z,  \cap_{i=2}^5 E_i^{(n)} \cap \{b_{\log n}=z+k\}\big).
%\nonumber\\
& = &
    \sum_{k=-\Gamma_n}^{\Gamma_n}
    \EE\big[\un_{
%    \cap_{i=2}^5 E_i^{(n)}(z)
%    E_3^{(n)} \cap E_4^{(n)}(z) \cap E_5^{(n)}  \cap E_6^{(n)}
    E_C^{(n)}(z)
    \cap \{b_{\log n}=z+k\}}\po(S_n=z)\big]
\nonumber\\
& = &
%    O\big((\log n)^{-3}\big)
    f_1(n,z)
    +
    \sum_{k=-\Gamma_n}^{\Gamma_n}
    J_0(k, n, z),
\label{Ineg_Proba_Zn_EL_1}
\end{eqnarray}
where
$|f_1(n,z)|\leq 5(\log n)^{-3}$ and (writing $E_i^{(n)}(z)$ instead of $E_i^{(n)}$ even for $i\neq 3$ for simplicity),
\begin{equation}\label{eq_def_J0}
    J_0(k, n, z)
:=
    \EE\big[\un_{E_C^{(n)}(z)\cap \{b_{\log n}=z+k\}}\widehat \nu_n(z)\big]
=
    \E\big[\un_{\cap_{i=3}^6 E_i^{(n)}(z)\cap \{b_{\log n}=z+k\}}\widehat \nu_n(z)\big].
\end{equation}
Notice that, using \eqref{Ineg_zn_Mplus} and in the remark below,
we have if $\omega \in E_C^{(n)}(z)\cap \{b_{\log n}=z+k\}$
with $|k|\leq \Gamma_n$
($M^\pm$ being defined in \eqref{eq_def_Mpm}),
\begin{equation}
    M^-
<
    z
=
    b_{\log n}-k
<
    M^+.
\label{Ineg_zn_dans_vallee_centrale}
\end{equation}
Hence, we have on $E_C^{(n)}(z)\cap \{b_{\log n}=z+k\}$ with $|k|\leq \Gamma_n$,
using the definitions of $\widehat \nu_n$ and $\widehat \mu_n$
(see \eqref{eq_def_nu_hat}  and \eqref{eq_def_mi_chapeau}),
$$
    \widehat \nu_n(z)
=
    \frac{\widehat \mu_n(z)}{\widehat\mu_n(2\Z+\un_{2\N+1}(n))}
%=
%    \frac{\widehat \mu_n(z)}{\sum_{\ell=x_{-1}}^{x_1-1}e^{-V(\ell)}}
=
    \frac{\widehat \mu_n(b_{\log n}-k)}
    {\sum_{i=M^-}^{M^+-1}e^{-V(i)}}
=
    \frac{e^{-V(b_{\log n}-k)}+e^{-V(b_{\log n}-k-1)}}
    {\sum_{i=M^-}^{M^+-1}e^{-V(i)}}
$$
since $z$ and $n$ have the same parity,
and $\widehat\mu_n(2\Z)=\widehat\mu_n(2\Z+1)=\sum_{i=M^-}^{M^+-1}e^{-V(i)}$,
and where we used the definition \eqref{eq_def_mi_chapeau} of $\widehat \mu_n$ on $]M^-, M^+[$
and \eqref{Ineg_zn_dans_vallee_centrale} in the last equality.

Now, we define for $j\in\{0,1\}$,
\begin{equation}\label{eq_def_J2pm}
    J_2^\pm(k, n, z,j)
:=
        \E\bigg[\un_{E_\pm^{(n)}
    %\cap E_3^{(n)} \cap E_5^{(n)}
    \cap \{b_{\log n}=z+k\}\cap E_6^{(n)}}
    \frac{e^{-V(b_{\log n}-k-j)}}
    {\sum_{i=M^-}^{M^+-1}e^{-V(i)}}
    \bigg].
\end{equation}
Notice that for $k\in\Z$ such that $|k|\leq \Gamma_n$, if $k\leq -z$ then
%on $\{b_{\log n}=z+k\}$,
%$b_{\log n}=z+k\leq 0$, so $\omega\in E_-^{(n)}$,
$
    \{b_{\log n}=z+k\}
\subset
    \{b_{\log n}\leq 0\}
=
E_-^{(n)}
$,
so
$    J_0(k, n, z)
\leq
    J_2^-(k, n, z,0)
+
    J_2^-(k, n, z,1),
$
whereas if $k> -z$, then
%on $\{b_{\log n}=z+k\}$,
%we have $b_{\log n}=z+k> 0$, so $\omega\in E_+^{(n)}$, so
$
    \{b_{\log n}=z+k\}
\subset
    \{b_{\log n}>0\}
=
    E_+^{(n)}
$,
so
$
    J_0(k, n, z)
\leq
    J_2^+(k, n, z,0)
+
    J_2^+(k, n, z,1)
$.
So we have, thanks to \eqref{Ineg_Proba_Zn_EL_1},
\begin{equation}\label{Ineg_Proba_J1_0_et_1_Plus_ET_Moins}
    \PP\big(S_n=z, E_C^{(n)}(z)\big)
\leq
    J_3(n,z,0)+J_3(n,z,1)
    +5(\log n)^{-3},
%    +O\big((\log n)^{-3}\big),
\end{equation}
where for $j\in\{0,1\}$,
\begin{equation}\label{eq_def_J3}
    J_3(n,z,j)
:=
    \sum_{k=-\Gamma_n}^{\Gamma_n}
        \Big[J_2^-(k, n, z,j) {\bf 1}_{\{k+z\leq 0\}}
        +
        J_2^+(k, n, z,j) {\bf 1}_{\{k+z> 0\}}
        \Big].
\end{equation}
We first consider $k\leq -z$, with $|k|\leq \Gamma_n$.
Hence on $\{b_{\log n}=z+k\}$, we have $b_{\log n}\leq 0$, so $\omega\in E_-^{(n)}$,
thus $M^-=x_{-1}$, $b_{\log n}=x_0$ and $M^+=x_1$ (recall that $x_i=x_i(V,\log n)$, $i\in\Z$  in this section). So
%\begin{equation}\label{Ineg_J0_J2}
%    J_0(k, n, z)
%\leq
%    J_2^-(k, n, z,0)
%+
%    J_2^-(k, n, z,1),
%\qquad
%    k\leq -z
%\end{equation}
%where
for $j\in\{0,1\}$,
recalling that for $u\in\Z$, $V_u(.):=V(u+.)-V(u)$,
and $V_u^-(.):=V(u-.)-V(u)$, we have
\begin{align}
&
    J_2^-(k, n, z,j) =
        \E\bigg[\un_{E_-^{(n)}
    %\cap E_3^{(n)} \cap E_5^{(n)}
    %\cap \{b_{\log n}=z+k\}\cap E_6^{(n)}}
    \cap \{x_0=z+k\}\cap E_6^{(n)}}
    %\frac{e^{-V(b_{\log n}-k-j)}}
    \frac{e^{-V(x_0-k-j)}}
    {\sum_{i=x_{-1}}^{x_1-1}e^{-V(i)}}
    \bigg]
\nonumber\\
& =
    \E\bigg[
    {\bf 1}_{\big\{V_{x_0}(x_1-x_0)\geq \log n,
    \ T_{V_{x_0}}(\log n)>\Gamma_n,
    \ T_{V_{x_0}^-}(\log n)>\Gamma_n
    \big\}}
    \frac{e^{-V_{x_0}(-k-j)}{\bf 1}_{\{x_0=z+k\}}}
    {\sum_{i=x_{-1}-x_0}^{x_1-x_0-1}e^{-V_{x_0}(i)}}
    \bigg].
\label{eq_J2_Indicatrice}
\end{align}
Notice that
$
    (V_{x_0}(i), \ 0\leq i \leq x_1-x_0)
=
    \theta(T_0(V,\log n))
$
and that
$
    (V_{x_0}^-(i), \ 0\leq i \leq x_0-x_{-1})
=
    (V(x_0-i)-V(x_0), \ 0\leq i \leq x_0-x_{-1})
=
    (V_{x_{-1}}(x_0-x_{-1}-i)-V_{x_{-1}}(x_0-x_{-1}), \ 0\leq i \leq x_0-x_{-1})
=
    \zeta[\theta(T_{-1}(V,\log n))]
$,
with $\zeta$ defined in \eqref{eq_def_zeta}.
Also, on the event in \eqref{eq_J2_Indicatrice},
$|k|\leq \Gamma_n$ implies that
$-k-j\leq \Gamma_n+1\leq T_{V_{x_0}}(\log n) \leq x_1-x_0=\ell[\theta(T_0(V,\log n))]$,
and similarly
$k+j\leq x_0-x_{-1}=\ell[\zeta(\theta(T_{-1}(V,\log n)))]$.
Hence, with the following notation for slopes $v$ and $t$,
%when $|k|\leq \Gamma_n$,
\begin{equation}\label{eq_def_phi_v_t}
    \varphi_v(t)
:=
    {\bf 1}_{\{t(\ell(t))\geq \log n,
    \ T_{t}(\log n)\wedge T_{v}(\log n)>\Gamma_n
%    \ >\lfloor (\log n)^{4/3+\delta_1}\rfloor
%    \ T_{t}(\log n)>\lfloor (\log n)^{4/3+\delta_1}\rfloor,
%    \ T_{v}(\log n)>\lfloor (\log n)^{4/3+\delta_1}\rfloor
    \}}
    \frac{e^{-t(-k-j)}{\bf 1}_{\{k+j\leq 0\}}+e^{-v(k+j)}{\bf 1}_{\{k+j> 0\}}}
    {\sum_{i=1}^{\ell(v)}e^{-v(i)}+\sum_{i=0}^{\ell(t)-1}e^{-t(i)}},
\end{equation}
in which we do not write the dependency on $n$, $k$, $j$ to simplify the notations,
we have for our fixed $n$, $k$ and $j$ since $|k|\leq \Gamma_n$,
\begin{equation}\label{eq_J2_Slopes}
    J_2^-(k, n, z,j)
=
    \E\big[\varphi_{\zeta[\theta(T_{-1}(V,\log n))]}[\theta(T_0(V,\log n))]{\bf 1}_{\{x_0=z+k\}}\big].
\end{equation}
In the rest of this section, all the slopes considered, such as $\mathcal T_{V,h}^{\uparrow}$,
$\mathcal T_{V,h}^{\downarrow}$, $\mathcal T_{V_-,h}^{\uparrow*}$, etc, are with $h=\log n$,
and we remove this subscript $h$ to simplify the notation.
That is, $\mathcal T_{V}^{\uparrow}$ denotes $\mathcal T_{V, \log n}^{\uparrow}$,
$\mathcal T_{V}^{\downarrow}$ denotes $\mathcal T_{V,\log n}^{\downarrow}$,
etc.
Due to Theorem \ref{Lemma_Independence_h_extrema} {\bf (i)}, conditionally on $E_-^{(n)}$,
$\zeta[\theta(T_{-1}(V,\log n))]$ is independent of
$(\theta[T_0(V,\log n)], x_0)$
%$T_0(V,\log n)$
%(which is a function of $\theta[T_0(V,\log n)]$ and $x_0$ by definition of $\theta$)
and
%$\zeta[\theta(T_{-1}(V,\log n))]$
has the same law as $\zeta(\mathcal T_V^\downarrow)$ (under $\p$)
and so as $\mathcal T_{V_-}^{\uparrow*}$ by Proposition \ref{Prop_Egalite_Loi_Zeta_Slopes}.
%(right $(\log n)$-upward slope for $V_-$ ) %{\bf (utile ?)}.
Hence, we get, since
$\varphi_{v}[\theta(T_0(V,\log n))]{\bf 1}_{E_+^{(n)}}=0$ for any $v$,
\begin{equation*}
    J_2^-(k, n, z,j)
%& = &
%    \E\big[\varphi_{\mathcal T_{V_-}^{\uparrow *}}[\theta(T_0(V,\log n))]{\bf 1}_{\{x_0=z+k\}}\big]
%\\
=
    \E\Big[\E\big(\varphi_{v}[\theta(T_0(V,\log n))]{\bf 1}_{\{x_0=z+k\}}
        \big)_{|v=\mathcal T_{V_-}^{\uparrow *}}\Big].
\end{equation*}
Thus, applying  the (renewal) Theorem \ref{Lemma_Central_Slope} eq. \eqref{eq_Central_Slope_General_phi}
with $h=\log n$, $\varphi=\varphi_v$, $\Delta_0=\{z+k\}$, $\Delta_1=\Z$
(notice that $\varphi_v(t)=0$ if $t$ is a downward slope
whereas ${\bf 1}_{\{t(\ell(t))\geq \log n\}}=1$ when $t$ is an upward $(\log n)$-slope),
we get, $\mathcal T_V^\uparrow$ and $\mathcal T_{V_-}^{\uparrow *}$ being here independent,
\begin{eqnarray}
    J_2^-(k, n, z,j)
& = &
    \E\Bigg[\frac{\E\big(
            \sharp\{0\leq i <\ell\big(\mathcal T_V^\uparrow\big), \ -i=z+k\}
            \varphi_{v}\big(\mathcal T_V^\uparrow\big)
        \big)_{|v=\mathcal T_{V_-}^{\uparrow *}}}
        {\E\big[\ell\big(\mathcal T_V^\uparrow\big)+\ell\big(\mathcal T_V^\downarrow\big)\big]}\Bigg]
\nonumber\\
& = &
    \E\Bigg(
        \un_{\big\{T_{\mathcal T_V^\uparrow}(\log n)\wedge T_{\mathcal T_{V_-}^{\uparrow *}}(\log n)
        >\Gamma_n\big\}}
%    \frac{\sharp\{0\leq i <\ell\big(\mathcal T_V^\uparrow\big), \ -i=z+k\}}
%    {\E\big(\ell\big(\mathcal T_V^\uparrow\big)+\ell\big(\mathcal T_V^\downarrow\big)\big)}
\nonumber\\
&&
    \frac{e^{-\mathcal T_V^\uparrow(-(k+j))}\un_{\{k+j\leq 0\}}+e^{-\mathcal T_{V_-}^{\uparrow *}(k+j)}\un_{\{k+j> 0\}}}
    {\sum_{i=1}^{\ell(\mathcal T_{V_-}^{\uparrow *})}e^{-\mathcal T_{V_-}^{\uparrow *}(i)}
    +\sum_{i=0}^{\ell(\mathcal T_V^\uparrow)-1}e^{-\mathcal T_V^\uparrow(i)}}
    \frac{\un_{\{-z-k <\ell(\mathcal T_V^\uparrow)\}}}
    {\E\big[\ell\big(\mathcal T_V^\uparrow\big)+\ell\big(\mathcal T_V^\downarrow\big)\big]}
    \Bigg),
~~~~
\label{eq_J2_Transforme}
\end{eqnarray}
where we used
$
    \sharp\{0\leq i <\ell\big(\mathcal T_V^\uparrow\big), \ -i=z+k\}
=
    \un_{\{-z-k <\ell(\mathcal T_V^\uparrow)\}}
$ when $z+k\leq 0$.

We now assume that $k> -z$, with $|k|\leq \Gamma_n$.
We have $b_{\log n}> 0$ on $\{b_{\log n}=z+k\}$, and so $\omega\in E_+^{(n)}$,
thus $b_{\log n}=x_1$, $M^-=x_0$ and  $M^+=x_2$.
So by \eqref{eq_def_J2pm}, for  $j\in\{0,1\}$,
%\begin{equation}\label{Ineg_J0_J2_bPlus}
%    J_0(k, n, z)
%\leq
%    J_2^+(k, n, z,0)
%+
%    J_2^+(k, n, z,1),
%\qquad
%    k> -z
%\end{equation}
%where for $j\in\{0,1\}$,
%\begin{align*}
\begin{equation*}
%&
    J_2^+(k, n, z,j)
%=
%    \E\bigg[\un_{E_+^{(n)}
%    \cap \{b_{\log n}=z+k\}\cap E_6^{(n)}}
%    \frac{e^{-V(x_1-k-j)}}
%    {\sum_{i=x_0}^{x_2-1}e^{-V(i)}}
%    \bigg]
%\\
%&
=
    \E\bigg[
    {\bf 1}_{\big\{V_{x_1}(x_0-x_1)\geq \log n,\ x_1=z+k,
    \ T_{V_{x_1}}(\log n)\wedge T_{V_{x_1}^-}(\log n)>\Gamma_n
    \big\}}
    \frac{e^{-V_{x_1}(-k-j)}}
    {\sum_{i=x_0-x_1}^{x_2-x_1-1}e^{-V_{x_1}(i)}}
    \bigg].
\end{equation*}
%\end{align*}
Notice that
$
    (V_{x_1}(i), \ 0\leq i \leq x_2-x_1)
=
    \theta(T_1(V,\log n))
$
and that
$
    (V_{x_1}(-i), \ 0\leq i \leq x_1-x_0)
=
    (V(x_1-i)-V(x_1), \ 0\leq i \leq x_1-x_0)
%=
%    (V_{x_0}(x_1-x_0-i)-V_{x_0}(x_1-x_0), \ 0\leq i \leq x_1-x_0)
=
    \zeta[\theta(T_0(V,\log n))]
$.
Hence, with
\begin{equation*}
    \varphi_v^+(t)
:=
    {\bf 1}_{\{t(\ell(t)) \geq \log n,
    \ T_{t}(\log n)\wedge T_{v}(\log n)>\Gamma_n
%    \ >\lfloor (\log n)^{4/3+\delta_1}\rfloor
%    \ T_{t}(\log n)>\lfloor (\log n)^{4/3+\delta_1}\rfloor,
%    \ T_{v}(\log n)>\lfloor (\log n)^{4/3+\delta_1}\rfloor
    \}}
    \frac{e^{-v(-k-j)}{\bf 1}_{\{k+j\leq 0\}}+e^{-t(k+j)}{\bf 1}_{\{k+j>  0\}}}
    {\sum_{i=1}^{\ell(t)}e^{-t(i)}+\sum_{i=0}^{\ell(v)-1}e^{-v(i)}},
\end{equation*}
in which we do not write the dependency on $n$, $k$, $j$ to simplify the notations,
we have
\begin{eqnarray*}
    J_2^+(k, n, z,j)
& = &
    \E\big[\varphi^+_{\theta[T_1(V,\log n)]}[\zeta(\theta(T_0(V,\log n)))]{\bf 1}_{\{x_1=z+k\}}\big].
\end{eqnarray*}
Since due to Theorem \ref{Lemma_Independence_h_extrema} {\bf (ii)}, conditionally on $E_+^{(n)}$,
$\theta(T_1(V,\log n))$
has the law $\mathscr{L}\big(\mathcal T_V^\uparrow\big)$,
and is independent of $(\theta(T_0(V,\log n)), x_1)$,
we have,
%with $\mathcal T_V^\uparrow$ denoting a slope independent of $T_0(V,\log n)$ (with a slight abuse of notation),
\begin{equation*}
    J_2^+(k, n, z,j)
=
%& = &
%    \E\big[\varphi^+_{\mathcal T_V^\uparrow}[\zeta(\theta(T_0(V,\log n)))]{\bf 1}_{\{x_1=z+k\}}\big]
%\\
%& = &
    \E\big[\E\big(\varphi^+_{v}[\zeta(\theta(T_0(V,\log n))   )]{\bf 1}_{\{x_1=z+k\}}
        \big)_{|v=\mathcal T_V^\uparrow}\big],
\end{equation*}
since $\varphi^+_{v}[\zeta(\theta(T_0(V,\log n))   )]{\bf 1}_{E_-^{(n)}}=0$ for any $v$.
Thus, applying  the (renewal) Theorem \ref{Lemma_Central_Slope}
with $h=\log n$,  $\varphi=\varphi_v^+\circ \zeta$, $\Delta_0=\Z$, $\Delta_1=\{z+k\}$
(we use once more that $\varphi_v^+ \circ \zeta (t)=0$ when $t$ is a (translated) upward slope,
since in this case $\zeta (t)$ is a downward slope),
we get
\begin{equation*}
    J_2^+(k, n, z,j)
%& = &
%    \E\big[\E\big(\varphi^+_{v}[\zeta(\theta(T_0(V,\log n))   )]{\bf 1}_{\{x_1=z+k\}}
%        \big)_{|v=\mathcal T_V^\uparrow}\big].
%\\
=
    \E\Bigg[\frac{\E\big(
            \sharp\{0\leq i <\ell\big(\mathcal T_V^\downarrow\big), \ \ell\big(T_V^\downarrow\big)-i=z+k\}
            \varphi_v^+\circ \zeta\big(\mathcal T_V^\downarrow\big)
        \big)_{|v=\mathcal T_V^\uparrow}}
        {\E\big[\ell\big(\mathcal T_V^\uparrow\big)+\ell\big(\mathcal T_V^\downarrow\big)\big]}\Bigg].
\end{equation*}
Recall that, by Proposition \ref{Prop_Egalite_Loi_Zeta_Slopes} {\bf (ii)},
$\zeta(\mathcal T_V^\downarrow)=_{law}\mathcal T_{V_-}^{\uparrow*}$.
Hence, $\mathcal T_V^\uparrow$ and $\mathcal T_{V_-}^{\uparrow*}$ being independent, and using
%With $\mathcal T_V^\uparrow$ being independent of
%$\mathcal T_{V_-}^{\uparrow*}=\zeta(\mathcal T_V^\downarrow)$,
%nd so
$
    \ell\big(\mathcal T_V^\downarrow\big)
=
    \ell\big(\zeta\big(\mathcal T_V^\downarrow\big)\big)
%=
%    \ell\big(\mathcal T_{V_-}^{\uparrow*}\big)
$,
we get
\begin{eqnarray}
    J_2^+(k, n, z,j)
& = &
    \E\Bigg(
        \un_{\big\{T_{\mathcal T_{V_-}^{\uparrow *}}(\log n)>\Gamma_n\big\}}
        \un_{\big\{T_{\mathcal T_V^\uparrow}(\log n)>\Gamma_n\big\}}
\nonumber\\
&&
    \frac{e^{-\mathcal T_V^\uparrow(-(k+j))}\un_{\{k+j\leq 0\}}+e^{-\mathcal T_{V_-}^{\uparrow *}(k+j)}\un_{\{k+j> 0\}}}
    {\sum_{i=1}^{\ell(\mathcal T_{V_-}^{\uparrow *})}e^{-\mathcal T_{V_-}^{\uparrow *}(i)}
    +\sum_{i=0}^{\ell(\mathcal T_V^\uparrow)-1}e^{-\mathcal T_V^\uparrow(i)}}
%    \frac{\sharp\{0\leq i <\ell\big(\mathcal T_{V_-}^{\uparrow *}\big), \ \ell\big(\mathcal T_{V_-}^{\uparrow *}\big)-i=z+k\}}
    \frac{\un_{\{z+k \leq \ell(\mathcal T_{V_-}^{\uparrow*})\}}}
    {\E\big[\ell\big(\mathcal T_V^\uparrow\big)+\ell\big(\mathcal T_V^\downarrow\big)\big]}
    \Bigg),
~~~~
\label{eq_J2_Transforme_Plus}
\end{eqnarray}
where we used
$
    \sharp\{0\leq i <\ell\big(\mathcal T_V^\downarrow\big), \ \ell\big(\mathcal T_V^\downarrow\big)-i=z+k\}
=
    \un_{\{z+k \leq \ell(\mathcal T_V^\downarrow)\}}
$
which becomes
$
    \un_{\{z+k \leq \ell(\mathcal T_{V_-}^{\uparrow*})\}}
$ since $z+k> 0$
and
${\bf 1}_{\{t(\ell(t))\geq \log n\}}=1$ for $t=\zeta\big(\mathcal T_V^\downarrow\big)$.
Notice that the only difference between this formula and \eqref{eq_J2_Transforme}
is that
%in $\sharp\{\dots\}$, $\ell\big(\mathcal T_V^\uparrow\big)$
%is replaced by $\ell\big(\mathcal T_{V_-}^{\uparrow *}\big)$
%and $i$ by $\ell\big(\mathcal T_{V_-}^{\uparrow *}\big)-i$.
$\un_{\{-z-k <\ell(\mathcal T_V^\uparrow)\}}$
is replaced by
$\un_{\{z+k \leq \ell(\mathcal T_{V_-}^{\uparrow*})\}}$.

We now define
\begin {eqnarray*}
    z_n^+
& := &
\left\{
\begin{array}{ll}
z+\Gamma_n & \text{ if } z\leq -\Gamma_n,\\
0 & \text{ if } -\Gamma_n< z \leq \Gamma_n,\\
z-\Gamma_n & \text{ if } z> \Gamma_n,
\end{array}
\right.
\\
    \psi_k\big(\mathcal T_V^\uparrow, \mathcal T_{V_-}^{\uparrow *}, z\big)
& := &
\left\{
\begin{array}{ll}
\un_{\{-z-k <\ell(\mathcal T_V^\uparrow)\}} & \text{ if } z\leq -\Gamma_n,\\
\un_{\{0 <\ell(\mathcal T_V^\uparrow)\}}
%=1
    & \text{ if } -\Gamma_n< z \leq \Gamma_n,\\
\un_{\{z+k \leq \ell(\mathcal T_{V_-}^{\uparrow *})\}} & \text{ if } z> \Gamma_n.
\end{array}
\right.
\end{eqnarray*}
%\noindent{\bf First case:}
Notice that in the case $z\leq -\Gamma_n$,
we have $z+k\leq 0$ for every $k$ in the sum in \eqref{eq_def_J3},
so, using \eqref{eq_J2_Transforme}, we have for each $j\in\{0,1\}$
(the inequality being an equality in this first case $z\leq -\Gamma_n$),
\begin{eqnarray}
    J_3(n,z,j)
& \leq &
    \E\Bigg(
        \un_{\big\{T_{\mathcal T_V^\uparrow}(\log n)>\Gamma_n\big\}}
        \un_{\big\{T_{\mathcal T_{V_-}^{\uparrow^*}}(\log n)>\Gamma_n\big\}}
\nonumber\\
&&
    \frac{
        \sum_{k=-\Gamma_n}^{\Gamma_n}
        \big(e^{-\mathcal T_V^\uparrow(-(k+j))}\un_{\{k+j\leq 0\}}
        +e^{-\mathcal T_{V_-}^{\uparrow *}(k+j)}\un_{\{k+j> 0\}}\big)
        \psi_k\big(\mathcal T_V^\uparrow, \mathcal T_{V_-}^{\uparrow *}, z\big)}
    {\Big(\sum_{i=0}^{\ell(\mathcal T_V^\uparrow)-1}e^{-\mathcal T_V^\uparrow(i)}
    +\sum_{i=1}^{\ell(\mathcal T_{V_-}^{\uparrow *})}e^{-\mathcal T_{V_-}^{\uparrow *}(i)}\Big)
    \E\big[\ell\big(\mathcal T_V^\uparrow\big)+\ell\big(\mathcal T_V^\downarrow\big)\big]}
    \Bigg).
~~~~~~
\label{Ineg_J3_a}
\end{eqnarray}
When $z> \Gamma_n$,
we have $z+k> 0$ for every $k$ in the sum in \eqref{eq_def_J3}.
So, combining \eqref{eq_def_J3} and  \eqref{eq_J2_Transforme_Plus}, inequality \eqref{Ineg_J3_a}
remains true in this case (and is actually an equality in this second case).

Finally,
assume that $-\Gamma_n<z\leq \Gamma_n$.
In this case,  notice that the quantity
$\un_{\{-z-k <\ell(\mathcal T_V^\uparrow)\}}$
which appears in \eqref{eq_J2_Transforme} for $k+z\leq 0$,
and the quantity $\un_{\{z+k \leq \ell(\mathcal T_{V_-}^{\uparrow*})\}}$
which appears in \eqref{eq_J2_Transforme_Plus} for $k+z> 0$
are both dominated by
$1=\un_{\{0 <\ell(\mathcal T_V^\uparrow)\}}=\psi_k\big(\mathcal T_V^\uparrow, \mathcal T_{V_-}^{\uparrow *}, z\big)$
$\p$-a.s.,
so
$J_2^-(k, n, z,j)$ and $J_2^+(k, n, z,j)$ are dominated by the same formula. So for $j\in\{0,1\}$,
\eqref{Ineg_J3_a} also remains true in this case.
So, \eqref{Ineg_J3_a} holds for every $z\in\Z$ and every $j\in\{0,1\}$.

Now, we notice that
for every
$
    -\Gamma_n
\leq
    k
\leq
    \Gamma_n
$,
we have
$
    \psi_k\big(\mathcal T_V^\uparrow, \mathcal T_{V_-}^{\uparrow *}, z\big)
=
    \un_{\{-z-k <\ell(\mathcal T_V^\uparrow)\}}
$
$
\leq
    \un_{\{-z-\Gamma_n<\ell(\mathcal T_V^\uparrow)\}}
=
    \un_{\{-z_n^+<\ell(\mathcal T_V^\uparrow)\}}
$
when $z\leq -\Gamma_n$,
also
$
    \psi_k\big(\mathcal T_V^\uparrow, \mathcal T_{V_-}^{\uparrow *}, z\big)
=
    \un_{\{-z_n^+<\ell(\mathcal T_V^\uparrow)\}}
$
when
$-\Gamma_n< z \leq \Gamma_n$,
whereas
$
    \psi_k\big(\mathcal T_V^\uparrow, \mathcal T_{V_-}^{\uparrow *}, z\big)
=
    \un_{\{z+k \leq \ell(\mathcal T_{V_-}^{\uparrow *})\}}
$
$
\leq
    \un_{\{z-\Gamma_n \leq \ell(\mathcal T_{V_-}^{\uparrow *}) \}}
$
$
=
    \un_{\{z_n^+ \leq \ell(\mathcal T_{V_-}^{\uparrow *}) \}}
$
when $z> \Gamma_n$.

Hence, \eqref{Ineg_J3_a} leads to, for every $j\in\{0,1\}$, $n\geq n_3$ and $z\in(2\Z+n)$, as explained below,
\begin{equation}
    J_3(n,z,j)
\leq
        \frac{\p\big[-z_n^+ <\ell\big(\mathcal T_V^\uparrow\big)\big]}
        {\E\big[\ell\big(\mathcal T_V^\uparrow\big)+\ell\big(\mathcal T_V^\downarrow\big)\big]}
    {\bf 1}_{\{z\leq \Gamma_n\}}
+
        \frac{\p\big[ z_n^+ \leq \ell(\mathcal T_{V_-}^{\uparrow *}) \big]}
        {\E\big[\ell\big(\mathcal T_V^\uparrow\big)+\ell\big(\mathcal T_V^\downarrow\big)\big]}
    {\bf 1}_{\{z> \Gamma_n\}}
\label{Ineg_Premiere}
=
    \p\big(b_{\log n}= z_n^+\big).
\end{equation}
Indeed, we first used
$
    \Gamma_n+1
\leq
    T_{\mathcal T_V^\uparrow}(\log n)
\leq
    \ell(\mathcal T_V^\uparrow)
$
and similarly
$
    \Gamma_n+1
%$
%$
%\leq
%    T_{\mathcal T_V^\uparrow}(\log n)
\leq
    \ell(\mathcal T_{V_-}^{\uparrow *})
$,
%and the same inequalities with $V$ replaced by $V_-$
%and $\mathcal T_V^\uparrow$ by $\mathcal T_{V_-}^{\uparrow *}$,
so that
$
\sum_{k=-\Gamma_n}^{\Gamma_n}(\dots+\dots)
\leq
\big(
\sum_{i=0}^{\ell(\mathcal T_V^\uparrow)-1}\dots
    +\sum_{i=1}^{\ell(\mathcal T_{V_-}^{\uparrow *})}\dots
\big)
$
in \eqref{Ineg_J3_a} to get the (first) inequality.
Then, to get the following equality, we used
%\eqref{eq_proba_bh_x_negatif}
%the second equality
eq. \eqref{eq_proba_bh_x_negatif} of Lemma \ref{Lemma_Proba_bh_egal}
when $z\leq \Gamma_n$,
and
$
    \ell\big(\mathcal T_{V_-}^{\uparrow*}\big)
=_{law}
    \ell\big(\zeta\big(\mathcal T_V^\downarrow\big)\big)
=
    \ell\big(\mathcal T_V^\downarrow\big)
$
by Proposition \ref{Prop_Egalite_Loi_Zeta_Slopes} {\bf (ii)}
and
%the first equality
eq. \eqref{eq_proba_bh_x_positif} of Lemma \ref{Lemma_Proba_bh_egal}
%\eqref{eq_proba_bh_x_positif}
when $z> \Gamma_n$.

Now, let $\e>0$.
By Theorem \ref{Th_Local_Limit_b_h}, there exists $n_4\geq n_3$ such that,
for every $j\in\{0,1\}$, $n\geq n_4$ and $z\in(2\Z+n)$,
$$
    J_3(n,z,j)
\leq
    \p\big(b_{\log n}= z_n^+\big)
\leq
    \frac{\sigma^2}{(\log n)^2}\varphi_\infty\bigg(\frac{\sigma^2 z_n^+}{
    (\log n)^2}\bigg)
    +\e(\log n)^{-2}.
$$
Now, recall that $\varphi_\infty$ is uniformly continuous on $\R$
since $\varphi_\infty$ is continuous on $\R$ and $\lim_{\pm\infty}\varphi_\infty=0$.
Also, $\sup_{z\in\Z}|\sigma^2 z_n^+(\log n)^{-2}-\sigma^2 z(\log n)^{-2}|\to 0$
as $n\to +\infty$ because $\delta_1<2/3$.
Thus, there exists $n_5\geq n_4$ such that for all $n\geq n_5$,
$\sup_{z\in\Z}|\varphi_\infty(\sigma^2 z_n^+(\log n)^{-2})-\varphi_\infty(\sigma^2 z(\log n)^{-2})|\leq \sigma^{-2}\e$.
Hence,
\begin{equation}
    \forall n\geq n_5, \, \forall z\in(2\Z+n), \forall j\in\{0,1\},
\quad
    J_3(n,z,j)
\leq
    \frac{\sigma^2}{(\log n)^2}\varphi_\infty\bigg(\frac{\sigma^2 z}{(\log n)^2}\bigg)
    +\frac{2\e}{(\log n)^2}.
\label{Ineg_J1j_1}
\end{equation}
Finally, \eqref{Ineg_Proba_J1_0_et_1_Plus_ET_Moins} and \eqref{Ineg_J1j_1}
lead to, for all $n\geq n_5$,
$$
    %\forall n\geq n_5, \,
    \forall z\in(2\Z+n),
\quad
    \PP\big(S_n=z, E_C^{(n)}(z)\big)
%\leq
%    J_3(n,z,0)+J_3(n,z,1)
%    +5(\log n)^{-3},
\leq
    \frac{2\sigma^2}{(\log n)^2}\varphi_\infty\bigg(\frac{\sigma^2 z}{(\log n)^2}\bigg)
    +\frac{4\e}{(\log n)^2}
    +\frac{5}{(\log n)^{3}}.
$$
This gives \eqref{Ineg_Minor_Sur_EL},
which proves the proposition.
\hfill$\Box$

%%%%%%%%%%%%%%%%%%%%%%%%%%%%%%%%%%%%%%%%%%%%%%%%%%%%%%%%%%%%%%%%%%%%%%%%%%%%%%%%%%%%%%%%
%                                                                                      %
%                               SECTION ENV NEGLIGEABLES                               %
%                                                                                      %
%%%%%%%%%%%%%%%%%%%%%%%%%%%%%%%%%%%%%%%%%%%%%%%%%%%%%%%%%%%%%%%%%%%%%%%%%%%%%%%%%%%%%%%%

\section{Proving that some environments or trajectories are negligible}
\label{Sect_Negligible}
The aim of this section is to prove that
$\sup_{z\in\Z}\PP\big(S_n=z, \big(E_C^{(n)}(z)\big)^c \big)$ is negligible
compared to $(\log n)^{-2}$ as $n\to+\infty$
(recall $E_C^{(n)}(z)$ from \eqref{eq_def_EL}).
To this aim we give upper bounds of the probabilities of different events, most of them depending both on the environment and on the walk,
except the event considered in Lemma \ref{Lemma_Proba_E3c}.

%%%%%%%%%%%%%%%%%%%%%%%%%%

\subsection{Contribution of $\big(E_4^{(n)}(z)\big)^c$}

As a warm up, we start with following estimate.

\begin{lem}\label{Lemma_Proba_E3c}
There exists $c_9>0$ such that
$$
    \forall n\geq n_5,\,
    \forall z\in\Z,
\qquad
%    \sup_{z\in\Z}
    \p\big[
        \big(E_4^{(n)}(z)\big)^c
        \cap E_3^{(n)}
        \cap E_6^{(n)}
        \cap E_7^{(n)}(z)
%        \cap \big\{|b_{\log n}-z|\leq (\log n)^{4/3+\delta_1}\big\}
    \big]
\leq
    c_9(\log_2 n)^3(\log n)^{-3}.
$$
\end{lem}

\begin{figure}[htbp]
\includegraphics[width=16.0cm,height=7.91cm]{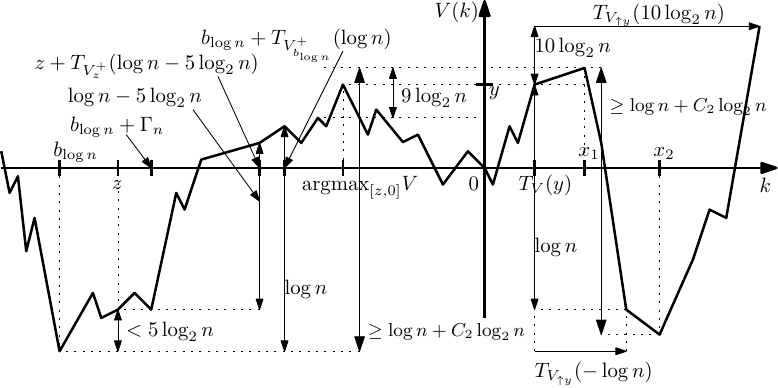}
\caption{Schema of the potential $V$ on $\omega\in E_{8, -}^{(n)}(z)$, with $x_i=x_i(V,\log n)$
and $y=\max_{[z,0]}V$.}
\label{figure_Proba_E3c}
\end{figure}

\noindent{\bf Proof:}
Let $n\geq n_5$ and $z\in\Z$. We introduce
\begin{equation*}
    E_8^{(n)}(z)
:=
    \big(E_4^{(n)}(z)\big)^c
    \cap
    E_3^{(n)}
    %\cap E_5^{(n)}
    \cap E_6^{(n)}
    \cap E_7^{(n)}(z),
%    \cap \{|b_{\log n}-z|\leq (\log n)^{4/3+\delta_1}\}
\qquad
    E_{8, \pm}^{(n)}(z)
:=
    E_{\pm}^{(n)}
    \cap
    E_8^{(n)}(z).
\end{equation*}
We first assume that $\omega\in E_{8, -}^{(n)}(z)$ (see Figure \ref{figure_Proba_E3c}).
Hence, $\omega \in E_-^{(n)}\cap \big(E_4^{(n)}(z)\big)^c$, so
$b_{\log n}=x_0(V,\log n)$,
\begin{equation}\label{Ineg_Vzn}
    V(z)-V(b_{\log n})
<
     5\log_2 n
\end{equation}
and
\begin{eqnarray}
    \max_{[b_{\log n}, 0]}V
& \geq &
    V[x_1(V,\log n)]-9\log_2 n
=
    V(b_{\log n})+H[T_0(V,\log n)]-9\log_2 n
\nonumber\\
& \geq &
    V(b_{\log n})+\log n +(C_2-9)\log_2 n
>
    V(b_{\log n})+\log n,
\label{Ineg_TVzn_0}
\end{eqnarray}
since
$
    H[T_0(V,\log n)]
=
    H[T_0(V,h_n-C_1\log_2 n)]
\geq
    \log n+C_2\log_2 n
$
by Remark \ref{Remark_x_i}
because $\omega\in E_3^{(n)}$, and where we used $C_2>9$.

Also, $\omega\in E_6^{(n)}\cap E_7^{(n)}(z)$, so as in \eqref{Ineg_zn_Mplus}, using \eqref{Ineg_TVzn_0} in the last inequality,
\begin{equation}\label{Ineg_zn_atteinte_et_0}
    b_{\log n}-T_{V_{b_{\log n}}^-}(\log n)
<
    z
%\leq
%    b_{\log n}+\Gamma_n
<
    b_{\log n}+T_{V_{b_{\log n}}^+}(\log n)
\leq 0,
%\leq
%    M^+.
\end{equation}
where for $x\in\Z$, $V_x^\pm(k)=V(x\pm k)-V(x)$, $k\in\N$, as before.
This and \eqref{Ineg_TVzn_0} also lead to
\begin{equation}\label{Ineg_max_zn0_x1}
    \max_{[z,0]}V
=
    \max_{[b_{\log n}, 0]}V
\geq
    V[x_1(V,\log n)]-9\log_2 n.
\end{equation}
We now introduce, for $y\geq 0$, $V_{\uparrow y}(k):=V[k+T_V(y)]-V[T_V(y)]$, $k\in\N$, and
\begin{eqnarray*}
    E_{9, \pm}^{(n)}(z)
& := &
    \big\{T_{V_{z}^\pm}(\log n-5\log_2 n)<T_{V_{z}^\pm}(-5\log_2 n)\big\},
\\
    E_{10}^{(n)}(y)
& := &
    \big\{T_{V_{\uparrow y}}(-\log n) < T_{V_{\uparrow y}}(10\log_2 n)\big\},
\end{eqnarray*}
where $\log n-5\log_2 n>h_n>0$ since $n\geq n_5\geq n_3$ and $C_1>20$.
Due to \eqref{Ineg_Vzn} and \eqref{Ineg_zn_atteinte_et_0} and since $b_{\log n}$ is a left $(\log n)$-minimum,
we have $\omega\in E_{9, -}^{(n)}(z)\cap E_{9, +}^{(n)}(z)$.

Also, notice that, using \eqref{Ineg_max_zn0_x1},
$
    V[x_1(V,\log n)]
=
    \max_{[0, x_2(V,\log n)]}V
=
    \max_{[b_{\log n}, x_2(V,\log n)]}V
\geq
    \max_{[b_{\log n}, 0]}V
=
    \max_{[z,0]}V
$
and
$
    H[T_1(V,\log n)]
=
    H[T_1(V,h_n-C_1\log_2 n)]
\geq
    \log n+C_2\log_2 n
$
with $C_2>9$
(by Remark \ref{Remark_x_i}
%similarly as for $T_0$
since $\omega\in E_3^{(n)}$).
So, after hitting $\big[\max_{[z,0]}V,+\infty\big[$, the potential $(V(u),\ u\geq 0)$ cannot take values larger than
$
    V[x_1(V,\log n)]
\leq
    \max_{[z,0]}V+9\log_2 n
$ (see \eqref{Ineg_max_zn0_x1})
before going (down) to $x_2(V,\log n)$ with
$
    V[x_2(V,\log n)]
=
    V[x_1(V,\log n)]
-
    H[T_1(V,\log n)]
\leq
    V[x_1(V,\log n)]-\log n-C_2\log_2 n
\leq
    \max_{[z,0]}V-\log n
$
by \eqref{Ineg_max_zn0_x1} and since $C_2>9$.
Hence, $\omega\in E_{10}^{(n)}\big(\max_{[z,0]}V\big)$.

Finally,
$
    z+T_{V_{z}^+}(\log n-5\log_2 n)
\leq
    b_{\log n}+T_{V_{b_{\log n}}^+}(\log n)
\leq
    0
$
by \eqref{Ineg_Vzn} and \eqref{Ineg_zn_atteinte_et_0},
and
$E_{9, -}^{(n)}(z)\cap E_{9, +}^{(n)}(z)\cap \{z+T_{V_{z}^+}(\log n-5\log_2 n)\leq 0\}$ depend only on
$V^-=(V(k),\ k\leq 0)$.
Hence, conditioning by $V^-$ to get the third line,
using \eqref{eqOptimalStopping2} for the forth,
the independence of $V_{z}^-$ and $V_{z}^+$ and $C_0<\log_2 n$ since $n\geq n_5\geq n_3$ for the fifth,
and once more \eqref{eqOptimalStopping2} for the sixth,
we get for every $n\geq n_5$ and $z\in\Z$,
\begin{eqnarray}
&&
    \p\big[E_{8, -}^{(n)}(z)\big]
\nonumber\\
& \leq &
    \p\big[
        E_{9, -}^{(n)}(z)\cap E_{9, +}^{(n)}(z)
        \cap \big\{z+T_{V_{z}^+}(\log n-5\log_2 n)\leq 0\big\}
        \cap E_{10}^{(n)}\big(\max\nolimits_{[z,0]}V\big)\big]
\nonumber\\
& = &
    \E\big[
            \un_{E_{9, -}^{(n)}(z)\cap E_{9, +}^{(n)}(z)\cap \{z+T_{V_{z}^+}(\log n-5\log_2 n)\leq 0\}}
            \p\big( E_{10}^{(n)}\big(\max\nolimits_{[z,0]}V\big) | V^-\big)
    \big]
\nonumber\\
& \leq &
    \E\big[
            \un_{E_{9, -}^{(n)}(z)\cap E_{9, +}^{(n)}(z)}
            (10\log_2 n+C_0)(\log n+10\log_2 n +C_0)^{-1}
    \big]
\nonumber\\
& \leq &
    \p\big[E_{9, -}^{(n)}(z)\big]
    \p\big[E_{9, +}^{(n)}(z)\big]
    (11 \log_2 n)(\log n)^{-1}
\nonumber\\
& \leq &
    (6\log_2 n)^2(11 \log_2 n)(\log n)^{-3}.
\label{Ineg_E27_Moins}
\end{eqnarray}
%Hence
%$
%    \p\big[E_{8, -}^{(n)}\big]
%=
%    O\big((\log_2 n)^3(\log n)^{-3}\big)
%$.

We show similarly that
$
    \p\big[E_{8, +}^{(n)}(z)\big]
%=
%    O\big((\log_2 n)^3
\leq
    396(\log_2 n)^3(\log n)^{-3}
%    (\log n)^{-3}\big)
$
for every $n\geq n_5$ and $z\in\Z$.
This, combined with \eqref{Ineg_E27_Moins},  ends the proof of the lemma.
\hfill$\Box$

%%%%%%%%%%%%%%%%%%%%%%%%%%

%\subsection{Contribution of events such that $b_{\log n}$ is far from $z$ when $\tau(b_{\log n})\leq n$}
\subsection{Case when $b_{\log n}$ is far from $z$ without subvalleys or small valleys}

In this subsection, we prove that the event constituted by environments and trajectories
such that $b_{\log n}$ is far from $z$ and $S_n=z$ while
$E_3^{(n)}\cap E_5^{(n)}$ holds is negligible.
More precisely, we prove the following proposition.

\begin{prop}\label{Lemma_Proba_z_b_far}
%We have for large $n$,
There exist $c_{10}>0$ and $n_6\geq n_5$ such that, for all $n\geq n_6$,
\begin{equation}\label{Ineg_Proba_z_b_far}
    \forall z\in\Z,
\qquad
    \PP\big(S_n=z,
    %|z-b_{\log n}|>(\log n)^{4/3+\delta_1},
    |z-b_{\log n}|>\Gamma_n,
    E_3^{(n)}, E_5^{(n)}\big)
%=
%    o\big((\log n)^{-2}\big).
\leq
    c_{10}(\log n)^{-2-\delta_1/2}.
\end{equation}
\end{prop}

Before giving a complete proof, we first introduce the different cases considered.

\noindent{\bf Organisation of the proof:}
We consider separately the case $\tau(b_{\log n})\leq n$ (see Lemma \ref{Lemma_Proba_z_far_b_atteint})
and the case $\tau(b_{\log n})> n$ (see Lemmas \ref{Lemma_Proba_z_far_b_non_atteint}, \ref{Lemma_Proba_Nulle} and \ref{Lemma_Proba_F1_E2})
since in this second case, we prove (see \eqref{Ineg_Proba_Sortie_x1x3})
that with large enough probability, $\tau[x_2(V,\log n)]\leq n$ on $E_-^{(n)}$
and similarly $\tau[x_{-1}(V,\log n)]\leq n$ on $E_+^{(n)}$.
So in the first case $\tau(b_{\log n})\leq n$,
$S$ goes before time $n$ to the bottom $b_{\log n}$ of the central valley of height at least $\log n$,
whereas in the second case $\tau(b_{\log n})> n$,
$S$ goes before time $n$ to the bottom of a neighbour valley of height at least $\log n$ with large probability.
Figure \ref{figure_Lemma_5_5} gives the schema of a potential for which $S$ can go
before time $n$, with relatively comparable quenched probability, to each of the bottoms
of the two valleys "surrounding" the origin, $x_0(V,\log n)$ and $x_2(V, \log n)$ in this figure.

%Assume for example that $b_{\log n}\leq 0$, so $b_{\log n}=x_0(V,\log n)$.
%We prove in \eqref{Ineg_Proba_Sortie_x1x3} that $\tau[x_0(V,\log n)]\wedge \tau[x_2(V,\log n)]\leq n$ %with large enough probability.

%If $\tau(b_{\log n})\leq n$, thanks to reversibility,
%we compare the quenched probability that $S_n=z$ with the invariant probability measure of $z$,
%which has small expectation since $z$ is far from $b_{\log n}$.

%Otherwise if $\tau(b_{\log n})> n$, then we have $\tau[x_2(V,\log n)]\leq n$,
%and  $z=S_n>b_{\log n}+(\log n)^{4/3+\delta_1}$.
%In this case, we show that, if $z$ is far from $x_2(V,\log n)$,
%then either $S$ has to cross $x_1(V,\log n)$ which has high potential before time $n$,
%either the potential of $z$ itself is high (compared to that of $x_2(V,\log n)$), so that the %probability that $S_n=z$ is low
%once more by reversibility.

%We now prove this rigourously.

\subsubsection{Case when $\tau(b_{\log n})\leq n$}

In this subsection, we
%prove that the contribution of environments such that $b_{\log n}$ is far from $z$, that is,
%$|b_{\log n}-z|> (\log n)^{4/3+\delta_1}$, to $\{S_n=z\}$ is negligible. We start with the case where $S$ hits $b_{\log n}$ before time $n$,
consider the case $\tau(b_{\log n})\leq n$ of Proposition \ref{Lemma_Proba_z_b_far},
since for this case we can use an inequality coming from the reversibility of $S$.
More precisely, we prove the following lemma.

\begin{lem} \label{Lemma_Proba_z_far_b_atteint}
There exists $c_{11}>0$ and $n_6\geq n_5$ such that for all $n\geq n_6$,
\begin{equation}
    \forall z\in\Z,
\qquad
    \PP\big(S_n=z,
    %|b_{\log n}-z|> (\log n)^{4/3+\delta_1},
    |b_{\log n}-z|> \Gamma_n,
    \tau(b_{\log n})\leq n, E_3^{(n)}, E_5  ^{(n)} \big)
\leq
    c_{11}(\log n)^{-2-\delta_1/2}.
\label{Ineg_Proba_z_far_b_atteint}
\end{equation}
\end{lem}

\noindent{\bf Proof:}
In this proof, $\mathcal T_V^\uparrow$ and $\mathcal T_V^\downarrow$
denote respectively $\mathcal T_{V,\log n}^\uparrow$ and $\mathcal T_{V,\log n}^\downarrow$.
By Lemma \ref{Lemma_Esperance_Longueur_Slope} applied with $h=\log n$, there exists $n_6\geq n_5$ such that
for all $n\geq n_6$,
$
    \E\big[\ell\big(\mathcal T_V^\uparrow\big)+\ell\big(\mathcal T_V^\downarrow\big)\big]
\geq
    c_7(\log n)^2
$.
Let $n\geq n_6$ and $z\in\Z$.
We separate the proof into different cases, first when $z\notin[M^-, M^+]$,
then when $z \in[M^-, M^+]-\big]\widehat L^-,\widehat L^+\big[$
and finally when $z\in \big]\widehat L^-,\widehat L^+\big[$, this last case being
cut into four subcases, depending on the signs of $b_{\log n}$ and of $z-b_{\log n}$.

{\bf First step:} we have,
% for $n\geq n_5$ and $z\in\Z$,
conditioning by $\omega$ and applying the strong Markov property at stopping time $\tau(\log n)$,
recalling $M^\pm$ from \eqref{eq_def_Mpm} (with $x_i=x_i(V,\log n)$, see Figure \ref{figure_Good_Env_Coupling}),
\begin{eqnarray}
&&
    \PP\big(S_n=z,
%        |b_{\log n}-z|> (\log n)^{4/3+\delta_1},
        |b_{\log n}-z|> \Gamma_n,
        \tau(b_{\log n})\leq n,
    E_3^{(n)}, E_5^{(n)}, z\notin[M^-,M^+] \big)
%    \PP\big(S_n=z,  |b_{\log n}-z|> (\log n)^{4/3+\delta_1}, \tau(b_{\log n})\leq n, E_3^{(n)}, z\notin[M^-,M^+] \big)
\nonumber\\
& \leq &
    \EE\big(
%    \un_{E_3^{(n)}\cap E_5^{(n)}}
    \un_{E_3^{(n)}}
    \un_{\{\tau(b_{\log n})\leq n\}}
                \po^{b_{\log n}}(S_k\notin[M^-,M^+])_{|k=n-\tau(b_{\log n)})}\big)
%    +\p\big[\big(E_5^{(n)}\big)^c\big]
\nonumber\\
& \leq &
    \E\big[
%    \un_{E_3^{(n)}\cap E_5^{(n)}}
    \un_{E_3^{(n)}}
    \po^{b_{\log n}}[\tau(M^-)\wedge \tau(M^+)\leq n]\big]
%    +\p\big[\big(E_5^{(n)}\big)^c\big]
\leq
%    2
    (\log n)^{-3},
\label{Ineg_Very_Very_Far}
\end{eqnarray}
where we used \eqref{Ineg_Sortie_M_pm},
which is still valid on $E_3^{(n)}$ for $n\geq n_3$ with $\widehat b(n)$ replaced by $b_{\log n}$, recalling that $n_6\geq n_3$.
%and eq. \eqref{Ineg_Proba_E4} of Lemma \ref{Lemma_Proba_E4}.

{\bf Second step:}
%Let $n\geq n_5$ and $z\in\Z$.
By reversibility (see \eqref{reversiblemeas}), we have for all $y\in\Z$, $k\in\N$
and a.s. every environment $\omega$,
$$
    \po^{b_{\log n}}(S_k=y)
=
    \po^y(S_k=b_{\log n})\frac{\mu_\omega(y)}{\mu_\omega(b_{\log n})}
\leq
    \frac{e^{-V(y)}+e^{-V(y-1)}}{e^{-V(b_{\log n})}+e^{-V(b_{\log n}-1)}}
\leq
    c_{12} e^{-[V(y)-V(b_{\log n})]}
$$
with $c_{12}:=(1+\e_0^{-1})$ by ellipticity.
Hence, recalling $M^\pm$ from \eqref{eq_def_Mpm}
and $\widehat L^\pm$ from \eqref{eq_def_L-_hat} and \eqref{eq_def_L+_hat},
conditioning by $\omega$ and applying the strong Markov property at time $\tau(b_{\log n})$,
\begin{align}
&
    \PP\big(S_n=z,
%        |b_{\log n}-z|> (\log n)^{4/3+\delta_1},
        |b_{\log n}-z|> \Gamma_n,
        \tau(b_{\log n})\leq n,
            E_3^{(n)}, E_5^{(n)}, z\in[M^-,M^+] \big)
\nonumber\\
& =
    \EE\Big[
%        \un_{\{|b_{\log n}-z|> (\log n)^{4/3+\delta_1}\}}
        \un_{\{|b_{\log n}-z|> \Gamma_n\}}
        \un_{\{\tau(b_{\log n})\leq n\}}
        \un_{E_3^{(n)}\cap E_5^{(n)}}\un_{\{z\in[M^-,M^+]\}}
        \po^{b_{\log n}}[S_k=z]_{|k=n-\tau(b_{\log n})}
    \Big]
\nonumber\\
& \leq
    \E\Big[
%        \un_{\{|b_{\log n}-z|> (\log n)^{4/3+\delta_1}\}}
        \un_{\{|b_{\log n}-z|> \Gamma_n\}}
        \un_{E_3^{(n)}\cap E_5^{(n)}}\un_{\{z\in[M^-,M^+]\}}
        c_{12} e^{-[V(z)-V(b_{\log n})]}
    \Big].
\label{Ineg_Very_Far}
\end{align}
We cut the expectation in \eqref{Ineg_Very_Far} into several parts.
We first notice that since $n\geq n_6\geq n_3$,
\begin{equation}
%&&
    \E\Big[
        \un_{\{z\in[M^-, \widehat L^-]\cup[\widehat L^+, M^+]\}}\un_{E_3^{(n)}\cap E_5^{(n)}}
        c_{12} e^{-[V(z)-V(b_{\log n})]}
    \Big]
\leq
        c_{12} (\log n)^{-C_1}
\leq
    (\log n)^{-3}
\label{Ineg_Very_Far_2}
\end{equation}
by \eqref{eq_Minimum_V_Couplage_Meeting} and \eqref{eq_Minimum_V_Couplage_Meeting_2},
and since $C_1>20$ and $\log n>2\e_0^{-1}$ because $n\geq n_6\geq n_3$.

{\bf Third step:}
Hence, there only remains to treat the case $z\in]\widehat L^-,\widehat L^+[$, which we divide into $4$
subcases, depending on the signs of $z$
and $z-b_{\log n}$.
In this step, we write $\mathcal T_i:=\theta(T_i(V,\log n))$ for $-1\leq i \leq 1$ to simplify the notation.
First, we have,
using $ T_{\mathcal T_1}(h_n)\leq  T_{\mathcal T_1}(\log n)$
and the fact that $\{z-b_{\log n}>\Gamma_n,\ b_{\log n}>0\}$ depends only on $b_{\log n}$
and so is measurable with respect to $\sigma(\mathcal T_0, x_0(V,\log n))$
in the first inequality,
using the law $\mathcal T_V^\uparrow$ of $\mathcal T_1$ and its independence with $(\mathcal T_0, x_0(V,\log n))$
conditionally on
$\{b_{\log n}> 0\}$
% or $\{b_{\log n}> 0\}$
(i.e. on $T_0(V,\log n)$ being downward)
%or downward)
by Theorem \ref{Lemma_Independence_h_extrema} {\bf (ii)} in the first equality,
then the law of $\mathcal T_V^\uparrow$ ($\mathcal T_{V,h}^\uparrow$ with $h=\log n$) by Theorem \ref{Lemma_Law_of_Slopes}
with $\Xi_{\log n}=\{ T_V(\log n)<T_V(\R_-^*)\}$ as defined in \eqref{eq_Def_A_n}
in the second inequality,
then Proposition \ref{Lemma_Laplace_V_Conditionne} in the third one
since $\Gamma_n\geq p_4$ because $n\geq n_6\geq n_3$ and $n_3\geq \exp(p_5)$, we get
\begin{eqnarray}
&&
    \E\Big[
        \un_{\{z>b_{\log n}+
        %(\log n)^{4/3+\delta_1}
        \Gamma_n
        \}}\un_{\{b_{\log n}>0\}}
        \un_{\{z\in]\widehat L^-,\widehat L^+[\}}\un_{E_3^{(n)}\cap E_5^{(n)}}
        e^{-[V(z)-V(b_{\log n})]}
    \Big]
\nonumber\\
& \leq &
    \E\Big[
    \un_{\{z-b_{\log n}>
    %(\log n)^{4/3+\delta_1}
    \Gamma_n
    \}}
    \un_{\{b_{\log n}>0\}}
    \E\Big[
        \un_{\{y_n< T_{\mathcal T_1}(\log n)\}}
        e^{-\mathcal T_1(y_n)}
    \big| \sigma(\mathcal T_0, x_0(V,\log n))
    \Big]_{y_n=z-b_{\log n}}
    \Big]
\nonumber\\
& = &
    \E\Big[
    \un_{\{z-b_{\log n}>
    %(\log n)^{4/3+\delta_1}
    \Gamma_n
    \}}
    \un_{\{b_{\log n}>0\}}
    \E\Big[
        \un_{\{y_n< T_{\mathcal T_V^\uparrow}(\log n)\}}
        e^{-\mathcal T_V^\uparrow(y_n)}
    \Big]_{y_n=z-b_{\log n}}
    \Big]
\nonumber\\
& \leq &
    \E\Big[
    \un_{\{z-b_{\log n}>
    %(\log n)^{4/3+\delta_1}
    \Gamma_n
    \}}
    \E\Big[
        \un_{\{y_n< T_{V}(\log n)\}}
        e^{-V(y_n)}
        |\Xi_{\log n}
    \Big]_{y_n=z-b_{\log n}}
%    \un_{\{b_{\log n}>0\}}
    \Big]
\nonumber\\
& \leq &
    \E\Big[
    \un_{\{z-b_{\log n}>
    (\log n)^{4/3+\delta_1}\}}
    c_{13} (z-b_{\log n})^{-3/2}
    \Big]
\leq
%    c_{13}\left((\log n)^{4/3+\delta_1}\right)^{-3/2}
%=
    c_{13}(\log n)^{-2-3\delta_1/2}.
\label{Ineg_106_Non_Central_1}
\end{eqnarray}
Also,
%writing here $\mathcal T_0:=\theta(T_0(V,\log n))$ and
using $x_0(V,\log n)=b_{\log n}<z< \widehat L^+=x_0(V,\log n)+T_{\mathcal T_0}(h_n)$ with $h_n\leq \log n$
(on the event of the second line below) and $E_5^{(n)}$
in the first inequality, we have
\begin{eqnarray}
\label{Esperance_Partie4_01}
&&
    \E\Big[
        \un_{\{z>b_{\log n}+
        %(\log n)^{4/3+\delta_1}
        \Gamma_n
        \}}\un_{\{b_{\log n}\leq0\}}
        \un_{\{z\in]\widehat L^-,\widehat L^+[\}}
        \un_{E_3^{(n)} \cap E_5^{(n)}}
        e^{-[V(z)-V(b_{\log n})]}
    \Big]
\\
& = &
    \sum_{y\leq 0}
    \E\Big[
        \un_{\{z>y+
        %(\log n)^{4/3+\delta_1}
        \Gamma_n
        \}}\un_{\{b_{\log n}=y\}}
        \un_{\{z\in]\widehat L^-,\widehat L^+[\}}\un_{E_3^{(n)} \cap E_5^{(n)}}
        e^{-[V(z)-V(y)]}
    \Big]
\nonumber\\
& \leq &
    \sum_{y\leq 0}
        \un_{\{z-y>
        %(\log n)^{4/3+\delta_1}
        \Gamma_n
        \}}
    \E\Big[
        \un_{\{b_{\log n}=y\}}
%     \un_{\{x_0(V, \log n)=y, \mathcal T_0(\ell(\mathcal T_0))>0\}}
        \un_{\{z-y<T_{\mathcal T_0}(\log n)\}}
        \un_{\{\ell(\mathcal T_0)\leq 2(\log n)^{2+\delta_1}\}}
        e^{-\mathcal T_0(z-y)}
    \Big]
%\nonumber\\
%& = &
%    \sum_{y\leq 0}
%        \un_{\{z-y> (\log n)^{4/3+\delta_1}\}}
%    \frac{
%        \E\big[
%        \sharp\big\{0\leq i < \ell(\mathcal T_V^\uparrow), \, -i=y\big\}
%        \un_{\{z-y<T_{\mathcal T_V^\uparrow}(\log n), \ell(\mathcal T_V^\uparrow)\leq 2(\log n)^{2+\delta_1}\}}
%        e^{-\mathcal T_V^\uparrow(z-y)}
%        \big]
%    }
%   {\E[\ell(\mathcal T_V^\uparrow)+\ell(\mathcal T_V^\downarrow)]}
\nonumber\\
& = &
    \sum_{y\leq 0}
        \un_{\{z-y>
        %(\log n)^{4/3+\delta_1}
        \Gamma_n
        \}}
    \frac{
        \E\Big[
        \un_{\{-\ell(\mathcal T_V^\uparrow) < y\big\}}
        \un_{\{z-y<T_{\mathcal T_V^\uparrow}(\log n)\}}\un_{\{\ell(\mathcal T_V^\uparrow)\leq 2(\log n)^{2+\delta_1}\}}
        e^{-\mathcal T_V^\uparrow(z-y)}
        \Big]
    }
    {\E\big[\ell\big(\mathcal T_V^\uparrow\big)+\ell\big(\mathcal T_V^\downarrow\big)\big]},
\nonumber
\end{eqnarray}
where we used, in the last equality, eq. \eqref{eq_Central_Slope_General_phi} of Theorem \ref{Lemma_Central_Slope}
with $\Delta_0=\{y\}$, $\Delta_1=\Z$ and $h=\log n$
and $\sharp\big\{0\leq i < \ell(\mathcal T_V^\uparrow), \, -i=y\big\}
%=\un_{\{z-y> (\log n)^{4/3+\delta_1}\}}
=\un_{\{-\ell(\mathcal T_V^\uparrow) < y\}}
$,
for which we recall that for $y\leq 0$, $b_{\log n}=y$ means that $x_0(V, \log n)=y$ and $\mathcal T_0(\ell(\mathcal T_0))>0$,
i.e. $\mathcal T_0$ is an upward slope.

Then,
%applying Lemma \ref{Lemma_Esperance_Longueur_Slope}
using the definition of $n_6$
and $y>-\ell(\mathcal T_V^\uparrow)\geq  -2(\log n)^{2+\delta_1}$
% in the first inequality,
and the law of slopes provided by {Theorem \ref{Lemma_Law_of_Slopes} \bf (i)}
%in the first equality,
in the first inequality,
and Proposition \ref{Lemma_Laplace_V_Conditionne} in the second inequality
since $\Gamma_n\geq p_4$ and $\log n \geq p_5$ because $n\geq n_6\geq n_3$,
we get,
with $c_{14}:=3c_7^{-1}c_{13}$,
%with $c$ a constant which value can change from line to line,
%{\bf (c'est ici qu'on a besoin du meme $\delta_1$ dans $E_5$ et $E_7$)}
\begin{eqnarray}
\eqref{Esperance_Partie4_01}
& \leq &
%    \frac{c_7^{-1}}{(\log n)^2}
%    \sum_{y=-\lceil 2(\log n)^{2+\delta_1}\rceil} ^0
%        \un_{\{z-y>
%        %(\log n)^{4/3+\delta_1}
%        \Gamma_n
%        \}}
%        \E\Big[
%        \un_{\{z-y<T_{\mathcal T_V^\uparrow}(\log n)\}}
%        e^{-\mathcal T_V^\uparrow(z-y)}
%        \Big]
%\nonumber\\
%& = &
    \frac{c_7^{-1}}{(\log n)^2}
    \sum_{y=-\lceil2(\log n)^{2+\delta_1}\rceil} ^0
        \un_{\{z-y>
        %(\log n)^{4/3+\delta_1}
        \Gamma_n
        \}}
        \E\big[
        \un_{\{z-y<T_V(\log n)\}}
        e^{-V(z-y)} | \Xi_{\log n}
        \big]
\nonumber\\
& \leq &
    \frac{c_7^{-1}}{(\log n)^2}
    \sum_{y=-\lceil 2(\log n)^{2+\delta_1}\rceil} ^0
        \un_{\{z-y> (\log n)^{4/3+\delta_1}\}}
        c_{13}(z-y)^{-3/2}
\nonumber\\
& \leq &
    c_7^{-1}(\log n)^{-2}
    \big[2(\log n)^{2+\delta_1}+1\big]
        c_{13}\big((\log n)^{4/3+\delta_1}\big)^{-3/2}
\leq
    c_{14}(\log n)^{-2-\delta_1/2}.~~~~
\label{Ineg_106_Central_1}
\end{eqnarray}
%We prove similarly as in \eqref{Ineg_106_Non_Central_1} and \eqref{Ineg_106_Central_1} that

%We write here $\mathcal T_{-1}:=\theta(T_{-1}(V,\log n))$,
Notice that
$
    \zeta(\mathcal T_{-1})
=
    \big(V(x_0-i)-V(x_0),\ 0\leq i \leq x_0-x_{-1}\big)
$,
with $x_j=x_j(V,\log n)$, $j\in\Z$.
Moreover,
by Theorem \ref{Lemma_Independence_h_extrema} {\bf (i)},
conditionally on $T_0(V,\log n)$ being upward, i.e. on $\{b_{\log n}\leq 0\}$,
$\zeta(\mathcal T_{-1})$ is independent of
$(\mathcal T_0, x_0(V,\log n))$
%$T_0(V,\log n)$
and has the same law as $\zeta(\mathcal T_V^\downarrow)$,
so is equal in law, by Proposition \ref{Prop_Egalite_Loi_Zeta_Slopes},  to $\mathcal T_{V^-, \log n}^{\uparrow*}$,
which law is given by Theorem \ref{Lemma_Law_of_Slopes_Right} {\bf (i)} applied to $V^-$
(with $\zeta$ defined in \eqref{eq_def_zeta}).
Using this in the second inequality,
% with $A_n$ defined in \eqref{eq_Def_A_n},
then Proposition \ref{Lemma_Laplace_V_Conditionne} in the third one, we get since $n\geq n_6$,
with
$
    \Xi_{\log n}^{*-}
:=
    \{T_{V^-}(\log n)<T_{V^-}^*(]-\infty,0])\}
$,
similarly as in \eqref{Ineg_106_Non_Central_1},
\begin{eqnarray}
&&
    \E\Big[
        \un_{\{z<b_{\log n}-
        %(\log n)^{4/3+\delta_1}
        \Gamma_n
        \}}\un_{\{b_{\log n}\leq 0\}}
        \un_{\{z\in]\widehat L^-,\widehat L^+[\}}\un_{E_3^{(n)}\cap E_5^{(n)}}
        e^{-[V(z)-V(b_{\log n})]}
    \Big]
\nonumber\\
& \leq &
    \E\Big[
        \un_{\{z<b_{\log n}-
        %(\log n)^{4/3+\delta_1}
        \Gamma_n, \,
        %\}}\un_{\{
        b_{\log n}\leq 0\}}
    \E\Big[
        \un_{\{
        %(\log n)^{4/3+\delta_1}<
        y_n< T_{\zeta(\mathcal T_{-1})}(\log n)\}}
        e^{-\zeta(\mathcal T_{-1})(y_n)}
    \big| \sigma(\mathcal T_0, x_0(V,\log n))
    \Big]_{y_n=b_{\log n}-z}
    \Big]
\nonumber\\
& \leq &
    \E\Big[
    \un_{\{b_{\log n}-z>
    %(\log n)^{4/3+\delta_1}
    \Gamma_n
    \}}
    \E\Big[
        \un_{\{y_n< T_{V^-}(\log n)\}}
        e^{-V^-(y_n)}
        |\Xi_{\log n}^{*-}
    \Big]_{y_n=b_{\log n}-z}
%    \un_{\{b_{\log n}>0\}}
    \Big]
\nonumber\\
& \leq &
    \E\Big[
    \un_{\{b_{\log n}-z>
    %(\log n)^{4/3+\delta_1}
    \Gamma_n
    \}}
    c_{13} (b_{\log n}-z)^{-3/2}
    \Big]
%\nonumber\\
%& \leq &
%    C\Big((\log n)^{4/3+\delta_1}\Big)^{-3/2}
%=
\leq
    c_{13}(\log n)^{-2-3\delta_1/2}.
\label{Ineg_106_Non_Central_1_Neg}
\end{eqnarray}
Also,
%writing here $\mathcal T_0:=\theta(T_0(V,\log n))$ and
using $x_1(V,\log n)=b_{\log n}>z> \widehat L^-=x_1(V,\log n)-T_{\zeta(\mathcal T_0)}(h_n)$ with $h_n\leq \log n$
%(on the event of the second line below)
and $E_5^{(n)}$
in the first inequality, we have
\begin{eqnarray}
\label{Esperance_Partie4_01_Bis}
&&
    \E\Big[
        \un_{\{z<b_{\log n}-
        %(\log n)^{4/3+\delta_1}
        \Gamma_n
        \}}\un_{\{b_{\log n} > 0\}}
        \un_{\{z\in]\widehat L^-,\widehat L^+[\}}
        \un_{E_3^{(n)} \cap E_5^{(n)}}
        e^{-[V(z)-V(b_{\log n})]}
    \Big]
\\
%& = &
%    \sum_{y> 0}
%    \E\Big[
%        \un_{\{z<y-
%        %(\log n)^{4/3+\delta_1}
%        \Gamma_n
%        \}}\un_{\{b_{\log n}=y\}}
%        \un_{\{z\in]\widehat L^-,\widehat L^+[\}}\un_{E_3^{(n)} \cap E_5^{(n)}}
%        e^{-[V(z)-V(y)]}
%    \Big]
%\nonumber\\
& \leq &
    \sum_{y> 0}
        \un_{\{z-y< -
        %(\log n)^{4/3+\delta_1}
        \Gamma_n
        \}}
    \E\Big[
        \un_{\{b_{\log n}=y\}}
%     \un_{\{x_0(V, \log n)=y, \mathcal T_0(\ell(\mathcal T_0))>0\}}
        \un_{\{y-z<T_{\zeta(\mathcal T_0)}(\log n)\}}
        \un_{\{\ell(\mathcal T_0)\leq 2(\log n)^{2+\delta_1}\}}
        e^{-\zeta(\mathcal T_0)(y-z)}
    \Big]
%\nonumber\\
%& = &
%    \sum_{y\leq 0}
%        \un_{\{z-y> (\log n)^{4/3+\delta_1}\}}
%    \frac{
%        \E\big[
%        \sharp\big\{0\leq i < \ell(\mathcal T_V^\uparrow), \, -i=y\big\}
%        \un_{\{z-y<T_{\mathcal T_V^\uparrow}(\log n), \ell(\mathcal T_V^\uparrow)\leq 2(\log n)^{2+\delta_1}\}}
%        e^{-\mathcal T_V^\uparrow(z-y)}
%        \big]
%    }
%   {\E[\ell(\mathcal T_V^\uparrow)+\ell(\mathcal T_V^\downarrow)]}
\nonumber\\
& = &
    \sum_{y> 0}
        \un_{\{z-y<    -
        %(\log n)^{4/3+\delta_1}
        \Gamma_n
        \}}
    \frac{
        \E\Big[
        \un_{\{y\leq \ell(\mathcal T_V^\downarrow) \big\}}
        \un_{\{y-z<T_{\zeta(\mathcal T_V^\downarrow)}(\log n),\ \ell(\mathcal T_V^\downarrow)\leq 2(\log n)^{2+\delta_1}\}}
        e^{-\zeta(\mathcal T_V^\downarrow)(y-z)}
        \Big]
    }
    {\E\big[\ell\big(\mathcal T_V^\uparrow\big)+\ell\big(\mathcal T_V^\downarrow\big)\big]},
\nonumber
\end{eqnarray}
where we used, in the last equality, eq. \eqref{eq_Central_Slope_General_phi} of Theorem \ref{Lemma_Central_Slope}
with $\Delta_0=\Z$, $\Delta_1=\{y\}$  and $h=\log n$
and $\sharp\big\{0\leq i < \ell(\mathcal T_V^\downarrow), \, \ell(\mathcal T_V^\downarrow)-i=y\big\}
%=\un_{\{z-y> (\log n)^{4/3+\delta_1}\}}
=\un_{\{y\leq \ell(\mathcal T_V^\downarrow)\}}
$ for $y>0$,
for which we recall that for $y> 0$,
$b_{\log n}=y$ means that $x_1(V, \log n)=y$ and $\mathcal T_0(\ell(\mathcal T_0))<0$.

Then, using the definition of $n_6$,
%applying Lemma \ref{Lemma_Esperance_Longueur_Slope} and
$y\leq \ell(\mathcal T_V^\downarrow)\leq 2(\log n)^{2+\delta_1}$
and Proposition \ref{Prop_Egalite_Loi_Zeta_Slopes}
in the first inequality,
then Theorem \ref{Lemma_Law_of_Slopes_Right} {\bf (i)}  in the equality,
and Proposition \ref{Lemma_Laplace_V_Conditionne} in the second inequality,
we get since
$\log n \geq p_5$ and $\Gamma_n\geq p_4$ because $n\geq n_6\geq n_3$,
%$n\geq n_6\geq \exp(p_5)$ and so $\Gamma_n\geq p_4$,
concluding as in \eqref{Ineg_106_Central_1},
\begin{align}
& \eqref{Esperance_Partie4_01_Bis}
 \leq
    \frac{c_7^{-1}}{(\log n)^2}
    \sum_{y=1}^{\lfloor 2(\log n)^{2+\delta_1}\rfloor }
        \un_{\{z-y<-
        %(\log n)^{4/3+\delta_1}
        \Gamma_n
        \}}
        \E\Big[
        \un_{\{y-z<T_{\mathcal T_{V^-}^{\uparrow*}}(\log n)\}}
        e^{-\mathcal T_{V^-}^{\uparrow*}(y-z)}
        \Big]
\nonumber\\
& =
    \frac{c_7^{-1}}{(\log n)^2}
    \sum_{y=1}^{\lfloor 2(\log n)^{2+\delta_1}\rfloor }
        \un_{\{y-z>
        %(\log n)^{4/3+\delta_1}
        \Gamma_n
        \}}
        \E\big[
        \un_{\{y-z<T_{V^-}(\log n)\}}
        e^{-V^-(y-z)} | \Xi_{\log n}^{*-}
        \big]
%\nonumber\\
%& \leq &
%    \frac{c_7^{-1}}{(\log n)^2}
%    \sum_{y=1}^{\lfloor  2(\log n)^{2+\delta_1}\rfloor }
%        \un_{\{y-z>
%        %(\log n)^{4/3+\delta_1}
%        \Gamma_n
%        \}}
%        c_{13}(y-z)^{-3/2}
%\nonumber\\
%& \leq &
%    \frac{c_7^{-1}}{(\log n)^2}
%    \big[2(\log n)^{2+\delta_1}+1\big]
%        c_{13}\big((\log n)^{4/3+\delta_1}\big)^{-3/2}
%\leq
%    c(\log n)^{-2-\delta_1/2}.
\leq
    \frac{c_{14}}{(\log n)^{2+\delta_1/2}}.
\label{Ineg_106_Central_1_Bis}
\end{align}
Combining %\eqref{Ineg_Very_Far},
\eqref{Ineg_106_Non_Central_1},
\eqref{Ineg_106_Central_1},
%\eqref{Ineg_106_Non_Central_2}
\eqref{Ineg_106_Non_Central_1_Neg}
and
\eqref{Ineg_106_Central_1_Bis}
%\eqref{Ineg_106_Central_2}
ensures that, with $c_{15}:=2c_{14}+2c_{13}$,
$$
    \E\Big[
        \un_{\{|b_{\log n}-z|>
        %(\log n)^{4/3+\delta_1}
        \Gamma_n\}}
        \un_{E_3^{(n)}\cap E_5^{(n)}}\un_{\{z\in]\widehat L^-,\widehat L^+[\}}
        e^{-[V(z)-V(b_{\log n})]}
    \Big]
\leq
    c_{15}(\log n)^{-2-\delta_1/2}.
$$
%{\bf (ou 2e methode : by symmetry (preciser ?), replacing left $(\log n)$-minima by right $(\log n)$-minima and consequently $A_n$ by $A_n^*$)}
%\begin{equation}
%    \E\Big[
%        \un_{\{z<b_{\log n}- (\log n)^{4/3+\delta_1}\}}\un_{\{b_{\log n}\leq 0\}}
%        \un_{\{z\in]\widehat L^-,\widehat L^+[\}}\un_{E_3^{(n)}\cap E_5^{(n)}}
%        e^{-[V(z)-V(b_{\log n})]}
%    \Big]
%\leq
%    \frac{C}{(\log n)^{2+3\delta_1/2}}.
%\label{Ineg_106_Non_Central_2}
%\end{equation}
%\begin{equation}
%    \E\Big[
%        \un_{\{z<b_{\log n}- (\log n)^{4/3+\delta_1}\}}\un_{\{b_{\log n} > 0\}}
%        \un_{\{z\in]\widehat L^-,\widehat L^+[\}}
%        \un_{E_3^{(n)} \cap E_5^{(n)}}
%        e^{-[V(z)-V(b_{\log n})]}
%    \Big]
%\leq
%    c(\log n)^{-2-\delta_1/2}.
%\label{Ineg_106_Central_2}
%\end{equation}
This, combined with \eqref{Ineg_Very_Far_2}, proves that the right hand side of
\eqref{Ineg_Very_Far} is $\leq c_{16}(\log n)^{-2-\delta_1/2}$ for all $n\geq n_6$ and $z\in\Z$
with $c_{16}:=(c_{15}c_{12}+1)$.
This together with \eqref{Ineg_Very_Very_Far} gives \eqref{Ineg_Proba_z_far_b_atteint}
since $\delta_1\in(0,2/3)$, with $c_{11}:= c_{16}+1$.
\hfill$\Box$

%%%%%%%%%%%%%%%%%%%%%%%%%%%%%%%%%%%%

%\subsubsection{Case with $z< b_{\log n}-(\log n)^{4/3+\delta_1}$, $b_{\log n}\leq 0$ and $\tau(b_{\log n})>n$}

\subsubsection{Case with $b_{\log n}$ far from $z$, without subvalleys and small valleys when $\tau(b_{\log n})>n$}
The aim of this subsection is to prove the following lemma.

\begin{lem} \label{Lemma_Proba_z_far_b_non_atteint}
There exists a constant $c_{17}>0$ such that, for all $n\geq n_6$ and all $z\in\Z$,
\begin{equation}
    \PP\big(S_n=z,  |b_{\log n}-z|>
    %(\log n)^{4/3+\delta_1},
    \Gamma_n,
    \tau(b_{\log n})> n, E_3^{(n)}, E_5^{(n)} \big)
\leq
    %c_{17}\frac{(\log_2 n)^3}{(\log n)^3}.
    c_{17}(\log_2 n)^3(\log n)^{-3}.
\label{Ineg_Proba_z_far_b_non_atteint}
\end{equation}
\end{lem}

We start with the case $b_{\log n}\leq 0$.
We first make the following simple remark.

\begin{lem}\label{Lemma_Proba_Nulle}
We have,
%for every $n\geq 3$,
\begin{equation}\label{eq_Proba_Nulle}
    \forall n\geq 3,\,
    \forall z\in\Z,
\qquad
    \PP\big(S_n=z,\, b_{\log n}\leq 0,\,
    z< b_{\log n}-
    %(\log n)^{4/3+\delta_1}
    \Gamma_n,\,
    \tau(b_{\log n})>n
    \big)
=
    0.
\end{equation}
\end{lem}

\noindent{\bf Proof:}
On $\{b_{\log n}\leq 0,\, z< b_{\log n}-
%(\log n)^{4/3+\delta_1}
\Gamma_n,\, \tau(b_{\log n})>n\}$,
we have $z<b_{\log n}\leq 0$, so for $S$ starting from $0$ (under $\po$ or $\PP$),
$\tau(z)>\tau(b_{\log n})>n$ and thus $S_n\neq z$.
This leads to \eqref{eq_Proba_Nulle}.
\hfill$\Box$

In order to prove Lemma \ref{Lemma_Proba_z_far_b_non_atteint},
we also have to give an upper bound for the probability of
$F_1^{(n)}(z)$, where
%{\bf (integrer $E_3$ a $F_1$ ? Voire $E_5$ ?)}
%{\bf (le renommer $F_1^-$ ?)}
$$
    F_1^{(n)}(z)
:=
    \big\{ S_n=z,\, b_{\log n}\leq 0,\, z > b_{\log n}+
    %(\log n)^{4/3+\delta_1}
    \Gamma_n,\, \tau(b_{\log n})>n\big\}\cap E_3^{(n)}\cap E_5^{(n)}.
$$

Loosely speaking, on $E_3^{(n)}$ by Remark \ref{Remark_x_i}, there are no subvalleys of height larger than $h_n-C_1\log_2 n$
in the $(\log n)$-central valley $[M^-, M^+]$ and in the two neighbor  valleys (of height at least $\log n$)
on its left and on its right, and the height of these three valleys is quite larger than $\log n$.
In particular, we prove:

\begin{lem} \label{Lemma_Proba_F1_E2}
There exists a constant $c_{18}>0$ such that
\begin{equation}\label{Ineg_Proba_F1_E2}
    \forall n\geq n_6, \,
    \forall z\in\Z,
\qquad
    \PP\big[F_1^{(n)}(z)\big]
\leq
    c_{18}(\log_2 n)^3(\log n)^{-3}.
\end{equation}
\end{lem}

\begin{figure}[htbp]
\includegraphics[width=15.95cm,height=6.83cm]{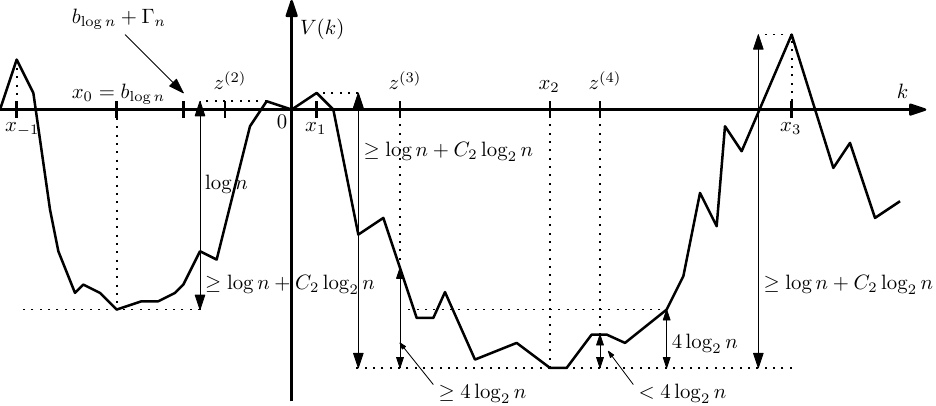}
\caption{Schema of the potential $V$  with $x_i=x_i(V,\log n)$,
and $z$ equal to $z^{(2)}$, $z^{(3)}$ and $z^{(4)}$ respectively for
step $2$, $3$ \big(on $F_1^{(n)}(z)$\big) and $4$ \big(on $F_2^{(n)}(z)$\big).
}
\label{figure_Lemma_5_5}
\end{figure}

\noindent {\bf Outline of the proof:}
See Figure \ref{figure_Lemma_5_5} for a schema of the potential.
Assume for example that $b_{\log n}\leq 0$, so $x_0=b_{\log n}$, with $x_i:=x_i(V,\log n)$, $i\in\Z$,
and that $F_1^{(n)}(z)$ holds.
Since $\tau(x_0)>n$, we first prove that, by Lemma \ref{Lemma_Proba_Descente}, with large probability, $\tau(x_2)\leq n$.
Second, if $z$
%($=S_n$)
is not in the valley $[x_1, x_3]$, then after first hitting $x_2$,
$S$ has to leave this valley before time $n$ (so that $S_n=z\notin[x_1,x_3]$),
which has negligible probability
since the height of this valley $[x_1,x_3]$ is quite larger than $\log n$ on $E_3^{(n)}$.
Third, if $z$ belongs to the valley $[x_1, x_3]$ with $V(z)\geq V(x_2)+4\log_2 n$,
then the probability that $S_n=z$ is negligible by reversibility, which we can
apply to $S$ started at $x_2$ by strong Markov property.
Finally, if $z$ belongs to the valley $[x_1, x_3]$ with $V(z)< V(x_2)+4\log_2 n$,
then $V(z+.)-V(z)$ goes up $\log n$ before going down $4\log_2 n $ on the left and on the right,
and conditionally on $(V(k), \ k\geq 0)$,
$
    \max_{[x_0,0]} V-\max_{[0, z]} V
=
    \max_{[x_0,0]} V-V(x_1)\in [-9\log_2 n, 0[
$ (otherwise $\tau(x_2)<\tau(x_0)$ would have small probability which would contradict our first step).
Since all these three conditions have probability less than $c(\log_2 n)(\log n)^{-1}$ for some $c>0$ with
some independence, this last case is also negligible compared to $(\log n)^{-2}$.
We now prove this rigorously.

\noindent {\bf Proof:}
Let $n\geq n_6$ and $z\in\Z$.
In all the proof, we write $x_i$ for $x_i(V,\log n)$ for every $i\in\Z$.

{\bf First step:}
Applying Lemma \ref{Lemma_Proba_Descente} with $h=\log n$, $\xi_2=1$,
$a=x_0<b=x_1<c=x_2$
(so that {\bf (i)} is satisfied for $\omega\in E_-^{(n)}$), $\xi_1=2C_1$
(so {\bf(ii)} and {\bf (iii)} are satisfied since there is
no left $(\log n-2C_1\log_2 n)$-extremum in $]x_0$, $x_1[$ nor in $]x_1, x_2[$
for $\omega\in E_-^{(n)}\cap E_3^{(n)}$ by Remark \ref{Remark_x_i}),
$\alpha=3$ (so {\bf (iv)} is satisfied for $\omega\in E_5^{(n)}$ since $0<\delta_1<2/3$) and $x=0$, we get
since $n\geq n_6\geq n_3$ and so $\log n\geq \widehat h_2(2C_1,1)$,
\begin{equation}\label{Ineg_Proba_Sortie_x1x3}
    \forall \omega \in E_-^{(n)}\cap E_3^{(n)}\cap E_5^{(n)},
\qquad
    \po\big[\tau(x_0)\wedge \tau(x_2) \geq n \big]
\leq
    (\log n)^{-4}.
\end{equation}
As a consequence, using
$
    \tau(x_0)
=
    \tau(b_{\log n})
>n
$
on $F_1^{(n)}(z)$,
we get
\begin{equation}
    \PP\big[F_1^{(n)}(z)\cap \{\tau(x_2) \geq n\} \big]
\leq
    \E\big[
        \un_{E_-^{(n)}\cap E_3^{(n)}\cap E_5^{(n)}}
    \po\big(\tau(x_0)\wedge \tau(x_2) \geq n \big)
    \big]
\leq
    (\log n)^{-4}.
\label{Ineg_Proba_F1E2_1}
\end{equation}

{\bf Second step:}
There only remains  to consider $F_1^{(n)}(z)\cap \{\tau(x_2) < n\}$.
This second step focuses on the case $z\notin]x_1, x_3[$.
We start with the case $z\leq x_1$
\big(see Figure \ref{figure_Lemma_5_5} with $z=z^{(2)}$\big).
In what follows we prove that in this case, the probability that, after hitting $x_2$, $S$
goes or goes back to $z\in]-\infty, x_1]$ before time $n$ is negligible.

To this aim, using $z\leq x_1< x_2$,
then \eqref{InegProba2} and ellipticity
%\eqref{eqEllipticity}
\eqref{eq_ellipticity_for_V}
in the second inequality,
we have %for large $n$,
for every $\omega\in E_-^{(n)}\cap E_3^{(n)}\cap\{z\leq x_1\}$ and  $k\in\{0,\dots, n\}$,
\begin{eqnarray*}
    \po^{x_2}[S_k=z]
& \leq &
    \po^{x_2}[\tau(x_1)\leq \tau(z)\leq k]
\\
& \leq &
    (k+1) \e_0^{-1} \exp(-H[T_1(V,\log n)])
\leq
    2\e_0^{-1} (\log n)^{-C_2}
\leq
    (\log n)^{-4}
\end{eqnarray*}
since
$
    V[x_1]-\min_{[x_1, x_2]}V
=
    H[T_1(V,\log n)]
\geq
    \log n +C_2 \log_2 n
$
on $E_3^{(n)}$, $C_2>9$ and $n\geq n_6\geq n_3$.
Hence, conditioning by $\omega$ then applying the strong Markov property at time $\tau(x_2)$,
\begin{eqnarray}
&&
    \PP\big[F_1^{(n)}(z)\cap \{\tau(x_2) < n\}\cap\{z\leq x_1\}\big]
\nonumber\\
& \leq &
    \EE\Big[\un_{E_-^{(n)}\cap E_3^{(n)}\cap\{z\leq x_1\}\cap \{\tau(x_2) < n\}}
        \po^{x_2}[S_k=z]_{|k=n-\tau(x_2)}
    \Big]
\nonumber\\
& \leq &
    (\log n)^{-4}.
\label{Ineg_Proba_F1_Gauche}
\end{eqnarray}
Similarly, using \eqref{InegProba1} instead of \eqref{InegProba2}, we have for large $n$,
\begin{equation}
    \PP\big[F_1^{(n)}(z)\cap \{\tau(x_2) < n\}\cap\{z\geq x_3\}\big]
\leq
    (\log n)^{-4}.
\label{Ineg_Proba_F1_Droite}
\end{equation}

\noindent{\bf Third step:}
Now,  on $\{x_1<z<x_3\}\cap \{V(z)\geq V(x_2)$
$+4\log_2 n\}\cap E_-^{(n)}$
\big(see Figure \ref{figure_Lemma_5_5} with $z=z^{(3)}$\big),
we have by reversibility (see \eqref{reversiblemeas}) and ellipticity
\eqref{eq_ellipticity_for_V},
%\eqref{eqEllipticity},
for $k\in\N$,
$$
    \po^{x_2}(S_k=z)
\leq
    \frac{\mu_\omega(z)}{\mu_\omega(x_2)}
\leq
    \e_0^{-1} \exp[V(x_2)-V(z)]
\leq
    \frac{\e_0^{-1}}{(\log n)^{4}}.
$$
As a consequence, once more conditioning by $\omega$ and applying the strong Markov property,
proceeding as in \eqref{Ineg_Proba_F1_Gauche},
\begin{eqnarray}
&&
    \PP\big[F_1^{(n)}(z)\cap \{\tau(x_2) < n\}
    \cap\{x_1<z<x_3\}\cap \{V(z)\geq V(x_2)+4\log_2 n\}\big]
\nonumber\\
& \leq &
    \EE\Big[\un_{E_-^{(n)}\cap\{x_1<z<x_3\}\cap \{V(z)\geq V(x_2)+4\log_2 n\}
    \cap \{\tau(x_2) < n\}}
        \po^{x_2}[S_k=z]_{|k=n-\tau(x_2)}
    \Big]
\nonumber\\
& \leq &
    \e_0^{-1} (\log n)^{-4}.
\label{Ineg_Proba_F1_Milieu_1}
\end{eqnarray}

{\bf Forth step:}
Finally, we study
\big(see Figure \ref{figure_Lemma_5_5} with $z=z^{(4)}$\big),
\begin{eqnarray*}
    F_2^{(n)}(z)
& := &
    F_1^{(n)}(z)\cap \{\tau(x_2) < n\}
    \cap
    \{x_1<z<x_3\}
    \cap \{V(z)< V(x_2)+4\log_2 n\}.
\end{eqnarray*}
This set is empty for $z<0$ because $x_1>0$, so we can assume that $z\geq 0$.

We once more define
$
    V_{z}^{\pm}(k)
:=
    V(z\pm k)-V(z)
$, $k\in\Z$, and notice that $V_{z}^{-}$ and $V_{z}^{+}$ are independent.
We also introduce
\begin{equation*}
    E_{11}^{(n)}
:=
    \{\tau(x_2)< \tau(x_0)\},
\qquad
    E_{12}^{(n)}
:=
    \Big\{\max_{[x_0, 0]} V \leq V(x_1)-9 \log_2 n\Big\}.
\end{equation*}
%We now write $x_i=x_i(V,\log n)$ in the rest of this proof.
We have by \eqref{probaatteinte}, for large $n$, for all $\omega\in E_-^{(n)}\cap E_5^{(n)}\cap E_{12}^{(n)}$,
\begin{equation*}
    \po\big(E_{11}^{(n)}\big)
=
%    \po[\tau(x_1)<\tau(x_0)]
%\\
%& = &
    \frac{\sum_{i=x_0}^{-1}e^{V(i)}}{\sum_{i=x_0}^{x_2-1}e^{V(i)}}
\leq
    \frac{(\log n)^{2+\delta_1}\exp[\max_{[x_0, 0]}V]}{\exp[V(x_1)]}
\leq
    (\log n)^{-6}.
\end{equation*}
%\begin{eqnarray*}
%    \po\big(E_{11}^{(n)}\big)
%\leq
%    x_0 \exp[\max\nolimits_{[x_0, 0]}V-V(x_1)]
%\leq
%    (\log n)^{-6}.
%\end{eqnarray*}
Consequently, since
$
%F_2^{(n)}(z) \subset
F_1^{(n)}(z)\cap \{\tau(x_2) < n\} \subset E_{11}^{(n)}$,
\begin{equation}
    \PP\big[F_2^{(n)}(z)\cap E_{12}^{(n)}\big]
\leq
    \PP\big[E_-^{(n)}\cap E_5^{(n)}\cap E_{12}^{(n)}\cap
    E_{11}^{(n)}
    %\big(E_{11}^{(n)}\big)^c
    \big]
\leq
    (\log n)^{-6}.
\label{Ineg_Proba_F2_E9}
\end{equation}
There remains to study $\PP\big[F_2^{(n)}(z)\cap \big(E_{12}^{(n)}\big)^c\big]$.
For a process $(v(k), \ k\in\N)$ and $y\in\N$,
we define $v_y(.):=v(y+.)-v(y)$
and
\begin{eqnarray*}
    E_{13}^{(n)}(v)
& := &
    \{T_v(]-\infty,-\log n])<T_v(9\log_2 n)\},
\\
    E_{14}^{(n)}(z)
& := &
    \big\{T_{V_{z}^+}(\log n) < T_{V_{z}^+}(-4\log_2 n)\big\}
\cap
    \big\{T_{V_{z}^-}(\log n) < T_{V_{z}^-}(-4\log_2 n)\big\}
\nonumber\\
&&
%\quad
    \cap
    \big\{T_{V_{z}^-}(\log n) \leq z
    %\wedge T_{V_{z}^-}(-4\log_2 n)
    \big\}.
\end{eqnarray*}
Also for $a\geq 0$, let
$
    V_{1,a}(k)
:=
    V^-\big[k+T_{V^-}(a)\big]-V^-\big[T_{V^-}(a)\big]
$, $k\in\N$.
We claim that
\begin{eqnarray}
&&
    F_2^{(n)}(z)\cap \big(E_{12}^{(n)}\big)^c
\subset
    E_{14}^{(n)}(z)\cap E_{13}^{(n)}\big(V_{1,\max(0,\max_{[0,z]}V-9\log_2 n)}\big).
\label{Inclusion_F2E9}
\end{eqnarray}
Indeed on $F_2^{(n)}(z)\cap \big(E_{12}^{(n)}\big)^c $, we have $b_{\log n}=x_0\leq 0$,
$z\in]x_1,x_3[$,
$\min_{[x_1,x_3]} V=V(x_2)> V(z)-4\log_2 n$,
and $V(x_1)\geq V(x_2)+\log n+C_2\log_2 n\geq V(z)+\log n$ due to $E_3^{(n)}$
and since $C_2>9$, the same being true also for $V(x_3)$ instead of $V(x_1)$.
So $V_{z}^\pm$ hits $[\log n,+\infty[$ before $]-\infty, -4\log_2 n]$,
%which leads to
so $F_2^{(n)}(z)\cap \big(E_{12}^{(n)}\big)^c $ is included in
the first two sets in
$E_{14}^{(n)}(z)$.

Also on $F_2^{(n)}(z)$, $x_0=b_{\log n}\leq 0$, thus $\max_{[x_0,x_2]}V=V(x_1)$, so $\max_{[0,z]}V=V(x_1)$ if $x_1<z\leq x_2$.
Assume now that $x_2<z<x_3$ and $F_2^{(n)}(z)$ holds. If $\max_{[x_2, z]}V>V(x_1)$, then
$\min\{u\in[x_2,z], \ V(u)=\max_{[x_2, z]}V\}$ would be a left $(\log n)$-maximum
(because its potential would be greater than $V(x_1)\geq V(x_2)+\log n+C_2\log_2 n\geq V(z)+\log n$
due to $E_3^{(n)}$ and $C_2>9$ as before, and greater than $V(x_2)+\log n$),
belonging to $]x_2,x_3[$, which is not possible, so $\max_{[x_2, z]}V\leq V(x_1)$.
Hence $\max_{[0,z]}V=V(x_1)\geq V(z)+\log n$ in both cases,
so $\max_{[0,z]} V_{z}^-\geq \log n$,
thus $F_2^{(n)}(z)\cap \big(E_{12}^{(n)}\big)^c $ is included in the third set in
$E_{14}^{(n)}(z)$.

Finally on $F_2^{(n)}(z)\cap \big(E_{12}^{(n)}\big)^c $,
we have $\max_{[x_0,0]}V<V(x_1)= \max_{[x_0,x_2]}V$ by definition of the $x_i$ and of $E_-^{(n)}$,
and $\max_{[x_0,0]}V> V(x_1)-9\log_2 n\geq V(x_0)+\log n$ by definition of $\big(E_{12}^{(n)}\big)^c$ and
since $H[T_0(V,\log n)]\geq \log n+C_2\log_2 n$ with $C_2>9$ on $E_3^{(n)}$.
Also, we just proved that $\max_{[0,z]}V=V(x_1)$.
Hence, starting from $0$, $V^-$ first hits $[\max_{[0,z]}V-9\log_2 n, +\infty[$, then goes down at least
$\log n$ before $|x_0|$ and so before going up $9\log_2 n$, so $\omega\in E_{13}^{(n)}\big(V_{1,a}\big)$
with $a=\max_{[0,z]}V-9\log_2 n$ if $\max_{[0,z]}V-9\log_2 n\geq 0$.
Otherwise, $\max_{[x_0,0]}V<V(x_1)=\max_{[0,z]}V< 9\log_2 n$,
with $V(x_0)= V(x_1)-H[T_0(V,\log n)]\leq -\log n +(9-C_2)\log_2 n<-\log n$ since $C_2>9$ due to $E_3^{(n)}$ and $E_5^{(n)}$,
so $\omega\in E_{13}^{(n)}\big(V^-\big)=E_{13}^{(n)}\big(V_{1,a}\big)$ with $a=0$.
So \eqref{Inclusion_F2E9} is proved in every case.

We have in particular, by \eqref{eqOptimalStopping2}, since $n\geq n_6\geq n_3$ so $\log_2 n>C_0$,
\begin{equation}
    \p[E_{13}^{(n)}(V^-)]
=
    \p\big(T_{V^-}(]-\infty,-\log n])<T_{V^-}([9\log_2 n,+\infty[)\big)
\leq
    10(\log_2 n)(\log n)^{-1}.
\label{Ineg_Proba_E10}
\end{equation}
%uniformly (in $y$) for large $n$ by \eqref{eqOptimalStopping2}, the same being true for $V_-$.
Also, using first the independence between $V_{z}^+$ and $V_{z}^-$, which have the same law as $V$ and $V^-$ respectively,
then applying \eqref{eqOptimalStopping2} again,
we have since $n\geq n_6\geq n_3$,
\begin{eqnarray}
    \p\big(E_{14}^{(n)}(z)\big)
& \leq &
    \p\big[T_{V_{z}^+}(\log n) < T_{V_{z}^+}(-4\log_2 n)\big]
    \p\big[
    T_{V_{z}^-}(\log n) < T_{V_{z}^-}(-4\log_2 n)\big]
\nonumber\\
& \leq &
    25(\log_2 n)^2(\log n)^{-2}.
\label{Ineg_Proba_E20}
\end{eqnarray}
Hence using \eqref{Inclusion_F2E9}, then conditioning by $V^+=(V(k),\ k\geq 0)$,
noting that $E_{14}^{(n)}(z)$ and $\max_{[0,z]}V$ depend only on
$V^+$ and for every $a\in\R_+$,
 $E_{13}^{(n)}\big(V_{1,a}\big)$ only on $V^-$, which is independent of $V^+$ and has the same law as $V_{1,a}$,
then applying \eqref{Ineg_Proba_E10} and \eqref{Ineg_Proba_E20}, we get
\begin{eqnarray*}
    \PP\big[F_2^{(n)}(z)\cap \big(E_{12}^{(n)}\big)^c\big]
& \leq &
    \p\big[E_{14}^{(n)}(z)\cap E_{13}^{(n)}\big(V_{1,\max(0,\max_{[0,z]}V-9\log_2 n)}\big)\big]
\\
& = &
    \E\big[\un_{E_{14}^{(n)}(z)}\p\big[E_{13}^{(n)}\big(V_{1,\max(0,\max_{[0,z]}V-9\log_2 n)}\big) | V^+\big]\big]
\\
& = &
    \E\big[\un_{E_{14}^{(n)}(z)}\p\big[E_{13}^{(n)}\big(V_{1,a}\big) \big]_{|a=\max(0,\max_{[0,z]}V-9\log_2 n)}\big]
\\
& \leq &
    10(\log_2 n)(\log n)^{-1}\p(E_{14}^{(n)}(z))
\leq
    250(\log_2 n)^3(\log n)^{-3}.
\end{eqnarray*}
This, together with \eqref{Ineg_Proba_F2_E9} gives
$
    \p\big(F_2^{(n)}(z)\big)
\leq
    251(\log_2 n)^3(\log n)^{-3}
$
for all $n\geq n_6$ and $z\in\Z$.
% since $\log n>1$ and $\log_2 n.

\noindent{\bf Conclusion:}
Combining this with \eqref{Ineg_Proba_F1E2_1}, \eqref{Ineg_Proba_F1_Gauche}, \eqref{Ineg_Proba_F1_Droite},
\eqref{Ineg_Proba_F1_Milieu_1} proves \eqref{Ineg_Proba_F1_E2}.
\hfill$\Box$

\noindent {\bf Proof of Lemma \ref{Lemma_Proba_z_far_b_non_atteint}:}
We prove, similarly as in Lemmas \ref{Lemma_Proba_Nulle} and \ref{Lemma_Proba_F1_E2}
(replacing in particular $x_0$, $x_1$, $x_2$ and $x_3$
respectively by $x_1$, $x_0$, $x_{-1}$ and $x_{-2}$ respectively in its proof,
nearly by symmetry)
that for every $n\geq n_6$ and every $z\in\Z$,
\begin{align*}
    \PP\big(S_n=z, b_{\log n}> 0, z> b_{\log n}+
    \Gamma_n,
    %(\log n)^{4/3+\delta_1},
    \tau(b_{\log n})>n
    \big)
& =
    0,
\\
    \PP\big(S_n=z, b_{\log n}> 0, z <
    b_{\log n}-
    %(\log n)^{4/3+\delta_1},
    \Gamma_n,
    \tau(b_{\log n})>n,E_3^{(n)}, E_5^{(n)}\big)
& \leq
    c_{18}(\log_2 n)^3(\log n)^{-3}.
\end{align*}
Combining this with Lemmas \ref{Lemma_Proba_Nulle} and \ref{Lemma_Proba_F1_E2}
proves Lemma \ref{Lemma_Proba_z_far_b_non_atteint} with $c_{17}:=2c_{18}$.
\hfill$\Box$

\noindent{\bf Proof of Proposition \ref{Lemma_Proba_z_b_far}:}
This proposition follows directly from Lemmas \ref{Lemma_Proba_z_far_b_atteint}
and \ref{Lemma_Proba_z_far_b_non_atteint}
with $c_{10}:=c_{11}+c_{17}$, since $(\log_2 n)^3\leq (\log n)^{1/2}$
for $n\geq n_6\geq n_3$ and $\delta_1\in]0, 2/3[$.
% and $n_6\geq n_3$.
\hfill$\Box$

%%%%%%%%%%%%%%%%%%%%%%%%%%%%%%%%%%%%%%%%%%%%%%%%%%%%%%%%%%%%%%%%%%%%%

%%%%%%%%%%%%%%%%%%%%%%%%%%%%%%%%%%%%%%%%%%%%%%%%%%%%%%%%%%%%%%%%%%%%%

%%%%%%%%%%%%%%%%%%%%%%%%%%%%%%%%%%%%%%%%%%%%%%%%%%%%%%%%%%%%%%%%%%%%%%%%%%%%%%%%%%%%%%

%%%%%%%%%%%%%%%%%%%%%%%%%%%%%%%%%%%%%%%%%%%%%%%%%%%%%%%%%%%%%%%%%%%%%%%%%%%%%%%%%%%%%%

%\subsection{Case with $b_{\log n}$ far from $z$, with subvalleys or small valleys}

\subsection{Case with at least one subvalley or small valley}

We now focus on the case where some of the valleys (of height $\geq \log n$)
close to the origin can be small (i.e. with height $<\log n +C_2\log_2 n$),
or can contain subvalleys of height less than but close to $\log n$.
%, which corresponds to the event $(E_3^{(n)})^c$.
More precisely, the aim of this subsection is to prove the following estimate.

\begin{prop}\label{Prop_Local_Limit_Si_Petite_Slope}
There exists $n_9\geq n_6$ and $c_{19}>0$ such that
$$
    \forall n\geq n_9,\, \forall z\in\Z,
\qquad
    \PP\big(S_n=z, (E_3^{(n)})^c\big)
\leq
    c_{19}(\log_2 n)^3 (\log n)^{-3}.
$$
\end{prop}

This case can be divided into many different subcases. For example, there can be, or not,  a subvalley of height close to $\log n$
inside the $(\log n)$-central valley, either at the right or at the left of $b_{\log n}$, or there can even be two such subvalleys.
There can also exist, close to the $(\log n)$-central valley,
one or two valleys with height close to $\log n$,
larger or smaller than $\log n$,
which can trap the random walk $(S_k)_k$ for some time.
Also, the height of the $(\log n)$-central valley can be close to $\log n$,
which can enable $S$ to escape it before time $n$ with not so small quenched probability.
Taking into account the indexes of the left $(h_n-C_1\log_2 n)$-slopes considered, i.e. with height less than $\log n+C_2\log_2 n$,
and their height, larger or smaller than $\log n$,
the indexes $i$ of the first left $h_n$-minimum $b_i(V,h_n)$ (defined in \eqref{eq_def_bi})
visited by $S$ before time $n$, of the second one etc,
the fact that $z$ is close or far from these left $h_n$-extrema,
this makes dozens of cases.
However we will combine together some of these cases,
for example with the help of Lemma \ref{Lem_Proba_zn_Proche_Fond} and of the notation $\I_k$
defined in \eqref{eq_def_Ik} below.

On $\big(E_3^{(n)}\big)^c$, there exists some $i\in\{-10,\dots 10\}$ such that
$H[T_i(V,h_n-C_1\log_2 n)]<\log n+C_2 \log_2 n$.
%First, if no such $j$ belongs to $\{-1,0,1\}$, we should be able to do the same proof as in the previous subsection
Also, we prove that with large probability, there are no more than two such $i$.
To this aim, we define
$$
    E_{15}^{(n)}
:=
    \big\{
        \sharp\{i\in\Z,\ -99\leq i \leq 99,\ H[T_i(V,h_n-C_1\log_2 n)]<\log n+C_2 \log_2 n \} \leq 2
    \big\}.
$$
More precisely, we prove the following estimate.

\begin{lem}\label{Lemma_Proba_E11_Petites_Slopes}
%As $n\to+\infty$,  {\bf (ceci ne depend pas de $z$)}
There exist $n_7\geq n_6$ and $c_{20}>0$ such that,
\begin{equation}\label{Ineg_Proba_E11_Petites_Slopes}
    \forall n\geq n_7,
\qquad
    \p\big[\big(E_{15}^{(n)}\big)^c\big]
\leq
    c_{20}(\log_2 n)^3 (\log n)^{-3}.
%    O\big( (\log_2 n)^3 (\log n)^{-3}\big).
\end{equation}
\end{lem}

\noindent {\bf Proof:}
Due to Lemma \ref{Lemma_Loi_excess_height_h_extrema}, we have
$
    \p\big[E_{16}^{(n)}(i) \mid b_{\widetilde h_n}\leq 0 \big]
=
    O\big((\log_2 n)(\log n)^{-1}\big)
$,
$i\in\Z$,
where
$
    \widetilde h_n
=
%    \log n-C_1\log_2 n
    h_n-C_1\log_2 n
$ as before and
$$
    E_{16}^{(n)}(i)
:=
    \big\{
        H[T_i(V,h_n-C_1\log_2 n)]<\log n+C_2 \log_2 n
    \big\},
\qquad
    i\in\Z.
$$
Hence, using the independence of the translated left $\widetilde h_n$-slopes conditionally on $\{b_{\widetilde h_n}\leq 0\}$
(see Theorem \ref{Lemma_Independence_h_extrema} {\bf (i)}),
we have
\begin{eqnarray*}
    \p\big[\big(E_{15}^{(n)}\big)^c \mid  b_{\widetilde h_n}\leq 0 \big]
& = &
    \p\big(\cup_{-99\leq i_1 < i_2 < i_3 \leq 99}
        E_{16}^{(n)}(i_1)\cap E_{16}^{(n)}(i_2)\cap E_{16}^{(n)}(i_3) \mid  b_{\widetilde h_n}\leq 0 \big)
\\
& \leq &
    \sum_{-99\leq i_1 < i_2 < i_3 \leq 99} \prod_{k=1}^3
    \p\big[E_{16}^{(n)}(i_k) \mid  b_{\widetilde h_n}\leq 0 \big]
=
    O\big((\log_2 n )^3(\log n)^{-3}\big)
\end{eqnarray*}
as $n\to+\infty$.
We prove similarly the same inequality with $b_{\widetilde h_n}\leq 0$ replaced by $b_{\widetilde h_n}> 0$,
which proves the lemma.
\hfill$\Box$

We define, for $h>0$ and $i\in\Z$
(this definition being different from that of \cite{DGP_Collision_Sinai}),
\begin{equation}\label{eq_def_bi}
    b_i(V,h)
:=
\left\{
\begin{array}{ll}
x_{2i}(V,h)  & \text{ if } x_0(V,h) \text{ is a left $h$-minimum},\\
x_{2i-1}(V,h) & \text{ otherwise}.\\
\end{array}
\right.
\end{equation}
So, the $b_i(V,h)$, $i\in\Z$, are the left $h$-minima for $V$, such that  $b_0(V,h)\leq 0<b_1(V,h)$
and $b_i(V,h)<b_{i+1}(V,h)$, $i\in\Z$.
We also denote by $M_i(V,h)$ the unique left $h$-maximum for $V$ between $b_i(V,h)$ and $b_{i+1}(V,h)$.
Hence, $M_i(V,h)=x_{j+1}(V,h)$ if $b_i(V,h)=x_j(V,h)$.

We now prove that the probability that $z$ is "close" (in terms of potential)
to the bottom $b_j(V,h_n)$ of a valley of height $h_n$
and that $\omega\in \big(E_3^{(n)}\big)^c$
is small. More precisely, we define, for $h>0$,
$$
    E_{17}^{(n)}(j, h, z)
:=
    \{M_{j-1}(V,h)\leq z\leq M_j(V,h),\ V(z)\leq V[b_j(V,h)]+4\log_2 n\},
\quad\
    j\in\Z.
$$

We now have the following lemma, which is useful to prove Lemma \ref{Lemma_Proba_Vallee_k_entreDi}
(in which we take $h_n'=h_n$)
and Lemma \ref{Lem_Fond_Velle_k_Droite}
(in which we take $h_n'=\widetilde h_n$)
and then Lemma \ref{Lemma_Valle_k_E2c}.

\begin{lem}\label{Lem_Proba_zn_Proche_Fond}
There exist $c_{21}>0$ and $n_8\geq n_7$ such that,
whether $h'_n=h_n$ or $h'_n=\widetilde h_n:=h_n-C_1\log_2 n$, we have
$$
    \forall n\geq n_8,\,
    \forall z\in\Z,
\qquad
    \p\big[\big( E_3^{(n)}\big)^c \cap \cup_{j=-8}^8 E_{17}^{(n)}(j, h'_n, z)\big]
\leq
    c_{21}(\log_2 n)^3(\log n)^{-3}.
$$
\end{lem}

Loosely speaking, in the case $h_n'=\widetilde{h}_n$,
on $E_{17}^{(n)}(j, \widetilde{h}_n, z)$, $V_{z}^+$ and $V_{z}^-$ go up $\widetilde{h}_n-4\log_2 n$ before going down $-4\log_2 n$,
which has probability $O\big (\log_2 n)^2(\log n)^{-2}\big)$.
Also, on $E_3^{(n)}$ one of the left $\widetilde h_n$-slopes around the origin has an excess height less than some $C\log_2 n$,
which has probability $O\big((\log_2 n)(\log n)^{-1}\big)$, with some independence, which leads to
Lemma \ref{Lem_Proba_zn_Proche_Fond} in the case $h_n'=\widetilde{h}_n$, the second case being
nearly a consequence of the first one.
We now prove this rigorously.

\noindent{\bf Proof of Lemma \ref{Lem_Proba_zn_Proche_Fond}:}
Let $n\geq n_7$ and $z\in\Z$.
We start with the case $h'_n=\widetilde h_n$.
On the one hand, we notice that for $-13\leq j \leq 13$,
on $\big(E_3^{(n)}\big)^c \cap E_{17}^{(n)}\big(j, \widetilde h_n, z\big)$,
$z$ belongs to the support $\big[x_k\big(V, \widetilde h_n\big), x_{k+1}(V, \widetilde h_n\big)\big]$
of a left $\widetilde h_n$-slope
$T_{k}':=T_k\big(V, \widetilde h_n\big)$ with $2j-2\leq k \leq 2j$,
the value of $k$ depending on $x_0\big(V, \widetilde h_n\big)$ being a left $\widetilde h_n$-maximum or minimum for $V$ and
on $z\leq b_j\big(V, \widetilde h_n\big)$
or $z> b_j\big(V, \widetilde h_n)$, with
$
    T_k'(z)-\inf_{y\in[x_k(V, \widetilde h_n), x_{k+1}(V, \widetilde h_n)]} T_k'(y)
\leq
    4\log_2 n
$.
%On the other hand, as in Remark \ref{Remark_x_i}, on $E_3^{(n)}$,
%$x_i(V, h_n')= x_i\big(V, \widetilde h_n\big)$ for every $-9\leq i \leq 10$,
%so $H[T_i(V, h_n')]\geq \log n +C_2\log_2 n$ for every $-10\leq i \leq 10$,
%since $H[T_i(V, h_n')]=H\big[T_i\big(V, \widetilde h_n\big)\big]\geq \log n +C_2\log_2 n$ if $|i|\leq 9$,
%whereas for $i=10$ (resp. $-10$), the support of $T_i\big(V, \widetilde h_n\big)$
%is included in the one of $T_i(V, h_n')$ because $\widetilde h_n \leq h_n'$
%and $x_i(V, h_n')=x_i\big(V, \widetilde h_n\big)$
%\big(resp. $x_{i+1}(V, h_n')=x_{i+1}\big(V, \widetilde h_n\big)$\big),
%so $H[T_i(V, h_n')]\geq H\big[T_i\big(V, \widetilde h_n\big)\big]\geq \log n +C_2\log_2 n$.
Hence, using $x_i\big(V, \widetilde h_n\big)=x_{i-k}\big(V_{z},\widetilde h_n\big)+z$, $i\in\Z$
on $\big\{x_k\big(V, \widetilde h_n\big) \leq z < x_{k+1}\big(V, \widetilde h_n\big) \big\}$
and the definition of $E_3^{(n)}$,
we get
\begin{eqnarray*}
&&
    \p\big[\big( E_3^{(n)}\big)^c \cap \cup_{j=-13}^{13} E_{17}^{(n)}\big(j, \widetilde h_n, z\big)\big]
\\
& \leq &
    \p\Big(
        \cup_{k=-28}^{27} \big\{x_k\big(V, \widetilde h_n\big) \leq z < x_{k+1}\big(V, \widetilde h_n \big) \big\}
            \cap\Big\{T_k'(z)-\inf_{[0, \ell(T_k')]} \theta(T_k')
                        \leq 4\log_2 n\Big\}
\\
    &&
\qquad
    \cap \cup_{i=-10}^{10} \big\{ H\big[T_i\big(V, \widetilde h_n\big)\big]< \log n +C_2\log_2 n\big\}
    \Big)
\\
& \leq &
    \p\Big(
%        \Big\{\inf_{[0, \ell(T_0)]} T_0 \geq -4\log_2 n\Big\}_{|T_0=T_0(V_{z}, h_n')}
        \Big\{\inf_{[x_0(V_{z}, \widetilde h_n), x_1(V_{z}, \widetilde h_n)]} V_{z} \geq -4\log_2 n\Big\}
    \cap \cup_{j=-37}^{38} \big\{ H\big[T_j\big(V_{z}, \widetilde h_n\big)\big]< \log n +C_2\log_2 n\big\}
    \Big),
\end{eqnarray*}
where $V_{z}$ has the same law as $V$, so the last probability does not depend on $z$.

Now, notice that, with $V^\pm=(V(\pm y), \ y \in\N)$ as before,
and
$
    \widetilde V_3(k)
:=
    V\big[k+T_V\big(\big[\widetilde h_n-4\log_2 n, +\infty\big[\big)\big]
%-
%    V[T_V([h_n'-4\log_2 n, +\infty[)]
$, $k\in\N$,
we have
\begin{equation}\label{Ineg_Probas_E18_0}
    E_{18, 0}^{(n)}
    %\big(\widetilde h_n\big)
    \cap\{b_{\widetilde h_n}\leq 0\}
\subset
    E_{19, +}^{(n)}
    %\big (\widetilde h_n\big)
\cap
    E_{19, -}^{(n)}
    %\big (\widetilde h_n\big)
\cap
    E_{20}^{(n)},
    %\big(\widetilde h_n \big),
\end{equation}
where for $i\in\Z$ and $h>0$,
\begin{eqnarray*}
    E_{18, i}^{(n)}
    %\big (\widetilde h_n\big )
& := &
    \big\{\inf\nolimits_{[x_0(V, \widetilde h_n), x_1(V, \widetilde h_n)]} V \geq -4\log_2 n,\,
        %b_{h_n'}\leq 0, \,
        H\big [T_i\big(V, \widetilde h_n\big )\big ]< \log n +C_2\log_2 n
    \big\},
\\
    E_{19, \pm}^{(n)}
    %\big (\widetilde h_n\big)
& := &
    \{
        T_{V^{\pm}}\big(\big[\widetilde h_n -4\log_2 n, +\infty\big[\big)
        <
        T_{V^{\pm}}(]-\infty, -4\log_2 n[)
    \},
\\
    E_{20}^{(n)}
    %\big(\widetilde h_n\big)
& := &
    \{
        T_{\widetilde V_3}\big(\big]-\infty, \log n+C_2\log_2 n-\widetilde h_n\big[\big)
        <
        T_{\widetilde V_3}\big(\big[\log n+C_2\log_2 n, +\infty\big[\big)
    \}.
\end{eqnarray*}
Using \eqref{eqOptimalStopping2} and $n\geq n_7\geq n_3$, we have
$
%    \p\big[E_{19, \pm}^{(n)}(h_n')\big]
    \p\big(E_{19, \pm}^{(n)}\big)
\leq
    10(\log_2 n)(\log n)^{-1}
$
% because $n\geq n_7\geq n_3$,
%which has, by \eqref{eqOptimalStopping2}, a probability
and
$
%    \p\big[E_{20}^{(n)}(h_n')\big]
    \p\big(E_{20}^{(n)}\mid V(k),\, k\leq T_V\big(\big[\widetilde h_n-4\log_2 n, +\infty\big[\big)\big)
$
$
\leq
    \big(\log n+C_2\log_2 n -\big(\widetilde h_n-4\log_2 n\big)+C_0\big)
    \big(\widetilde h_n+C_0\big)^{-1}
\leq
    2(2C_1+C_2+5)(\log_2 n)(\log n)^{-1}
$.
%because $n\geq n_7\geq n_3$.
Hence, using \eqref{Ineg_Probas_E18_0},
%using first the strong Markov property at time $T_V\big(\big[\widetilde h_n-4\log_2 n, +\infty\big[\big)$
conditioning by $\sigma\big(V(k),\, k\leq T_V\big(\big[\widetilde h_n-4\log_2 n, +\infty\big[\big)$
then using the independence of $V^+$ and $V^-$, we
have, with $c_{22}:=200(2C_1+C_2+5)$,
\begin{equation*}
    \p\big(
        E_{18, 0}^{(n)},
%        (h_n'),
        \, b_{\widetilde h_n}\leq 0
    \big)
\leq
%%    \p\big[E_{19, -}^{(n)}(h_n')\big]
%    \p\big(E_{19, -}^{(n)}\big)
%%    \p\big[E_{19, +}^{(n)}(h_n')\big]
%    \p\big(E_{19, +}^{(n)}\big)
%%    \p\big[E_{20}^{(n)}(h_n')\big]
%    \p\big(E_{20}^{(n)}\big)
%\leq
    c_{22}(\log_2 n)^3(\log n)^{-3}.
\end{equation*}
We get similarly  the same result with $b_{\widetilde h_n}\leq 0$ replaced by $b_{\widetilde h_n}> 0$.

Finally, for $i\neq 0$, using Theorem \ref{Lemma_Independence_h_extrema} {\bf (i)}
since
$
    H\big[\theta\big(T_i\big(V, \widetilde h_n\big)\big)\big]
=
    H\big[T_i\big(V, \widetilde h_n\big)\big]
$,
\begin{align*}
&
    \p\big(
        E_{18, i}^{(n)},\,
        %(h_n'),\,
        b_{\widetilde h_n}\leq 0
    \big)
\\
& =
    \p\Big(
        \inf_{[x_0(V, \widetilde  h_n), x_1(V, \widetilde h_n)]} V \geq -4\log_2 n,\, b_{\widetilde h_n}\leq 0
    \Big)
    \p\big(
        H\big[T_i\big(V, \widetilde h_n\big)\big]< \log n +C_2\log_2 n \mid b_{\widetilde h_n}\leq 0
%        H[\theta(T_i(V, h_n'))]< \log n +C_2\log_2 n, \, b_{h_n'}\leq 0
    \big)
\\
%& \leq &
%    25(\log_2 n)^2(\log n)^{-2}
%    .12(\log n +C_2\log_2 n -h_n'+C_0)(h_n')^{-1}
%\\
& \leq
    200c_8(2C_1+C_2+C_0)(\log_2 n)^3(\log n)^{-3}
\end{align*}
for large $n$
since the first probability in the second line is
$
\leq
%    \p\big[E_{19, -}^{(n)}(h_n')\big]
%    \p\big[E_{19, +}^{(n)}(h_n')\big]
    \p\big(E_{19, -}^{(n)}\big)
    \p\big(E_{19, +}^{(n)}\big)
$
and the second one is
$\leq c_8\big(\log n +C_2\log_2 n -\widetilde h_n\big)\big(\widetilde h_n\big)^{-1}$
for large $n$ by Lemma \ref{Lemma_Loi_excess_height_h_extrema}.
We get similarly  the same result with $b_{\widetilde  h_n}\leq 0$ replaced by $b_{\widetilde  h_n}> 0$,
using Theorem \ref{Lemma_Independence_h_extrema} {\bf (ii)} instead of {\bf (i)}.
Thus,
there exists some $c_{23}>0$ and some $n_8\geq n_7$ such that
$    \p\big(
        E_{18, i}^{(n)}
        %(h_n')
    \big)
\leq
    c_{23} (\log_2 n)^3(\log n)^{-3}
$
for all $n\geq n_8$ and all $-37\leq i \leq 38$.
%whereas $h'_n=h_n$ or $h'_n=\widetilde h_n=h_n-C_1\log_2 n$.

Finally, for all $n\geq n_8$ for all $z\in\Z$,
\begin{equation}\label{Ineg_Proba_E17_cas_1}
    \p\big[\big( E_3^{(n)}\big)^c \cap \cup_{j=-13}^{13} E_{17}^{(n)}\big(j, \widetilde  h_n, z\big)\big]
\leq
    \sum_{i=-37}^{38}
    \p\big(
        E_{18, i}^{(n)}
        %(h_n')
%        \inf_{[x_0(V, h_n'), x_1(V, h_n')]} V \geq -4\log_2 n, H[T_i(V, h_n')]< \log n +C_2\log_2 n
    \big)
\leq
    76c_{23} (\log_2 n)^3(\log n)^{-3},
%    10^5(3C_1+C_2+C_0)(\log_2 n)^3(\log n)^{-3}
\end{equation}
which proves the lemma in the case $h'_n=\widetilde h_n$.

We now turn to the case $h_n'=h_n$. Let $z\in\Z$.
To this aim, we introduce some notation, which will also be useful in the proof of Lemma \ref{Lemma_Nombre_Vallees_Visitees} below.
For $j\in\Z$, let
\begin{equation}
\label{eq_def_lambda_j}
    \Lambda_j
:=
    \sharp\big\{k\in\Z, \ x_k\big(V, \widetilde h_n\big)\in[x_j(V, h_n), x_{j+1}(V, h_n)[\big\},
\end{equation}
which belongs to $(2\N+1)$ since left $\widetilde h_n$-maxima and minima alternate and $\widetilde h_n<h_n$.
If for $j\in\Z$, $\Lambda_j=2k+1$ with $k>1$, then
$[x_j(V, h_n), x_{j+1}(V, h_n)[
$
$
=
\big[x_{\ell}\big(V, \widetilde h_n\big),$ $x_{\ell+2k+1}\big(V, \widetilde h_n\big)\big[
$
for some $\ell\in\Z$.
Also for each $0\leq i < k$,
$H\big[T_{\ell+2i+1}\big(V, \widetilde h_n\big)\big]<h_n$,
otherwise, if moreover $x_{\ell}\big(V, \widetilde h_n\big)$ is a left $h_n$-minimum (resp. maximum), then
$
    \widetilde u
:=
    \min\big\{u\in\big[x_\ell\big(V, \widetilde h_n\big), x_{\ell+2i+1}\big(V, \widetilde h_n\big)\big],
        \, V(u)=\max_{[x_\ell(V, \widetilde h_n)\leq u \leq x_{\ell+2i+1}(V, \widetilde h_n)]} V
    \big\}
$
would be a left $h_n$-extremum
\big(since
$
    V\big(\widetilde u\big)
%\geq
%    V\big[x_{\ell+2i+1}\big(V, \widetilde h_n\big)\big]
\geq
    V\big[x_{\ell+2i+2}\big(V, \widetilde h_n\big)\big]
$
$
    +h_n
\geq
    V\big[x_{\ell}\big(V, \widetilde h_n\big)\big]+h_n
$
\big),
belonging to $]x_j(V, h_n), x_{j+1}(V, h_n)[$,
which is not possible
(resp. similar argument with $\max$ replaced by $\min$).

Hence on $E_{15}^{(n)}$, $\Lambda_0\leq 5$, otherwise the support of $T_0(V, h_n)$ would contain the support of at least
$(\Lambda_0-1)/2\geq 3$ slopes $T_p\big(V, \widetilde h_n\big)$ with height
$H\big[T_p\big(V, \widetilde h_n\big)\big]< h_n < \log n+C_2\log_2 n$, with at least three of them such that
$|p|\leq 5$, which is not possible on $E_{15}^{(n)}$.
Also, notice that %by induction,
for $j\geq 1$,
$
    x_j(V, h_n)
=
x_{\ell}\big(V, \widetilde h_n\big)
$
with
$1\leq \ell\leq \Lambda_0+\dots+\Lambda_{j-1}$.
Thus by induction,
$\Lambda_j\leq 5$ for every $0\leq j \leq 17$,
for which we use for $0<j\leq 17$ the same argument as for $\Lambda_0$
with $1\leq p\leq \Lambda_0+\dots+\Lambda_{j-1}+5$ ($\leq 5(j+1)\leq 90$ by hypothesis of induction).
Similarly on $E_{15}^{(n)}$,
$\Lambda_j\leq 5$ for every $-17\leq j \leq 0$,
and so for  every $-17\leq j \leq 17$.

Consequently, for the same reasons, on $E_{15}^{(n)}$, if for $-17\leq j \leq 17$, $\Lambda_j=3$
(resp. $\Lambda_j=5$), then the support of $T_j(V,h_n)$ contains the support of at least one (resp. at least two)
slope(s) $T_p\big(V, \widetilde h_n\big)$ with height
$H\big[T_p\big(V, \widetilde h_n\big)\big]< h_n < \log n+C_2\log_2 n$ with $|p|<99$.
Thus, $\Lambda_0+\dots +\Lambda_j\leq j+5$ for every $0\leq j \leq 17$
and $\Lambda_j+\dots +\Lambda_0\leq |j|+5$ for every $-17\leq j \leq 0$.

Notice that for $-8\leq j \leq 8$,
on $E_{15}^{(n)}\cap E_{17}^{(n)}(j, h_n, z)$, we have
$b_j(V,h_n)=x_k(V,h_n)$ with $k\in\{2j-1,2j\}
%\subset\{-17,\dots,16\}
$ by \eqref{eq_def_bi},
so $b_j(V,h_n)=x_\ell\big(V,\widetilde h_n\big)$ with
$1\leq\ell\leq \Lambda_0+\dots+\Lambda_{2j-1}\leq (2j-1)+5\leq 20$ if $1\leq j\leq 8$,
and $|\ell|\leq \Lambda_{2j-1}+\dots+\Lambda_0\leq |2j-1|+5\leq 22$ if $-8\leq j\leq 0$, using the previous paragraph.
So $b_j(V,h_n)=x_\ell\big(V,\widetilde h_n\big)=b_{j_0}\big(V,\widetilde h_n\big)$
with $\ell\in\{2j_0-1, 2j_0\}$, thus $-11\leq j_0\leq 10$.
Consequently, there exists $j_1\in\Z$ such that
$
    z
\in
    \big[M_{{j_1}-1}\big(V,\widetilde h_n\big),M_{j_1}\big(V,\widetilde h_n\big)\big]
\subset
    [M_{j-1}(V,h_n),M_j(V,h_n)]
\subset
    \big[x_{\ell-5}\big(V,\widetilde h_n\big),x_{\ell+5}\big(V,\widetilde h_n\big)\big]
\subset
    \big[x_{-27}\big(V,\widetilde h_n\big),x_{25}\big(V,\widetilde h_n\big)\big]
$
(so that $-13\leq j_1\leq 13$),
with
$
    V(z)
\leq
    V[b_j(V,h_n)]+4\log_2 n
$
$
\leq
    V\big[b_{j_1}\big(V,\widetilde h_n\big)\big]+4\log_2 n
$,
so the conditions defining $E_{17}^{(n)}\big(j_1, \widetilde h_n, z\big)$ are satisfied.
Hence,
\begin{align*}
&
    \p\big[\big( E_3^{(n)}\big)^c \cap \cup_{j=-8}^8 E_{17}^{(n)}(j, h_n, z)\big]
\\
\leq &
    \p\big[\big(E_{15}^{(n)}\big)^c\big]
    +
    \p\big[\big( E_3^{(n)}\big)^c \cap E_{15}^{(n)}
    \cap \cup_{j_1=-13}^{13} E_{17}^{(n)}\big(j_1, \widetilde  h_n, z\big)\big]
\leq
    (c_{20}+76c_{23})(\log_2 n)^3 (\log n)^{-3}
\end{align*}
by Lemma \ref{Lemma_Proba_E11_Petites_Slopes} and \eqref{Ineg_Proba_E17_cas_1} since $n\geq n_8\geq n_7$,
which proves the lemma when $h'_n=h_n$.
\hfill$\Box$

%\begin{equation}\label{Ineg_Proba_zn_E2_Compl}
%    \PP\big[S_n=z,  \big(E_3^{(n)}\big)^c\big]
%\leq
%    o\big((\log n)^{-2}\big).
%\end{equation}

We now introduce some notation. Recall that $\tau[b_i(V,h_n)]<\infty$ $\PP$-a.s. for every $i\in\Z$
since $S=(S_k)_k$ is $\p$-almost surely recurrent.
We define by induction
\begin{eqnarray}
    \I_1
& := &
    {\bf 1}_{\{\tau[b_1(V,h_n)]<\tau[b_0(V,h_n)]\}},
\nonumber\\
    \I_k
& := &
    \sum_{\ell\in\Z-\{\I_j,\, 1\leq j<k\}} \ell
    \prod_{i\in\Z,\, i\notin\{\I_j,\, 1\leq j<k\}\cup\{\ell\}} {\bf 1}_{\{\tau[b_\ell(V,h_n)] <\tau[b_i(V,h_n)]\}},
\qquad
    k\geq 2.
\label{eq_def_Ik}
\end{eqnarray}
In words, $\I_1$ is the index $\ell$ of the first $b_\ell[V,h_n]$ visited by $S$, so that
$\I_1=0$ if $\tau[b_0(V,h_n)]<\tau[b_1(V,h_n)]$ and $\I_1=1$ if $\tau[b_1(V,h_n)]<\tau[b_0(V,h_n)]$,
which are the only possible cases since $b_0(V,h_n)\leq 0=S_0<b_1(V,h_n)$ $\PP$-a.s.
Similarly, $\I_2$ is the index $\ell$ of the second $b_\ell(V,h_n)$ visited by $S$, so $\I_2\neq \I_1$, and more generally
$\I_k$ is for $k\in\N^*$ the index $\ell$ of the $k$-th $b_\ell(V,h_n)$ visited by $S$, so that
$\I_k\notin \{\I_1,\I_2, \dots, \I_{k-1}\}$.
Notice that $\tau[b_{\I_1}(V,h_n)]=\tau[b_0(V,h_n)]\wedge \tau[b_1(V,h_n)]$ is a stopping time under $\po$ with the natural filtration of $S$,
and more generally $\tau[b_{\I_k}(V,h_n)]$ is a stopping time for every $k\geq 1$.

Recall  that $0\in[b_0(V,h_n), b_1(V,h_n)[$, that $b_0(V,h_n)$ and $b_1(V,h_n)$ are consecutive left $h_n$-minima,
and $M_0(V,h_n)$ is the only left $h_n$-maximum between them.
So, applying Lemma \ref{Lemma_Proba_Descente} with $h=\log n$,  $\xi_2=1/10$,
$a=b_0(V,h_n)<b=M_0(V,h_n)<c=b_1(V,h_n)$ which satisfy {\bf (i)} due to the previous remark,
$\xi_1=C_1$ so that {\bf (ii)} and {\bf (iii)} are satisfied since there is no left $(h_n=\log n-C_1\log_2 n)$-extremum
in $]M_0(V,h_n), b_1(V,h_n)[$ nor in $]b_0(V,h_n), M_0(V,h_n)[$,
$\alpha=3$ (so that {\bf (iv)} is satisfied for $\omega \in E_5^{(n)}$,
since $|x_i(V,h_n)|\leq |x_i(V,\log n)|$ for every $i\in\Z$ and $\delta_1<2/3$) and $x=0$,
we get for $n\geq n_8$ (which implies that $n\geq n_3$ so $\log n\geq \widehat h_2(C_1,1/10)$),
for almost all $\omega\in E_5^{(n)}$,
$$
    \po\big[\tau(b_{\I_1}(V,h_n))\geq n/10\big]
=
    \po\big[\tau(b_0(V,h_n))\wedge \tau(b_1(V,h_n))\geq n/10\big]
\leq
    (\log n)^{-4}.
$$
Consequently, using Lemma \ref{Lemma_Proba_E4}, for $n\geq \max(n_8, p_3)=:n_9$,
\begin{equation}\label{Ineg_Proba_Atteinte_bI1}
    \PP\big[\tau(b_{\I_1}(V,h_n))\geq n/10\big]
\leq
    \PP\big[\tau(b_{\I_1}(V,h_n))\geq n/10, E_5^{(n)}\big]
+
    \p\big[\big(E_5^{(n)}\big)^c\big]
\leq
    2(\log n)^{-3}.
\end{equation}

%The aim of the present subsection is to prove the following proposition, which says, loosely speaking, that
%if there is a small slope, among the slopes around the origin, then the probability that $S_n=z$ is negligible.

We now prove several lemmas which are useful to prove Proposition \ref{Prop_Local_Limit_Si_Petite_Slope}.
In what follows, for $i\in\Z$, we write
%$b_i$
%$\textsl{b}_i$
%$\textsf{b}_i$
$\bb_i$
%{\bf (conflit notation avec $b_h$ donc remplace par $\bb_i$)}
and $M_i$ respectively for $b_i(V,h_n)$ and $M_i(V,h_n)$ (which are defined in and after \eqref{eq_def_bi}).
We first prove that, with large enough probability,
$S$ only visits up to $3$ different $\bb_i$ before time $n$:

\begin{lem}\label{Lemma_Nombre_Vallees_Visitees}
There exists $c_{24}>0$ such that,
%{\bf (ceci ne depend pas de $z$)}
$$
    \forall n\geq n_9,
\qquad
    \PP\big[\tau(\bb_{\I_4})\leq n\big]
\leq
    c_{24}(\log_2 n)^3(\log n)^{-3}.
$$
\end{lem}

The main idea is that, loosely speaking, on $E_{15}^{(n)}$,
$S$ has to cross, before $\tau(\bb_{\I_4})$, at least one slope with height at least $\log n +C_2\log_2 n$,
which takes more than $n$ units of time with large probability.
We now prove this rigorously.

\noindent{\bf Proof of Lemma \ref{Lemma_Nombre_Vallees_Visitees}:}
Let $n\geq n_9$.
First, for every $1\leq k \leq 3$, using $\bb_{\I_k}<M_{\I_k}<\bb_{\I_{k}+1}   \leq \bb_{\I_{k+1}}$ when $\I_k<\I_{k+1}$ in the first inequality,
then conditioning by $\omega$ then applying the strong Markov property at time $\tau(\bb_{\I_k})$ in the following line,
and finally \eqref{InegProba1} and ellipticity in the last line, we have
\begin{eqnarray}
    p_{1,k,n}
& := &
    \PP\big[\tau(\bb_{\I_k}, \bb_{\I_{k+1}})\leq n,
    V(M_{\I_k})-V(\bb_{\I_k})\geq \log n+C_2\log_2 n, I_{k+1}>\I_k
    \big]
\nonumber\\
& \leq &
    \PP\big[\tau(\bb_{\I_k}, M_{\I_k})\leq n,
    V(M_{\I_k})-V(\bb_{\I_k})\geq \log n+C_2\log_2 n\big]
\nonumber\\
& = &
    \EE\big[{\bf 1}_{\{V(M_{\I_k})-V(\bb_{\I_k})\geq \log n+C_2\log_2 n\}}
    \po^{\bb_{\I_k}}[\tau(M_{\I_k})\leq n]
    \big]
\nonumber\\
& \leq &
    2\e_0^{-1}(\log n)^{-C_2}
\leq
    (\log n)^{-4}
\label{Ineg_p1kn}
\end{eqnarray}
since $C_2>9$ and $n\geq n_9\geq n_3$.
Similarly, using \eqref{InegProba2} instead of \eqref{InegProba1}
and $\bb_{\I_{k+1}}\leq \bb_{\I_{k}-1}<M_{\I_k-1}<\bb_{\I_k}$ when $I_{k+1}<\I_k$, we have
\begin{eqnarray}
    p_{2,k,n}
& := &
    \PP\big[\tau(\bb_{\I_k}, \bb_{\I_{k+1}})\leq n,
    V(M_{\I_k-1})-V(\bb_{\I_k})\geq \log n+C_2\log_2 n, I_{k+1}<\I_k
    \big]
\nonumber\\
& \leq &
    (\log n)^{-4}
\label{Ineg_p2kn}
\end{eqnarray}
for every $1\leq k \leq 3$
since $C_2>9$ and $n\geq n_9\geq n_3$.

We now prove that on $E_{15}^{(n)}$,
\begin{equation}\label{ineq_majoration_nombre_petits_extrema}
    \sharp\{-6\leq j \leq 6, \ H[T_j(V,h_n)]< \log n +C_2\log_2 n \}
\leq
    2.
\end{equation}
%the number of $-6\leq j \leq 6$ such that
%$H[T_j(V,h_n)]< \log n +C_2\log_2 n $ is less than or equal to $2$.
To this aim, we use \eqref{eq_def_lambda_j} and the following paragraphs.
We claim that on $E_{15}^{(n)}$, if for some $-6\leq j \leq 6$, $H[T_j(V,h_n)]< \log n +C_2\log_2 n $,
then the support of $T_j(V,h_n)$
contains at least the support of one $T_k\big(V, \widetilde h_n\big)$ with
$H\big[T_k\big(V, \widetilde h_n\big)\big]<\log n +C_2\log_2 n$
with $|k|<99$.
Indeed, on $E_{15}^{(n)}$, $x_j(V,h_n)=x_k\big(V, \widetilde h_n\big)$
with $|k|\leq \Lambda_0+\dots+\Lambda_j\leq 35$
and the support of $T_j(V,h_n)$ contains
the support of $T_{k}\big(V, \widetilde h_n\big)$,
so $H\big[T_k\big(V, \widetilde h_n\big)\big]\leq H[T_j(V,h_n)]<\log n+C_2\log_2n$.
Since there are at most two slopes $H\big[T_k\big(V, \widetilde h_n\big)\big]$, $|k|\leq 35$  with height $<\log n +C_2\log_2 n$
on $E_{15}^{(n)}$, there are at most two $j\in\{-6,\dots, 6\}$
such that $H[T_j(V,h_n)]< \log n +C_2\log_2 n $,
which proves \eqref{ineq_majoration_nombre_petits_extrema}.

Also,
%$\bb_{\I_1}$ is equal to some $x_i(V,h_n)$
%with $i\in\{-1,0,1,2\}$ (see \eqref{eq_def_bi})
$\{\I_j,\ 1\leq j \leq k\}\subset \{1-k,\dots, k\}$
%$1-k\leq \I_k\leq k $
for every $k\in\N^*$
by induction,
since for $k\geq 2$, $\min\{\I_j,\, 1\leq j<k\}-1 \leq \I_k\leq \max\{\I_j,\, 1\leq j<k\}+1$,
because $S$ only makes $\pm 1$ jumps.
So by \eqref{eq_def_bi},
$\bb_{\I_k}=x_{i_k}(V,h_n)$ with $i_k\in\{-5,\dots, 6\}$ when $1\leq k\leq 3$.
Hence, each height $V(M_{\I_k})-V(\bb_{\I_k})$ or $V(M_{\I_k-1})-V(\bb_{\I_k})$ with $k\in\{1,2,3\}$ is equal to
some $H[T_j(V,h_n)]$ with $|j|\leq 6$, so at most two of them
are less than $\log n+C_2\log_2 n$ on $E_{15}^{(n)}$ by \eqref{ineq_majoration_nombre_petits_extrema}.

Hence, for $n\geq n_9$, using \eqref{Ineg_p1kn} and \eqref{Ineg_p2kn} in the last inequality,
\begin{align}
&
    \PP\big[\tau(\bb_{\I_4})\leq n, E_{15}^{(n)}\big]
\nonumber\\
& \leq
    \PP\big[\cap_{j=1}^3 \{\tau(\bb_{\I_j}, \bb_{\I_{j+1}})\leq n\} \cap
        \cup_{k=1}^3\big(\{V(M_{\I_k})-V(\bb_{\I_k})\geq \log n+C_2\log_2 n, \I_{k+1}>\I_k\}
\nonumber\\
&
\qquad\qquad\qquad\qquad\qquad\qquad\qquad
                 \cup\{V(M_{\I_k-1})-V(\bb_{\I_k})\geq \log n+C_2\log_2 n, \I_{k+1}<\I_k\}\big)
    \big]
\nonumber\\
& \leq
    \sum_{k=1}^3 (p_{1,k,n}+p_{2,k,n})
\leq
    6(\log n)^{-4}.
\label{eq_avant_numero_bIk}
\end{align}
This together with Lemma \ref{Lemma_Proba_E11_Petites_Slopes}
proves Lemma \ref{Lemma_Nombre_Vallees_Visitees} since $n_9\geq n_7\geq n_3$.
\hfill$\Box$

In the following lemma, we study separately the cases in which $z\in[\bb_{\I_k-1}, \bb_{\I_k+1}]$ for $1\leq k \leq 3$
(in view of Lemma \ref{Lemma_Nombre_Vallees_Visitees} since $S_i\in\cup_{k=1}^3 [\bb_{\I_k-1}, \bb_{\I_k+1}]$ for $i\leq\tau(\bb_{\I_4})$).

\begin{lem}\label{Lemma_Valle_k_E2c}
There exists $c_{25}>0$ such that, for all $n\geq n_9$, for all $z\in\Z$,
for all $1\leq k\leq 3$,
$$
    \PP\big(S_n=z, \bb_{\I_k-1}\leq z \leq \bb_{\I_k+1}, \tau(\bb_{\I_k})\leq n ,E_5^{(n)},
    \big(E_3^{(n)}\big)^c\big)
\leq
    c_{25}(\log_2 n)^3 (\log n)^{-3}.
$$
\end{lem}

\begin{figure}[htbp]
\includegraphics[width=15.98cm,height=6.67cm]{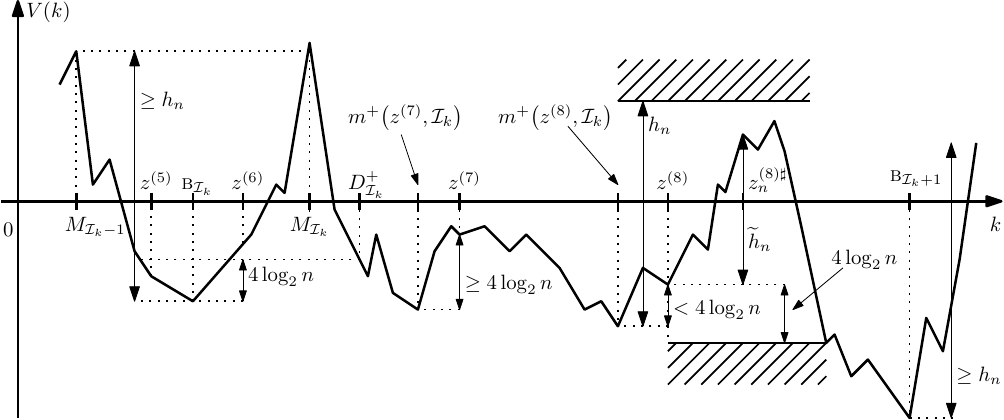}
\caption{Schema of the potential $V$,
with $z$ represented as
$z^{(5)}$ in the first case of the proof of Lemma \ref{Lemma_Proba_Vallee_k_entreDi},
$z^{(6)}$ in the second one,
and as
$z^{(7)}$ on $E_{21,k}^{(n, z)}$ and
$z^{(8)}$ on $E_{22,k}^{(n, z)}$ for the proof of Lemma \ref{Lem_Fond_Velle_k_Droite}.
}
\label{figure_Ik_cas_1_et_2}
\end{figure}

Before proving Lemma \ref{Lemma_Valle_k_E2c}, we introduce some notation.
For $i\in\Z$, let (see Figure \ref{figure_Ik_cas_1_et_2}),
\begin{eqnarray}\label{eq_def_Di+}
    D_  i^+
& := &
    \min\{j\geq M_i,\ V(j)\leq V(\bb_i)+4\log_2 n\},
\\
    D_i^-
& := &
    \max\{j\leq M_{i-1},\ V(j)\leq V(\bb_i)+4\log_2 n\},
\nonumber
\end{eqnarray}
so that, by ellipticity, $V(j)\geq V(\bb_i)+4\log_2 n+\log \e_0$ for each  $j\in([D_i^-, M_{i-1}]\cup [M_i,D_i^+])$.

We cut the proof of Lemma \ref{Lemma_Valle_k_E2c} into two main parts. First we consider the case
$z\in\big[D_{\I_k}^-,D_{\I_k}^+\big]$
in Lemma \ref{Lemma_Proba_Vallee_k_entreDi}, then
$z\in\big]D_{\I_k}^+, \bb_{\I_k+1}\big]$ in Lemma \ref{Lem_Fond_Velle_k_Droite},
the case $z\in\big[\bb_{\I_k-1}, D_{\I_k}^-\big[$ being obtained by symmetry in \eqref{Ineg_Fond_Velle_k_Gauche}.

\begin{lem}\label{Lemma_Proba_Vallee_k_entreDi}
There exists $c_{26}>0$ such that, for all $n\geq n_9$, for all $z\in\Z$,
for all $1\leq k\leq 3$,
\begin{equation}\label{Ineg_Proba_Vallee_k_entreDi}
    \PP\big[S_n=z, D_{\I_k}^-\leq z \leq D_{\I_k}^+, \tau(\bb_{\I_k})\leq n, \big( E_3^{(n)}\big)^c \big]
\leq
    c_{26}(\log_2 n)^3(\log n)^{-3}.
\end{equation}
\end{lem}

\noindent{\bf Proof:}
The proof is divided into two cases, one for which we use Lemma \ref{Lem_Proba_zn_Proche_Fond}
if $V(z)-V(\bb_{\I_k})$ is small enough ($\leq 4\log_2 n$), and one for which we use reversibility
if it is larger.
More precisely, let $n\geq n_9$ and $z\in\Z$.
First,
recall that $\{\I_j,\ 1\leq j \leq k\}\subset \{1-k,\dots, k\}$
%$1-k\leq \I_k\leq k $
for every $k\in\N^*$.
%notice that $\I_1\in\{0,1\}$ and that for $k\geq 2$, $\min\{\I_j, j<k\}-1 \leq \I_k\leq \max\{\I_j, j<k\}+1$,
%and so $|\I_k|\leq k$ for every $k\geq 1$.
So by Lemma \ref{Lem_Proba_zn_Proche_Fond} with $h_n'=h_n$,
since $n\geq n_9\geq n_8$,
we have for every $1\leq k \leq 3$,
taking into account all the possible values $j$ of $\I_k$
(see Figure \ref{figure_Ik_cas_1_et_2} with $z=z^{(5)}$),
\begin{eqnarray}
&&
    \PP\big[S_n=z, M_{\I_k-1}\leq z \leq M_{\I_k}, V(z)\leq V(\bb_{\I_k})+4\log_2 n, \big( E_3^{(n)}\big)^c\big]
\nonumber\\
& \leq &
    \p\big[\big( E_3^{(n)}\big)^c \cap \cup_{j=-2}^3 \{ M_{j-1}\leq z \leq M_j, V(z)\leq V(\bb_j)+4\log_2 n\}\big]
\nonumber\\
& \leq &
    c_{21}(\log_2 n)^3(\log n)^{-3}.
\label{Ineg_Proba_zn_bas}
\end{eqnarray}
Second, conditioning by $\omega$, then applying the strong Markov property
at stopping time $\tau(\bb_{\I_k})$ in the first equality,
we get (see Figure \ref{figure_Ik_cas_1_et_2} with $z=z^{(6)}$),
\begin{eqnarray}
&&
    \PP\big[S_n=z, V(z)\geq V(\bb_{\I_k})+4\log_2 n+\log \e_0, \tau(\bb_{\I_k})\leq n \big]
\nonumber\\
& = &
    \EE\big[{\bf 1}_{\{V(z)\geq V(\bb_{\I_k})+4\log_2 n+\log \e_0\}}
    {\bf 1}_{\{\tau(\bb_{\I_k})\leq n\}}
    \po^{\bb_{\I_k}}(S_\ell=z)_{|\ell=n-\tau(\bb_{\I_k})} \big]
\nonumber\\
& \leq &
    \big(1+e^{C_0}\big)\e_0^{-1}
    (\log n)^{-4},
\label{Ineg_Proba_zn_haut}
\end{eqnarray}
since
$
    \po^{\bb_{\I_k}}(S_\ell=z)
\leq
    \frac{\mu_\o(z)}{\mu_\o(\bb_{\I_k})}
\leq
    \big(1+e^{C_0}\big) \exp(-[V(z)-V(\bb_{\I_k})])$ for all $\ell\in\N$
by reversibility and ellipticity (see \eqref{reversiblemeas} and \eqref{eq_ellipticity_for_V}).

Finally, notice that if $D_{\I_k}^-\leq z \leq D_{\I_k}^+$,
then either $V(z)\geq V(\bb_{\I_k})+4\log_2 n+\log \e_0$,
either $M_{\I_k-1}\leq z \leq M_{\I_k}$ and $V(z)\leq V(\bb_{\I_k})+4\log_2 n$
(by the remark after \eqref{eq_def_Di+} and since $\log \e_0\leq 0$).
Hence, combining \eqref{Ineg_Proba_zn_bas} and \eqref{Ineg_Proba_zn_haut}, we get \eqref{Ineg_Proba_Vallee_k_entreDi},
since $n\geq n_9\geq n_3$.
\hfill $\Box$

We now consider the case $z\in\big]D_{\I_k}^+, \bb_{\I_k+1}\big]$ (notice that this interval may be empty).
We prove the following lemma.

\begin{lem}\label{Lem_Fond_Velle_k_Droite}
There exists $c_{27}>0$ such that, for all $n\geq n_9$, for all $z\in\Z$,
for all $1\leq k\leq 3$,
\begin{equation}\label{Ineg_Fond_Velle_k_Droite}
%    \PP\big[\big\{D_{\I_k}^+<z\leq \bb_{\I_k+1}\big\}\cap \{S_n=z, \tau(\bb_{\I_k})\leq n\}\cap E_5^{(n)}\cap \big( %E_3^{(n)}\big)^c\big]
    \PP\big[S_n=z, D_{\I_k}^+<z\leq \bb_{\I_k+1},  \tau(\bb_{\I_k})\leq n, E_5^{(n)},  \big( E_3^{(n)}\big)^c\big]
\leq
    c_{27} (\log_2 n)^3 (\log n)^{-3}.
\end{equation}
\end{lem}

Before giving the proof, we introduce some notation.
Let $n\geq n_9$ and $z\in\Z$.
We define for $i\in\Z$ (see Figure \ref{figure_Ik_cas_1_et_2}),
\begin{equation}\label{eq_def_mi_plus}
    m^+(z,i)
:=
    \min\Big\{D_i^+\leq j \leq z,\quad V(j)=\min\nolimits_{[D_i^+, z]}V\Big\},
\end{equation}
with by convention, $\min\emptyset=+\infty$, so $m^+(z,i)$ is defined in every case, even if we use it only when $z\geq D_i^+$.

\begin{figure}[htbp]
\includegraphics[width=15.99cm,height=7.42cm]{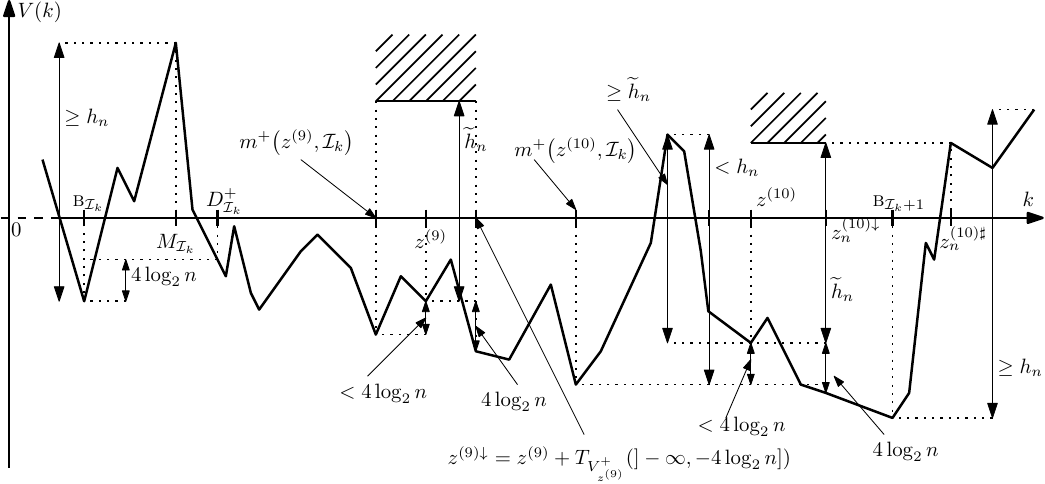}
\caption{Schema of the potential $V$, with $z$ equal to
$z^{(9)}$ on $E_{23,k}^{(n, z)}\cap E_{27,k}^{(n)}$,
and $z^{(10)}$ on $E_{23,k}^{(n, z)}\cap \big(E_{27,k}^{(n)}\big)^c$.
}
\label{figure_Ik}
\end{figure}

\noindent{\bf Idea of the proof:} (see Figures \ref{figure_Ik_cas_1_et_2} and \ref{figure_Ik} for the different cases).
First, loosely speaking,
if $V(z)$ is quite larger than the minimum of $V$ in $[D_{\I_k}^+,z]$
\big(see $E_{21,k}^{(n, z)}$ and Figure \ref{figure_Ik_cas_1_et_2} below\big)
and $n\geq \tau(\bb_{\I_k})$, then by reversibility the probability that $S_n=z$ is negligible.
So we can assume that $V(z)$ is just slightly higher than $\min\nolimits_{[D_{\I_k}^+, z]}V$.
If moreover on the right of $z$, the potential $V$ goes up $\widetilde h_n$ before going down $4\log_2 n$
\big(see $E_{22,k}^{(n, z)}$ and Figure \ref{figure_Ik_cas_1_et_2} below\big),
we prove that we are in $\cup_{q=-8}^8 E_{17}^{(n)}\big(q, \widetilde h_n, z\big)$,
so, applying Lemma \ref{Lem_Proba_zn_Proche_Fond}, the probability of this case is also negligible.
Thus we can also assume that
on the right of $z$, the potential $V$ does not go   up $\widetilde h_n$ before going down $4\log_2 n$
\big(see $E_{23,k}^{(n, z)}$ below and Figure \ref{figure_Ik}\big).
In this case, if $\tau(\bb_{\I_k})+\tau[\bb_{\I_k}, m^+(z,\I_k)] >n$, then $S_n\neq z$.
Also, we can choose some constant $c_{28}$ such that, applying \eqref{ineg_Proba_Atteinte_egalite_DGP},
$\tau(\bb_{\I_k})+\tau[\bb_{\I_k}, m^+(z,\I_k)] \in\big[n-n(\log n)^{-c_{28}},n\big]$
has a negligible probability.
Finally, if $\tau(\bb_{\I_k})+\tau[\bb_{\I_k}, m^+(z,\I_k)] < n-n(\log n)^{-c_{28}} $,
then we prove that quite quickly and in particular before time $n$ (if some very probable additional condition
is satisfied, see \eqref{eq_def_E30} and \eqref{Ineg_Proba_E27c}), $S$
goes to some place $z_n^\downarrow$ with $V(z_n^\downarrow)\leq V(z)-4\log_2 n$,
and then the probability that $S_n=z$ is negligible, once more by reversibility.
We now prove this rigorously.

\noindent{\bf Proof of Lemma \ref{Lem_Fond_Velle_k_Droite}:}
Let $n\geq n_9$, $z\in\Z$ and $1\leq k \leq 3$.
The proof is divided into three main cases, corresponding to the following events,
the last one being itself divided into four subcases
(which are defined around \eqref{eq_def_E16_a_18} and \eqref{eq_def_E30}):
\begin{eqnarray*}
    E_{21,k}^{(n, z)}
& := &
    \big\{D_{\I_k}^+<z\leq \bb_{\I_k+1}\big\}
    \cap \Big\{V(z)\geq \min\nolimits_{[D_{\I_k}^+, z]}V+4\log_2 n\Big\},
\\
    E_{22,k}^{(n, z)}
& := &
   \big\{D_{\I_k}^+<z\leq \bb_{\I_k+1}\big\}
   \cap \Big\{V(z)< \min\nolimits_{[D_{\I_k}^+, z]}V+4\log_2 n\Big\}
\\
&&
\quad\qquad\qquad\qquad\qquad
    \cap \big\{ T_{V_{z}^+}\big(\big[\widetilde h_n,+\infty\big[\big)<T_{V_{z}^+}(]-\infty, -4\log_2 n])   \big\},
\\
    E_{23,k}^{(n, z)}
& := &
    \big\{D_{\I_k}^+<z\leq \bb_{\I_k+1}\big\}
    \cap \Big\{V(z)< \min\nolimits_{[D_{\I_k}^+, z]}V+4\log_2 n\Big\}
\\
&&
\quad\qquad\qquad\qquad\qquad
    \cap \big\{ T_{V_{z}^+}(]-\infty, -4\log_2 n]) < T_{V_{z}^+}\big(\big[\widetilde{h}_n,+\infty\big[\big)  \big\}.
\end{eqnarray*}
where
$
    V_{z}^+(\ell)
=
    V(z+ \ell)-V(z)
$, $\ell\in\N$
as before and $\widetilde h_n:=h_n-C_1\log_2 n=\log n -2C_1\log_2 n$.
See figures \ref{figure_Ik_cas_1_et_2} and \ref{figure_Ik}.

\noindent {\bf First case:}
We consider the event $E_{21,k}^{(n, z)}$.

We have, once more conditioning by $\omega$ then applying
the strong Markov property at stopping time $\tau(\bb_{\I_k})$ in the first equality,
then using $\bb_{\I_k}\leq D_{\I_k}^+\leq m^+(z,\I_k)\leq z$ on $E_{21,k}^{(n, z)}$ in the second equality,
then the strong Markov property at time $\tau[m^+(z,\I_k)]$,
\begin{eqnarray}
&&
    \PP\big[S_n=z, \tau(\bb_{\I_k})\leq n, E_{21,k}^{(n, z)}\big]
\nonumber\\
& = &
    \EE\big[{\bf 1}_{\{\tau(\bb_{\I_k})\leq n\}}
    {\bf 1}_{E_{21,k}^{(n, z)}}\po^{\bb_{\I_k}}(S_\ell=z)_{|\ell=n-\tau(\bb_{\I_k})}\big]
\nonumber\\
& = &
    \EE\big[{\bf 1}_{\{\tau(\bb_{\I_k})\leq n\}} {\bf 1}_{E_{21,k}^{(n, z)}}
        \po^{\bb_{\I_k}}(S_\ell=z, \tau[m^+(z,i)]\leq \ell)_{|i=\I_k,\, \ell=n-\tau(\bb_{\I_k})}\big]
\nonumber\\
& = &
    \EE\big[{\bf 1}_{\{\tau(\bb_{\I_k})\leq n\}} {\bf 1}_{E_{21,k}^{(n, z)}}
        \E_\omega^{\bb_{\I_k}}\big({\bf 1}_{\{\tau[m^+(z,i)]\leq \ell\}}
                    \po^{m^+(z,i)}(S_t=z)_{|t=\ell-\tau[m^+(z,i)]} \big)_{|i=\I_k,\,  \ell=n-\tau(\bb_{\I_k})}
    \big]
\nonumber\\
& \leq &
    \big(1+e^{C_0}\big)(\log n)^{-4}
\label{Ineg_Proba_E14}
\end{eqnarray}
since
$
    \po^{m^+(z,\I_k)}(S_t=z)
\leq
    \frac{\mu_\o(z)}{\mu_\o[m^+(z,\I_k)]}
\leq
    \big(1+e^{C_0}\big)\exp(-[V(z) -V(m^+(z,\I_k))])
\leq
    \big(1+e^{C_0}\big)(\log n)^{-4}
$
for all $t\in\N$ on $E_{21,k}^{(n, z)}$
by reversibility and ellipticity (see \eqref{reversiblemeas} and \eqref{eq_ellipticity_for_V}).

\noindent{\bf Second case:}
We now focus on $ E_{22,k}^{(n, z)}$.
Notice in particular that $E_{22,k}^{(n, z)}$ includes the case where the potential of $z$ is "close" to the one of $\bb_{\I_k+1}$
(with a difference of potential lower than $4\log_2 n$).

We now assume that we are on $E_{22,k}^{(n, z)}$.
% \cap \big\{D_{\I_k}^+<z\leq \bb_{\I_k+1}\big\}$.
Hence we have, by definition \eqref{eq_def_Di+} of $D_{\I_k}^+$,
\begin{equation}
    \min\nolimits_{[M_{\I_k}, z]}V
=
    \min\nolimits_{[D_{\I_k}^+, z]}V
>
    V(z)-4\log_2 n.
\label{Ineg_Min_V_zn_gauche}
\end{equation}
Also, $V(M_{\I_k})=\max_{[M_{\I_k}, \bb_{\I_k+1}]}V$ and $[M_{\I_k},z] \subset [M_{\I_k}, \bb_{\I_k+1}]$, so
\begin{equation}
    \max_{[M_{\I_k}, z]}V
=
    V(M_{\I_k})
%=
%   V(\bb_{\I_k})
%    +
%    H[T_{\I_k}(V,h_n)]
\geq
    V(\bb_{\I_k})+h_n
\geq
    V(D_{\I_k}^+)-4\log_2 n +h_n,
\label{Ineg_Max_V_zn_gauche}
\end{equation}
since
$
%    H[T_{\I_k}(V,h_n)]=
    V(M_{\I_k})-V(\bb_{\I_k})
\geq
    h_n
$
and
once more by definition of $D_{\I_k}^+$.
%since $D_{\I_k}^+<\bb_{\I_k+1}$.

Now, let $z_n^\sharp:=z+T_{V_{z}^+}\big(\big[\widetilde h_n,+\infty\big[\big)$.
%{\bf (remplace $z_n^+$ conflit)}
%and
%$$
%    m_0:=m_0(\I_k, z)
%=
%    \min\big\{\ell\in[D_{\I_k}^+, z_n^\sharp],\ V(\ell)=\min\nolimits_{[D_{\I_k}^+, z_n^\sharp]} V \big\}.
%$$
%{\bf (which is a $\widetilde h_n$ minimum and can be, or not, equal to $\bb_{\I_k+1})$}
By definition of $D_{\I_k}^+$
and due to the first event defining $E_{22,k}^{(n, z)}$,
then due to the last two events defining $E_{22,k}^{(n, z)}$, we have
\begin{equation}\label{Ineg_m0}
%\label{eq_m0_min_Mz}
    \min\nolimits_{[M_{\I_k}, z_n^\sharp]}V
%=
%    V(m_0)
=
    \min\nolimits_{[D_{\I_k}^+, z_n^\sharp]}V
\geq
    V(z)-4\log_2 n.
\end{equation}

%Also, by definition of $D_{\I_k}^+$,
%\begin{equation}\label{eq_m0_min_Mz}
%    \min\nolimits_{[M_{\I_k}, z_n^\sharp]}V
%=
%    \min\nolimits_{[D_{\I_k}^+, z_n^\sharp]}V
%\geq
%    V(z)-4\log_2 n.
%\end{equation}
%So $    m_0
%=
%    \min\big\{\ell\in[M_{\I_k}, z_n^\sharp],\ V(\ell)=\min\nolimits_{[M_{\I_k}, z_n^\sharp]} V \big\}
%$, and by definition of $M_{\I_k}$ and of $m_0$,

There exists a unique index $p\in\Z$ such that
$
    M_{p-1}\big(V, \widetilde h_n\big)
\leq
    z
<
    M_p\big(V, \widetilde h_n\big)
$.
So $M_{p-1}\big(V, \widetilde h_n\big)$ is the largest
left $\widetilde h_n$-maximum less than or equal to $z$.
Since $M_{\I_k}$ is a left $h_n$ and then left $\widetilde h_n$-maximum and is $\leq z$,
we have
$
    M_{\I_k}
\leq
    M_{p-1}\big(V, \widetilde h_n\big)
$.

Assume that $z_n^\sharp<b_p\big(V, \widetilde h_n\big)$.
We define
$$
    b_n^\sharp
:=
    \inf\big\{q\in\Z,\ q\geq M_{p-1}\big(V, \widetilde h_n\big),\ V(q)= \min\nolimits_{[M_{p-1}(V, \widetilde h_n),\, z_n^\sharp]}V \big\}.
$$
We would have
$
    M_{p-1}\big(V, \widetilde h_n\big)
\leq
    z
<
    z_n^\sharp
<
    b_p\big(V, \widetilde h_n\big)
$
and so
$
    V(z_n^\sharp)
\geq
    V(z)+\widetilde h_n
\geq
    V(b_n^\sharp)+\widetilde h_n
$
by definition of $z_n^\sharp$ and $b_n^\sharp$,
and
$$
    V\big[M_{p-1}\big(V, \widetilde h_n\big)\big]
=
    \max_{[M_{p-1}(V, \widetilde h_n),\, b_p(V, \widetilde h_n)]}V
\geq
    V(z_n^\sharp)
\geq
    V(b_n^\sharp)+\widetilde h_n.
$$
Hence, $b_n^\sharp$ would be a left $\widetilde h_n$-minimum of $V$, strictly between
$M_{p-1}\big(V, \widetilde h_n\big)$ and $b_p\big(V, \widetilde h_n\big)$,
which is not possible
because $M_{p-1}\big(V, \widetilde h_n\big)$ and $b_p\big(V, \widetilde h_n\big)$
 are consecutive left $\widetilde h_n$-extrema
(see \eqref{eq_def_bi} and the comments below).
So, $b_p\big(V, \widetilde h_n\big)\leq z_n^\sharp$.

Thus,
$
    M_{\I_k}
\leq
    M_{p-1}\big(V, \widetilde h_n\big)
\leq
    b_p\big(V, \widetilde h_n\big)
\leq
    z_n^\sharp
$,
Hence,
$
    V\big[b_p\big(V, \widetilde h_n\big)\big]
\geq
     \min_{[M_{\I_k}, z_n^\sharp]}V
\geq
    V(z)-4\log_2 n
$
by \eqref{Ineg_m0} for the last inequality,
and thus
$
V(z) \leq V\big[b_p\big(V, \widetilde h_n\big)\big]+4\log_2 n
$.
Hence, using the definition of $p$, we are in $E_{17}^{(n)}\big(p, \widetilde h_n, z\big)$
(defined after \eqref{eq_def_bi}).

Also,
%$M_{\I_k}^+<z\leq \bb_{\I_k+1}$
$M_{\I_k}<z\leq \bb_{\I_k+1}$
on $E_{22,k}^{(n, z)}$.
So, either
$
    z_n^\sharp
<
    \bb_{\I_k+1}
<
    M_{\I_k+1}
$,
either $z_n^\sharp\geq \bb_{\I_k+1}$.
In the second case,
$
    V(\bb_{\I_k+1})
\geq
    \min\nolimits_{[M_{\I_k}, z_n^\sharp]}V
\geq
    V(z)-4\log_2 n
$
by \eqref{Ineg_m0},
so
$
    V(z)+\widetilde h_n
\leq
    V(\bb_{\I_k+1})+\widetilde h_n+4\log_2 n
\leq
    V(\bb_{\I_k+1})+h_n
$
since $C_1>4$,
thus
$
    z_n^\sharp
%M=
%   z_n+T_{V_{z_n}^+}\big(\big[\widetilde h_n,+\infty\big[\big)
\leq
    \bb_{\I_k+1}
    +T_{V_{\bb_{\I_k+1}}^+}([h_n,+\infty[)
\leq
    M_{\I_k+1}
$
by definition of $z_n^\sharp$ and $M_{\I_k+1}$.
Hence in every case,
$
    \bb_{\I_k}
<
    M_{\I_k}
\leq
    b_p\big(V, \widetilde h_n\big)
\leq
    z_n^\sharp
\leq
    M_{\I_k+1}
<
    \bb_{\I_k+2}
$,
and so
$
    b_p\big(V, \widetilde h_n\big)
\in
    ]\bb_{\I_k},\bb_{\I_k+2}[
$.

We now also assume that $\omega\in E_{15}^{(n)}$.
We recall that since $1\leq k\leq 3$, there exists $i_k\in\{-5,\dots, 6\}$
such that $\bb_{\I_k}=x_{i_k}(V,h_n)$
(as proved before \eqref{eq_avant_numero_bIk}).
%$\bb_{\I_k}$ is by definition (see \eqref{eq_def_bi}) equal to some $x_i(V,h_n)$ with
%%$|i|\leq 2|\I_k|+1\leq 2|k|+1\leq 9$.
%$|i|\leq 2|\I_k|\leq 2|k|\leq 6$.

Also $\bb_{\I_k}$ is a left $h_n$-minimum and since $\widetilde h_n<h_n$,  it is
a fortiori a left $\widetilde h_n$-minimum, so is equal to a $b_j\big(V,\widetilde h_n\big)$,
with $-4\leq j\leq 5$ since
$-2\leq \I_k\leq 3$ (see before \eqref{eq_avant_numero_bIk})
and on $E_{15}^{(n)}$, as already proved,
all the left $\widetilde h_n$-minima $b_\ell\big(V, \widetilde h_n\big)$ with
$|\ell|\leq 8$
%$|\ell|\leq 17$
%$|\ell|\leq 40$
are also left $h_n$-minima except at most two of them because $h_n<\log n+C_2\log_2 n$,
%each one
%making appear at most one left $\widetilde h_n$-minimum which is not a left $h_n$-minimum,
thus the number of left $\widetilde h_n$-minima in $\big]0, \bb_{\I_k}\big]$ if $\bb_{\I_k}>0$
(resp. $\big[\bb_{\I_k}, 0\big]$ if $\bb_{\I_k}\leq 0$) is at most $|\I_k|+2$ (resp. $|\I_k|+3$).
Also for this last reason, there are no more than three left $\widetilde h_n$-minima
in
$]\bb_{\I_k},\bb_{\I_k+2}[$,
%$]\bb_{\I_k}, \bb_{\I_k+1}[$,
interval to which $b_p\big(V, \widetilde h_n\big)$ belongs as proved previously on $ E_{22,k}^{(n, z)}$,
so $|p|\leq |j|+3\leq 8$.
Since we already proved that we are on $E_{17}^{(n)}\big(p, \widetilde h_n, z\big)$, this gives
$
    E_{22,k}^{(n, z)}\cap E_{15}^{(n)}
\subset
    \cup_{q=-8}^8 E_{17}^{(n)}\big(q, \widetilde h_n, z\big)
$.

Finally, by Lemmas \ref{Lemma_Proba_E11_Petites_Slopes} and \ref{Lem_Proba_zn_Proche_Fond},
we have since $n\geq n_9\geq n_8\geq n_7$,
%$$
%    \PP\big[E_{22,k}^{(n, z)}\cap E_{15}^{(n)}\cap \big( E_3^{(n)}\big)^c \big]
%\leq
%    \p\big[\big( E_3^{(n)}\big)^c \cap \cup_{q=-8}^8 E_{17}^{(n)}\big(q, \widetilde h_n, z\big)\big]
%\leq
%    c_{21}(\log_2 n)^3(\log n)^{-3}.
%$$
%This and
%%\eqref{Ineg_Proba_E11_Petites_Slopes}
%Lemma \ref{Lemma_Proba_E11_Petites_Slopes} give, since $n\geq n_9\geq n_7$,
\begin{eqnarray}
    \PP\big[E_{22,k}^{(n, z)}\cap \big( E_3^{(n)}\big)^c \big]
& \leq &
    \p\big[\big(E_{15}^{(n)}\big)^c\big]
+
    \PP\big[E_{22,k}^{(n, z)}\cap E_{15}^{(n)}\cap \big( E_3^{(n)}\big)^c \big]
\nonumber\\
& \leq &
    \p\big[\big(E_{15}^{(n)}\big)^c\big]
+
    \p\big[\big( E_3^{(n)}\big)^c \cap \cup_{q=-8}^8 E_{17}^{(n)}\big(q, \widetilde h_n, z\big)\big]
\nonumber\\
& \leq &
    (c_{20}+c_{21})(\log_2 n)^3 (\log n)^{-3}.
\label{Ineg_Proba_E13}
\end{eqnarray}

\noindent{\bf Third case:}
There remains to consider $ E_{23,k}^{(n, z)}$.
We recall $m^+(z,i)$ from \eqref{eq_def_mi_plus},
and the definition of the return time
$
    \tau^*(y)
:=
    \inf\{k\geq 1\ :\ S_k=y\}
$ for $y\in\Z$.

Using \eqref{ineg_Proba_Atteinte_egalite_DGP} (with $\bb_{\I_k}<M_{\I_k}<m^+(z,\I_k)$)  in the first line,
the Markov property in the second one,
\eqref{probaatteinte}
in the third one, we have on $ E_{23,k}^{(n, z)}$ for every $\ell\in\N$,
\begin{eqnarray}
    \po^{\bb_{\I_k}}[\tau(m^+(z,i))=\ell]_{|i=\I_k}
& \leq &
    \po^{\bb_{\I_k}}[\tau(m^+(z,i))<\tau^*(\bb_{i})]_{|i=\I_k}
\nonumber\\
& = &
    \o_{\bb_{\I_k}}\po^{\bb_{\I_k}+1}[\tau(m^+(z,i))<\tau(\bb_{i})]_{|i=\I_k}
\nonumber\\
& \leq &
    \exp[V(\bb_{\I_k})-V(M_{\I_k})]
\nonumber\\
%& = &
%    \exp[-H(T_{\I_k}(V,h_n))]
&  \leq &
   \exp(-h_n)=(\log n)^{C_1}/n
\label{Ineg_Proba_Temps_Atteinte_m_plus}
\end{eqnarray}
since $V(M_{\I_k})-V(\bb_{\I_k})=H[T_q(V, h_n)]\geq h_n$ with $q$ such that $\bb_{\I_k}=b_q(V,h_n)$.

Let $c_{28}:=C_1+4$. The next step is to prove that
%for $n\geq n_9$,
\begin{equation}\label{Ineg_Proba_E16_E4}
    \PP\big[S_n=z,\, \tau(\bb_{\I_k})\leq n,\, E_{23,k}^{(n, z)}\cap E_5^{(n)}\big]
\leq
    c_{29}(\log_2 n)^3(\log n)^{-3}
\end{equation}
for some constant $c_{29}>0$.
To this aim, we consider the three following events, defined as
\begin{eqnarray}
    E_{24,k}^{(n, z)}
& := &
    \big\{\tau(\bb_{\I_k})+\tau[\bb_{\I_k}, m^+(z,\I_k)] < n-n(\log n)^{-c_{28}} \big\},
\nonumber\\
    E_{25,k}^{(n, z)}
& := &
    \big\{\tau(\bb_{\I_k})+\tau[\bb_{\I_k}, m^+(z,\I_k)] \in\big[n-n(\log n)^{-c_{28}},n\big] \big\},
\label{eq_def_E16_a_18}
\\
    E_{26,k}^{(n, z)}
& := &
    \big\{\tau(\bb_{\I_k})+\tau[\bb_{\I_k}, m^+(z,\I_k)] >n \big\}.
\nonumber
\end{eqnarray}

First, we have, conditioning by $\omega$ then applying the strong Markov property at stopping time
$\tau(\bb_{\I_k})$, then summing \eqref{Ineg_Proba_Temps_Atteinte_m_plus} for all the integers $\ell$
in $\big[t-n(\log n)^{-c_{28}},t\big]\cap\N$,
\begin{eqnarray}
&&
    \PP\big[\tau(\bb_{\I_k})\leq n, E_{23,k}^{(n, z)}\cap E_{25,k}^{(n, z)}\big]
\nonumber\\
& = &
    \EE\Big[{\bf 1}_{\{\tau(\bb_{\I_k})\leq n\}\cap E_{23,k}^{(n, z)}}
        \po^{\bb_{\I_k}}\big(\tau(m^+(z,i))\in\big[t-n(\log n)^{-c_{28}},t\big]\big)_{|i=\I_k,\, t=n-\tau(\bb_{\I_k})}
    \Big]
\nonumber\\
& \leq &
%    2(\log n)^{C_1-c_{28}}
    [n(\log n)^{-c_{28}}+1](\log n)^{C_1}/n
\leq
    2(\log n)^{-c_{28}+C_1}
\leq
    (\log n)^{-3}
\label{Ineg_E18}
\end{eqnarray}
since $C_1-c_{28}=-4$ and $n\geq n_9\geq n_3$.

%{\bf (attention $\I_k$ est aleatoire ; confus dans \eqref{Ineg_Proba_Temps_Atteinte_m_plus} ?)}

Also, on $E_{23,k}^{(n, z)}\cap \{\tau(\bb_{\I_k})\leq n\}\cap E_{26,k}^{(n, z)}$,
$m^+(z,\I_k)\leq z$, and
after hitting $\bb_{\I_k}$,
$S$ does not hit $m^+(z,\I_k)> \bb_{\I_k}$ before time $n$,
so $S_n< m^+(z,\I_k) \leq z$ thus $S_n\neq z$. Hence,
\begin{equation}
    \PP\big[S_n=z, \tau(\bb_{\I_k})\leq n, E_{23,k}^{(n, z)}\cap E_{26,k}^{(n, z)}\big]
\leq
    \PP\big[S_n=z, S_n<m^+(z,\I_k)\leq z\big]
=
    0.
\label{Ineg_E19}
\end{equation}

There only remains to consider $E_{24,k}^{(n, z)}$.
To this aim, we introduce
\begin{equation}\label{eq_def_E30}
    E_{27,k}^{(n,z)}
:=
    \big\{\max\nolimits_{[m^+(z,\I_k),\, z]}V\leq V(z)+\widetilde h_n\big\}.
\end{equation}

We have,
conditioning by $\omega$ then applying the strong Markov property at stopping time
$\tau(\bb_{\I_k})+\tau[\bb_{\I_k}, m^+(z,\I_k)]$,
\begin{align}
&
    \PP\big[S_n=z, \tau(\bb_{\I_k})\leq n, E_{23,k}^{(n, z)}\cap E_{24,k}^{(n, z)}\cap E_{27,k}^{(n,z)}\cap E_5^{(n)}\cap E_{15}^{(n)}\big]
\nonumber\\
& =
    \EE\Big[{\bf 1}_{\{\tau(\bb_{\I_k})\leq n\}\cap E_{23,k}^{(n, z)}\cap E_{24,k}^{(n, z)}\cap E_{27,k}^{(n,z)}\cap E_5^{(n)}\cap E_{15}^{(n)}}
        \po^{m^+(z,\I_k)}\big(S_t=z\big)_{|t=n-\tau(\bb_{\I_k})-\tau[\bb_{\I_k}, m^+(z,\I_k)]}
    \Big].
\label{Ineg123456}
\end{align}

We introduce
$
    z_n^\downarrow
:=
    z
    +
    T_{V_{z}^+}(]-\infty, -4\log_2 n])
$.
Assume that $E_{23,k}^{(n, z)}$ holds and that $\bb_{\I_k+1}<z_n^\downarrow$.
So we would have
$z \leq \bb_{\I_k+1}<z_n^\downarrow< z_n^\sharp$, and then
$V(\bb_{\I_k+1}) > V(z)-4\log_2 n$,
so
$
    \max_{[\bb_{\I_k+1},z_n^\sharp]}V
=
    V(z_n^\sharp)
\leq
    V(z)+\widetilde h_n+C_0
<
    V(\bb_{\I_k+1})+4\log_2 n + \widetilde h_n+C_0
<
    V(\bb_{\I_k+1})+h_n
$
since $n\geq n_9\geq n_3$ and $C_1>20$,
thus
$
    z_n^\sharp
<
    M_{\I_k+1}
$.
So we would have
$z_n^\downarrow\in[\bb_{\I_k+1}, M_{\I_k+1}]$ with
$
    V(z_n^\downarrow)
\leq
    V(z)-4\log_2 n
<
    V(\bb_{\I_k+1})
=
    \min_{[\bb_{\I_k+1}, M_{\I_k+1}]} V
\leq
    V(z_n^\downarrow)
$
which is not possible.
So, $z_n^\downarrow\leq \bb_{\I_k+1}$  on $E_{23,k}^{(n, z)}$.

Also on $E_{23,k}^{(n, z)}$,
$
    \min_{[M_{\I_k}, z]}V
=
    V[m^+(z,\I_k)]
>
    V(z)-4\log_2 n
$
as in
\eqref{Ineg_Min_V_zn_gauche},
and
$
    \min_{[z, z_n^\downarrow]}V
\geq
    V(z)-4\log_2 n+\log(\e_0)
$
by ellipticity. So we have on $E_{23,k}^{(n, z)}$,
\begin{equation}\label{Ineg_Min_et_Max}
    \min\nolimits_{[M_{\I_k}, z_n^\downarrow]}V
\geq
    V(z)-4\log_2 n+\log(\e_0).
\end{equation}
Notice that
$
    \max_{[z, z_n^\downarrow]}V
<
    V(z)+\widetilde h_n
$
on $E_{23,k}^{(n, z)}$.
So we have on $E_{23,k}^{(n, z)} \cap E_{27,k}^{(n,z)}$,
\begin{equation}\label{Ineg_Min_et_Max_2}
    \max\nolimits_{[m^+(z,\I_k), z_n^\downarrow]}V
\leq
    V(z)+\widetilde h_n.
\end{equation}
Now on $E_{23,k}^{(n, z)}\cap E_{27,k}^{(n,z)}\cap E_5^{(n)}\cap E_{15}^{(n)}$,
by Markov inequality and \eqref{InegEsperance1},
then by \eqref{Ineg_Min_et_Max}, \eqref{Ineg_Min_et_Max_2} and
$
-(\log n)^3\leq x_{-10}(V, \log n)
\leq x_{-10}\big(V, \widetilde{h}_n\big) \leq  M_{-5}\big(V, \widetilde{h}_n\big) \leq M_{-3}\leq M_{\I_k}< D_{\I_k}^+
\leq m^+(z,\I_k)
\leq z< z_n^\downarrow\leq \bb_{\I_k+1}\leq \bb_4\leq b_6(V, \widetilde{h}_n\big)
\leq x_{12}\big(V, \widetilde{h}_n\big)\leq x_{12}(V, \log n)\leq (\log n)^3$
(because on $E_{15}^{(n)}$ there are similarly as previously, in $[M_{-3}, 0]$,  at most two $M_j\big(V,\widetilde h_n\big)$ which are not equal to some $M_\ell(V, h_n)$ so $M_{-5}\big(V, \widetilde{h}_n\big) \leq M_{-3}$ and similarly $\bb_4\leq b_6(V, \widetilde{h}_n\big)$),
we get
\begin{eqnarray}
&&
    \po^{m^+(z,\I_k)}\big[\tau(M_{i})\wedge \tau(z_n^\downarrow)\geq 2^{-1}n(\log n)^{-c_{28}}
    \big]_{|i=\I_k}
\label{Ineg_Esperance_Atteinte_M_zdown}
\\
& \leq &
    2n^{-1}(\log n)^{c_{28}}\e_0^{-1}(z_n^\downarrow-M_{\I_k})^2
    \exp\Big(\max\nolimits_{[m^+(z,\I_k), z_n^\downarrow]}V-\min\nolimits_{[M_{\I_k}, z_n^\downarrow]}V\Big)
\nonumber\\
& \leq &
    8(\log n)^{c_{28}+6}n^{-1}\e_0^{-2}
    \exp\big(\widetilde h_n+4\log_2 n\big)
=
    8(\log n)^{c_{28}-2C_1+10}\e_0^{-2}
\leq
    (\log n)^{-3}
\nonumber
\end{eqnarray}
since
$
    c_{28}-2C_1+10
%=
%    4-C_1+10
=
    14-C_1
<
    -6
$
and $n\geq n_9\geq n_3$.
Moreover on $E_{23,k}^{(n, z)}$,
we get by definition of $E_{23,k}^{(n, z)}$,
$$
    V(z)
<
    \min\nolimits_{[D_{\I_k}^+, z]}V+4\log_2 n
\leq
    V\big[D_{\I_k}^+\big]+4\log_2 n,
$$
and as a consequence, using \eqref{Ineg_Max_V_zn_gauche} which remains true on $E_{23,k}^{(n, z)}$,
\begin{equation}\label{Ineg_Mino_VDIk_Vzn}
    V(M_{\I_k})
\geq
    V(D_{\I_k}^+)-4\log_2 n +h_n
\geq
    V(z)-8\log_2 n +h_n.
\end{equation}
Hence on $E_{23,k}^{(n, z)}\cap E_{27,k}^{(n,z)}\cap E_5^{(n)}\cap E_{15}^{(n)}$, using \eqref{probaatteinte},
then \eqref{Ineg_Min_et_Max_2} and \eqref{Ineg_Mino_VDIk_Vzn},
\begin{eqnarray}
&&
    \po^{m^+(z,\I_k)}\big[\tau(M_i)< \tau(z_n^\downarrow)\big]_{|i=\I_k}
\label{Ineg_Proba_Atteinte_M_zdown}
\\
& \leq &
    \big(z_n^\downarrow-m^+(z,\I_k)\big)
    \exp\Big[\max\nolimits_{[m^+(z,\I_k),\, z_n^\downarrow]}V-V(M_{\I_k})\Big]
%\\
%& \leq &
%    2(\log n)^3
%    \exp\Big[V(z)+\widetilde h_n - [V(z)-8\log_2 n +h_n]\Big]
\nonumber\\
& \leq &
    2(\log n)^3
    \exp\big[\widetilde h_n+8\log_2 n -h_n\big]
%\\
%& = &
%    2(\log n)^3
%    \exp\big[\log n -2C_1\log_2 n +8\log_2 n -(\log n - C_1\log_2 n )\big]
%\\
=
    2(\log n)^{11-C_1}
\leq
    (\log n)^{-3}
\nonumber
\end{eqnarray}
since $11-C_1<-9$ and $n\geq n_9\geq n_3$.
Consequently on $E_{23,k}^{(n, z)}\cap E_{27,k}^{(n,z)} \cap E_5^{(n)}\cap E_{15}^{(n)}$,
\begin{equation}\label{Ineg_Proba_Atteinte_zndown}
    \po^{m^+(z,\I_k)}\big[\tau(z_n^\downarrow)\geq 2^{-1}n(\log n)^{-c_{28}}\big]
\leq
    \eqref{Ineg_Esperance_Atteinte_M_zdown} + \eqref{Ineg_Proba_Atteinte_M_zdown}
\leq
    2(\log n)^{-3}.
\end{equation}
Notice that for every $\ell\in\N$, by reversibility and ellipticity (see \eqref{reversiblemeas} and \eqref{eq_ellipticity_for_V}),
\begin{equation}\label{Ineg_Proba_reversibility_zndown}
    \po^{z_n^\downarrow}(S_\ell=z)
\leq
    \mu_\o(z)
    /\mu_\o(z_n^\downarrow)
\leq
    \big(1+e^{C_0}\big)\exp[-V(z)+V(z_n^\downarrow)]
\leq
    \e_0^{-1}(\log n)^{-4}.
\end{equation}
On $E_{23,k}^{(n, z)}\cap E_{27,k}^{(n,z)} \cap E_5^{(n)}\cap E_{15}^{(n)}$, for every $t\geq n(\log n)^{-c_{28}}$,
by \eqref{Ineg_Proba_Atteinte_zndown}, the strong Markov property and \eqref{Ineg_Proba_reversibility_zndown},
since $n\geq n_9\geq n_3$,
\begin{eqnarray*}
&&
    \po^{m^+(z,\I_k)}\big(S_t=z\big)
\\
& \leq &
    \po^{m^+(z,\I_k)}\big[\tau(z_n^\downarrow)\geq 2^{-1}n(\log n)^{-c_{28}}\big]
+
    \po^{m^+(z,\I_k)}\big[S_t=z, \tau(z_n^\downarrow)< 2^{-1}n(\log n)^{-c_{28}}\big]
\\
& \leq &
    2(\log n)^{-3}
+
    \eo^{m^+(z,\I_k)}\big({\bf 1}_{\{\tau(z_n^\downarrow)< 2^{-1}n(\log n)^{-c_{28}}\}}
        \po^{z_n^\downarrow}(S_\ell=z)_{|\ell=t-\tau(z_n^\downarrow)}
    \big)
\\
& \leq &
    2(\log n)^{-3}+\e_0^{-1}(\log n)^{-4}
\leq
    3(\log n)^{-3}.
%    \big(2+\e_0^{-1}\big)(\log n)^{-3}.
\end{eqnarray*}
Finally, this and \eqref{Ineg123456}
(on which
$
    t
\geq
    n(\log n)^{-c_{28}}
$
thanks to $E_{24,k}^{(n, z)}$) give
\begin{equation}\label{Ineg_Proba_sur_E17_E30}
    \PP\big[S_n=z, \tau(b_{\I_k})\leq n, E_{23,k}^{(n, z)}\cap E_{24,k}^{(n, z)}\cap E_{27,k}^{(n,z)}\cap E_5^{(n)}\cap E_{15}^{(n)}\big]
\leq
    3(\log n)^{-3}.
\end{equation}

There only remains to estimate $\p\big[\big(E_{27,k}^{(n,z)}\big)^c\cap E_{23,k}^{(n, z)}\big]$.
We define (see Figure \ref{figure_Ik} with $z=z^{(10)}$),
\begin{eqnarray*}
    V_{2,n}^-
& := &
    V_{z}^-[.+T_{V_{z}^-}([\widetilde h_n, +\infty[)]-V_{z}^-[T_{V_{z}^-}([\widetilde h_n, +\infty[)],
\\
    V_{3,n}^-
& := &
    V_{2,n}^-[.+T_{V_{2,n}^-}(]-\infty, -\widetilde h_n])]
-
    V_{2,n}^-[T_{V_{2,n}^-}(]-\infty, -\widetilde h_n])],
\\
    E_{28}^{(n)}
& := &
    \big\{T_{V_{z}^-}([\widetilde h_n,+\infty[)<T_{V_{z}^-}(]-\infty, -4\log_2 n[)\big\},
\\
    E_{29}^{(n)}
& := &
    \big\{T_{V_{2,n}^-}(]-\infty, -\widetilde h_n])
            <
            T_{V_{2,n}^-}([C_1\log_2 n, +\infty[)
    \big\},
\\
    E_{30}^{(n)}
& := &
    \big\{T_{V_{3,n}^-}([\widetilde h_n, +\infty[)
            <
          T_{V_{3,n}^-}(]-\infty, -4\log_2 n-C_0[)
    \big\}.
\end{eqnarray*}
Notice that $E_{23,k}^{(n, z)}\cap \big(E_{27,k}^{(n,z)}\big)^c$
is included in $E_{28}^{(n)}$
because $\max_{[0, z-m^+(z,\I_k) ]} V_{z}^-> \widetilde h_n$
by \eqref{eq_def_E30}
and
$
    \min_{[0, z-m^+(z,\I_k) ]} V_{z}^-
=
    V_{z}^-(z-m^+(z,\I_k))
>
    -4\log_2 n
$
by definitions  of $m^+(z,\I_k)$ (see \eqref{eq_def_mi_plus}) and of $E_{23,k}^{(n)}$.
It is also included in  $E_{29}^{(n)}$,
otherwise there would be a left $h_n$-maximum of $V$ in $]m^+(z,\I_k), z[$
and so in $]M_{\I_k}, \bb_{\I_k+1}[$ which is not possible.
Finally, $E_{23,k}^{(n, z)}\cap \big(E_{27,k}^{(n,z)}\big)^c$
is also included in $E_{30}^{(n)}$
because
$
    \min\nolimits_{[M_{\I_k}, z]}V
>
    V(z)-4\log_2 n
$
as in \eqref{Ineg_Min_V_zn_gauche}
and
$
    V(M_{\I_k})
\geq
    V(z)-8\log_2 n +h_n
\geq
    V(z)+\widetilde h_n
+
    C_0
$
by \eqref{Ineg_Mino_VDIk_Vzn} and since $C_1>20$ and $n\geq n_9\geq n_3$.
Using the independence of $E_{28}^{(n)}$, $E_{29}^{(n)}$ and  $E_{30}^{(n)}$,
provided by the strong Markov property,
then applying \eqref{eqOptimalStopping2},
we get
\begin{equation}\label{Ineg_Proba_E27c}
    \p\big[\big(E_{27,k}^{(n,z)}\big)^c\cap E_{23,k}^{(n, z)}\big]
\leq
    \p\big[E_{28}^{(n)}\big]
    \p\big[E_{29}^{(n)}\big]
    \p\big[E_{30}^{(n)}\big]
%\\
%& \leq &
%    \frac{(4\log_2 n+C_0)}{(\widetilde h_n+4\log_2 n+C_0)}
%    \frac{C_1\log_2 n+C_0}{(\widetilde h_n+C_1\log_2 n+C_0)}
%    \frac{4\log_2 n+2C_0}{(\widetilde h_n+4\log_2 n+2C_0)}
%\\
%& \leq &
%    \frac{(5\log_2 n)}{(\log n)/2}
%    \frac{(C_1+1)\log_2 n}{(\log n)/2}
%    \frac{6\log_2 n}{(\log n)/2}
%\\
%& \leq &
%    240(C_1+1)\frac{(\log_2 n)}{(\log n)}
%    \frac{\log_2 n}{(\log n)}
%    \frac{\log_2 n}{(\log n)}
%\\
\leq
    c_{30}(\log_2 n)^3(\log n)^{-3},
\end{equation}
with $c_{30}:=10\times 2(C_1+1)\times 12$
since $n\geq n_9\geq n_3$.
%Notice that on $E_{23,k}^{(n, z)}\cap \big(E_{27,k}^{(n,z)}\big)^c$ (see Figure \ref{figure_Ik} with $z=z^{(10)}$),
%we have
%$T_{V_{z}^-}([\widetilde h_n -4\log_2 n[)<T_{V_{z}^-}(]-\infty, -4\log_2 n[)$,
%%$V_{z}^-$ goes up $\widetilde h_n -4\log_2 n$
%%before going down $-4\log_2 n$,
%that $V_{2,n}^-:=V_{z}^-[.+T_{V_{z}^-}([\widetilde h_n -4\log_2 n[)]-V_{z}^-[T_{V_{z}^-}([\widetilde h_n -4\log_2 n[)] $
%hits  $]-\infty, -(\widetilde h_n -4\log_2 n)]$ before $[(C_1+4)\log_2 n, +\infty[$
%%then goes down $\widetilde h_n-4\log_2 n$ before going up $h_n$
%(otherwise there would be a left $h_n$-maximum of $V$ in $]M_{\I_k}, b_{\I_k+1}[$ which is not possible),
%and
%$
%    V_{2,n}^-[.+T_{V_{2,n}^-}(]-\infty, -(\widetilde h_n -4\log_2 n)])]
%-
%    V_{2,n}^-[T_{V_{2,n}^-}(]-\infty, -(\widetilde h_n -4\log_2 n)])]
%$
%hits $[\widetilde h_n, +\infty[$ before $]-\infty, -4\log_2 n-C_0[$
%by definitions \eqref{eq_def_Di+} of $D_i^+$ and $E_{23,k}^{(n, z)}$.
%%then goes up  $\widetilde h_n$ before going down $C_1\log_2 n$.
%These three events are independent by strong Markov property, and each have probability
%$O((\log_2n)(\log n)^{-1})$ by \eqref{eqOptimalStopping2},
%so
%$$
%    \p\big[\big(E_{27,k}^{(n)}\big)^c\cap E_{23,k}^{(n, z)}\big]=O\big((\log_2 n)^3(\log n)^{-3}\big).
%   $$
This, combined with \eqref{Ineg_Proba_sur_E17_E30} and Lemma \ref{Lemma_Proba_E11_Petites_Slopes}, gives,
where LHS means left hand side and since $n\geq n_9\geq n_7$,
\begin{eqnarray}
&&
    \PP\big[S_n=z, \tau(b_{\I_k})\leq n, E_{23,k}^{(n, z)}\cap E_{24,k}^{(n, z)}\cap E_5^{(n)}\big]
\nonumber\\
& \leq &
    %\PP\big[S_n=z, \tau(b_{\I_k})\leq n, E_{23,k}^{(n, z)}\cap E_{24,k}^{(n, z)}\cap E_{27,k}^{(n,z)}\cap E_5^{(n)}\cap %E_{15}^{(n)}\big]
    \text{LHS of } \eqref{Ineg_Proba_sur_E17_E30}
    +
    \p\big[\big(E_{15}^{(n)}\big)^c\big]
    +
    \p\big[\big(E_{27,k}^{(n,z)}\big)^c\cap E_{23,k}^{(n, z)}\big]
\nonumber\\
& \leq &
    (3+c_{20}+c_{30})(\log_2 n)^3(\log n)^{-3}.
%    O\big((\log_2 n)^3(\log n)^{-3}\big).
\label{Ineg_Proba_sur_E17_seul}
\end{eqnarray}
Combining this, \eqref{Ineg_E18} and \eqref{Ineg_E19} proves \eqref{Ineg_Proba_E16_E4}
with $c_{29}:=c_{20}+c_{30}+4$ since $n\geq n_3$.
Finally, \eqref{Ineg_Proba_E14}, \eqref{Ineg_Proba_E13} and \eqref{Ineg_Proba_E16_E4} prove \eqref{Ineg_Fond_Velle_k_Droite}
with $c_{27}:=\e_0^{-1}+c_{20}+c_{21}+c_{29}$
for every $n\geq n_9$, $z\in\Z$ and $1\leq k \leq 3$,
which ends the proof of Lemma \ref{Lem_Fond_Velle_k_Droite}.
\hfill$\Box$

\noindent {\bf Proof of Lemma \ref{Lemma_Valle_k_E2c}:}
We prove similarly as in the proof of Lemma \ref{Lem_Fond_Velle_k_Droite}
that for all $n\geq n_9$, $z\in\Z$ and $1\leq k \leq 3$,
%that for every $1\leq k \leq 3$, as $n\to+\infty$,
\begin{equation}\label{Ineg_Fond_Velle_k_Gauche}
    \PP\big[S_n=z, \bb_{\I_k-1} \leq z< D_{\I_k}^-,  \tau(\bb_{\I_k})\leq n,
    E_5^{(n)}, (E_3^{(n)})^c\big]
\leq
    c_{27}(\log_2 n)^3 (\log n)^{-3}.
\end{equation}
Combining \eqref{Ineg_Fond_Velle_k_Gauche}, \eqref{Ineg_Proba_Vallee_k_entreDi} and \eqref{Ineg_Fond_Velle_k_Droite}
proves Lemma \ref{Lemma_Valle_k_E2c} with $c_{25}:=c_{26}+2c_{27}$.
\hfill$\Box$

\noindent {\bf Proof of Proposition \ref{Prop_Local_Limit_Si_Petite_Slope}:}
%Let $n\geq n_9$ and $z\in\Z$.
Notice that for $n\geq n_9$ and $k\geq 1$, on $\{\tau(\bb_{I_k})\leq n <\tau(\bb_{\I_{k+1}})\}$,
the random walk $S$ does not reach the $\bb_i$ with $i\in\Z\setminus\{\I_1, \dots, \I_k\}$
before time $n$, and so $S_n$ belongs to
$\big]\min\{\bb_{{\I_i}-1}, 1\leq i \leq k\}, \max\{\bb_{{\I_i}+1}, 1\leq i \leq k\}\big[$,
which is equal to $\cup_{i=1}^k ]\bb_{\I_i-1}, \bb_{\I_i+1}[$.
Consequently, using \eqref{Ineg_Proba_Atteinte_bI1} and Lemma \ref{Lemma_Nombre_Vallees_Visitees} in the second inequality,
for all $n\geq n_9$ and all $z\in\Z$, with $c_{31}:=c_{24}+2$,
\begin{eqnarray*}
&&
    \PP\big(S_n=z, (E_3^{(n)})^c\big)
\\
& \leq &
    \PP\big[\tau(\bb_{\I_1})>n\big]
+
    \PP\big[S_n=z, (E_3^{(n)})^c\cap\cup_{k=1}^3\{\tau(\bb_{I_k})\leq n <\tau(\bb_{\I_{k+1}})\}\big]
+
    \PP\big[\tau(\bb_{\I_4})\leq n\big]
\\
& \leq &
    \sum_{k=1}^3
    \PP\big(\tau(\bb_{I_k})\leq n <\tau(\bb_{\I_{k+1}}), S_n=z\in \cup_{i=1}^k ]\bb_{\I_i-1}, \bb_{\I_i+1}[,(E_3^{(n)})^c\big)
+
    c_{31}\frac{(\log_2 n)^3}{(\log n)^{3}}
\\
& \leq &
%    \sum_{k=1}^3\sum_{i=1}^k
%    %{\bf 1}_{\{i\leq k\}}
%    \PP\big(\tau(\bb_{I_k})\leq n <\tau(\bb_{\I_{k+1}}), S_n=z\in ]\bb_{\I_i-1}, \bb_{\I_i+1}[,(E_3^{(n)})^c\big)
%+
%    O\bigg( \frac{(\log_2 n)^3}{(\log n)^{3}}\bigg)
%    C \frac{(\log_2 n)^3}{(\log n)^{3}}
%\\
%& = &
%    \sum_{i=1}^3\sum_{k=1}^3{\bf 1}_{\{i\leq k\}}
%    \PP\big(\tau(b_{I_k})\leq n <\tau(b_{\I_{k+1}}), S_n=z\in ]b_{\I_i-1}, b_{\I_i+1}[,(E_3^{(n)})^c\big)
%+
%    C \frac{(\log_2 n)^3}{(\log n)^{3}}
%\\
%& = &
    \sum_{i=1}^3\sum_{k=i}^3
    \PP\big(\tau(\bb_{I_k})\leq n <\tau(\bb_{\I_{k+1}}), S_n=z\in ]\bb_{\I_i-1}, \bb_{\I_i+1}[,(E_3^{(n)})^c\big)
+
    c_{31}\frac{(\log_2 n)^3}{(\log n)^{3}}
%    O\bigg( \frac{(\log_2 n)^3}{(\log n)^{3}}\bigg)
\\
& = &
    \sum_{i=1}^3
%    \Big(
    \PP\big[\tau(\bb_{I_i})\leq n < \tau(\bb_{I_4}), S_n=z\in ]\bb_{\I_i-1}, \bb_{\I_i+1}[,(E_3^{(n)})^c\big]
%    +\p\big[\big(E_5^{(n)}\big)^c\big]
%    \Big)
+
    c_{31}\frac{(\log_2 n)^3}{(\log n)^{3}}
%    O\bigg( \frac{(\log_2 n)^3}{(\log n)^{3}}\bigg)
\\
& \leq  &
    c_{19}(\log_2 n)^3(\log n)^{-3},
\end{eqnarray*}
with $c_{19}:=3c_{25}+3+ c_{31}$ and where we used Lemmas \ref{Lemma_Valle_k_E2c}
and \ref{Lemma_Proba_E4} in the last line since $n\geq n_9\geq \max(n_3, p_3)$.
This proves Proposition \ref{Prop_Local_Limit_Si_Petite_Slope}.
\hfill$\Box$

\subsection{Proof of the upper bound in Theorem \ref{Th_Local_Limit_Sinai}}
Recall $E_C^{(n)}(z)$ from \eqref{eq_def_EL}.
We have, for all $n\geq \max(n_9, p_2)$ and all $z\in\Z$,
\begin{eqnarray*}
&&
    \PP\big(S_n=z, \big(E_C^{(n)}(z)\big)^c\big)
\\
& \leq &
%    \PP(S_n=z)
%& \leq &
%    \PP(S_n=z, E_C^{(n)}(z))
%+
    \PP\big(S_n=z, \big(E_3^{(n)}\big)^c\big)
+
    \p\big[\big(E_5^{(n)}\big)^c\big]
+
    \p\big[\big(E_6^{(n)}\big)^c\cap E_5^{(n)}\big]
\\
&&
+
    \PP\big[S_n=z, \big(E_7^{(n)}(z)\big)^c, E_3^{(n)}, E_5^{(n)}\big]
%    \PP\big(S_n=z, |z-b_{\log n}|>(\log n)^{4/3+\delta_1}, E_3^{(n)}, E_5^{(n)}\big)
%\\
%&&
+
    \p\big[\big(E_4^{(n)}(z)\big)^c \cap E_3^{(n)} \cap E_6^{(n)}\cap E_7^{(n)}(z) \big]
%    \p\big[\big(E_4^{(n)}\big)^c \cap E_3^{(n)} \cap E_6^{(n)}\cap\{|b_{\log n}-z|\leq (\log n)^{4/3+\delta_1}\}\big]
\\
& \leq &
    (c_{19}+2+c_9)(\log_2 n)^3(\log n)^{-3}
    +c_{10}(\log n)^{-2-\delta_1/2}
\end{eqnarray*}
by %\eqref{Ineg_Minor_Sur_EL},
Proposition \ref{Prop_Local_Limit_Si_Petite_Slope},
Lemmas \ref{Lemma_Proba_E4} and \ref{Lemma_Proba_46},
Proposition \ref{Lemma_Proba_z_b_far} and Lemma \ref{Lemma_Proba_E3c}.
This and Proposition \ref{Ineg_Upper_Bound_ELT}
give,
since $\delta_1\in]0,2/3[$,
\begin{equation}\label{Ineg_Minor_Sur_EL_Bis}
    \sup_{z\in(2\Z+n)}
    \bigg[
    \PP\big(S_n=z\big)
-
    \frac{2\sigma^2}{(\log n)^2}\varphi_\infty\bigg(\frac{\sigma^2 z}{(\log n)^2}\bigg)
    \bigg]
\leq
    o\big((\log n)^{-2}\big),
\end{equation}
as $n\to+\infty$,
which proves the upper bound in Theorem \ref{Th_Local_Limit_Sinai}.
\hfill$\Box$

%%%%%%%%%%%%%%%%%%%%%%%%%%%%%%%%%%%%%%%%%%%%%%%%%%%%%%%%%%%%%%%%%%%%%%%%%%%%%%%%%%%%%%%%
%                                                                                      %
%                                 SECTION LOWER BOUND                                  %
%                                                                                      %
%%%%%%%%%%%%%%%%%%%%%%%%%%%%%%%%%%%%%%%%%%%%%%%%%%%%%%%%%%%%%%%%%%%%%%%%%%%%%%%%%%%%%%%%

\section{Proof of the lower bound of Theorem \ref{Th_Local_Limit_Sinai}}\label{Sect_lower_Bonud}

%Let $A>0$.
%We assume that $n$ and $z$ have the same parity, and that
%$z=O\big((\log n)^2\big)$.
%$|z|\leq A(\log n)^2$.
Let $\e>0$.
Since $\lim_{\pm\infty}\varphi_\infty=0$,
we can fix some $A_0>0$ such that
$\sup_{|x|\geq A_0}|\varphi_\infty(\sigma^2 x)|<\sigma^{-2}\e$.
%We define, for $z\in\Z$ and $n\geq n_9$,
%$$
%    z^-
%:=
%\left\{
%\begin{array}{ll}
%z-\Gamma_n & \text{ if } z\leq -\Gamma_n,\\
%0 & \text{ if } -\Gamma_n< z \leq \Gamma_n,\\
%z+\Gamma_n & \text{ if } z> \Gamma_n.
%\end{array}
%\right.
%$$

In this section, $\mathcal T_V^\uparrow$ and $\mathcal T_V^\downarrow$ always denote $\widetilde h_n$-slopes,
that is, $\mathcal T_V^\uparrow=\mathcal T_{V,\widetilde h_n}^\uparrow$
and $\mathcal T_V^\downarrow=\mathcal T_{V,\widetilde h_n}^\downarrow$,
where $\widetilde h_n=\log n-2C_1\log_2 n=h_n-C_1\log_2 n$ as before.
In what follows, we consider independent slopes
$Z_{2k}^\uparrow$, $-9\leq k \leq 9$
and $Z_{2k+1}^\downarrow$, $-9\leq k \leq 9$,
%being independent,
%(non central) left $\widetilde h_n$-slopes of $V$,
%these slopes being upward ones for $Z_{2k}^\uparrow$
%(having the same law as $\mathcal T_V^\uparrow$ when $h=\widetilde h_n$),
each $Z_{2k}^\uparrow$ having the same law as $\mathcal T_V^\uparrow$
%(with $h=\widetilde h_n$),
\big(i.e. $\mathcal T_{V,\widetilde h_n}^\uparrow$\big),
and each $Z_{2k+1}^\downarrow$ having the same law as $\mathcal T_V^\downarrow$
\big(i.e. $\mathcal T_{V,\widetilde h_n}^\downarrow$\big).
%(with $h=\widetilde h_n$).
%and being downward ones for  $Z_{2k+1}^\downarrow$
%(having the same law as $\mathcal T_V^\downarrow$ when $h=\widetilde h_n$),

%In what follows, $\mathcal T_V^\uparrow$ and $\mathcal T_V^\downarrow$
%are independent
%and have the law of a (non central) left $\widetilde h_n$-upward (resp. downward) slope of $V$,
%%and respectively of $V_-$,
%where $\widetilde h_n=\log n-2C_1\log_2 n=h_n-C_1\log_2 n$ as before.

Recall that $\zeta$ is defined in \eqref{eq_def_zeta}.
We also introduce
$Y_{-1}^\uparrow:=\zeta\big(Z_{-1}^\downarrow\big)$,
%$\mathcal T_{V_-}^\uparrow=\zeta\big(Z_{-1}^\downarrow\big)$,
%$\mathcal T_{V_-}^\uparrow=\zeta\big(\mathcal T_V^\downarrow\big)$,
which is independent of $Z_0^\uparrow$,
with
$
    Y_{-1}^\uparrow
=_{law}
    \zeta\big(\mathcal T_V^\downarrow\big)
=_{law}
    \mathcal T_{V^-, \widetilde h_n}^{\uparrow*}
=:
    \mathcal T_{V^-}^{\uparrow*}
$
by Proposition \ref{Prop_Egalite_Loi_Zeta_Slopes},
%and has the law of a (non central) upward right $\widetilde h_n$-slope for $V_-$,
%%independent of $\mathcal T_V^\uparrow$,
and
$\ell\big(Y_{-1}^\uparrow\big)=\ell\big(Z_{-1}^\downarrow\big)$.
%$\ell\big(\mathcal T_{V_-}^\uparrow\big)=\ell\big(Z_{-1}^\downarrow\big)$.
%with $\ell\big(\mathcal T_{V_-}^\uparrow\big)=\ell\big(\mathcal T_V^\downarrow\big)$.

\noindent{\bf First case:}
We start with the case $z\leq -\Gamma_n$.
%so that $z+k\leq 0$ for every $k$ in the sum in \eqref{Ineg_Proba_Zn_E\Gamma_1}.

%In what follows, $\mathcal T_V^\uparrow$ and $\mathcal T_{V_-}^\uparrow=\zeta\big(\mathcal T_V^\downarrow\big)$
%are independent
%($\zeta$ being defined after \eqref{eq_def_J3}),
%and have the law of a (non central) left $\widetilde h_n$-upward (resp. downward) slope of $V$,
%%and respectively of $V_-$,
%where $\widetilde h_n=\log n-2C_1\log_2 n=h_n-C_1\log_2 n$ as before.

Using Lemma \ref{Lemma_Proba_bh_egal} eq. \eqref{eq_proba_bh_x_negatif},
we have for each $z\in\Z$ such that  $z-\Gamma_n\leq 0$,
\begin{equation}\label{eq_def_J6}
    J_6(n,z)
:=
    \p\big(b_{\widetilde h_n}= z-\Gamma_n\big)
=
        \frac{\p\big[-z+\Gamma_n <\ell\big(Z_0^\uparrow\big)\big]}
        {\E\big(\ell\big(Z_0^\uparrow\big)+\ell\big(Z_1^\downarrow\big)\big)}.
\end{equation}
Using the uniform continuity of $\varphi_\infty$ on $\R$
and
$
    \sup_{z\in[-A_0(\log n)^2, A_0(\log n)^2]}
    \bigg|\frac{\sigma^2 z}{(\log n)^2}
            - \frac{\sigma^2 \big(z-\Gamma_n\big)}{ \big(\widetilde h_n\big)^2}
    \bigg|
$
$
=
    o(1)
$
as $n\to+\infty$
since $\delta_1<2/3$ and $\widetilde h_n\sim_{n\to+\infty} \log n$
in the first inequality,
then $\|\varphi_\infty\|_\infty=:\sup_{\R}|\varphi_\infty|<\infty$
and $\widetilde h_n\sim_{n\to+\infty} \log n$ in the second one,
and finally Theorem \ref{Th_Local_Limit_b_h} in the last one,
there exists $n_{10}\geq \max(n_9, p_2)$
(with $p_2$ defined in Lemma \ref{Lemma_Proba_46})
such that for all $n\geq n_{10}$,
for all $z\in[-A_0(\log n)^2, A_0(\log n)^2]$ such that $z-\Gamma_n\leq 0$,
\begin{eqnarray}
    \frac{\sigma^2}{(\log n)^2}\varphi_\infty\bigg(\frac{\sigma^2 z}{(\log n)^2}\bigg)
%& \leq &
%    \frac{\sigma^2}{\big(\widetilde h_n\big)^2}\varphi_\infty\bigg(\frac{\sigma^2 z}{(\log n)^2}\bigg)
%    +o\big((\log n)^{-2}\big)
%\nonumber\\
& \leq &
    \frac{\sigma^2}{(\log n)^2}\varphi_\infty\bigg(\frac{\sigma^2 \big(z-\Gamma_n\big)}{
    \big(\widetilde h_n\big)^2}\bigg)
    +\e(\log n)^{-2}
\nonumber\\
& \leq &
    \frac{\sigma^2}{\big(\widetilde h_n\big)^2}\varphi_\infty\bigg(\frac{\sigma^2 \big(z-\Gamma_n\big)}{
    \big(\widetilde h_n\big)^2}\bigg)
    +2\e(\log n)^{-2}
\nonumber\\
& \leq &
%    \p\big(b_{\log n}= z-\Gamma_n\big)
%    +o\big((\log n)^{-2}\big)
%\\
%& = &
    J_6(n,z)
    +3\e(\log n)^{-2}.
\label{Ineg_phi_infty_J6}
\end{eqnarray}
Also for $n\geq n_{10}$, if $|z|>A_0(\log n)^2$, then by definition of $A_0$,
\begin{equation} \label{Ineg_phi_infty_J6_2}
    \frac{\sigma^2}{(\log n)^2}\varphi_\infty\bigg(\frac{\sigma^2 z}{(\log n)^2}\bigg)
\leq
    \frac{\e}{(\log n)^2}
\leq
    J_6(n,z)
    +3\e(\log n)^{-2}.
\end{equation}
The objective is to approximate progressively this quantity $J_6(n,z)$ by $\PP(S_n=z)$,
by using Theorems \ref{Lemma_Independence_h_extrema} {\bf (i)}
and \ref{Lemma_Central_Slope} eq. \eqref{eq_Central_Slope_General_phi} (see \eqref{eq_Application_Renewal_et_independence} below)
and Lemma \ref{Lemma_Approx_PQuenched_Nu}.
To this aim, we introduce the following events.
\begin{eqnarray*}
    E_{31}^{(n)}
& := &
    \big\{T_{Z_0^\uparrow}\big(\widetilde h_n\big)>\Gamma_n\big\}
\cap
    \big\{T_{Y_{-1}^\uparrow}\big(\widetilde h_n\big)>\Gamma_n\big\},
\\
    E_{32}^{(n)}
& := &
    \big\{\forall k\in\big[T_{Z_0^\uparrow}\big(\widetilde h_n\big),
    \ell\big( Z_0^\uparrow\big)\big], \ Z_0^\uparrow(k)\geq 9\log_2 n \big\},
\\
    E_{33}^{(n)}
& := &
    \big\{\forall k\in\big[T_{Y_{-1}^\uparrow}\big(\widetilde h_n\big),
    \ell\big(Y_{-1}^\uparrow\big)\big], \ Y_{-1}^\uparrow(k)\geq 9\log_2 n \big\},
\\
    E_{34}^{(n)}
& := &
    \cap_{k=-9}^{9}\big\{H\big(Z_{2k}^\uparrow\big)\geq \log n+C_2\log_2 n  \big\}
\cap
    \cap_{k=-9}^{9}\big\{H\big(Z_{2k+1}^\downarrow\big)\geq \log n+C_2\log_2 n  \big\},
\\
    E_{35}^{(n)}
& := &
    \bigg\{
    \sum_{k=-9}^{9} \ell\big(Z_{2k}^\uparrow\big)
    +
    \sum_{k=-9}^{9} \ell\big(Z_{2k+1}^\downarrow\big)
    \leq
    (\log n)^{2+\delta_1}
    \bigg\}.
\end{eqnarray*}
Recall that $Z_0^\uparrow=_{law}\mathcal T_V^\uparrow$
and
$
    Y_{-1}^\uparrow
%=_{law}
%    \zeta\big(\mathcal T_V^\downarrow\big)
=_{law}
    \mathcal T_{V^-}^{\uparrow*}
$
(with $h=\widetilde h_n$).
So by Lemma \ref{Lemma_Proba_46} eq. \eqref{Ineg_avec_E20_Lemma},
there exists $n_{11}\geq n_{10}$ such that for all $n\geq n_{11}$,
\begin{eqnarray}
    \p\big[\big(E_{31}^{(n)}\big)^c\big]
\leq
    \p\big[T_{\mathcal T_V^\uparrow}\big(\widetilde h_n\big) \leq \Gamma_n\big]
+
    \p\big[T_{\mathcal T_{V^-}^{\uparrow*}}\big(\widetilde h_n\big) \leq \Gamma_n\big]
\leq
    \e(\log n)^{-2}.
\label{Ineg_avec_E20}
\end{eqnarray}
%We now prove that
%\begin{equation}\label{Ineg_Proba_E21}
%    \p\big[\big(E_{32}^{(n)}\big)^c\big]
%=
%    O\big((\log_2 n)(\log n)^{-1}\big).
%\end{equation}
Notice that,  using Theorem \ref{Lemma_Law_of_Slopes} {\bf (i)} with its notation and $h=\widetilde h_n$,
since $Z_0^\uparrow\big(T_{Z_0^\uparrow}\big(\widetilde h_n\big)\big)\in\big[\widetilde h_n, \widetilde h_n+C_0\big[$
by ellipticity \eqref{eq_ellipticity_for_V},
there exists $n_{12}\geq n_{11}$ such that for all $n\geq n_{12}$,
%$\widetilde h_n\geq 10\log_2 n$ and
\begin{eqnarray}
    \p\big[\big(E_{32}^{(n)}\big)^c\big]
& = &
    \p\big[
    \exists k\in\big[T_{Z_0^\uparrow}\big(\widetilde h_n\big),
    \ell\big(Z_0^\uparrow\big)\big], \ Z_0^\uparrow(k)< 9\log_2 n
    \big]
\nonumber\\
& \leq &
    \p\big[T_V\big(9\log_2 n-\widetilde h_n\big)<M_{\widetilde h_n}^\sharp,\
    \min\nolimits_{[0, M_{\widetilde h_n}^\sharp]}V>-\widetilde h_n-C_0,\
    V\big(M_{\widetilde h_n}^\sharp\big)\geq 0
    \big]
\nonumber\\
& \leq &
    \p\big[T_V\big(\widetilde h_n-9\log_2 n\big)<T_V(-9\log_2 n-C_0)\big]
\leq
    22(\log_2 n)(\log n)^{-1},
~~~~
\label{Ineg_Proba_E21}
\end{eqnarray}
%where we used (see Theorem \ref{Lemma_Law_of_Slopes} {\bf (i)}) the fact that after hitting $\big[\widetilde h_n,+\infty\big[$,
%$\mathcal T_V^\uparrow$ has the law of $V$ up to
%its maximum (at $\ell(\mathcal T_V^\uparrow)$)
%before its first decrease of $\widetilde h_n$,
by the strong Markov property at $T_V\big(9\log_2 n-\widetilde h_n\big)$,
%its first hitting time of $]-\infty, 9\log_2 n]$ after $T_{\mathcal T_V^\uparrow}\big(\widetilde h_n\big)$,
%$\mathcal T_V^\uparrow(\ell(\mathcal T_V^\uparrow))\geq \widetilde h_n$,
\eqref{eq_ellipticity_for_V}
and \eqref{eqOptimalStopping2} since $n\geq n_{12}\geq n_3$.
%and so $\widetilde h_n-9\log_2 n>\log_2 n>C_0$ and .
We prove similarly that
$
    \p\big[\big(E_{33}^{(n)}\big)^c\big]
\leq
    22(\log_2 n)(\log n)^{-1}
$ for all $n\geq n_{12}$,
using Theorem \ref{Lemma_Law_of_Slopes_Right} {\bf (i)} instead of Theorem \ref{Lemma_Law_of_Slopes}.

Also, using \eqref{ineg_Proba_Excess_1} with $h=\widetilde h_n$,
there exists $n_{13}\geq n_{12}$ such that, for all $n\geq n_{13}$,
%Theorem \ref{Lemma_Law_of_Slopes} {\bf (i)} with its notation and \eqref{eqOptimalStopping2},
and all $-9\leq k \leq 9$, since $n_{13}\geq n_3$,
\begin{eqnarray*}
&&
    \p\big[H(Z_{2k}^\uparrow)< \log n+C_2\log_2 n\big]
=
    \p\big[H\big(\mathcal T_V^\uparrow\big) - \widetilde h_n< (2C_1+C_2)\log_2 n\big]
\\
& \leq &
    4(2C_1+C_2)(\log_2 n)(\log n)^{-1}.
\end{eqnarray*}
%\begin{eqnarray*}
%&&
%    \p\big[H(Z_{2k}^\uparrow)< \log n+C_2\log_2 n\big]
%%=
%%    \p\big[H\big(\mathcal T_V^\uparrow\big) < \log n+C_2\log_2 n\big]
%\\
%& \leq &
%    \p\big[V(\tilde \tau_1)< \log n+C_2\log_2 n-\widetilde h_n\big]
%\\
%& \leq &
%    \p\big[T_V\big((2C_1+C_2)\log_2 n-\widetilde h_n\big) < T_V((2C_1+C_2)\log_2 n)\big]
%\\
%& \leq &
%(2C_1+2C_2)(\log_2 n)(\log n)^{-1}.
%\end{eqnarray*}
This remains true for $H\big(Z_{2k}^\uparrow\big)$ replaced by
%$H\big(\mathcal T_{V_-}^\uparrow\big)$ and for
$
    H\big(Z_{2k+1}^\downarrow\big)
=_{law}
    H\big(\mathcal T_V^\downarrow\big)
=_{law}
    H\big(\mathcal T_{-V}^\uparrow\big)
$,
$-9\leq k \leq 9$ by Theorem \ref{Lemma_Law_of_Slopes} {\bf (ii)} (or by the inequality after \eqref{ineg_Proba_Excess_1}).
Consequently, we get
$
    \p\big[\big(E_{34}^{(n)}\big)^c\big]
\leq
    152(2C_1+C_2)(\log_2 n)(\log n)^{-1}
%=O\big((\log_2 n)(\log n)^{-1}\big)
$
for all $n\geq n_{13}$.
%$\p\big[\big(E_{34}^{(n)}\cap E_{35}^{(n)} \cap E_{625}^{(n)}\big)^c\big]=O\big((\log_2 n)(\log n)^{-1}\big)$.

Moreover, we have
$
    \p\big[\big(E_{35}^{(n)}\big)^c\big]
\leq
    19\p\big[\ell\big(\mathcal T_V^\uparrow\big)>(\log n)^{2+\delta_1}/50\big]
    +
    19\p\big[\ell\big(\mathcal T_V^\downarrow\big)>(\log n)^{2+\delta_1}/50\big]
$
$
\leq
$
$
    38(\log n)^{-8}
%=O\big((\log n)^{-4}\big)
$
for all $n\geq n_{13}$
by Lemma \ref{Lemma_Proba_E4} eq. \eqref{Ineg_Proba_Longueur_slopes_1}
and \eqref{Ineg_Proba_Longueur_slopes_2}, since $n_{13}\geq n_9\geq p_3$.

Also,
using \eqref{eq_def_J6}  then
$\p\big[\big(E_i^{(n)}\big)^c\big]=o(1)$ for $31\leq i \leq 35$,
there exists $n_{14}\geq n_{13}$ such that, for all $n\geq n_{14}$ and all $z\leq \Gamma_n$,
$
\E\big(\ell\big(Z_0^\uparrow\big)+\ell\big(Z_1^\downarrow\big)\big)
\geq c_7(\log n)^2
$
by Lemma \ref{Lemma_Esperance_Longueur_Slope}
and
\begin{equation}
    J_6(n,z)
\leq
    \E\Bigg(
        \frac{\un_{\{-z+\Gamma_n <\ell(Z_0^\uparrow)\}}
              \un_{\cap_{i=31}^{35} E_i^{(n)}}
        }
        {\E\big(\ell\big(Z_0^\uparrow\big)+\ell\big(Z_1^\downarrow\big)\big)}
    \frac
    {
        \Big(\sum_{i=0}^{\ell(Z_0^\uparrow)-1}e^{-Z_0^\uparrow(i)}
             +\sum_{i=1}^{\ell(Y_{-1}^\uparrow)}e^{-Y_{-1}^\uparrow(i)}
%             +\sum_{i=1}^{\ell(\mathcal T_{V_-}^\uparrow)}e^{-\mathcal T_{V_-}^\uparrow(i)}
        \Big)
    }
    {
        \Big(\sum_{i=0}^{\ell(Z_0^\uparrow)-1}e^{-Z_0^\uparrow(i)}
             +\sum_{i=1}^{\ell(Y_{-1}^\uparrow)}e^{-Y_{-1}^\uparrow(i)}
%             +\sum_{i=1}^{\ell(\mathcal T_{V_-}^\uparrow)}e^{-\mathcal T_{V_-}^\uparrow(i)}
        \Big)
    }
    \Bigg)
    +
    \frac{\e}{(\log n)^2}.
\label{Ineg_J6_quotient_sum}
\end{equation}
The next step is to deal with the sums in numerator in the previous expectation.
Notice that on $E_{32}^{(n)}\cap E_{35}^{(n)}$, we have $\ell\big(Z_0^\uparrow\big)\leq (\log n)^3-1$ and so
\begin{equation}\label{Ineq_apres_TTV_avantlTV}
    \sum_{i=T_{Z_0^\uparrow}(\widetilde h_n)}^{\ell(Z_0^\uparrow)-1}e^{-Z_0^\uparrow(i)}
\leq
    \big[\ell(Z_0^\uparrow)-T_{Z_0^\uparrow}\big(\widetilde h_n\big)\big](\log n)^{-9}
\leq
    (\log n)^{-6}.
\end{equation}
Similarly,
$
    \sum_{i=T_{Y_{-1}^\uparrow}(\widetilde h_n)}^{\ell(Y_{-1}^\uparrow)}e^{-Y_{-1}^\uparrow(i)}
%    \sum_{i=T_{\mathcal T_{V_-}^\uparrow}(\widetilde h_n)}^{\ell(\mathcal T_{V_-}^\uparrow)}e^{-\mathcal T_{V_-}^\uparrow(i)}
\leq
    (\log n)^{-6}
$ on $E_{33}^{(n)}\cap E_{35}^{(n)}$
since $\ell\big(Y_{-1}^\uparrow\big)=\ell\big(Z_{-1}^\downarrow\big)$.

Also
using Theorem \ref{Lemma_Law_of_Slopes} {\bf (i)} since $Z_0^\uparrow=_{law}\mathcal T_V^\uparrow$,
then applying Proposition \ref{Lemma_Laplace_V_Conditionne},
for large $n$, for all $i\geq \Gamma_n$,
$$
    \E\Big[\exp\big(-Z_0^\uparrow(i)\big){\bf 1}_{\{i < T_{Z_0^\uparrow}(\widetilde h_n)\}}\Big]
=
    \E\Big[ e^{-V(i)}\un_{\{i < T_V(\widetilde h_n)\}} \mid T_V(\widetilde h_n)<T_V(\R_-^*)\Big]
\leq
    c_{13} i^{-3/2}.
$$
This remains true with $Z_0^\uparrow$ and $V$ replaced
by $Y_{-1}^\uparrow=_{law} \mathcal T_{V^-}^{\uparrow*}$ and $V_-$,
and $T_V(\R_-^*)$ by $T_{V^-}^*(\R_-)$ by Theorem \ref{Lemma_Law_of_Slopes_Right} {\bf (i)}
and Proposition \ref{Lemma_Laplace_V_Conditionne}.
So there exists $c_{32}>0$ and $n_{15}\geq n_{14}$ such that, for all $n\geq n_{15}$,
\begin{equation}
    \E\Bigg(
        %\sum_{i=\lfloor (\log n)^{4/3+\delta_1}\rfloor}^{T_{\mathcal T_{V_{\pm}}^\uparrow}(\widetilde h_n)}
        \sum\limits_{\Gamma_n\leq i < T_{Z_0^\uparrow}(\widetilde h_n)}
        e^{-Z_0^\uparrow(i)}
    \Bigg)
=
    \sum_{i=\Gamma_n}^\infty
    \E\bigg[
        e^{-Z_0^\uparrow(i)}
        {\bf 1}_{\big\{i< T_{Z_0^\uparrow}(\widetilde h_n)\big\}}
    \bigg]
    %\E\left[\frac{{\bf 1}_{\big\{i\leq T_{\mathcal T_{V_{\pm}}^\uparrow}(\widetilde h_n)\big\}}}{e^{\mathcal T_{V_{\pm}}^\uparrow(i)}}\right]
\leq
    \sum_{i=\Gamma_n}^\infty \frac{c_{13}}{i^{3/2}}
\leq
    \frac{c_{32}}{(\log n)^{\frac{2}{3}+\frac{\delta_1}{2}}},
\label{ineq_middle_of_slope}
\end{equation}
since $\Gamma_n=\big\lfloor (\log n)^{4/3+\delta_1}\big\rfloor$.
This remains true with $Z_0^\uparrow$ replaced by
$Y_{-1}^\uparrow$.

Combining \eqref{Ineg_J6_quotient_sum} with \eqref{Ineq_apres_TTV_avantlTV}, \eqref{ineq_middle_of_slope},
and the corresponding inequalities for $V_-$ and
$Y_{-1}^\uparrow$,
%$\mathcal T_{V_-}^\uparrow$,
$
    \ell\big(Z_0^\uparrow\big)
\geq
    T_{Z_0^\uparrow}\big(\widetilde h_n\big)
>
    \Gamma_n
$
and
$
    \ell\big(Y_{-1}^\uparrow\big)
%    \ell\big(\mathcal T_{V_-}^\uparrow\big)
%=
%    \ell\big(Z_{-1}^\downarrow\big)
\geq
    T_{Y_{-1}^\uparrow}\big(\widetilde h_n\big)
>
    \Gamma_n
$
on $E_{31}^{(n)}$,
$\sum_{i=0}^{\ell(Z_0^\uparrow)-1}e^{-Z_0^\uparrow(i)}\geq 1$
%$
%    e^{-Y_{-1}^\uparrow(\ell(Y_{-1}^\uparrow))}
%\leq
%    \exp\big(-\widetilde h_n\big)
%$
and %once more
again
$
\E\big(\ell\big(Z_0^\uparrow\big)+\ell\big(Z_1^\downarrow\big)\big)
%\sim_{n\to+\infty} c\big(\widetilde h_n\big)^2\sim_{n\to+\infty}
\geq c_7(\log n)^2
$,
there exists $n_{16}\geq n_{15}$ such that, for all $n\geq n_{16}$,
for every $j\in\{0,1\}$, for all $z\leq \Gamma_n$ (although $J_6$ does not depend on $j$),
\begin{eqnarray}
    J_6(n,z)
& \leq &
    \E\Bigg(
        \frac{\un_{\{-z+\Gamma_n <\ell(Z_0^\uparrow)\}}
              \un_{\cap_{i=31}^{35} E_i^{(n)}}
        }
        {\E\big(\ell\big(Z_0^\uparrow\big)+\ell\big(Z_1^\downarrow\big)\big)}
\nonumber\\
&&
    \frac
    {
        \Big(\sum_{i=0}^{\Gamma_n-1}e^{-Z_0^\uparrow(i)}
             +\sum_{i=1}^{\Gamma_n-1}e^{-Y_{-1}^\uparrow(i)}
%             +\sum_{i=1}^{\lfloor (\log n)^{4/3+\delta_1}\rfloor-1}e^{-\mathcal T_{V_-}^\uparrow(i)}
        \Big)
    }
    {
        \Big(\sum_{i=0}^{\ell(Z_0^\uparrow)-1}e^{-Z_0^\uparrow(i)}
             +\sum_{i=1}^{\ell(Y_{-1}^\uparrow)}e^{-Y_{-1}^\uparrow(i)}
%             +\sum_{i=1}^{\ell(\mathcal T_{V_-}^\uparrow)}e^{-\mathcal T_{V_-}^\uparrow(i)}
        \Big)
    }
    \Bigg)
    +
    2\e(\log n)^{-2}
\nonumber\\
& \leq &
    J_7(j, n, z)
    +
    2\e(\log n)^{-2},
\label{Ineg_J6_quotient_sum_J7}
\end{eqnarray}
where for $j\in\{0,1\}$,
\begin{align}
&
    J_7(j, n, z)
\\
& :=
    \E\Bigg(
        \frac{\un_{\{-z+\Gamma_n <\ell(Z_0^\uparrow)\}}
              \un_{\cap_{i=31}^{35} E_i^{(n)}}
        }
        {\E\big(\ell\big(Z_0^\uparrow\big)+\ell\big(Z_1^\downarrow\big)\big)}
\nonumber\\
&
    \quad
    \frac
    {
        \Big(\sum_{k=-\Gamma_n}^{\Gamma_n}
                \big[
                    e^{-Z_0^\uparrow(-(k+j))}{\bf 1}_{\{k+j\leq 0\}}
                    +
                    e^{-Y_{-1}^\uparrow(k+j)}{\bf 1}_{\{k+j>0\}}
%                    e^{-\mathcal T_{V_-}^\uparrow(k+j)}{\bf 1}_{\{k+j>0\}}
                \big]
        \Big)
    }
    {
        \Big(\sum_{i=0}^{\ell(Z_0^\uparrow)-1}e^{-Z_0^\uparrow(i)}
             +\sum_{i=1}^{\ell(Y_{-1}^\uparrow)}e^{-Y_{-1}^\uparrow(i)}
%             +\sum_{i=1}^{\ell(\mathcal T_{V_-}^\uparrow)}e^{-\mathcal T_{V_-}^\uparrow(i)}
        \Big)
    }
    \Bigg)
\label{eq_def_J7}
\\
& \leq
    \sum_{k=-\Gamma_n}^{\Gamma_n}
    \E\Bigg(
        \frac{\un_{\{-z-k<\ell(Z_0^\uparrow)\}\cap \cap_{i=31}^{35} E_i^{(n)}}
        \Big(
                    e^{-Z_0^\uparrow(-(k+j))}{\bf 1}_{\{k+j\leq 0\}}
                    +
                    e^{-Y_{-1}^\uparrow(k+j)}{\bf 1}_{\{k+j>0\}}
%                    e^{-\mathcal T_{V_-}^\uparrow(k+j)}{\bf 1}_{\{k+j>0\}}
        \Big)
        }
        {\E\big(\ell\big(Z_0^\uparrow\big)+\ell\big(Z_1^\downarrow\big)\big)
        \Big(\sum_{i=0}^{\ell(Z_0^\uparrow)-1}e^{-Z_0^\uparrow(i)}
             +\sum_{i=1}^{\ell(Y_{-1}^\uparrow)}e^{-Y_{-1}^\uparrow(i)}
%             +\sum_{i=1}^{\ell(\mathcal T_{V_-}^\uparrow)}e^{-\mathcal T_{V_-}^\uparrow(i)}
        \Big)
        }
    \Bigg).
\nonumber
\end{align}
%Recall that $z\leq -\Gamma_n$ in this first case.
Now, for $-\Gamma_n\leq k \leq \Gamma_n$,
applying Theorems \ref{Lemma_Independence_h_extrema} {\bf (i)}
and \ref{Lemma_Central_Slope} eq. \eqref{eq_Central_Slope_General_phi}, we have for every nonnegative measurable function $\varphi$,
since
$
    \{b_{\widetilde h_n}=z+k\}
=
    \big\{x_0\big(V ,\widetilde h_n\big)=z+k\big\}
    \cap
    \big\{\theta\big(T_0\big(V,\widetilde h_n\big)\big)\in \bigsqcup_{t\in\mathbb N^*}\mathbb R_+^t\big\}
=
    \big\{x_0\big(V ,\widetilde h_n\big)=z+k\big\}
    \cap
    \big\{  V\big(x_0\big(V, \widetilde h_n\big)\big)
            <
            V\big(x_1\big(V, \widetilde h_n\big)\big)
    \big\}
$
and
$
    \sharp\big\{0\leq i < \ell(\mathcal T_{V}^\uparrow), \, -i=z+k\big\}
=
    {\bf 1}_{\{-z-k<\ell(\mathcal T_{V}^\uparrow)\}}
$
for each $k\in\Z$ such that $z+k\leq 0$,
\begin{align}
&
    \E\big[
        \varphi\big(\theta\big(T_i\big(V,\widetilde h_n\big)\big), -18\leq i \leq 19\big)
        {\bf 1}_{\{b_{\widetilde h_n}=z+k\}}
    \big]
\nonumber\\
= &
    \E\bigg[
        \varphi\big(Z_{-18}^\uparrow, Z_{-17}^\downarrow, Z_{-16}^\uparrow, \dots, Z_{-2}^\uparrow,
            Z_{-1}^\downarrow, Z_0^\uparrow,
            Z_1^\downarrow, \dots, Z_{18}^\uparrow, Z_{19}^\downarrow
        \big)
        \frac{{\bf 1}_{\{-z-k<\ell(Z_0^\uparrow)\}}}
        {\E\big(\ell\big(Z_0^\uparrow\big)+\ell\big(Z_1^\downarrow\big)\big)}
    \bigg].
\label{eq_Application_Renewal_et_independence}
\end{align}
In the previous equality, $\theta(T_i(V,\widetilde h_n))$ becomes $Z_{i}^\uparrow$ or
$Z_{i}^\downarrow$ depending on the parity of $i$.

So, since
$Y_{-1}^\uparrow=\zeta\big(Z_{-1}^\downarrow\big)$
and $z\leq -\Gamma_n$ in this first case,
%$\mathcal T_{V_-}^\uparrow=\zeta\big(Z_{-1}^\downarrow\big)$,
we get, as explained below,
\begin{eqnarray}
    J_7(j, n, z)
& \leq &
    \sum_{k=-\Gamma_n}^{\Gamma_n}
    \E\Bigg(
    \frac{
%        \Big(
                    e^{-\theta[T_0(V, \tilde h_n)](-(k+j))}{\bf 1}_{\{k+j\leq 0\}}
                    +
                    e^{-\zeta[\theta(T_{-1}(V,\tilde h_n))](k+j)}{\bf 1}_{\{k+j>0\}}
%        \Big)
        }
        {
%        \Big(
        \sum_{i=0}^{\ell(T_0(V, \tilde h_n))-1}e^{-\theta[T_0(V, \tilde h_n)](i)}
             +\sum_{i=1}^{\ell(T_{-1}(V,\tilde h_n))}e^{-\zeta[\theta(T_{-1}(V,\tilde h_n))](i)}
%        \Big)
        }
\nonumber\\
&&
\qquad \qquad \qquad \qquad
        \un_{\{b_{\widetilde h_n}=z+k\} \cap E_3^{(n)} \cap E_5^{(n)} \cap E_6^{(n)} }
    \Bigg)
\nonumber\\
& = &
%    \sum_{k=-\lfloor (\log n)^{4/3+\delta_1}\rfloor}^{\lfloor (\log n)^{4/3+\delta_1}\rfloor}
%    \E\Bigg(
%        \frac{e^{-V(z-j)}\un_{\{b_{\widetilde h_n}=z+k\}\cap E_3^{(n)} \cap E_5^{(n)} \cap E_6^{(n)} }
%        }
%        {
%        \sum_{i=x_{-1}(V,\widetilde h_n)}^{x_1(V, \widetilde h_n)-1}e^{-V(i)}
%        }
%    \Bigg).
    \E\Bigg(
        \frac{
            \sum_{k=-\Gamma_n}^{\Gamma_n}
            e^{-V(b_{\log n}-k-j)}\un_{\{b_{\log n}=z+k\}\cap E_3^{(n)} \cap E_5^{(n)} \cap E_6^{(n)} }
        }
        {
        \sum_{i=x_{-1}(V,\log n)}^{x_1(V, \log n)-1}e^{-V(i)}
        }
    \Bigg),
\label{Ineg_J7_transforme}
\end{eqnarray}
for all $n\geq n_{16}$, $j\in\{0,1\}$ and $z\leq -\Gamma_n$.
Indeed, when applying \eqref{eq_Application_Renewal_et_independence} to the quantity after \eqref{eq_def_J7},
$E_{34}^{(n)}$ corresponds to (i.e. becomes)  a set $\widetilde E_{34}^{(n)}$ included in $E_3^{(n)}$,
on which we have in particular
$H\big[T_i\big(V, \widetilde h_n\big)\big]\geq \log n+C_2\log_2 n$ for all $-13\leq i \leq 13$
and so
$x_j\big(V,\widetilde h_n\big)=x_j(V, \log n)$ for $j\in\{-12,\dots, 12\}$
and so $b_{\tilde h_n}=b_{\log n}=x_0(V,\log n)$
(when $b_{\tilde h_n}=z+k\leq 0$);
$Z_0^\uparrow$ corresponds to
$
    \theta\big(T_0\big(V,\widetilde h_n\big)\big)
=
    \theta\big(T_0\big(V,\log n\big)\big)
=
    \big(V(b_{\log n}+i)-V(b_{\log n})$,
$0\leq i \leq x_1(V, \log n)-b_{\log n}\big)$
and
$Y_{-1}^\uparrow$ to
$
    \zeta\big(\theta\big[T_{-1}\big(V,\widetilde h_n\big)\big]\big)
=
    (V(b_{\log n}-i)-V(b_{\log n}),\ 0\leq i \leq b_{\log n}-x_{-1}(V,\log n))
$
so
$E_{31}^{(n)}$ corresponds to a set included in $E_6^{(n)}$
since $\widetilde h_n<\log n$,
$E_{35}^{(n)}$ corresponds to a set included in
$\big\{\big|x_{-12}\big(V, \widetilde h_n\big)-x_{12}\big(V,\widetilde h_n\big)\big|\leq (\log n)^{2+\delta_1}\big\}$,
 and the intersection of this and $\widetilde E_{34}^{(n)}$ is itself included in
 $\big\{\big|x_{-12}\big(V, \log n\big)-x_{12}\big(V,\log n\big)\big|\leq (\log n)^{2+\delta_1}\big\}$,
and so in
 $E_5^{(n)}$,
whereas $E_{32}^{(n)}$ and $E_{33}^{(n)}$ are not necessary anymore.

Notice that
$
    \sum_{k=-\Gamma_n}^{\Gamma_n}
    e^{-V(b_{\log n}-k-j)}
\leq
    \sum_{i=x_{-1}(V,\log n)}^{x_1(V, \log n)-1}e^{-V(i)}
=
    \sum_{i=M^-}^{M^+-1}e^{-V(i)}
$
on $E_6^{(n)}\cap\{b_{\log n}\leq 0\}\cap E_3^{(n)}$
with $M^\pm$ defined in \eqref{eq_def_Mpm}
since $V(x_{\pm 1})-V(x_0)\geq \log n+C_0$ for $n\geq n_{16}\geq n_3$.
Thus, using
Lemma \ref{Lemma_Proba_E3c},
there exists $n_{17}\geq n_{16}$ such that,
for all $n\geq n_{17}$, $j\in\{0,1\}$ and $z\leq -\Gamma_n$
(writing $E_i^{(n)}(z)$ instead of $E_i^{(n)}$ for $i\neq 3$),
\begin{equation}\label{Ineg_J7_fin}
    J_7(j, n, z)
\leq
    \E\Bigg(
        \frac{
            \sum_{k=-\Gamma_n}^{\Gamma_n}
            e^{-V(z-j)}\un_{\{z=b_{\log n}-k\}
%            e^{-V(b_{\log n}-k-j)}\un_{\{b_{\log n}=z+k\}
            \cap \cap_{\ell=3}^6 E_\ell^{(n)}(z)  }
        }
        {
        \sum_{i=M^-}^{M^+-1}e^{-V(i)}
        }
    \Bigg)
+
    \e(\log n)^{-2},
%    o\big((\log n)^{-2}\big).
\end{equation}
where we write $E_\ell^{(n)}(z)$ for $E_\ell^{(n)}$ for $\ell\in\{3,5,6\}$ for convenience.
Hence, using \eqref{Ineg_J6_quotient_sum_J7},
then \eqref{Ineg_J7_fin},
$M^-< z<M^+$ on $E_6^{(n)}\cap{\{z=b_{\log n}-k\}}$ for $|k|\leq  \Gamma_n$
%with $M^\pm$ defined in \eqref{eq_def_Mpm},
and \eqref{eq_def_mi_chapeau}
gives for all $n\geq n_{17}$ and $z\in(2\Z+n)$ such that $z\leq -\Gamma_n$,
\begin{eqnarray*}
    2J_6(n,z)
& \leq &
    J_7(1, n, z)+J_7(0, n, z)+4\e(\log n)^{-2}
\\
& \leq &
    \E\Bigg(
        \frac{
            \sum_{k=-\Gamma_n}^{\Gamma_n}
            \widehat \mu_n(z)\un_{\{z=b_{\log n}-k\}\cap \cap_{\ell=3}^6 E_\ell^{(n)}(z)  }
        }
        {
        \sum_{i=M^-}^{M^+-1}e^{-V(i)}
        }
    \Bigg)
+
    6\e(\log n)^{-2}
\\
& = &
    \E\Big(
            \widehat \nu_n(z)\un_{\{|b_{\log n}-z|\leq \Gamma_n\}\cap \cap_{\ell=3}^6 E_\ell^{(n)}(z)  }
    \Big)
+
    6\e(\log n)^{-2},
\end{eqnarray*}
where we used $\widehat \mu_n(2\Z)=\widehat \mu_n(2\Z+1)=\sum_{i=M^-}^{M^+-1}e^{-V(i)}$
and the definition \eqref{eq_def_nu_hat} of $\widehat \nu_n$
since $n$ and $z$ have the same parity.

Applying Lemma \ref{Lemma_Approx_PQuenched_Nu},
there exists $n_{18}\geq n_{17}$ such that,
for all $n\geq n_{18}$ and all $z\in(2\Z+n)$ such that $z\leq -\Gamma_n$,
$$
    2J_6(n,z)
\leq
    \E\big(
            P_\omega[S_n=z]
            +5(\log n)^{-3}
    \big)
+
    6\e(\log n)^{-2}
\leq
    \PP(S_n=z)
+
    7\e(\log n)^{-2}.
$$
This, \eqref{Ineg_phi_infty_J6} and \eqref{Ineg_phi_infty_J6_2} lead to
$
    \frac{2\sigma^2}{(\log n)^2}\varphi_\infty\big(\frac{\sigma^2 z}{(\log n)^2}\big)
\leq
    \PP(S_n=z)
+
    13\e(\log n)^{-2}
%    o\big((\log n)^{-2}\big)
$
for all $n\geq n_{18}$ and all $z\in(2\Z+n)$ such that $z\leq -\Gamma_n$.
%which
%proves the lower bound in Theorem \ref{Th_Local_Limit_Sinai} in
%concludes
%this first case.

\noindent{\bf Second case:}
We now consider the case $z> \Gamma_n$.
%so that $z+k\leq 0$ for every $k$ in the sum in \eqref{Ineg_Proba_Zn_E\Gamma_1}.
We use the same
$Z_{2k}^\uparrow$, $Z_{2k+1}^\downarrow$,
$Y_{-1}^\uparrow=\zeta\big(Z_{-1}^\downarrow\big)$ and $E_i^{(n)}$
as in the first case.
%{\bf (mettre avant separation des cas ?)}

Using Lemma \ref{Lemma_Proba_bh_egal} with $x=z+\Gamma_n$,
%\eqref{eq_proba_bh_x_positif}, let
we have when $z+\Gamma_n>0$,
\begin{equation}\label{eq_def_J6c}
    J_6^+(n,z)
:=
    \p\big(b_{\widetilde h_n}= z+\Gamma_n\big)
=
    \frac{
        \p\big(z+\Gamma_n \leq \ell(Z_{-1}^\downarrow)\big)}
        {\E\big(\ell\big(Z_0^\uparrow\big)+\ell\big(Z_1^\downarrow\big)\big)}.
\end{equation}
Similarly as in \eqref{Ineg_phi_infty_J6},
for all $n\geq n_{18}$, for all $z\in]\Gamma_n, A_0(\log n)^2]$,
%Using the uniform continuity of $\varphi_\infty$ on $\R$, $\|\varphi_\infty\|_\infty<\infty$,
%$\widetilde h_n\sim_{n\to+\infty} \log n$, $z=O\big((\log n)^2\big)$
% and $\delta_1<2/3$ in the first inequality,
%then Theorem \ref{Th_Local_Limit_b_h} in the second one, we get
\begin{eqnarray}
    \frac{\sigma^2}{(\log n)^2}\varphi_\infty\bigg(\frac{\sigma^2 z}{(\log n)^2}\bigg)
\leq
    \frac{\sigma^2}{\big(\widetilde h_n\big)^2}\varphi_\infty\bigg(\frac{\sigma^2 \big(z+\Gamma_n\big)}{
    \big(\widetilde h_n\big)^2}\bigg)
    +\frac{2\e}{(\log n)^2}
\leq
%    \p\big(b_{\log n}= z-\Gamma_n\big)
%    +o\big((\log n)^{-2}\big)
%\\
%& = &
    J_6^+(n,z)
    +\frac{3\e}{(\log n)^2}.
\label{Ineg_phi_infty_J6_c}
\end{eqnarray}
Also
$
    \frac{\sigma^2}{(\log n)^2}\varphi_\infty\bigg(\frac{\sigma^2 z}{(\log n)^2}\bigg)
\leq
    J_6^+(n,z)
    +\frac{3\e}{(\log n)^2}
$
for all $n\geq n_{18}$ and all $z\geq \max(\Gamma_n, $ $A_0(\log n)^2)$ as in \eqref{Ineg_phi_infty_J6_2},
and so for all $z> \Gamma_n$.

Similarly as in \eqref{Ineg_J6_quotient_sum} and \eqref{Ineg_J6_quotient_sum_J7},
using
$
    \un_{\{z+\Gamma_n \leq \ell(Z_{-1}^\downarrow)\}}
\leq
    \un_{\{z+k\leq \ell(Z_{-1}^\downarrow)\}}
$
instead of
%using
$
    \un_{\{-z+\Gamma_n <\ell(Z_0^\uparrow)\}}
$
$
\leq
$
$
    \un_{\{-z-k<\ell(Z_0^\uparrow)\}}
$,
we get for all $n\geq n_{18}$ and all $z> \Gamma_n$,
$
    J_6^+(n,z)
\leq
    J_7^+(j, n, z)
    +
    2\e(\log n)^{-2}
%\label{Ineg_J6_quotient_sum_J7_c}
$ for each $j\in\{0,1\}$,
where for $j\in\{0,1\}$,
\begin{align}
&
    J_7^+(j, n, z)
\label{eq_def_J7plus}
\\
& :=
    \sum_{k=-\Gamma_n}^{\Gamma_n}
    \E\Bigg(
        \frac{\un_{\{z+k\leq \ell(Z_{-1}^\downarrow)\}\cap \cap_{i=31}^{35} E_i^{(n)}}
        \Big(
                    e^{-Z_0^\uparrow(-(k+j))}{\bf 1}_{\{k+j\leq 0\}}
                    +
                    e^{-Y_{-1}^\uparrow(k+j)}{\bf 1}_{\{k+j>0\}}
%                    e^{-\mathcal T_{V_-}^\uparrow(k+j)}{\bf 1}_{\{k+j>0\}}
        \Big)
        }
        {\E\big(\ell\big(Z_0^\uparrow\big)+\ell\big(Z_1^\downarrow\big)\big)
        \Big(\sum_{i=0}^{\ell(Z_0^\uparrow)-1}e^{-Z_0^\uparrow(i)}
             +\sum_{i=1}^{\ell(Y_{-1}^\uparrow)}e^{-Y_{-1}^\uparrow(i)}
%             +\sum_{i=1}^{\ell(\mathcal T_{V_-}^\uparrow)}e^{-\mathcal T_{V_-}^\uparrow(i)}
        \Big)
        }
    \Bigg).
\nonumber
\end{align}
Now, applying Theorems \ref{Lemma_Independence_h_extrema} {\bf (ii)}
and \ref{Lemma_Central_Slope} eq. \eqref{eq_Central_Slope_General_phi}, we have for every nonnegative measurable function $\varphi$,
since
$
    \{b_{\widetilde h_n}=z+k\}
=
    \big\{x_1\big(V ,\widetilde h_n\big)=z+k\big\}
    \cap
    \big\{\theta\big(T_0\big(V,\widetilde h_n\big)\big)\in \bigsqcup_{t\in\mathbb N^*}\mathbb R_-^t\big\}
$
when $z+k> 0$,
\begin{align}
&
    \E\big[
        \varphi\big(\theta(T_i(V,\widetilde h_n)), -17\leq i \leq 20\big)
        {\bf 1}_{\{b_{\widetilde h_n}=z+k\}}
    \big]
\nonumber\\
= &
    \E\bigg[
        \varphi\big(Z_{-18}^\uparrow, Z_{-17}^\downarrow, Z_{-16}^\uparrow, \dots,
            Z_{-2}^\uparrow,
            Z_{-1}^\downarrow, Z_0^\uparrow, \dots, Z_{18}^\uparrow, Z_{19}^\downarrow
        \big)
        \frac{{\bf 1}_{\{z+k\leq \ell(Z_{-1}^\downarrow)\}}}
        {\E\big(\ell\big(Z_0^\uparrow\big)+\ell\big(Z_1^\downarrow\big)\big)}
    \bigg].
\label{eq_Application_Renewal_et_independence_c}
\end{align}
In the previous equality, $\theta(T_i(V,\widetilde h_n))$ becomes $Z_{i-1}^\uparrow$ or
$Z_{i-1}^\downarrow$ depending on the parity of $i$.

So, since
$Y_{-1}^\uparrow=\zeta\big(Z_{-1}^\downarrow\big)$ and $z> \Gamma_n$,
%$\mathcal T_{V_-}^\uparrow=\zeta\big(Z_{-1}^\downarrow\big)$,
we get, similarly as in \eqref{Ineg_J7_transforme},
% and \eqref{Ineg_J7_fin},
 with $b_{\tilde h_n}=b_{\log n}=x_1(V,\log n)$,
and using the definition \eqref{eq_def_Mpm} of $M^\pm$ on $\{b_{\log n}>0\}$.
\begin{eqnarray}
    J_7^+(j, n, z)
& \leq &
    \sum_{k=-\Gamma_n}^{\Gamma_n}
    \E\Bigg(
    \frac{
%        \Big(
                    e^{-\theta[T_1(V, \tilde h_n)](-(k+j))}{\bf 1}_{\{k+j\leq 0\}}
                    +
                    e^{-\zeta[\theta(T_0(V,\tilde h_n))](k+j)}{\bf 1}_{\{k+j>0\}}
%        \Big)
        }
        {
%        \Big(
        \sum_{i=0}^{\ell(T_1(V, \tilde h_n))-1}e^{-\theta[T_1(V, \tilde h_n)](i)}
             +\sum_{i=1}^{\ell(T_0(V,\tilde h_n))}e^{-\zeta[\theta(T_0(V,\tilde h_n))](i)}
%        \Big)
        }
\nonumber\\
&&
\qquad \qquad \qquad \qquad
        \un_{\{b_{\widetilde h_n}=z+k\} \cap E_3^{(n)} \cap E_5^{(n)} \cap E_6^{(n)} }
    \Bigg)
\nonumber\\
& = &
    \E\Bigg(
        \frac{
            \sum_{k=-\Gamma_n}^{\Gamma_n}
            e^{-V(b_{\log n}-k-j)}\un_{\{b_{\log n}=z+k\}\cap E_3^{(n)} \cap E_5^{(n)} \cap E_6^{(n)} }
        }
        {
%        \sum_{i=x_0(V,\log n)}^{x_2(V, \log n)-1}e^{-V(i)}
        \sum_{i=M^-}^{M^+-1}e^{-V(i)}
        }
    \Bigg),
\label{Ineg_J7plus_bis}
\end{eqnarray}
for all $n\geq n_{18}$, $j\in\{0,1\}$ and $z> \Gamma_n$.

We conclude as in the first case that
$
    \frac{2\sigma^2}{(\log n)^2}\varphi_\infty\big(\frac{\sigma^2 z}{(\log n)^2}\big)
\leq
    \PP(S_n=z)
+
    13\e(\log n)^{-2}
%    o\big((\log n)^{-2}\big)
$
for all $n\geq n_{18}$ and all $z\in(2\Z+n)$ such that $z > \Gamma_n$.

%which proves the lower bound in Theorem \ref{Th_Local_Limit_Sinai} in this second case.

\noindent{\bf Third case:}
We finally consider the case
$-\Gamma_n<z \leq \Gamma_n$.

We use the same notation as in the first case.
Notice that \eqref{eq_def_J6},
\eqref{Ineg_phi_infty_J6}, \eqref{Ineg_phi_infty_J6_2},
    \eqref{Ineg_J6_quotient_sum},  \eqref{Ineg_J6_quotient_sum_J7}
and \eqref{eq_def_J7}
remain valid when $n\geq n_{18}$ and $-\Gamma_n<z \leq \Gamma_n$,
with the same definitions of $J_6$ and $J_7$.
However in this third case, that is, for every $n\geq n_{18}$ and $-\Gamma_n<z \leq \Gamma_n$,
for  $j\in\{0,1\}$,
\begin{equation}
    J_7(j, n, z)
\leq
    J_8(j, n, z)+J_9(j, n, z)+J_{10}(j, n, z),
\label{Ineg_J7_8_9_10}
\end{equation}
where
\begin{align*}
& J_8(j, n, z)
\\
& :=
    \sum_{k=-\Gamma_n}^{-z}
    \E\Bigg(
        \frac{\un_{\{-z-k<\ell(Z_0^\uparrow)\}\cap \cap_{i=31}^{35} E_i^{(n)}}
        \Big(
                    e^{-Z_0^\uparrow(-(k+j))}{\bf 1}_{\{k+j\leq 0\}}
                    +
                    e^{-Y_{-1}^\uparrow(k+j)}{\bf 1}_{\{k+j>0\}}
%                    e^{-\mathcal T_{V_-}^\uparrow(k+j)}{\bf 1}_{\{k+j>0\}}
        \Big)
        }
        {\E\big(\ell\big(Z_0^\uparrow\big)+\ell\big(Z_1^\downarrow\big)\big)
        \Big(\sum_{i=0}^{\ell(Z_0^\uparrow)-1}e^{-Z_0^\uparrow(i)}
             +\sum_{i=1}^{\ell(Y_{-1}^\uparrow)}e^{-Y_{-1}^\uparrow(i)}
%             +\sum_{i=1}^{\ell(\mathcal T_{V_-}^\uparrow)}e^{-\mathcal T_{V_-}^\uparrow(i)}
        \Big)
        }
    \Bigg),
\\
& J_9(j, n, z)
\\
& :=
    \sum_{k=-z+1}^{\Gamma_n}
    \E\Bigg(
        \frac{\big(1- \un_{\{z+k\leq \ell(Z_{-1}^\downarrow)\}} \big)
        \un_{\cap_{i=31}^{35} E_i^{(n)}}
        \Big(
                    e^{-Z_0^\uparrow(-(k+j))}{\bf 1}_{\{k+j\leq 0\}}
                    +
                    e^{-Y_{-1}^\uparrow(k+j)}{\bf 1}_{\{k+j>0\}}
%                    e^{-\mathcal T_{V_-}^\uparrow(k+j)}{\bf 1}_{\{k+j>0\}}
        \Big)
        }
        {\E\big(\ell\big(Z_0^\uparrow\big)+\ell\big(Z_1^\downarrow\big)\big)
        \Big(\sum_{i=0}^{\ell(Z_0^\uparrow)-1}e^{-Z_0^\uparrow(i)}
             +\sum_{i=1}^{\ell(Y_{-1}^\uparrow)}e^{-Y_{-1}^\uparrow(i)}
%             +\sum_{i=1}^{\ell(\mathcal T_{V_-}^\uparrow)}e^{-\mathcal T_{V_-}^\uparrow(i)}
        \Big)
        }
%\\
%&
%\qquad\qquad\qquad\qquad
%\frac{
%    }
%    {
%    }
    \Bigg),
\\
& J_{10}(j, n, z)
\\
&
 :=
    \sum_{k=-z+1}^{\Gamma_n}
    \E\Bigg(
        \frac{\un_{\{z+k\leq \ell(Z_{-1}^\downarrow)\}}\un_{\cap_{i=31}^{35} E_i^{(n)}}
        \Big(
                    e^{-Z_0^\uparrow(-(k+j))}{\bf 1}_{\{k+j\leq 0\}}
                    +
                    e^{-Y_{-1}^\uparrow(k+j)}{\bf 1}_{\{k+j>0\}}
%                    e^{-\mathcal T_{V_-}^\uparrow(k+j)}{\bf 1}_{\{k+j>0\}}
        \Big)
        }
        {\E\big(\ell\big(Z_0^\uparrow\big)+\ell\big(Z_1^\downarrow\big)\big)
        \Big(\sum_{i=0}^{\ell(Z_0^\uparrow)-1}e^{-Z_0^\uparrow(i)}
             +\sum_{i=1}^{\ell(Y_{-1}^\uparrow)}e^{-Y_{-1}^\uparrow(i)}
%             +\sum_{i=1}^{\ell(\mathcal T_{V_-}^\uparrow)}e^{-\mathcal T_{V_-}^\uparrow(i)}
        \Big)
        }
%\\
%&
%\qquad\qquad\qquad\qquad
%\frac{
%    }
%    {
%    }
    \Bigg).
\end{align*}
We first notice that, since
$
    -\ell(Z_0^\uparrow)+2
\leq
    -\Gamma_n+1
\leq
    -z+1
\leq
    \Gamma_n
\leq
    \ell(Y_{-1}^\uparrow)-1
$
on $E_{31}^{(n)}$
for $-\Gamma_n<z \leq \Gamma_n$,
and using
$
    \un_{\{z+\Gamma_n\leq \ell(Z_{-1}^\downarrow)\}}
\leq
    \un_{\{z+k\leq \ell(Z_{-1}^\downarrow)\}}
$,
there exists $n_{19}\geq n_{18}$ such that, for all $n\geq n_{19}$, all $-\Gamma_n<z \leq \Gamma_n$,
and all $j\in\{0,1\}$,
\begin{eqnarray}
    0
\leq
    J_9(j, n, z)
& \leq &
    \E\Bigg(
        \frac{1- \un_{\{z+\Gamma_n\leq \ell(Z_{-1}^\downarrow)\}}
        }
        {\E\big(\ell\big(Z_0^\uparrow\big)+\ell\big(Z_1^\downarrow\big)\big)
        }
    \Bigg)
=
    \p\big(b_{\tilde h_n}=0\big)
-
    \p\big(b_{\tilde h_n}=z+\Gamma_n\big)
\nonumber\\
& \leq &
    \frac{\sigma^2}{\big(\widetilde h_n\big)^2}
    \bigg(\varphi_\infty(0)
        -
        \varphi_\infty\bigg(\frac{\sigma^2(z+\Gamma_n)}{\big(\widetilde h_n\big)^2}\bigg)
    \bigg)
    +\e(\log n)^{-2}/4
\nonumber\\
& \leq &
    \e(\log n)^{-2}/2
\label{Ineg_J9}
\end{eqnarray}
by Lemma \ref{Lemma_Proba_bh_egal}, then Theorem \ref{Th_Local_Limit_b_h},
and finally by continuity of $\varphi_\infty$ since $\delta_1<2/3$ and $|z|\leq \Gamma_n$.

In order to prove an inequality for $J_8(j, n, z)$, we can do the same proof as in the first case
from the line following \eqref{eq_def_J7} to \eqref{Ineg_J7_fin},
replacing
$\sum_{k=-\Gamma_n}^{\Gamma_n}$
by
$\sum_{k=-\Gamma_n}^{-z}$
since $|z|\leq \Gamma_n$ (so $z+k\leq 0$), which gives,
for all $n\geq n_{19}$, all $-\Gamma_n<z \leq \Gamma_n$ and all $j\in\{0,1\}$,
\begin{equation}\label{ineg_J8}
    J_8(j, n, z)
\leq
    \E\Bigg(
        \frac{
            \sum_{k=-\Gamma_n}^{-z}
            e^{-V(z-j)}\un_{\{z=b_{\log n}-k\}
            \cap \cap_{\ell=3}^6 E_\ell^{(n)}(z)  }
        }
        {
        \sum_{i=M^-}^{M^+-1}e^{-V(i)}
        }
    \Bigg)
+
    \e(\log n)^{-2}/4.
\end{equation}
In order to prove an inequality for $J_{10}(j, n, z)$,
we can do the same proof as in the second case, between the definition \eqref{eq_def_J7plus} of $J_7^+$
and \eqref{Ineg_J7plus_bis},
replacing
$\sum_{k=-\Gamma_n}^{\Gamma_n}$
by
$\sum_{k=-z+1}^{\Gamma_n}$
since $|z|\leq \Gamma_n$ (so $z+k>0$),
then using once more Lemma \ref{Lemma_Proba_E3c} as in \eqref{Ineg_J7_fin},
we get
for all $n\geq n_{19}$, all $-\Gamma_n<z \leq \Gamma_n$ and all $j\in\{0,1\}$,
$$
    J_{10}(j, n, z)
\leq
    \E\Bigg(
        \frac{
            \sum_{k=-z+1}^{\Gamma_n}
            e^{-V(z-j)}\un_{\{z=b_{\log n}-k\}
%            e^{-V(b_{\log n}-k-j)}\un_{\{b_{\log n}=z+k\}
            \cap \cap_{\ell=3}^6 E_\ell^{(n)}(z)  }
        }
        {
        \sum_{i=M^-}^{M^+-1}e^{-V(i)}
        }
    \Bigg)
+
    \e(\log n)^{-2}/4.
$$
This, \eqref{Ineg_J7_8_9_10}, \eqref{Ineg_J9}  and \eqref{ineg_J8}
prove that  \eqref{Ineg_J7_fin} remains true
for all $n\geq n_{19}$, all $-\Gamma_n<z \leq \Gamma_n$ and all $j\in\{0,1\}$.

Since \eqref{Ineg_J6_quotient_sum_J7}, \eqref{Ineg_phi_infty_J6} and \eqref{Ineg_phi_infty_J6_2}
also remain true, we conclude as in the first case that
for all $n\geq n_{19}$ and all $z\in(2\Z+n)$ such that $-\Gamma_n<z \leq \Gamma_n$,
we have
$
    \frac{2\sigma^2}{(\log n)^2}\varphi_\infty\big(\frac{\sigma^2 z}{(\log n)^2}\big)
\leq
    \PP(S_n=z)
+
    13\e(\log n)^{-2}
%    o\big((\log n)^{-2}\big)
$.

Finally, combining the conclusions of the three cases gives for all $n\geq n_{19}$,
\begin{equation*}
    \sup_{z\in(2\Z+n)}
    \bigg[
    \frac{2\sigma^2}{(\log n)^2}\varphi_\infty\bigg(\frac{\sigma^2 z}{(\log n)^2}\bigg)
-
    \PP\big(S_n=z\big)
    \bigg]
\leq
    13\e(\log n)^{-2},
\end{equation*}
which proves the lower bound in Theorem \ref{Th_Local_Limit_Sinai}.
This and \eqref{Ineg_Minor_Sur_EL_Bis}
prove Theorem \ref{Th_Local_Limit_Sinai}.
\hfill $\Box$

%%%%%%%%%%%%%%%%%%%%%%%%%%%%%%%%%%%%%%%%%%%%%%%%%%%%%%%%%%%%%%%%%%%%%%%%%%%%%%%%%%%%%%%

%%%%%%%%%%%%%%%%%%%%%%%%%%%%%%%%%%%%%%%%%%%%%%%%%%%%%%%%%%%%%%%%%%%%%%%%%%%%%%%%%%%%%%%

\section{Some estimates concerning the environment}\label{Sect_Technical}

\subsection{Probabilities of  $\big(E_5^{(n)}\big)^c$ and $(E_6^{(n)})^c$}
The aim of this subsection is to give upper bounds of some probabilities related to the events $E_i^{(n)}$,
which are defined between equations \eqref{eq_def_E1} and \eqref{eq_def_E5}.

\begin{lem}\label{Lemma_Proba_E4}
There exists $p_3\geq 2$ such that,
\begin{equation}\label{Ineg_Proba_E4}
    \forall n\geq p_3,
\qquad
    \p\big[(E_5^{(n)})^c\big]
\leq
    (\log n)^{-7}.
\end{equation}
%Also, let $\mathcal T_V^\uparrow$, $\mathcal T_V^\downarrow$ and $\mathcal T_{V_-}^\uparrow$
%having respectively the law of (non central) left $\widetilde h_n$-upward slope of $V$ ,
%left $\widetilde h_n$-downward slope of $V$,
%and left $\widetilde h_n$-upward slope $V_-$,
%where $\widetilde h_n=\log n-2C_1\log_2 n$ as before.
Also, we have for $n\geq p_3$, with $\widetilde h_n=\log n-2C_1\log_2 n$ as before,
    \begin{equation}\label{Ineg_Proba_Longueur_slopes_1}
    \p\big[\ell\big(\mathcal T_{V,\widetilde h_n}^\uparrow\big)> (\log n)^{2+\delta_1}/50\big]
\leq
    (\log n)^{-8},
\end{equation}
\begin{equation}\label{Ineg_Proba_Longueur_slopes_2}
    \p\big[\ell\big(\mathcal T_{V_-, \widetilde h_n}^\uparrow\big)> (\log n)^{2+\delta_1}/50\big]
=
    \p\big[\ell\big(\mathcal T_{V,\widetilde h_n}^\downarrow\big)> (\log n)^{2+\delta_1}/50\big]
\leq
    (\log n)^{-8}.
\end{equation}
\end{lem}

\noindent{\bf Proof:}
The idea is to approximate $V$ by a two-sided Brownian motion,
in order to transfer to $V$ some results already known for Brownian motions.

To this aim, we recall the definition of $h$-extrema introduced by Neveu et al. \cite{NP} for continuous functions.
If $w$ is a continuous function $\R\to\R$, $h>0$, and $y\in\R$, it is said that $w$ admits an {\it $h$-minimum at $y$}
if there exists real numbers $u$ and $v$ such that $u<y<v$,
$w(y)=\inf\{w(z),\ z\in[u,v]\}$, $w(u)\geq w(y)+h$ and $w(v)\geq w(y)+h$.
It is said that $w$ admits an {\it $h$-maximum} at $y$ if $-w$ admits an $h$-minimum at $y$.
In these two cases we say that $w$ admits an {\it $h$-extremum} at $y$.
Notice that contrary to Definition \ref{def_left_extrema}, all the inequalities are large.

It is known (see \cite{Cheliotis}, Lemma 8) that, when $w=W$ or $w=\sigma W$, almost surely,
{\bf (a)} $w$ is continuous on $\R$;
{\bf (b)} for every $h>0$, the set of $h$-extrema of $w$ can  be written
$\{x_k(w,h),\ k\in\Z\}$, where $(x_k(w,h))_{k\in\Z}$ is strictly increasing, unbounded from above and below, with $x_0(w,h)\leq 0<x_1(w,h)$;
{\bf (c)} for all $h>0$ and $k\in\Z$, $x_{k+1}(w,h)$ is an $h$-maximum if and only if $x_{k}(w,h)$ is an $h$-minimum
(we use the same notation as for left extrema of $V$ since no confusion is possible).

According to a slightly modified version
(see e.g. \cite{DGP_Collision_Sinai}, Lemma 4.3, with $(\log n)^\alpha$ replaced by $K$ and a single potential $V$
instead of two)
of the Koml\'os--Major--Tusn\'ady strong approximation
theorem (see Koml\'os et al. \cite{KMT}), there exist (strictly) positive constants  $C_3$
%$c_{10}$
and
$C_4$, independent of $K\in\N^*$, such that, possibly in an enlarged probability
space, there exists a two-sided standard Brownian
motion $(W(t),\ t\in\R)$, such that
\begin{equation*}
%    \mathcal{B}_1(K)
    E_{36}(K)
:=
    \left\{\sup_{-K\leq t\leq K} \Big|V(\lfloor t\rfloor)-\sigma W(t)\Big|\leq C_3\log
K\right\}
\end{equation*}
satisfies
$\p([E_{36}(K)]^c)\leq K^{-C_4}$
%$\p([\mathcal{B}_1(K)]^c)\leq K^{-C_4}$
%$\p([\mathcal{B}_1(K)]^c)\leq c_{10} K^{-C_4}$
for large $K$.

Let $n\geq n_3$ and $\alpha>0$, and recall that $0<\delta_1<2/3$.
We define $h_n':=\log n+3 C_3(3+8/C_4)\log_2 n$.
On $E_{36}\big((\log n)^{3+8/C_4}\big)$, consider,
if they exist,  two consecutive $h'_n$-minima for $\sigma W$,
denoted by $y_i:=x_i(\sigma W, h'_n)$ and $y_{i+2}:=x_{i+2}(\sigma W, h'_n)$,
such that
$|y_i|\leq \alpha(\log n)^{2+\delta_1}$ and $|y_{i+2}|\leq \alpha(\log n)^{2+\delta_1}$.
Let $z_{i+1}:=\min\big\{k\in[\lfloor y_i\rfloor , \lfloor y_{i+2}\rfloor]\cap\Z,\ V(k)=\max_{[\lfloor y_i\rfloor , \lfloor y_{i+2}\rfloor]} V\big\}$.
We have, using $\omega\in E_{36}\big((\log n)^{3+8/C_4}\big)$ in the second and forth inequalities,
for $n$ large enough so that $(\log n)^{3+8/C_4}> \alpha(\log n)^{2+\delta_1}$,
%uniformly for large $n$,
\begin{align*}
    V(z_{i+1})
& =
    \max_{[\lfloor y_i\rfloor , \lfloor y_{i+2}\rfloor]} V
\geq
    V(\lfloor x_{i+1}(\sigma W, h'_n)\rfloor)
\geq
    \sigma W[x_{i+1}(\sigma W, h'_n)]-C_3(3+8/C_4)\log_2 n
%    \sigma W[x_{i+1}(\sigma W, h'_n)]-C_3(9+3/C_4)\log_2 n
\\
& \geq
    \sigma W[x_i(\sigma W, h'_n)]+h'_n
    -C_3(3+8/C_4)\log_2 n
%    -C_3(9+3/C_4)\log_2 n
\\
& \geq
    V(\lfloor y_i\rfloor)+h'_n-2C_3(3+8/C_4)\log_2 n
\geq
    V(\lfloor y_i\rfloor)+\log n.
\end{align*}
We prove similarly that
$
    V(z_{i+1})
\geq
    V(\lfloor y_{i+2}\rfloor)+\log n
$, and so $\lfloor y_i\rfloor  < z_{i+1} < \lfloor y_{i+2}\rfloor$.
Since $\max_{[\lfloor y_i\rfloor, z_{i+1}[}V<V(z_{i+1})$
and $\max_{]z_{i+1}, \lfloor y_{i+2}\rfloor]}V\leq V(z_{i+1})$, $z_{i+1}$ is a left $(\log n)$-maximum for $V$.

So we have proved that for large $n$ on $E_{36}\big((\log n)^{3+8/C_4}\big)$,
between two consecutive $h'_n$-minima for $\sigma W$
belonging to the interval $\big[-\alpha(\log n)^{2+\delta_1}, \alpha(\log n)^{2+\delta_1}\big]$,
there is at least one left $(\log n)$-maximum for $V$.
Notice in particular that for such $n$, on $E_{36}\big((\log n)^{3+8/C_4}\big)$,
if $x_{17}(\sigma W, h'_n)$ $\leq \alpha(\log n)^{2+\delta_1}$,
then in $[x_1(\sigma W, h'_n),$ $x_{17}(\sigma W, h'_n)]$,
there are at least eight consecutive $h'_n$-minima for $\sigma W$,
and then at least seven left $(\log n)$-maxima for $V$,
and so $x_{13}(V, \log n)\leq x_{17}(\sigma W, h'_n)\leq \alpha(\log n)^{2+\delta_1}$.
Hence for large $n$,
\begin{eqnarray}
&&
    \p\big[x_{13}(V,\log n)> \alpha(\log n)^{2+\delta_1}, E_{36}\big((\log n)^{3+8/C_4}\big)\big]
\nonumber\\
& \leq &
    \p\big[x_{17}(\sigma W, h'_n) > \alpha(\log n)^{2+\delta_1}\big]
\nonumber\\
& \leq &
    \sum_{i=0}^{16} \p\bigg[\ell(T_i(\sigma W, h'_n))> \frac{\alpha(\log n)^{2+\delta_1}}{17}\bigg],
~~~~~~
\label{Ineg_Proba_Large_Slope_1}
\end{eqnarray}
where $\ell(T_i(w, h)):=x_{i+1}(w, h)-x_i(w, h)$ for $i\in\Z$, $h>0$ and any continuous function $w$,    is the
{\it length of the $i$-th $h$-slope of $w$}.

The length of a non central $1$-slope of $W$,
that is, $\ell(T_i(W, 1))$ for $i\neq 0$, has a density, which is (see \cite{Cheliotis}, eq. (7))
$
	f_\ell(x)
:=
%	\frac{\pi}{2} \sum_{k\in\Z} (-1)^k (k+1/2) \exp\big( -\pi^2(k+1/2)^2 x/2 \big){\bf 1}_{\R_+^*}(x)
	\pi \sum_{k\in\N} (-1)^k (k+1/2) \exp\big( -\pi^2(k+1/2)^2 x/2 \big){\bf 1}_{\R_+^*}(x)
$.
Also, the length of the central $1$-slope $\ell(T_0(W, 1))$ has a density, which is (see \cite{Cheliotis}, eq. (10)) equal to
$
	f_{\ell(T_0)}(x)
:=
	x f_\ell(x)
$.
Notice that
$
	f_\ell(x)
\leq
    (\pi/2) \exp[-\pi^2 x/8]
$
for large $x$. Hence for large $x$,
$
	f_{\ell(T_0)}(x)
\leq
    \exp[-\pi^2 x/10]
$
and
$
	f_\ell(x)
\leq
    \exp[-\pi^2 x/10]
$.
Thus,
 $ \p\big[\ell(T_i(W, 1))>u) = O( \exp(-\pi^2 u/10) )$ as $u\to+\infty$ for any $i\in\Z$,
so for large $n$,
\begin{align*}
&
    \p\big[\ell(T_i(\sigma W, h'_n))> \alpha (\log n)^{2+\delta_1}/17\big]
 =
    \p\big[\ell(T_i(W, 1))> \sigma^2 \alpha (\log n)^{2+\delta_1}/(17(h'_n)^2)\big]
\\
&
\leq
    \p\big[\ell(T_i(W, 1))>  \sigma^2 \alpha (\log n)^{\delta_1}/20\big]
 =
    O\big( \exp(-\pi^2 \sigma^2 \alpha (\log n)^{\delta_1}/200)\big),
%    O\big( n^{-\sigma^2 \pi^2/200}\big),
\end{align*}
as $n\to+\infty$,
where we used
$
    \ell(T_i(\sigma W, h'_n))
=
    \ell(T_i(W, h'_n/\sigma))
=_{law}
    (h'_n/\sigma)^2\ell(T_i(W, 1))
$
by scaling.
This, \eqref{Ineg_Proba_Large_Slope_1} and $\p([E_{36}(K)]^c)\leq K^{-C_4}$ for large $K$    lead to
\begin{eqnarray}
&&
    \p\big[x_{13}(V,\log n)> \alpha(\log n)^{2+\delta_1}\big]
\nonumber\\
& \leq &
    O\big( \exp(-\pi^2 \sigma^2 \alpha (\log n)^{\delta_1}/200)\big)
%    O\big(n^{-\sigma^2 \pi^2/200}\big)
    +
    \p\big[\big(E_{36}\big((\log n)^{3+8/C_4}\big)\big)^c\big]
\leq
    (\log n)^{-8}
\label{Ineg_Proba_x13}
\end{eqnarray}
for large $n$.
We prove similarly that
$
    \p\big[x_{-12}(V,\log n)< -\alpha(\log n)^{2+\delta_1}\big]
\leq
    (\log n)^{-8}
$.
Finally,
$$
    \p\big[(E_5^{(n)})^c\big]
\leq
    \p\big[x_{12}(V,\log n)> (\log n)^{2+\delta_1}\big]
+
    \p\big[x_{-12}(V,\log n)< -(\log n)^{2+\delta_1}\big]
\leq
    (\log n)^{-7}
$$
for large $n$, which proves \eqref{Ineg_Proba_E4}.

Since $x_3(V,\widetilde h_n) \leq x_3(V, \log n) <x_{13}(V,\log n)$, we get
$$
    \p\big[x_3(V,\widetilde h_n)> (\log n)^{2+\delta_1}/50\big]
\leq
    \p\big[x_{13}(V,\log n)> (\log n)^{2+\delta_1}/50\big]
\leq
    (\log n)^{-8}
$$
for large $n$ by \eqref{Ineg_Proba_x13}.
Since $x_3\big(V,\widetilde h_n\big)>x_3\big(V,\widetilde h_n\big) - x_1\big(V,\widetilde h_n\big)$, which
has the same law as
$\ell\big(\mathcal T_{V,\widetilde h_n}^\uparrow\big)+\ell\big(\mathcal T_{V,\widetilde h_n}^\downarrow\big)$,
by Theorem \ref{Lemma_Independence_h_extrema},
this gives
$$
    \p\big[\ell\big(\mathcal T_{V,\widetilde h_n}^\uparrow\big)> (\log n)^{2+\delta_1}/50\big]
\leq
    \p\big[x_3(V,\widetilde h_n)> (\log n)^{2+\delta_1}/50\big]
\leq
    (\log n)^{-8}
$$
and similarly
$
    \p\big[\ell\big(\mathcal T_{V,\widetilde h_n}^\downarrow\big)> (\log n)^{2+\delta_1}/50\big]
\leq
    (\log n)^{-8}
$ for large $n$.
Since
$
    \ell\big(\mathcal T_{V_-,\widetilde h_n}^\uparrow\big)
=_{law}
    \ell\big(\mathcal T_{-V,\widetilde h_n}^\uparrow\big)
=_{law}
  \ell\big(\mathcal T_{V,\widetilde h_n}^\downarrow\big),
$
by Theorem \ref{Lemma_Law_of_Slopes} {\bf (ii)},

this proves \eqref{Ineg_Proba_Longueur_slopes_1} and \eqref{Ineg_Proba_Longueur_slopes_2}
up to a change of $p_3$,
which ends the proof of the lemma.
\hfill$\Box$

We now turn to the probability of $(E_6^{(n)})^c\cap E_5^{(n)}$.

\begin{lem}\label{Lemma_Proba_46}
Recall that  $\delta_1\in]0, 2/3[$.
There exist $c_{33}>0$ and $p_2\in\N$ such that
\begin{equation}\label{Ineg_Proba_E5cE4}
    \forall n\geq p_2,
\qquad
    \p\big[(E_6^{(n)})^c\cap E_5^{(n)}\big]
\leq
    \exp\big[-c_{33}(\log n)^{2/3-\delta_1}\big]
\leq
    (\log n)^{-3}.
\end{equation}
We now consider left and right $\widetilde h_n$-slopes.
As $n\to+\infty$,
\begin{eqnarray}
%    \forall n\geq p_2,
%\qquad
    \p\Big[T_{\mathcal T_{V_\pm, \widetilde h_n}^\uparrow}\big(\widetilde h_n\big)\leq \Gamma_n\Big]
=
    o\big((\log n)^{-2}\big),
\qquad
    \p\Big[T_{\mathcal T_{V_\pm, \widetilde h_n}^{\uparrow*}}\big(\widetilde h_n\big)\leq \Gamma_n\Big]
=
    o\big((\log n)^{-2}\big),
\label{Ineg_avec_E20_Lemma}
\end{eqnarray}
recalling that
$T$,
$T_{\mathcal T_{V, h}^{\uparrow}}$,
$T_{\mathcal T_{V, h}^{\downarrow}}$,
$T_{\mathcal T_{V, h}^{\uparrow*}}$ and
$T_{\mathcal T_{V, h}^{\downarrow*}}$
are defined in \eqref{eq_def_TY},
Definition \ref{deff_law_slopes},
\eqref{eq_def_slope_etoile_1}
and
\eqref{eq_def_slope_etoile_2},
and that $V_{\pm}(.)=V(\pm .)$.
%where $T_{\mathcal T_{V_+}^\uparrow}(x)$
%\big(resp. $T_{\mathcal T_{V_-}^\uparrow}(x)$\big) is for $x>0$ the hitting time of $[x,+\infty[$
%by a left upward $\widetilde h_n$-slope $\mathcal T_{V_+}^\uparrow$ for $V_+=V$
%\big(resp. $\mathcal T_{V_-}^\uparrow$ for $V_-(.)=V(-.)$, and similarly for $\mathcal T_{V_\pm}^{\uparrow*}$\big).
\end{lem}

\noindent{\bf Proof:}
First, for $n\geq n_3$, $b\in\Z$ and $0<|i|\leq \Gamma_n$,
we have by Hoeffding's inequality (see \cite{Hoeffding63}, Theorem 2),
\begin{eqnarray}
    \p\big[V(b+i)-V(b)\geq \log n\big]
& = &
    \p\big[V(i)\geq \log n\big]
\leq
    \exp\big[-2(\log n)^2/\big(|i| \big(2C_0)^2\big)\big]
\nonumber\\
& \leq &
    \exp\big[-c_{34}(\log n)^2/|i|\big]
\leq
    \exp\big[-c_{34}(\log n)^{2/3-\delta_1}\big]~~~~~~~
\label{eq_Hoeffding}
\end{eqnarray}
with $c_{34}:=1/(2 C_0^2)>0$, since $V(i)$ is the sum of $|i|$
independent random variables with zero mean, bounded by $\pm C_0$
by ellipticity (see \eqref{eq_ellipticity_for_V}).

Notice that on $\big(E_6^{(n)}\big)^c\cap E_5^{(n)}$, there exists
$b=b_{\log n}\in\Z$ and $i\in\Z$ such that
$V(b+i)-V(b)\geq \log n$, $|i|\leq \Gamma_n$
and $|b|\leq (\log n)^{2+\delta_1}$ since $\omega\in E_5^{(n)}$.
Thus by \eqref{eq_Hoeffding},
\begin{eqnarray*}
    \p\big[(E_6^{(n)})^c\cap E_5^{(n)}\big]
& \leq &
    \sum_{|b|\leq \lfloor (\log n)^{2+\delta_1}\rfloor }\sum_{|i|\leq \lfloor (\log n)^{4/3+\delta_1}\rfloor}
    \p\big[V(b+i)-V(b)\geq \log n\big]
\\
& \leq &
    9(\log n)^5 \exp\big[-c_{34}(\log n)^{2/3-\delta_1}\big],
\end{eqnarray*}
since $0<\delta_1<2/3$.
This proves \eqref{Ineg_Proba_E5cE4}, e.g. with $c_{33}:=c_{34}/2$.

Now, notice that, using the law of $\mathcal T_{V,\widetilde h_n}^\uparrow$
provided by Theorem \ref{Lemma_Law_of_Slopes} {\bf (i)},
then \eqref{eq_Proba_Atteinte_logn_avant0} and once more
Hoeffding's inequality and $\widetilde h_n\sim_{n\to+\infty}\log n$,
\begin{eqnarray*}
&&
    \p\big[T_{\mathcal T_{V, \widetilde h_n}^\uparrow}\big(\widetilde h_n\big)\leq \Gamma_n\big]
\\
& = &
    \p\big[T_V\big(\widetilde h_n\big)\leq \Gamma_n, T_V\big(\widetilde h_n\big)<T_V(\R_-^*)\ \big]
    /\p\big[T_V\big(\widetilde h_n\big)<T_V(\R_-^*)\ \big]
\\
& \leq &
%    \frac{2\log n}{c^*}
    \frac{2\log n}{c_1}
    \sum_{i=1}^{\Gamma_n} \p\big[V(i)\geq \widetilde h_n\big]
\leq
%    \frac{2\log n}{c^*}
    \frac{2\log n}{c_1}
    \sum_{i=1}^{\lfloor (\log n)^{4/3+\delta_1}\rfloor}
    \exp\big[-c_{34}\big(\widetilde h_n\big)^2/i\big]
\\
& \leq &
    %(2/c^*)
    (2/c_1)
    (\log n)^3\exp[-c_{34}(\log n)^{2/3-\delta_1}/2]
\end{eqnarray*}
for large $n$.
This proves \eqref{Ineg_avec_E20_Lemma} for $V_+$ since $0<\delta_1<2/3$.
The proof for $V_-$ is similar.
The proof for $T_{\mathcal T_{V_\pm, \widetilde h_n}^{\uparrow*}}$ is the similar, with
Theorem \ref{Lemma_Law_of_Slopes_Right} and $c_1^*$ instead of
Theorem \ref{Lemma_Law_of_Slopes} and $c_1$.\hfill$\Box$

%%%%%%%%%%%%%%%%%%%%%%%%%%%%%%%%%%%%%%%%%%%%%%%%%%%%%%%%%%%%%%%%%%%%%%%%%

%%%%%%%%%%%%%%%%%%%%%%%%%%%%%%%%%%%%%%%%%%%%%%%%%%%%%%%%%%%%%%%%%%%%%%%%%%%

\subsection{Laplace transform of $V$ conditioned to stay positive or nonnegative}

The main tools of this subsection are local limit theorems for random walks conditioned to stay positive,
by Vatutin and Wachtel (\cite{Vatutin_Wachtel}, Theorems 4 and 6 and Lemma 12 with $\alpha=2$ and $\rho=1/2$).

%We first assume \eqref{eq_Hpothesis_4}. In particular,  conditioning $V$ to be $>0$ and $\geq 0$ are the same.
%We define for $n\geq 3$,
We define for $h\geq 0$,
with $T_V$ and $T_V^*$ defined in \eqref{eq_def_TY},
 and \eqref{eq_def_TY*},
\begin{eqnarray}\label{eq_Def_A_n}
    \Xi_h
&  := &
    \Big\{\inf_{[1, T_V([h,+\infty[)]}V\geq 0\Big\}
=
    \{T_V(h)<T_V(\R_-^*)\},
\\
    \Xi_h^*
&  := &
    \Big\{\inf_{[1, T_V([h,+\infty[)]}V> 0\Big\}
=
    \{T_V(h)<T_V^*(]-\infty,0])\}.
\nonumber
\end{eqnarray}
%\begin{eqnarray}\label{eq_Def_A_n}
%    A_n
%&  := &
%    \Big\{\inf_{[1, T_V([\log n,+\infty[)]}V\geq 0\Big\}
%=
%    \{T_V(\log n)<T_V(\R_-^*)\},
%\\
%    A_n^*
%&  := &
%    \Big\{\inf_{[1, T_V([\log n,+\infty[)]}V> 0\Big\}
%=
%    \{T_V(\log n)<T_V^*(]-\infty,0])\}.
%\nonumber
%\end{eqnarray}

The aim of this subsection is to prove the following uniform upper bound:

%{\bf (dans le lemme qui suit, remplace $\log n$ par $h$, et $\Xi_n$
%par un analogue avec $T_V(h)$ ; il faudra donc adapter le $p_8$ utilise dans $n_3$ etc ; OK pour $p_8$ ie $p_5$)}

\begin{prop}\label{Lemma_Laplace_V_Conditionne}
There exist $c_{13}>0$, $p_4>0$ and $p_5>0$ such that
\begin{equation*}
    \forall x\geq p_4, \forall h\geq p_5,
\qquad
    \E\big[ e^{-V(x)}\un_{\{x< T_V(h)\}} \mid T_V(h)<T_V(\R_-^*)\big]
\leq
    c_{13} x^{-3/2}.
\end{equation*}
This remains true when $T_V(\R_-^*)$ is replaced by $T_V^*(]-\infty,0])$.
\end{prop}

Before proving this lemma, we introduce some notation and some technical lemmas.
First, let
%{\bf (remplacement de $C_x$ par $G_x$ )}
\begin{equation}\label{eq_Def_CxPrime}
    G_x
:=
    \{\forall 1\leq k\leq x,\ V(k)\geq 0\},
\qquad
    G_x^*
:=
    \{\forall 1\leq k\leq x,\ V(k)> 0\},
\qquad
    x> 0.
\end{equation}
We know (due to the Spitzer and R\`osen theorem, see
Vatutin and Wachtel \cite{Vatutin_Wachtel} eq. (18),
or \cite{Bingham_Goldie_Teugels} Theorem 8.9.23, p. 382) that
\begin{equation}\label{eq_Proba_CxPrime}
    \p[G_x]
\sim_{x\to+\infty}
    c_{35} x^{-1/2},
\qquad
    \p[G_x^*]
\sim_{x\to+\infty}
    c_{35}^* x^{-1/2},
\end{equation}
where $c_{35}>0$ and $c_{35}^*>0$.

The following (uniform) estimates are maybe already known.
However we did not find them in the literature, so we provide their proof.

\begin{lem}\label{Lem_Proba_An_Ordre2}
For large $h>0$, for every $0\leq z <h$,
\begin{eqnarray}\label{Ineg_Pz_An}
    \frac{z-\E^z[V(T_V(\R_-^*))]}{h}
    -\frac{3C_0(z+C_0)}{h^2}
& \leq &
    \p^z(\Xi_{h})
\leq
    \frac{z-\E^z[V(T_V(\R_-^*))]}{h},
\\
    \frac{z-\E^z[V(T_V^*(\R_-))]}{h}
    -\frac{3C_0(z+C_0)}{h^2}
& \leq &
    \p^z(\Xi_{h}^*)
\leq
    \frac{z-\E^z[V(T_V^*(\R_-))]}{h}.
~~~~
\label{Ineg_Pz_An_Etoile}
\end{eqnarray}
Also, for $z=0$,
\begin{equation}\label{eq_Proba_An}
    h\p[\Xi_h]
\to_{h\to+\infty}
%\sim_{h\to+\infty}
    -\E[V(T_V(\R_-^*))]=:c_1>0,
\end{equation}
\begin{equation}\label{eq_Proba_An_Etoile}
    h\p[\Xi_h^*]
\to_{h\to+\infty}
%\sim_{h\to+\infty}
    -\E[V(T_V^*(\R_-))]
=:
    c_1^*>0.
\end{equation}
\end{lem}

\noindent{\bf Proof:}
Let $h>0$, $U_h:=T_V([h,+\infty[)\wedge T_V(\R_-^*)$, and $0\leq z <h$.
Since $(V(k),\ k\geq 0)$ is under $\p^z$ a martingale starting at $z$ for its natural filtration due to \eqref{eqRecurrence},
and $|V(k\wedge U_h)|\leq h +C_0$ a.s. for every $k\in\N$ thanks to ellipticity \eqref{eq_ellipticity_for_V},
the optimal stopping theorem gives
\begin{equation}
    z
=
    \E^z[V(U_h)]
=
%    \E^z[V(U_h){\bf 1}_{A_n}]
%    +
%    \E^z[V(U_h){\bf 1}_{A_n^c}]
%\nonumber\\
%& = &
    \E^z[V(T_V([h,+\infty[)){\bf 1}_{\Xi_{h}}]
    +
    \E^z[V(T_V(\R_-^*)){\bf 1}_{(\Xi_{h})^c}].
\label{Ineg_Unif_PAn_1}
\end{equation}
Since $h \leq V(T_V([h,+\infty[)) \leq h+C_0$ a.s. by ellipticity, we have
\begin{equation}
    h\p^z[\Xi_h]
\leq
    \E^z[V(T_V([h,+\infty[)){\bf 1}_{\Xi_h}]
\leq
    (h+C_0)\p^z[\Xi_h].
\label{Ineg_Unif_PAn_2}
\end{equation}
Also, $-C_0 \leq V[T_V(\R_-^*)]\leq 0$ a.s. by ellipticity, thus
\begin{eqnarray}
    \E^z[V(T_V(\R_-^*))]
& \leq &
    \E^z[V(T_V(\R_-^*)){\bf 1}_{(\Xi_h)^c}]
=
    \E^z[V(T_V(\R_-^*))]-    \E^z[V(T_V(\R_-^*)){\bf 1}_{\Xi_h}]
\nonumber\\
& \leq &
    \E^z[V(T_V(\R_-^*))]
    +C_0\p^z(\Xi_h).
\label{Ineg_Unif_PAn_3}
\end{eqnarray}
Hence, using first \eqref{Ineg_Unif_PAn_2}  and \eqref{Ineg_Unif_PAn_1}
and then the first inequality in \eqref{Ineg_Unif_PAn_3},
\begin{equation}
    h\p^z[\Xi_h]
\leq
    z-
    \E^z[V(T_V(\R_-^*)){\bf 1}_{(\Xi_h)^c}]
\leq
    z-
    \E^z[V(T_V(\R_-^*))].
\label{Ineg_Unif_PAn_4}
\end{equation}
Similarly,
$$
    (h+C_0)
    \p^z[\Xi_h]
\geq
        z-
        \E^z[V(T_V(\R_-^*)){\bf 1}_{(\Xi_h)^c}]
\geq
        z
        - \E^z[V(T_V(\R_-^*))]-C_0\p^z(\Xi_h),
$$
and so for large $h$ for every $0\leq z <h$, since $z+C_0\geq z-\E^z[V(T_V(\R_-^*))]\geq z \geq 0$,
\begin{eqnarray*}
    \p^z[\Xi_h]
& \geq &
    \frac{z-\E^z[V(T_V(\R_-^*))]}{h +2C_0}
\geq
    \frac{z-\E^z[V(T_V(\R_-^*))]}{h}
    \bigg(1-\frac{3C_0}{h}\bigg)
\\
& \geq &
    \frac{z-\E^z[V(T_V(\R_-^*))]}{h}
    -\frac{3C_0(z+C_0)}{h^2}.
\end{eqnarray*}
This and \eqref{Ineg_Unif_PAn_4} prove \eqref{Ineg_Pz_An}.
The proof of \eqref{Ineg_Pz_An_Etoile} is similar.
We get \eqref{eq_Proba_An} and \eqref{eq_Proba_An_Etoile} as a consequence.
\hfill$\Box$

In order to apply the results of Vatutin et al. (\cite{Vatutin_Wachtel}, thm. 4 and 6),
we introduce some of its notation (see its pages 177 and 179).
Let $\chi^+:=V(\tau^+)$, where $\tau^+:=\min\{k\geq 1,\ V(k)>0\}=T_V(\R_+^*)$, and
$\chi_k^+$, $k\geq 1$ be independent copies of $\chi^+$. We can now define the (left-continuous) renewal function
$$
    H(u)
:=
    \un_{\{u>0\}}+\sum_{k=1}^\infty \p\big(\chi^+_1+\dots+ \chi^+_k<u\big),
\qquad
    u\in\R.
$$
Also it is well known that $H(x)<\infty$ for every $x\in\R$
(see e.g. \cite{Vatutin_Wachtel} Lem. 13).

As in \cite{Vatutin_Wachtel} (page 180), we say that
the random variable $\log\frac{1-\o_0}{\o_0}$ is $(\ell,a)$-lattice if its distribution is lattice with span $\ell>0$ and shift $a\in[0,\ell[$,
that is, if $\ell$ is the maximal real number such that
the support of the distribution of $\log\frac{1-\o_0}{\o_0}$ is included in the set $(a+\ell\Z)=\{a+k\ell, k\in\Z\}$.
We say that the random variable $\log\frac{1-\o_0}{\o_0}$ is {\it non-lattice}
if its distribution is not supported in $(a+\ell\Z)$ for any $a\in\R$, $\ell>0$.
The two following lemmas are a bit more precise that what is needed in the present paper.
They may be of independent interest and will be useful in a work in progress \cite{Devulder_Rates_CV}.

\begin{lem} \label{LemLaplaceSachantC}
Assume
%\eqref{eq_Hpothesis_4}and
that $\log\frac{1-\o_0}{\o_0}$ is  non-lattice.
We have for $p\geq 0$,
%$p\in\{0,1\}$,
\begin{equation}\label{eqEquivalentLaplaceV}
    \E\Big[\big(V(x)\big)^p e^{-V(x)} | G_x^*\Big]
\sim_{x\to+\infty}
    \frac{f_2(p)}{x},
\qquad
    f_2(p)
:=
    %\frac{g_{2,0}(0)}{\sigma c_{35}^*}
    \frac{1}{c_{35}^*\sigma \sqrt{2\pi}}
    \int_0^\infty u^p e^{-u} H(u)\dd u
    \in]0,\infty[.
\end{equation}
\end{lem}

The case $p=1$ was already proved in Afanasyev et al. (\cite{Afanasyev_et_al_2012}, Prop. 2.1)
and Hirano (\cite{Hirano_98} Lemma 5)
with different methods.

\noindent{\bf Proof of Lemma \ref{LemLaplaceSachantC}:}
We fix $p\geq 0$, and define $\beta_p:=\sup_{y\geq 0}(y^p e^{-y/9})\in]0,\infty[$.
%$p\in\{0,1\}$.
%Let $\Delta>0$.
We first observe that for large $x$,
\begin{equation}\label{eqEsperanceLaplaceVxGrand}
    \E\big[\big(V(x)\big)^p e^{-V(x)}\un_{\{V(x)\geq 9\log x\}} | G_x^*\big]
\leq
%    9e^{-1}
    \beta_p
    \E\big[e^{-8V(x)/9}\un_{\{V(x)\geq 9\log x\}} | G_x^*\big]
\leq
%    9e^{-1}
    \beta_p
    x^{-8}.
\end{equation}
%since $y^p e^{-y/9}\leq 9 e^{-1}$ for all $y\geq 0$.
%In order to deal with $\E\big[V(x) e^{-V(x)} \un_{\{V(x)< 9\log x\}}| G_x'\big]$, we introduce the following notation.

Our potential $V$ is a random walk with i.i.d. bounded, non constant and zero mean jumps $\rho_x$, $x\in\Z$
by \eqref{eqEllipticity}, \eqref{eqRecurrence} \eqref{eq_def_sigma} and \eqref{eqDefPotentialV},
and by hypotheses, its jumps have a non lattice distribution.
%and satisfying \eqref{eq_Hpothesis_4}.
%Also due to the same hypothesis, $\{V(k)>0\  \forall 1\leq k \leq x\}=G_x'$.
So we can use the following result
(\cite{Vatutin_Wachtel}, Theorem 4 with $\alpha=2$, $\beta=0$ and $c_x\sim_{x\to+\infty}\sigma\sqrt{x}$,
as seen  in the line after its eq. (3))
 and with $g_{2,0}(0)=1/\sqrt{2\pi}$: for $\Delta>0$,
    \begin{equation}\label{eqVatutinWachtelNonLattice}
    \sigma\sqrt{x} \p\Big[V(x)\in[y,y+\Delta[\, \big|\, G_x^*\Big]
\sim_{x\to+\infty}
    \frac{1}{x\p[G_x^*]\sqrt{2\pi}}\int_{y}^{y+\Delta}H(u)\dd u
%    \frac{g_{2,0}(0)}{x\p[G_x^*]}\int_{y}^{y+\Delta}H(u)\dd u
\end{equation}
uniformly in $y\in]0, \delta_x \sqrt{x}]$, where $\delta_x\to 0$ as $x\to +\infty$.

We prove that this convergence is in fact uniform in $y\in[0, \delta_x \sqrt{x}]$ as $x\to+\infty$.
To this aim, notice that for fixed $x>0$ and $\Delta>0$, $\p\big[V(x)\in[y,y+\Delta[\, \big|\, G_x^*\big]$
tends to
$
    \p\big[V(x)\in]0,\Delta]\, \big|\, G_x^*\big]
=
    \p\big[V(x)\in[0,\Delta]\, \big|\, G_x^*\big]
$
as
%$y\downarrow 0$,
$y\to 0$ with $y>0$,
since
$\p\big[V(x)=0\, \big|\, G_x^*\big]=0$ by definition \eqref{eq_Def_CxPrime} of $G_x^*$.
Now, fix some $\e>0$.
Using the uniformity in $y\in]0, x^{1/4}]$ in  \eqref{eqVatutinWachtelNonLattice},
 there exists $A_\e>0$ such that for all $x>A_\e$, for all $y\in]0, x^{1/4}]$,
\begin{equation}\label{eq_pour_uniformite}
    1-\e
\leq
    \frac{x\p[G_x^*]\sigma\sqrt{x} \p\big[V(x)\in[y,y+\Delta[\, \big|\, G_x^*\big]}
    {    \frac{1}{\sqrt{2\pi}}\int_{y}^{y+\Delta}H(u)\dd u}
\leq
    1+\e.
\end{equation}
Letting $y\downarrow 0$ in \eqref{eq_pour_uniformite} for fixed $x>A_\e$ and using the convergence before \eqref{eq_pour_uniformite},
\eqref{eq_pour_uniformite} remains true with $[y,y+\Delta[$ and $\int_y^{y+\Delta}$ replaced respectively by
$[0,\Delta]$ and $\int_0^\Delta$.
Hence,
\begin{equation}\label{Ineg_Delta_Intermediaire_en_0}
    x\p[G_x^*]\sigma\sqrt{x} \p\Big[V(x)\in[0,\Delta]\, \big|\, G_x^*\Big]
\to_{x\to+\infty}
    %g_{2,0}(0)
    \frac{1}{\sqrt{2\pi}}
    \int_{0}^{\Delta}H(u)\dd u.
\end{equation}
Moreover, applying once more \eqref{eqVatutinWachtelNonLattice}
with $[y, y+\Delta[$ replaced by $[\Delta-\eta,\Delta+\eta[$
%for $0<\eta<\Delta$
for fixed $\Delta$ and $0<\eta<\Delta$ gives,
for large $x$, $
    x\p[G_x^*]\sigma \sqrt{x} \p\big[V(x)=\Delta\, \big|\, G_x^*\big]
\leq
    (2/\sqrt{2\pi})\int_{\Delta-\eta}^{\Delta+\eta}H(u)\dd u\
\leq
    (4/\sqrt{2\pi})H(2\Delta)\eta
    %to_{\eta\to 0} 0
$.
Since this is true for any $\eta>0$, we get
$
    x\p[G_x^*]\sigma \sqrt{x} \p\big[V(x)=\Delta\, \big|\, G_x^*\big]
\to
    0
$
as $x\to+\infty$.
So, \eqref{Ineg_Delta_Intermediaire_en_0} remains true with $[0,\Delta]$ replaced by $[0,\Delta[$.
This and \eqref{eqVatutinWachtelNonLattice}
prove that the convergence in \eqref{eqVatutinWachtelNonLattice} is in fact uniform in $y\in[0, \delta_x \sqrt{x}]$ as $x\to+\infty$,
where $\delta_x\to 0$ as $x\to +\infty$.

So, we have for any $\e>0$ and $\Delta>0$, for large $x$,
\begin{eqnarray*}
    \E\bigg[\frac{(V(x))^p}{e^{V(x)}} \un_{\{V(x)< 9\log x\}}\, \big|\, G_x^*\bigg]
& = &
    \sum_{k=0}^\infty
    \E\bigg[\frac{(V(x))^p}{e^{V(x)}} \un_{\{V(x)< 9\log x\}}\un_{\{V(x)\in[k\Delta, (k+1)\Delta[\}} \mid G_x^*\bigg]
%\\
%& \leq &
%    \sum_{k=0}^\infty
%    \E\bigg[\frac{((k+1)\Delta)^p}{e^{k\Delta}} \un_{\{V(x)< 9\log x\}}\un_{\{V(x)\in[k\Delta, (k+1)\Delta[]\}} \mid G_x^*\bigg]
%\\
%& = &
%    \sum_{k=0}^\infty
%    \int_{k\Delta}^{(k+1)\Delta} ye^{-y}\un_{\{y<9\log x\}} \p\big[V(x)\in \dd y|G_x'\big]
\\
& \leq &
    \sum_{k=0}^{\lfloor 9\Delta^{-1}\log x\rfloor } \frac{((k+1)\Delta)^p}{e^{k \Delta}}
    \p\big[V(x)\in [k\Delta, (k+1)\Delta[\, \mid\, G_x^*\big]
\\
& \leq &
%    \frac{(1+\e)}{\sigma\sqrt{x}}
    %\frac{g_{2,0}(0)}{x\p[G_x^*]}
    \frac{(1+\e)}{\sigma x^{3/2}\p[G_x^*]\sqrt{2\pi}}
    \sum_{k=0}^{\lfloor 9\Delta^{-1}\log x\rfloor } \frac{((k+1)\Delta)^P}{e^{k \Delta}}
    \int_{k\Delta}^{(k+1)\Delta}H(u)\dd u
\\
& \leq &
    \frac{(1+\e)^2}{c_{35}^* \sigma \sqrt{2\pi}  x}
    %\frac{(1+\e)^2g_{2,0}(0)}{\sigma c_{35}^* x}
    \sum_{k=0}^{\lfloor 9\Delta^{-1}\log x\rfloor } \frac{((k+1)\Delta)^p}{e^{k \Delta}}
    \Delta H\big[(k+1)\Delta\big],
\end{eqnarray*}
where we used \eqref{eq_Proba_CxPrime} and since $H$ is nondecreasing.
So,
\begin{eqnarray*}
&&
    \limsup_{x\to+\infty}
    \Big(
    x\E\big[\big(V(x)\big)^p e^{-V(x)} \un_{\{V(x)< 9\log x\}} \mid G_x^*\big]
    \Big)
\\
& \leq &
%    \frac{(1+\e)^2g_{2,0}(0)}{\sigma c_{35}^*}
    \frac{(1+\e)^2}{c_{35}^*\sigma \sqrt{2\pi}}
    e^\Delta
    \sum_{k=0}^{\infty} \Delta
    \frac{((k+1)\Delta)^p}{e^{(k+1) \Delta}} H[(k+1)\Delta]
%\\
\,
\to_{\Delta\to 0}\,
    %\frac{(1+\e)^2 g_{2,0}(0)}{\sigma c_{35}^*}
    \frac{(1+\e)^2 }{c_{35}^*\sigma \sqrt{2\pi}}
    \int_0^\infty u^p e^{-u} H(u)\dd u
    <\infty,
\end{eqnarray*}
since $H$ is a nondecreasing function and $H(x)=O(x^2)$ as $x\to +\infty$
e.g. by (\cite{Vatutin_Wachtel} Lem. 13 with $\alpha=2$ and $\rho=1/2$ as explained at the end of its p. 181,
following from Rogozin \cite{Rogozin}
% and with its $I_3$ being slowly varying at infinity)
and from the Spitzer-R\`ozen theorem).
%, Theorem 2).

This, combined with \eqref{eqEsperanceLaplaceVxGrand} gives
$$
    \limsup_{x\to+\infty}
    \Big(x\E\big[\big(V(x)\big)^p e^{-V(x)} | G_x^*\big]\Big)
\leq
    %\frac{g_{2,0}(0)}{\sigma c_{35}^*}
    \frac{1}{c_{35}^*\sigma \sqrt{2\pi}}
    \int_0^\infty u^p e^{-u} H(u)\dd u.
$$
Since we get a similar inequality for $\liminf$ and $H\geq 1$ on $]0,\infty[$,
this proves \eqref{eqEquivalentLaplaceV} and the lemma.
\hfill$\Box$

\begin{lem} \label{LemLaplaceSachantC_Lattice}
Assume that $\log\frac{1-\o_0}{\o_0}$ is $(h,a)$-lattice for some $h>0$ and $a\in[0,h[$.
We have for $p\geq 0$,
\begin{equation}\label{eqEquivalentLaplaceV_Lattice}
    \E\big[\big(V(x)\big)^p e^{-V(x)} | G_x^*\big]
\sim_{x\to+\infty}
    %\frac{h g_{2,0}(0) }{c_{35}^* \sigma x}
    \frac{h }{c_{35}^* \sqrt{2\pi} \sigma x}
    \psi_p[(a x) \textnormal{ mod } h],
\end{equation}
where
$
    \psi_p(y)
:=
    \sum_{k\in\N} (y +k h)^p e^{-(y +k h)} H(y +k h)
$, $y\in[0,h]$,
is a function bounded on $[0,h]$ between two (strictly) positive constants.
\end{lem}

\noindent{\bf Proof:}
Let $p\geq 0$, $h>0$ and $a\in[0,h[$,
%$p\in\{0,1\}$,
and assume that $\log\frac{1-\o_0}{\o_0}$ is $(h,a)$-lattice.
First, notice that for every $y\in[0,h]$,
$
    \psi_p(y)
\leq
    \sum_{k\in\N}(h+k h)^p e^{-k h }H(h+kh )
    =e^h\psi_p(h)
<
    \infty
$
since
$H$ is nondecreasing and
$H(x) = O(x^2)$ as $x\to+\infty$
as in the previous lemma.
%$H(x)\sim_{x\to+\infty} C x$ for some $C>0$
%by (\cite{Vatutin_Wachtel}, Lem. 13 with $\alpha=2$,  and with $\rho=1/2$ as explained in its page 181).
Moreover, taking into account only $k=1$, we have
$
    \psi_p(y)
\geq
    h^p e^{-2h}H(h)
>
    0
$
for every $y\in[0,h]$,
so $\psi_p$ is bounded on $[0,h]$ between two (strictly) positive constants.

Let $\e>0$.
Applying (\cite{Vatutin_Wachtel}, Theorem 6, extending previous results obtained when $a=0$ by Alili and Doney \cite{Alili_Doney}),
again with $\alpha=2$, $\beta=0$, $c_x\sim_{x\to+\infty}\sigma\sqrt{x}$,
%as seen  in the line after its eq. (3))
 and $g_{2,0}(0)=1/\sqrt{2\pi}$:
\begin{equation}\label{Eq_VW_lattice}
    \sigma\sqrt{x} \p[V(x)=a x+y \mid G_x^*]
\sim_{x\to+\infty}
    %h g_{2,0}(0)
    \frac{h H(a x +y)}{\sqrt{2\pi}x\p[G_x^*]}
\end{equation}
uniformly in $y\in]-a x, -a x +\delta_x \sqrt{x}]\cap (h\Z)$, where $\delta_x\to 0$ as $ x\to+\infty$.
Also, notice that for $y=-ax$ when $x>0$, we have
$
%    \p[V(x)=a x+y \mid C_x^*]
%=
    \p[V(x)=0 \mid G_x^*]
=
    0
=
h H(0)/[\sqrt{2\pi}x\p(G_x^*)\sigma\sqrt{x}]
    %h g_{2,0}(0)  H(0)/[x\p(G_x^*)\sigma\sqrt{x}]
$
by definitions of $G_x^*$ and $H$.
Hence for large $x$,
\begin{eqnarray*}
&&
    \E\big[(V(x))^p e^{-V(x)} \un_{\{V(x)< 9\log x\}}\, \big|\, G_x^*\big]
\\
& = &
    \sum_{k\in\Z,\ a x+k h\geq 0} (a x+k h)^p e^{-(a x+k h)}  \un_{\{a x + k h< 9\log x\}} \p[V(x)=a x+k h \mid G_x^*]
\\
& \leq &
    \sum_{k\in\Z,\ 0\leq a x+k h <9\log x } (a x+k h)^p e^{-(a x+k h)}
    %(1+\e)h g_{2,0}(0)
    (1+\e)
    \frac{h H(a x +k h)}{\sqrt{2\pi}\sigma x^{3/2}\p[G_x^*]}
\\
& \leq &
    %\frac{(1+\e)h g_{2,0}(0) }
    \frac{(1+\e)h }
    {\sqrt{2\pi} \sigma x^{3/2}\p[G_x^*]}
    \psi_p[(a x) \text{ mod } h]
\leq
    %\frac{(1+2\e)h g_{2,0}(0) }
    \frac{(1+2\e)h }
    {c_{35}^*  \sqrt{2\pi}\sigma x}
    \psi_p[(a x) \text{ mod } h]
\end{eqnarray*}
by \eqref{Eq_VW_lattice} applied with $\delta_x=9(\log x)/\sqrt{x}$ and \eqref{eq_Proba_CxPrime}.
This and \eqref{eqEsperanceLaplaceVxGrand}
 give for large $x$,
\begin{equation}\label{Ineg_Laplace_Cas_Arithmetique1}
     \E\big[\big(V(x)\big)^p e^{-V(x)} | G_x^*\big]
\leq
    %\frac{(1+2\e)h g_{2,0}(0) }
    \frac{(1+2\e)h }
    {c_{35}^* \sqrt{2\pi}\sigma x}
    \psi_p[(a x) \text{ mod } h]
    +\beta_p x^{-8}.
\end{equation}
Similarly as in \eqref{Ineg_Laplace_Cas_Arithmetique1},  for large $x$,
$$
     \E\big[\big(V(x)\big)^p e^{-V(x)} | G_x^*\big]
\geq
    \frac{(1-2\e)h}{c_{35}^* \sqrt{2\pi} \sigma x}
    \Big(\psi_p[(a x) \text{ mod } h]
    -O(x^{-8})
    \Big),
$$
since
$
    \sum_{k\in\Z,\ a x+k h \geq 9\log x } (a x+k h)^p e^{-(a x+k h)}
    H(a x +k h)
=
    O(x^{-8})
$
as $x\to+\infty$
because
$H(x)=O(x^2)$ as in the previous lemma.
This and \eqref{Ineg_Laplace_Cas_Arithmetique1} prove \eqref{eqEquivalentLaplaceV_Lattice}
since $x^{-8}=o\big(\psi_p[(a x) \text{ mod } h]/x\big)$ as $x\to+\infty$
because $\inf_{[0,h]} \psi_p>0$.
\hfill$\Box$

{\bf Proof of Proposition \ref{Lemma_Laplace_V_Conditionne}:}
Let $h>0$ and $x\in\N^*$.
We first provide a relation between conditioning by $\Xi_h^*$ and by $G_x^*$.
We have,
%for $x> 0$ and $n\in\N^*$,
due to the Markov property,
\begin{eqnarray}
&&
    \E\big[ e^{-V(x)}\un_{\{x< T_V(h)\}} \big| \Xi_h^*\big]
%\nonumber\\
%& = &
=
    \E\big[ e^{-V(x)}\un_{\{x< T_V(h)\}} \un_{\{T_V(h)<T_V^*(\R_-)\}}\big]/\p[\Xi_h^*]
\nonumber\\
& = &
    \E\big[ e^{-V(x)}\un_{\{\forall 0< k\leq x,\ 0< V(k)<h\}}
        \un_{\forall k\in[x,T_V(h)],\ V(k)> 0\}}\big]/\p[\Xi_h^*]
\nonumber\\
& = &
    \E\big[ e^{-V(x)}\un_{\{\forall 0< k\leq x,\ 0< V(k)<h\}}
        \p^{V(x)}\big(T_V(h)<T_V^*(\R_-)\big)\big]/\p[\Xi_h^*].
\label{Ineg_Laplace_V_sachant_An_1_etoile}
\end{eqnarray}
Hence for large $h>0$, for every $x\in\N^*$, by \eqref{Ineg_Laplace_V_sachant_An_1_etoile} and Lemma \ref{Lem_Proba_An_Ordre2} eq. \eqref{Ineg_Pz_An_Etoile},
\begin{align}
    \E\big[ e^{-V(x)}\un_{\{x< T_V(h)\}} | \Xi_h^*\big]
& \leq
    \E\bigg[ e^{-V(x)}\un_{\{\forall 0< k\leq x,\ 0< V(k)\}}
    \frac{V(x)-\E^{V(x)}[V(T_V^*(\R_-))]}{h\p[\Xi_h^*]}\bigg]
\nonumber\\
& =
    \frac{\p[G_x^*]}{h\p[\Xi_h^*]}
    \E\Big[  \big[V(x)-\E^{V(x)}[V(T_V^*(\R_-))]\big] e^{-V(x)} \, \big| \, G_x^*\Big].
~~~~
\label{eqEsperanceSachantCPrime_Etoile}
\end{align}
Also,
\begin{equation*}
    \E\big[ e^{-V(x)}\un_{\{x= T_V(h)\}} | \Xi_h^*\big]
\leq
    \E\big[ e^{-V(x)}\un_{\{V(x)\geq h\}} | \Xi_h^*\big]
\leq
    e^{-h}.
\end{equation*}
Let $\e>0$.
By \eqref{eqEsperanceSachantCPrime_Etoile}, \eqref{eq_Proba_CxPrime}  and \eqref{eq_Proba_An_Etoile},
then by ellipticity \eqref{eqEllipticity},
there exists $p_6>0$ and $p_7>0$ such that for $x\geq p_6$ and $h\geq p_7$,
\begin{eqnarray}\label{InegLaplaceSachantAn_Etolie}
%    \E\bigg[ \frac{\un_{\{x< T_V(h)\}}}{e^{V(x)}} | \Xi_h^*\bigg]
    \E\big[ e^{-V(x)} \un_{\{x< T_V(h)\}} | \Xi_h^*\big]
& \leq &
    \E\Big[ \big[V(x)-\E^{V(x)}\big[V(T_V^*(\R_-))\big]\big] e^{-V(x)} | G_x^*\Big]
    \frac{(1+\e)c_{35}^*}{c_1^* x^{1/2}}
    ~~~~~~~
\\
& \leq &
    (1+\e)(c_{35}^*/c_1^*) x^{-1/2}
    \big[\E\big(V(x)e^{-V(x)} | G_x^*\big) + C_0 \E\big(e^{-V(x)} | G_x^*\big) \big].
\nonumber
\end{eqnarray}
Thanks to Lemmas \ref{LemLaplaceSachantC} and \ref{LemLaplaceSachantC_Lattice},
there exists $p_4>p_6$ such  that, for $x\geq p_4$, for each $p\in\{0,1\}$,
$
    \E\big[\big(V(x)\big)^p e^{-V(x)} | G_x^*\big]
\leq
    \frac{f_3(p)}{x}
$,
with $f_3(p):=2f_2(p)$ when $\log\frac{1-\o_0}{\o_0}$ is  non lattice,
and $f_3(p):=2h \sup_{[0,h]}\psi_p/(c_{35}^*\sqrt{2\pi}\sigma)$
if $\log\frac{1-\o_0}{\o_0}$ is $(h,a)$-lattice for some $h>0$ and $a\in[0,h[$.
This together with \eqref{InegLaplaceSachantAn_Etolie}
gives
for $x\geq p_4$ and $h\geq p_7$,
\begin{eqnarray}
    \E\big[ e^{-V(x)}\un_{\{x< T_V(h)\}} | \Xi_h^*\big]
& \leq &
    (1+\e)(c_{35}^*/c_1^*)
    \big[ f_3(1)+ C_0 f_3(0) \big]
    x^{-3/2}.
\label{Ineg_Laplace_An_Etoile}
\end{eqnarray}

We now aim to prove a similar inequality, conditioning by $\Xi_h$ instead of $\Xi_h^*$.
There exists $c>0$ such that $\p[V(1)\in[c,2c]   ]>0$,
thanks to \eqref{eqRecurrence} and \eqref{eq_def_sigma}.
For such a (fixed) $c$, there exists $p_8\geq p_7$ such that for all $h\geq p_8$,  we have
$h/10>2c$,
$\p(\Xi_{h+2c}^*)/\p(\Xi_{h})\leq 2c_1^*/c_1$ (by Lemma \ref{Lem_Proba_An_Ordre2})
and
$
    \p\big[T_V(h/10)
        <
        T_V( ]-\infty,-h/10])
    \big]
\geq
    1/3
$
(e.g. by \eqref{eqOptimalStopping2}).
So with $\widetilde V_1(k):=V(k+1)-V(1)$, $k\geq 0$,
using the independence of $V(1)$ and $\widetilde V_1$,
then the independence of $(V(u), \ u\leq T_V(h))$
and $\widetilde V_2$,
defined by $\widetilde V_2(k):=V[T_V(h)+k]-V[T_V(h)]$, $k\geq 0$,
we have for $h\geq p_8$ and for $x\geq p_4$,
\begin{eqnarray*}
&&
    \E\big[ e^{-V(x+1)}\un_{\{x+1< T_V(h+2c)\}} | \Xi_{h+2c}^*\big]
=
    \E\big[ e^{-V(x+1)}\un_{\{x+1< T_V(h+2c)\}} {\bf 1}_{\Xi_{h+2c}^*}\big]/\p(\Xi_{h+2c}^*)
\\
& \geq &
    \E\big[ e^{-V(1)-\widetilde V_1(x)}\un_{\{V(1)\in [c,2c]\}}
        \un_{\{x< T_{\widetilde V_1}(h)\}}
        {\bf 1}_{\forall y\in[1,T_{\widetilde V_1}(h+2c)],\ \widetilde V_1(y)\geq 0}\big]/\p(\Xi_{h+2c}^*)
\\
& \geq &
    \frac{e^{-2c}\p[V(1)\in[c,2c]]}{\p(\Xi_{h+2c}^*)}
    \E\big[ e^{-V(x)}\un_{\{x< T_{V}(h)\}}
        {\bf 1}_{\forall y\in[1,T_{V}(h+2c)],\ V(y)\geq 0}\big]
\\
& \geq &
    \frac{e^{-2c}\p[V(1)\in[c,2c]]}{\p(\Xi_{h+2c}^*)}
    \E\big[ e^{-V(x)}\un_{\{x< T_V(h)\}}
    %{\bf 1}_{\forall y\in[1,T_{V}(\log m)],\ V(y)\geq 0}
    \un_{\Xi_{h}}
    \un_{T_{\widetilde V_2}(h/10)<T_{\widetilde V_2}( ]-\infty,-h/10])}\big]
\\
& \geq &
    \frac{\p[V(1)\in[c,2c]]\p(\Xi_{h})}{e^{2c}\p(\Xi_{h+2c}^*)}
    \E\bigg[ \frac{\un_{\{x< T_V(h)\}}}{e^{V(x)}} | \Xi_{h}\bigg]
    \p\bigg[T_V\bigg(\frac{h}{10}\bigg)<T_V\bigg( \bigg]-\infty,-\frac{h}{10}\bigg]\bigg)\bigg].
%\\
%& \geq &
%    \p[V(1)\in[c,2c]]
%    \E\big[ e^{-V(x)}\un_{\{x< T_V(h)\}} | \Xi_{h}\big]
%    c_1/(6c_1^*e^{2c})
\end{eqnarray*}
%Since the last probability is greater than $1/3$ for large $h$, e.g. by \eqref{eqOptimalStopping2},
%and using Lemma \ref{Lem_Proba_An_Ordre2}
%the definition of $n(m)$
So, using the definition of $p_8$ then \eqref{Ineg_Laplace_An_Etoile}, we get
with $c_{36}:=6 e^{2c}c_1^*/(c_1\p[V(1)\in[c,2c]])$,
for every $x\geq p_4$ and $h\geq p_8$,
\begin{equation*}
    \E\big[ e^{-V(x)}\un_{\{x< T_V(h)\}} \mid \Xi_{h}\big]
\leq
    c_{36}
    \E\big[ e^{-V(x+1)}\un_{\{x+1< T_V(h+2c)\}} | \Xi_{h+2c}^*\big]
\leq
    c_{13} x^{-3/2}
\end{equation*}
for some constant $c_{13}>0$.
This and \eqref{Ineg_Laplace_An_Etoile} prove Proposition \ref{Lemma_Laplace_V_Conditionne},
up to a change of $c_{13}$.
\hfill$\Box$

%We define, for $z\in\R$, $z \textnormal{ mod}^* h$ as the only element of $]0, h]\cap (z+h\Z)$
%(beware that we replace $0$ by $h$ in the standard definition of $z$ mod $h$).

\subsection{Two lemmas about left $h$-extrema}
For the sake of completeness, we prove the two following lemmas.
We recall that $\VVV$ is defined before \eqref{eqDefbh}.

\begin{lem}\label{Lemma_Alternate}
Let
$v\in\VVV$,
%$v$ be a function : $\Z\to \R$,
and let $h>0$.
The  left (resp. right) $h$-minima and left (resp. right) $h$-maxima for $v$ alternate.
\end{lem}

\noindent{\bf Proof:} Assume that $y_1$ and $y_2$ are two left $h$-minima for $v$, with $y_1<y_2$.
It is enough to prove that there exists at least a left $h$-maximum for $v$ between $y_1$ and $y_2$.
By Definition \ref{def_left_extrema},
for each $j\in\{1,2\}$, there exists $\alpha_j<y_j<\beta_j$ such that
$\min_{[\alpha_j,y_j-1]}v>v(y_j)$,
$\min_{[y_j+1,\beta_j]}v\geq v(y_j)$,
$v(\alpha_j)\geq v(y_j)+h$ and
$v(\beta_j)\geq v(y_j)+h$.
We define $x:=\min\{u\geq y_1,\ v(u)=\max_{[y_1,\  y_2]}v\}$.
The goal is to prove that $x$ is a left $h$-maximum for $v$.

Assume that $y_2\leq \beta_1$.
If $\alpha_2\leq y_1$, then
$\alpha_2\leq y_1<y_2\leq \beta_1$, so
$v(y_2)\geq \min_{[y_1+1,\beta_1]}v\geq v(y_1)$
and
$v(y_1)\geq \min_{[\alpha_2,y_2-1]}v>v(y_2)$,
which contradicts $v(y_2)\geq  v(y_1)$.
So $\alpha_2> y_1$, thus
$y_1<\alpha_2<y_2\leq \beta_1$.
%Let $x:=\min\{u\geq y_1,\ v(u)=\max_{[y_1,\  y_2]}v\}$.
We have
$v(x)=\max_{[y_1,\  y_2]}v \geq v(\alpha_2)\geq v(y_2)+h$
and
$v(x)\geq v(y_2)+h\geq \min_{[y_1+1, \beta_1]} v+h\geq v(y_1)+h$.
%Also by definition of $x$,
%$y_1<x<y_2$,
%$\max_{[y_1,x-1]}v<v(x)$
%and
%$\max_{[x+1,y_2]}v\leq v(x)$,
%so
%$x$ is a left $h$-maximum such that $y_1<x<y_2$.

Now, assume that $y_2> \beta_1$ and $\alpha_2\leq y_1$.
Thus,
$\alpha_2\leq y_1<\beta_1<y_2$, so
$v(y_1)\geq \min_{[\alpha_2,y_2-1]}v>v(y_2)$.
%Once more, let $x:=\min\{u\geq y_1,\ v(u)=\max_{[y_1,\  y_2]}v\}$.
We have $v(x)=\max_{[y_1,y_2]}v\geq v(\beta_1) \geq v(y_1)+h \geq v(y_2)+h$.
%Also, as before,
%$y_1<x<y_2$,
%$\max_{[y_1,x-1]}v<v(x)$
%and
%$\max_{[x+1,y_2]}v\leq v(x)$,
%so
%$x$ is a left $h$-maximum such that $y_1<x<y_2$.

Finally, assume that $y_2> \beta_1$ and $\alpha_2 > y_1$.
Hence,
$y_1<\beta_1<y_2$, so
$v(x)=\max_{[y_1,y_2]}v \geq v(\beta_1) \geq v(y_1)+h$.
Also,
$y_1<\alpha_2<y_2$, so
$v(x)=\max_{[y_1,y_2]}v \geq v(\alpha_2) \geq v(y_2)+h$.

So in every case, we have
$v(x)\geq v(y_1)+h$ and
$v(x)\geq v(y_2)+h$, with $h>0$, thus  by definition of $x$,
$y_1<x<y_2$,
$\max_{[y_1,x-1]}v<v(x)$
and
$\max_{[x+1,y_2]}v\leq v(x)$,
so
$x$ is a left $h$-maximum for $v$ such that $y_1<x<y_2$.
%, which concludes the proof for left $h$-minima and left $h$-maxima.

Applying this to $-v$ proves that, if $y_1$ and $y_2$ are two left $h$-maxima for $v$ with $y_1<y_2$,
there exists at least a left $h$-minimum for $v$ between $y_1$ and $y_2$, which concludes the proof of the lemma for left $h$-extrema.
The proof is similar for right ones by symmetry.
\hfill$\Box$

For the following lemma, see definitions \eqref{eqDef_tau_1_V}--\eqref{eqDef_m_2_V}, represented in Figure \ref{figure_tau_i_m_i_slopes}.

\begin{lem}\label{Lemma_Only_h_extrema}
Assume that $V\in\VVV$ (which has probability one if \eqref{eqEllipticity}, \eqref{eqRecurrence} and \eqref{eq_def_sigma} are satisfied).
{\bf (i)}
For $i\geq 1$, $m^{(V)}_{2i+1}(h)$ is a left $h$-minimum for $V$,
and there is no other left $h$-extremum for $V$ in $\big[\tau^{(V)}_{2i}(h),\tau^{(V)}_{2i+1}(h)\big[$.
{\bf (ii)}
For $i\geq 0$, $m^{(V)}_{2i+2}(h)$ is a left $h$-maximum for $V$,
and there is no other left $h$-extremum for $V$ in $\big[\tau^{(V)}_{2i+1}(h),\tau^{(V)}_{2i+2}(h)\big[$.
\end{lem}

\noindent{\bf Proof:}
Let $i\geq 1$.
First,
$
    m^{(V)}_{2i}(h)
<
    m^{(V)}_{2i+1}(h)
<
    \tau^{(V)}_{2i+1}(h)
$
by definition.
We also have
$
    V\big( \tau^{(V)}_{2i+1}(h)\big)
\geq
    V\big( m^{(V)}_{2i+1}(h)\big)
    +h
$
by \eqref{eqDef_tau_1_V} and \eqref{eqDef_m_1_V}
and
$
    V\big( m^{(V)}_{2i}(h)\big)
\geq
    V\big( \tau^{(V)}_{2i}(h)\big)
    +h
\geq
    V\big( m^{(V)}_{2i+1}(h)\big)
    +h
$
since $i\geq 1$
by \eqref{eqDef_m_2_V}, \eqref{eqDef_tau_2_V} and \eqref{eqDef_m_1_V}.
Also,
$
    \min_{[m^{(V)}_{2i+1}(h)+1,\tau^{(V)}_{2i+1}(h)]}V
\geq
    V\big(m^{(V)}_{2i+1}(h)\big)
$
by \eqref{eqDef_m_1_V},
$
    \min_{[m^{(V)}_{2i}(h), \tau^{(V)}_{2i}(h)-1]}V
>
    V\big(\tau^{(V)}_{2i}(h)\big)
\geq
    V\big(m^{(V)}_{2i+1}(h)\big)
$
by \eqref{eqDef_tau_2_V}, \eqref{eqDef_m_2_V} and \eqref{eqDef_m_1_V},
and
$
    \min_{[\tau^{(V)}_{2i}(h),m^{(V)}_{2i+1}(h)-1]}V
>
    V\big(m^{(V)}_{2i+1}(h)\big)
$
by \eqref{eqDef_m_1_V}.
So, $m^{(V)}_{2i+1}(h)$ is a left $h$-minimum for $V$.

\noindent{\bf First case:}
Assume that there exists a left $h$-minimum $y\neq m^{(V)}_{2i+1}(h)$ for $V$
in $\big[\tau^{(V)}_{2i}(h),$ $\tau^{(V)}_{2i+1}(h)\big[$,
and let $\alpha<y$ and $\beta>y$ be as in Definition \ref{def_left_extrema} with $v=V$.
Assume first that $y\in\big[\tau^{(V)}_{2i}(h), m^{(V)}_{2i+1}(h)\big[$.
If $\beta<m^{(V)}_{2i+1}(h)$, then
$
    V(\beta)
\geq
    V(y)+h
%>
%    V\big(m^{(V)}_{2i+1}(h)\big)
%    +h
$
with $\tau^{(V)}_{2i}(h)\leq y <\beta<\tau^{(V)}_{2i+1}(h)$,
which contradicts the definition of $\tau^{(V)}_{2i+1}(h)$.
If $\beta\geq m^{(V)}_{2i+1}(h)$, then
$y+1\leq m^{(V)}_{2i+1}(h)\leq \beta$
so
$
    V\big(m^{(V)}_{2i+1}(h)\big)
\geq
    \min_{[y+1,\beta]}V
\geq
    V(y)
$,
which contradicts
$
    V(y)
\geq
    \min_{[\tau^{(V)}_{2i}(h),m^{(V)}_{2i+1}(h)-1]}V
>
    V\big(m^{(V)}_{2i+1}(h)\big)
$
by \eqref{eqDef_m_1_V}.

So $y\in\big] m^{(V)}_{2i+1}(h), \tau^{(V)}_{2i+1}(h)\big[$.
If $\alpha>m^{(V)}_{2i+1}(h)$, then
$
    V(\alpha)
\geq
    V(y)+h
\geq
    V\big(m^{(V)}_{2i+1}(h)\big)+h
$
with $\tau^{(V)}_{2i}(h)\leq m^{(V)}_{2i+1}(h) < \alpha <\tau^{(V)}_{2i+1}(h)$,
which contradicts the definition of $\tau^{(V)}_{2i+1}(h)$.
If $\alpha \leq m^{(V)}_{2i+1}(h)$, then, since $y>m^{(V)}_{2i+1}(h)$, we have
$
    V\big(m^{(V)}_{2i+1}(h)\big)
\geq
    \min_{[\alpha,y-1]}V
>
    V(y)
$
by definition of $\alpha$,
which contradicts
$
    V(y)
\geq
    \min_{[\tau^{(V)}_{2i}(h),\tau^{(V)}_{2i+1}(h)]}V
=
    V\big(m^{(V)}_{2i+1}(h)\big)
$
by \eqref{eqDef_m_1_V}.
So there is no left $h$-minimum for $V$ in $\big[\tau^{(V)}_{2i}(h),\tau^{(V)}_{2i+1}(h)\big[
-\big\{m^{(V)}_{2i+1}(h)\big\}$.

\noindent{\bf Second case:}
Now, we assume that there exists a left $h$-maximum $y$
%$y\neq m^{(V)}_{2i+1}(h)$
for $V$ in $\big[\tau^{(V)}_{2i}(h)$, $\tau^{(V)}_{2i+1}(h)\big[$,
and let $\alpha<y$ and $\beta>y$ be as in Definition \ref{def_left_extrema} for left $h$-maxima.
%Assume first that $y\in\big[\tau^{(V)}_{2i}(h), m^{(V)}_{2i+1}(h)\big[$.
If $\alpha\geq \tau^{(V)}_{2i}(h)$, then
$
    V(y)
\geq
    V(\alpha)+h
$
with
$\tau^{(V)}_{2i}(h)\leq \alpha < y < \tau^{(V)}_{2i+1}(h)$,
which contradicts the definition of $\tau^{(V)}_{2i+1}(h)$.
If $\alpha< \tau^{(V)}_{2i}(h)$, then
$
    V(\alpha)
\leq
    V(y)-h
<
    V\big(\tau^{(V)}_{2i}(h)\big)
$
by definition of $\tau^{(V)}_{2i+1}(h)$
since $\tau^{(V)}_{2i}(h)\leq y <\tau^{(V)}_{2i+1}(h)$.
So if $m^{(V)}_{2i}(h) \leq \alpha< \tau^{(V)}_{2i}(h)$, then
$
    V(\alpha)
<
    V\big(\tau^{(V)}_{2i}(h)\big)
$
contradicts $V\big(\tau^{(V)}_{2i}(h)\big)<\min_{[m^{(V)}_{2i}(h),\tau^{(V)}_{2i}(h)[}V\leq V(\alpha)$,
coming from \eqref{eqDef_tau_2_V} and \eqref{eqDef_m_2_V} since $i\geq 1$.
Finally if $\alpha<m^{(V)}_{2i}(h)$, then
$\alpha<m^{(V)}_{2i}(h)<\tau^{(V)}_{2i}(h)\leq y$
by \eqref{eqDef_tau_2_V} and \eqref{eqDef_m_2_V}
since $i\geq 1$,
so
$
    V\big(m^{(V)}_{2i}(h)\big)
\leq
    \max_{[\alpha, y-1]}V
<
    V(y)
<
    V\big(\tau^{(V)}_{2i}(h)\big)+h
$
by definition of $\alpha$ and \eqref{eqDef_tau_1_V} since $y\in\big[\tau^{(V)}_{2i}(h),\tau^{(V)}_{2i+1}(h)\big[$,
which contradicts
$
    V\big(\tau^{(V)}_{2i}(h)\big)
\leq
    V\big(m^{(V)}_{2i}(h)\big)-h
$
coming from \eqref{eqDef_tau_2_V} and \eqref{eqDef_m_2_V} since $i\geq 1$.
So there is no left $h$-maximum for $V$ in $\big[\tau^{(V)}_{2i}(h),\tau^{(V)}_{2i+1}(h)\big[
%-\big\{m^{(V)}_{2i+1}(h)\big\}
$ for $i\geq 1$.

Thus {\bf (i)} is proved.
The proof of {\bf (ii)} is similar.
\hfill$\Box$

%%%%%%%%%%%%%%%%%%%%%%%%%%%%%%%%%%%%%%%%%%%%%%%%%%%%%%%%%%%%%%%%%%%%%%%%%%%%%%%%%%%%%%

%%%%%%%%%%%%%%%%%%%%%%%%%%%%%%%%%%%%%%%%%%%%%%%%%%%%%%%%%%%%%%%%%%%%%%%%%%%%%%%%%%%%%%%%
%                                                                                      %
%                                   BIBLIOGRAPHY                                       %
%                                                                                      %
%%%%%%%%%%%%%%%%%%%%%%%%%%%%%%%%%%%%%%%%%%%%%%%%%%%%%%%%%%%%%%%%%%%%%%%%%%%%%%%%%%%%%%%%

\bibliographystyle{alpha}

\begin{thebibliography}{07}


\bibitem{Afanasyev_et_al_2012}
    \textsc{Afanasyev, V. I., B\"{o}inghoff, C., Kersting, G. and Vatutin, V. A.}:
     Limit theorems for weakly subcritical branching processes in random environment.
    {\it J. Theoret. Probab.} {\bf 25} (2012), 703--732.

\bibitem{Alili_Doney}%[AD99]
    \textsc{Alili, L. and  Doney, R.A.}:
    Wiener-Hopf factorization revisited and some applications.
    {\it Stoc. Stoc. Rep.} {\bf 66} (1999), 87-102.




\bibitem{AndreolettiDevulder}%[AD15]
    \textsc{Andreoletti, P. and Devulder, A.}:
    Localization and number of visited valleys for a transient diffusion in random environment.
    {\it Electron. J. Probab.} {\bf 20}, no 56 (2015), 1--58.

\bibitem{AndreolettiDevulderVechambre}%[ADV16]
    \textsc{Andreoletti, P., Devulder, A. and V\'echambre, G.}:
    Renewal structure and local time for diffusions in random environment.
    {\it ALEA, Lat. Am. J. Probab. Math. Stat.} {\bf 13} (2016), 863--923.


\bibitem{Andres_Taylor}
    \textsc{Andres, S. and Taylor, P. A.}:
    Local limit theorems for the random conductance model and applications to the
    {G}inzburg-{L}andau {$\nabla\phi$} interface model.
    {\it J. Stat. Phys.} {\bf 182} (2021), paper no. 35.

\bibitem{AurzadaDevulderGuillotinPene}%[ADGP17]
    \textsc{Aurzada, F., Devulder, A., Guillotin-Plantard, N. and P\`ene, F.}:
    Random walks and branching processes in correlated Gaussian environment.
    {\it J. Stat. Phys.} {\bf 166} (2017), 1--23.


\bibitem{Barbu_Limnios}%[BL08]
    \textsc{Barbu, V. S. and Limnios, N.}:
     Semi-{M}arkov chains and hidden semi-{M}arkov models toward applications, Their use in reliability and DNA analysis.
    Lecture Notes in Statistics {\bf 191},
    Springer, New York, 2008.

\bibitem{Berger_CR}
    \textsc{Berger, N. , Cohen, M.  and Rosenthal, R. }:
    Local limit theorem and equivalence of dynamic and static points of view for certain ballistic random walks in i.i.d.
    environments.
    {\it Ann. Probab.} {\bf 44} (2016),  2889--2979.

\bibitem{Bertoin_Split}%[B93]
    \textsc{Bertoin, J.}:
    Splitting at the infimum and excursions in half-lines for random walks and L\'evy processes.
    {\it Stochastic Process. Appl.} {\bf 47} (1993), %no. 1,
    17--35.
    %MR-1232850


\bibitem{Bingham_Goldie_Teugels}%[BGT87]
    \textsc{Bingham N.H., Goldie C.M. and    Teugels J.L.}:
    Regular variation. Cambridge: Cambridge University Press, 1987, 494 pp.


\bibitem{Bovier_Faggionato}%[BF08]
    \textsc{Bovier, A. and Faggionato, A.}:
    Spectral analysis of {S}inais walk for small eigenvalues.
   {\it Ann. Probab.} {\bf 36} (2008), 198--254.


\bibitem{Brox}%[B86]
    \textsc{Brox, Th.}:  A one-dimensional diffusion process in a {W}iener medium.
    {\it Ann. Probab.} {\bf 14} (1986),  1206--1218.

\bibitem{Buraczewski_Dyszewski}
    \textsc{Buraczewski, D. and Dyszewski, P.}:
    Precise large deviations for random walk in random environment.
    {\it Electron. J. Probab.}  {\bf 23}, no 114 (2018), 1--26.

\bibitem{Cheliotis}%[C05]
    \textsc{Cheliotis, D.}:  Diffusion in random environment and the renewal theorem.
              {\it Ann. Probab.}  {\bf 33} (2005), 1760--1781.


\bibitem{Cheliotis_Favorite}%[C08]
    \textsc{Cheliotis, D.}:  Localization of favorite points for diffusion in a random environment.
    {\it Stoch. Proc. Appl.} {\bf 118} (2008), 1159--1189.

\bibitem{Chiarini}
    \textsc{Chiarini, A. and Deuschel, J.-D.}:
    Local central limit theorem for diffusions in a degenerate and unbounded random medium.
    {\it Electron. J. Probab.} {\bf 20}, no. 112 (2015), 1--30.



\bibitem{CoccoMonasson}%[CM08]
    \textsc{Cocco, S. and Monasson, R.}: Reconstructing a random potential from its random walks.
              {\it Europhysics Letters}  {\bf 81} (2008), 20002.

\bibitem{CometsGantertZeitouni}%[CGZ00]
    \textsc{Comets, F., Gantert, N. and Zeitouni, O.}:
    Quenched, annealed and functional large deviations for one-dimensional random walk in random environment.
    {\it Probab. Theory Related Fields} {\bf 118} (2000), 65--114.

\bibitem{CometsPopov}%[CP03]
    \textsc{Comets, F. and Popov, S.}:
    Limit law for transition probabilities and moderate deviations
    for Sinai's random walk in random environment.
    {\it Probab. Theory Related Fields} {\bf 126} (2003), 571--609.


\bibitem{Dembo_Gantert_Peres_Shi}%[DGPS07]
    \textsc{Dembo, A., Gantert, N.,  Peres, Y. and Shi, Z.}:
    Valleys and the maximum local time for random walk in random environment.
    \textit{Probab. Theory Related Fields} \textbf{137} (2007), 443--473.



\bibitem{Devulder_Persistence}%[D16]
    \textsc{Devulder, A.}:
    Persistence of some additive functionals of Sinai's walk.
    {\it Ann. Inst. H. Poincar\'e Probab. Stat.} {\bf 52},  No. 3 (2016), 1076--1105.


\bibitem{Devulder_SPL}%[D07]
    \textsc{Devulder, A.}: The speed of a branching system of random walks in random environment. {\it Statist. Probab. Lett.} {\bf 77} (2007), 1712--1721.

\bibitem{Devulder_Rates_CV}%[D22+]
    \textsc{Devulder, A.}:
    Rates of convergence in Sinai and Golosov localization theorems for random walks in random environments.
    Work in progress, (2023+).

\bibitem{DGP_Collision_Sinai}%[DGP18]
    \textsc{Devulder, A., Gantert N. and Pene F.}:
    Collisions of several walkers in recurrent random environments.
    {\it Electron. J. Probab.} {\bf 23}, no. 90 (2018), 1--34.

\bibitem{DGP_Collision_Transient}%[DGP19]
    \textsc{Devulder, A., Gantert N. and Pene F.}:
    Arbitrary many walkers meet infinitely often in a subballistic random environment.
    {\it Electron. J. Probab.} {\bf 24}, no. 100 (2019), 1--25.


\bibitem{Diel_Stat}%[D11]
    \textsc{Diel, R. and Lerasle, M.}:
    Non parametric estimation for random walks in random environment.
    {\it Stochastic Process. Appl.}  {\bf 128}, no 1 (2018),  132--155.


\bibitem{Dolgo_gold_13}%[DG13]
    \textsc{Dolgopyat, D. and Goldsheid, I.}:
    Local Limit Theorems for Random Walks in a 1D Random Environment.
    {\it Archiv der Mathmathematik} {\bf 101} (2013), 191--200.

\bibitem{Dolgo_gold_19}%[DG19]
    \textsc{Dolgopyat, D. and Goldsheid, I.}:
    Local Limit Theorems for Random Walks in a Random Environment on a Strip.
    {\it Pure Appl. Funct. Anal.} {\bf 5} (2020), 1297--1318.

\bibitem{Dolgo_gold_21}%[DG21]
    \textsc{Dolgopyat, D. and Goldsheid, I.}:
    Constructive approach to limit theorems for recurrent diffusive random walks on a strip.
    {\it Asymptot. Anal.} {\bf 122} (2021), 271--325.

\bibitem{Doyle_Snell}%[DS84]
    \textsc{Doyle, P. G. and Snell, E. J.}:
    Probability: Random walks and Electrical Networks.
    Carus Math. Monographs {\bf 22}, Math. Assoc. Amer., Washington DC, 1984.


\bibitem{ESZ2}%[ESZ09a]
    \textsc{Enriquez, N., Sabot, C. and Zindy, O.}:
    Limit laws for transient random walks in random environment on Z.
    {\it Ann. Inst. Fourier (Grenoble)}  {\bf 59} (2009), 2469--2508.


\bibitem{ESZ3}%[ESZ09b]
    \textsc{Enriquez, N., Sabot, C. and Zindy, O.}:
    Aging and quenched localization for one dimensional random walks in random
    environment in the sub-ballistic regime.
    {\it Bull. Soc. Math. France}  {\bf 137} (2009), 423-452.

\bibitem{Freire}%[F15]
    \textsc{Freire, M.V.}:
    Application of Moderate Deviation Techniques to Prove Sinai Theorem on RWRE
    {\it J. Stat. Phys.}
    {\bf 160} (2) (2015), 357--370.


\bibitem{NMFa}%[GKP14]
    \textsc{Gantert, N., Kochler M. and P\`ene, F.}:
    On the recurrence of some random walks in random environment.
    \textit{ALEA} \textbf{11} (2014), 483--502.



\bibitem{GPS_2010}
    \textsc{Gantert, N., Peres, Y. and Shi, Z.}:
    The infinite valley for a recurrent random walk in random environment.
    {\it Ann. Inst. Henri Poincar\'e Probab. Stat.} {\bf 46} (2010), 525--536.


\bibitem{Gantert_Peterson}%[GP11]
    \textsc{Gantert, N. and Peterson, J.}:
    Maximal displacement for bridges of random walks in a random environment.
   {\it Ann. Inst. Henri Poincar\'{e} Probab. Stat.} {\bf 47} (2011), 663--678.


\bibitem{Gantert_Sshi_2002}
    \textsc{Gantert N.  and Shi Z.}:
    Many visits to a single site by a transient random walk in
    random environment.
    {\it Stoch. Proc. Appl.} {\bf 99} (2002), 159--176.

\bibitem{Gantert_Zeitouni}%[GP11]
    \textsc{Gantert, N. and Zeitouni, O.}:
    Large deviations for one-dimensional random walk in a random environment---a survey.
    Random walks ({B}udapest, 1998), 127--165,
    Bolyai Soc. Math. Stud. {\bf 9},
    J\'{a}nos Bolyai Math. Soc., Budapest, 1999.

\bibitem{Golosov84} %[G84]
    \textsc{Golosov, A. O.}:
    Localization of random walks in one-dimensional random environments.
    \textit{Commun. Math. Phys.} \textbf{92} (1984), 491--506.

\bibitem{Golosov86}%[G86]
    \textsc{Golosov, A. O.}:
    Limit distributions for random walks in random environments. {\it Soviet Math. Dokl.} {\bf 28} (1986), 18--22.


\bibitem{Geoffrey_Stirzaker_3}
    \textsc{Grimmett, G. R. and Stirzaker, D. R.}:
    {\it Probability and random processes}.
    Oxford University Press, New York, third edition, 2001.

\bibitem{Hirano_98}
    \textsc{Hirano, K.}:
    Determination of the limiting coefficient for exponential functionals of random walks with positive drift.
    {\it J. Math. Sci. Univ. Tokyo} {\bf 5} (1998), 299--332.


\bibitem{Hoeffding63}%[H63]
    \textsc{Hoeffding, W.}:
    Probability inequalities for sums of bounded random variables.
    {\it J. Amer. Statist. Assoc.} {\bf 58} (1963), 13--30.


\bibitem{HuLocal}%[H00]
    \textsc{Hu Y.}:
    Tightness of localization and return time in random environment.
    {\it Stoch. Proc. Appl.} {\bf 86} (2000), 81--101.



\bibitem{HuShiLimits}
    \textsc{Hu, Y. and Shi, Z.}:
    The limits of {S}inai's simple random walk in random environment.
    {\it Ann. Probab.} {\bf 26} (1998), 1477--1521.

\bibitem{Hu_Shi_Problem}
    \textsc{Hu, Y. and Shi, Z.}:
    The problem of the most visited site in random environment.
    {\it Probab. Theory Related Fields} {\bf 116} (2000), 273--302.


\bibitem{Hug}%[H96]
    \textsc{Hughes, B.D.}: {\it Random Walks and Random Environment,
    vol. II: Random Environments}. Oxford Science Publications,
    Oxford, 1996.


\bibitem{Kesten}%[K86]
    \textsc{Kesten, H.}:  The limit distribution of Sinai's random walk in random environment.
    {\it Physica} 138A (1986), 299--309.

\bibitem{These_Jochler}%[K12]
        \textsc{Kochler M.}:
    Random Walks in Random Environment, Random Orientations and Branching.
    Doktors der Naturwissenschaften genehmigten Dissertation,
    Technische Universit{\"a}t M{\"u}nchen (2012).


\bibitem{KMT}%[KMT75]
    \textsc{Koml\'os, J., Major, P. and Tusn\'ady, G.}:
    An approximation of partial sums of independent rv's and the sample df. I.
    \textit{Wahrsch verw Gebiete/Probability Theory and Related Fields} \textbf{32} (1975), 111--131



\bibitem{Levin_Peres}
    \textsc{Levin, D. A., Peres, Y. and Wilmer, E. L.}:
    Markov chains and mixing times.
    American Mathematical Society, Providence, RI, with a chapter by James G. Propp and David B.Wilson, 2009.


\bibitem{Leskela_Stenlund}%[LS11]
    \textsc{Leskela L. and Stenlund M.}:
    A local limit theorem for a transient chaotic walk in a frozen environment.
    {\it Stochastic Process. Appl.} {\bf 121} (2011),  2818--2838.


\bibitem{Nauenberg}%[N85]
    \textsc{Nauenberg, M.}:
    Random walk in a random medium in one dimension.
    {\it J Stat Phys} {\bf 41} (1985), 803--810.

\bibitem{NP}%[NP89]
    \textsc{Neveu J. and Pitman J.}:
    Renewal property of the extrema and tree property of the
              excursion of a one-dimensional {B}rownian motion.
              {\it S\'eminaire de {P}robabilit\'es XXIII, Lecture Notes in Math.}   {\bf 1372}, 239--247,
              Springer, Berlin, 1989.


\bibitem{Padash}
    \textsc{Padash, A., Aghion, E., Schulz, A., Barkai, E., Chechkin, A. V., Metzler, R. and Kantz, H.}:
    Local equilibrium properties of ultraslow diffusion in the Sinai model.
    {\it New J. Phys.} {\bf 24} (2022), July, Paper No. 073026.

\bibitem{Revesz}%[R05]
    \textsc{R{\'e}v{\'e}sz, P.}:
    {\it Random walk in random and non-random environments}, second edition.
    World Scientific, Singapore, 2005.

\bibitem{Rogozin}%[R71]
    \textsc{Rogozin, B.A.}:
    On the distrbution of the first ladder moment and height and fluctustions of a random walk.
    {\it Theory Probab. Appl.} {\bf 16} (1971), 575--595.


\bibitem{Shi_98}
    \textsc{Shi, Z.}:
    A local time curiosity in random environment.
    {\it Stochastic Process. Appl.} {\bf 76} (1998), 231--250.


\bibitem{S2}%[S01]
     \textsc{Shi, Z.}:
     Sinai's walk via stochastic calculus. {\it Panoramas et Synth\`eses}
     {\bf 12} (2001), 53--74,
     Soci\'et\'e math\'ematique de France.



\bibitem{ShiZindy}%[SZ07]
    \textsc{Shi, Z. and Zindy, O.}:
    A weakness in strong localization for {S}inai's walk.
    \textit{Ann. Probab.}, \textbf{35} (2007), 1118--1140.


\bibitem{S82}%[S82]
    \textsc{Sinai, Ya. G.}:
    The limiting behavior of a one-dimensional random walk in a random
    medium.
    {\it Th. Probab. Appl.}, {\bf 27} (1982), 256--268.



\bibitem{S75}%[S75]
    \textsc{Solomon, F.}:
    Random walks in a random environment.
    {\it Ann. Probab.}, {\bf 3} (1975), 1--31.

\bibitem{Takenami}
    \textsc{Takenami, T.}:
    Local limit theorem for random walk in periodic environment.
    {\it Osaka J. Math.} {\bf 39} (4) (2002), 867--895.



\bibitem{Vatutin_Wachtel}%[VW09]
    \textsc{Vatutin, V. A. and Wachtel, V.}:
    Local probabilities for random walks conditioned to stay positive.
    {\it Probab. Theory Related Fields} {\bf 143} (2009), 177--217.


\bibitem{Z01}%[Z01]
    \textsc{Zeitouni, O.}:
    Lecture notes on random walks in random environment.
    {\it \'Ecole d'\'et\'e de probabilit\'es de
    Saint-Flour 2001}. Lecture Notes in
    Math. {\bf 1837},
    %pp.~
    189--312. Springer, Berlin, 2004.

\bibitem{Zindy}%[Z08]
    \textsc{Zindy, O.}:
    Upper limits of Sinai's walk in random scenery.
    {\it Stoch. Proc. Appl.} {\bf 118} (2008), 981--1003.


\end{thebibliography}

\end{document}